\documentclass[11pt,twoside]{article}
\newcommand{\ifims}[2]{#1} 
\newcommand{\ifAMS}[2]{#1}   
\newcommand{\ifau}[4]{#1}  
\newcommand{\ifbook}[2]{#1}   
\newcommand{\ifunivariate}[2]{#1} 
\newcommand{\ifapp}[2]{#2}  
\newcommand{\ifKL}[2]{#2}  
\newcommand{\ifinexact}[2]{#2}  

\usepackage{ifthen}
\usepackage[ampersand]{easylist}
\usepackage{amsmath,amssymb,amsthm}
\usepackage{mathalpha}
\usepackage{mathtools}
\usepackage{natbib}
\usepackage{epsfig,graphicx}
\usepackage{comment}
\usepackage{color}
\usepackage{srcltx}
\usepackage[mathscr]{eucal}
\usepackage[math]{easyeqn}
\usepackage{etoolbox}
\usepackage{hyperref}
\hypersetup{
            colorlinks,
            linkcolor=hookersgreen,
            linktoc=blue,
            citecolor=ultramarine,
            urlcolor=black,
            filecolor=black
            }
\usepackage{nameref}
\usepackage{multirow}

\usepackage{accents}
\usepackage{mathrsfs}
\usepackage{bbm}
\usepackage{algorithm}
\usepackage{algpseudocode}

\ifims{
\textheight=22cm
\textwidth=14.8cm
\topmargin=0pt
\oddsidemargin=1.0cm
\evensidemargin=1.0cm
\linespread{1.3}

}{ 
}

\numberwithin{equation}{section}
\numberwithin{figure}{section}
\newcounter{example}[section]
\numberwithin{example}{section}
\newcounter{remark}[section]
\numberwithin{remark}{section}
\newtheorem{theorem}{Theorem}[section]
\newtheorem{proposition}[theorem]{Proposition}
\newtheorem{lemma}[theorem]{Lemma}
\newtheorem{corollary}[theorem]{Corollary}

\newtheorem{exmp}[example]{Example}
\newtheorem{rmrk}[remark]{Remark}
\newenvironment{example}{\begin{exmp}\rm}{\end{exmp}}
\newenvironment{remark}{\begin{rmrk}\rm}{\end{rmrk}}

\ifbook{
    \newcommand{\Chapter}[1]{\section{#1}}
    \newcommand{\Section}[1]{\subsection{#1}}
    \newcommand{\Subsection}[1]{\subsubsection{#1}}
    \def\Chname{Section }

  }
  {
    \newcommand{\Chapter}[1]{\chapter{#1}}
    \newcommand{\Section}[1]{\section{#1}}
    \newcommand{\Subsection}[1]{\subsection{#1}}
    \def\Chname{Chapter}
    
  }

\renewcommand{\(}{$\,}
\renewcommand{\)}{\,$}

\def\nquad{\hspace{-1cm}}
\def\eqdef{\stackrel{\operatorname{def}}{=}}
\def\eqd{\stackrel{\operatorname{d}}{=}}

\DeclareMathAlphabet{\mathbbmsl}{U}{bbm}{bx}{sl}

\DeclareMathSymbol{\Alpha}{\mathalpha}{operators}{"41}
\DeclareMathSymbol{\Beta}{\mathalpha}{operators}{"42}
\DeclareMathSymbol{\Epsilon}{\mathalpha}{operators}{"45}
\DeclareMathSymbol{\Zeta}{\mathalpha}{operators}{"5A}
\DeclareMathSymbol{\Eta}{\mathalpha}{operators}{"48}
\DeclareMathSymbol{\Iota}{\mathalpha}{operators}{"49}
\DeclareMathSymbol{\Kappa}{\mathalpha}{operators}{"4B}
\DeclareMathSymbol{\Mu}{\mathalpha}{operators}{"4D}
\DeclareMathSymbol{\Nu}{\mathalpha}{operators}{"4E}
\DeclareMathSymbol{\Omicron}{\mathalpha}{operators}{"4F}
\DeclareMathSymbol{\Rho}{\mathalpha}{operators}{"50}
\DeclareMathSymbol{\Tau}{\mathalpha}{operators}{"54}
\DeclareMathSymbol{\Chi}{\mathalpha}{operators}{"58}
\DeclareMathSymbol{\omicron}{\mathord}{letters}{"6F}

\newcommand{\cc}[1]{\mathscr{#1}}
\newcommand{\bb}[1]{\boldsymbol{#1}}

\renewcommand{\bar}[1]{%
  \hbox{%
    \vbox{%
      \hrule height 0.5pt 
      \kern0.35ex
      \hbox{%
        \kern-0.05em
        \ensuremath{#1}%
        \kern-0.05em
      }%
    }%
  }%
} 

\makeatletter
\newcommand*\rel@kern[1]{\kern#1\dimexpr\macc@kerna}
\newcommand*\widebar[1]{%
  \begingroup
  \def\mathaccent##1##2{%
    \rel@kern{0.8}%
    \overline{\rel@kern{-0.8}\macc@nucleus\rel@kern{0.2}}%
    \rel@kern{-0.2}%
  }%
  \macc@depth\@ne
  \let\math@bgroup\@empty \let\math@egroup\macc@set@skewchar
  \mathsurround\z@ \frozen@everymath{\mathgroup\macc@group\relax}%
  \macc@set@skewchar\relax
  \let\mathaccentV\macc@nested@a
  \macc@nested@a\relax111{#1}%
  \endgroup
}
\makeatother

\renewcommand{\hat}[1]{\widehat{#1}}
\renewcommand{\tilde}[1]{\widetilde{#1}}

\makeatletter
\def\mathcenterto#1#2{\mathclap{\phantom{#1}\mathclap{#2}}\phantom{#1}}
\let\old@widetilde\widetilde
\def\widetildeto#1#2{\mathcenterto{#2}{\old@widetilde{\mathcenterto{#1}{#2\,}}}}
\let\old@widehat\widehat
\def\widehatto#1#2{\mathcenterto{#2}{\old@widehat{\mathcenterto{#1}{#2\,}}}}
\makeatother

\newcommand{\thankstitle}[1]{\ifthenelse{\equal{#1}{}}{}{\thanks{#1}}}
\newcommand{\thanksau}[1]{\ifthenelse{\equal{#1}{}}{}{\thanks{#1}}}

\newcommand{\aua}[6]
{\def\authora{#1}
\def\runauthora{#2}
\def\addressa{#3}
\def\emaila{#4}
\def\affiliationa{#5}
\def\thanksa{#6}}

\def\theauthors{
\ifau{ 
  \author{
    \authora
    \thanksau{\thanksa}
    \\[5.pt]
    \addressa \\
    \texttt{ \emaila}
  }
}
{  
  \author{
    \authora
    \thanksau{\thanksa}
    \\[5.pt]
    \addressa \\
    \texttt{ \emaila}
    \and
    \authorb
    \thanksau{\thanksb}
    \\[5.pt]
    \addressb \\
    \texttt{ \emailb}
  }
}
{   
  \author{
    \authora
    \thanksau{\thanksa}
    \\[5.pt]
    \addressa \\
    \texttt{ \emaila}
    \and
    \authorb
    \thanksau{\thanksb}
    \\[5.pt]
    \addressb \\
    \texttt{ \emailb}
    \and
    \authorc
    \thanksau{\thanksc}
    \\[5.pt]
    \addressc \\
    \texttt{ \emailc}
  }
} {   
  \author{
    \authora
    \thanksau{\thanksa}
    \\[5.pt]
    \addressa \\
    \texttt{ \emaila}
    \and
    \authorb
    \thanksau{\thanksb}
    \\[5.pt]
    \addressb \\
    \texttt{ \emailb}
    \and
    \authorc
    \thanksau{\thanksc}
    \\[5.pt]
    \addressc \\
    \texttt{ \emailc}
    \and
    \authord
    \thanksau{\thanksd}
    \\[5.pt]
    \addressd \\
    \texttt{ \emaild}
  }
}
}

\renewcommand{\Gamma}{\varGamma}
\renewcommand{\Pi}{\varPi}
\renewcommand{\Sigma}{\varSigma}
\renewcommand{\Delta}{\varDelta}
\renewcommand{\Lambda}{\varLambda}
\renewcommand{\Psi}{\varPsi}
\renewcommand{\Phi}{\varPhi}
\renewcommand{\Theta}{\varTheta}
\renewcommand{\Omega}{\varOmega}
\renewcommand{\Xi}{\varXi}
\renewcommand{\Upsilon}{\varUpsilon}

\def\argmax{\operatornamewithlimits{argmax}}
\def\argmin{\operatornamewithlimits{argmin}}

\def\av{\bb{a}}

\def\ev{\bb{e}}
\def\fv{\bb{f}}

\def\uv{\bb{u}}

\def\wv{\bb{w}}
\def\xv{\bb{x}}

\def\zv{\bb{z}}

\def\Av{\bb{A}}
\def\Bv{\bb{B}}

\def\Uv{\bb{U}}

\def\Xv{\bb{X}}
\def\Yv{\bb{Y}}
\def\Zv{\bb{Z}}

\def\deltav{\bb{\delta}}

\def\epsv{\bb{\varepsilon}}
\def\etav{\bb{\eta}}
\def\gammav{\bb{\gamma}}

\def\lambdav{\bb{\lambda}}

\def\thetav{\bb{\theta}}

\def\xiv{\bb{\xi}}

\def\Psiv{\bb{\Psi}}

\def\sumi{\sum_{i=1}^{n}}

\usepackage{color}

\definecolor{blue(pigment)}{rgb}{0.2, 0.2, 0.6}
\definecolor{ultramarine}{rgb}{0.07, 0.04, 0.56}
\definecolor{darkspringgreen}{rgb}{0.09, 0.45, 0.27}
\definecolor{hookersgreen}{rgb}{0.0, 0.44, 0.0}
\definecolor{plum(traditional)}{rgb}{0.56, 0.27, 0.52}
\definecolor{purple(html/css)}{rgb}{0.5, 0.0, 0.5}
\definecolor{magenta(dye)}{rgb}{0.79, 0.08, 0.48}

\def\supH{\lambda}
\def\rexH{\epsilon}

\def\dCross{\mathcal{G}}
\def\dPartial{\mathcal{H}}

\def\discr{\delta}

\def\prmt{\ups}

\def\prmtv{\bb{\prmt}}

\def\prmtvs{\prmtv^{*}}

\def\targ{x}
\def\targv{\bb{\targ}}
\def\targvs{\targv^{*}}
\def\targvb{\breve{\targv}}

\def\XXL{\mathcal{X}}

\def\CFT{\cc{C}}

\def\CDG{\CONSTi_{\GP}}

\def\hmax{\mathsf{c}}

\def\hL{h}

\def\elll{\ell}

\def\lgd{f}
\def\lgdb{\breve{\lgd}}

\def\lgdL{\elll}

\def\PfL{\P_{\lgd}}

\def\dagg{\prime}

\def\proj{P}

\def\blk{\operatorname{block}}

\def\amax{\nu}

\def\Matr{\mathfrak{M}}

\def\vH{\vA}

\def\Idd{\Gamma}

\def\WV{\mathcal{W}}

\def\feta{\phi}

\def\Rem{\mathcal{R}}

\def\Eta{\mathcal{H}}

\def\HS{\cc{H}}

\def\HV{\mathsf{H}}
\def\HVB{\mathbbmsl{V}} 

\def\HL{\mathbb{m}}

\def\smp{s}

\def\dltw{\delta}

\def\dltwb{\omega}

\def\dltwbs{\dltwb^{*}}
\def\dltwm{\kappa}

\def\dltwu{\tau}
\def\dltwd{\dltw^{\dagg}}
\def\dltwbd{\dltwb^{\dagg}}

\def\dltwbss{\dltwb^{+}}

\def\II{\mathcal{I}}

\def\gblk{H}
\def\fgblk{J}

\def\AFblk{\fgblk}

\def\AFblkn{\mathtt{j}}

\def\HFblk{\gblk}

\def\R{\mathbbmsl{R}}
\def\E{\mathbbmsl{E}}

\def\P{\mathbbmsl{P}}

\def\kappa{\varkappa}

\def\Frobg{\Lambda}

\def\blk{\operatorname{block}}

\def\diag{\operatorname{diag}}

\def\Fr{\operatorname{Fr}}

\def\ND{\mathcal{N}}

\def\oper{\operatorname{op}}

\def\Var{\operatorname{Var}}

\def\T{\top}
\def\tr{\operatorname{tr}}

\def\TV{\operatorname{TV}}

\def\cond{\, \big| \,}

\def\nsize{{n}}

\def\sumi{\sum_{i=1}^{\nsize}}

\def\ex{\mathrm{e}}

\def\Id{I\!\!\!I}
\def\Ind{\operatorname{1}\hspace{-4.3pt}\operatorname{I}}

\def\bias{\mathsf{b}}

\def\biasv{\bb{b}}

\def\biasQ{\bias}

\def\biase{\mathbb{b}}

\def\biasev{\biase}

\def\BB{\mathbbmsl{B}}

\def\BBB{\cc{B}}

\def\BBH{B}

\def\Bv{\bb{\BB}}

\def\CA{\mathcal{A}}

\def\CGP{w}

\def\CONST{\mathtt{C} \hspace{0.1em}}
\def\CONSTi{\mathtt{C}}

\def\CONSTIFT{\CONSTi_{\IFT}}
\def\CONSTgp{\CONSTi_{\gp}}

\def\CONSTV{\CONSTi_{\VP}}

\def\CONSTdv{\CONSTi_{0}}
\def\CONSTbias{\CONSTi_{\operatorname{bias}}}
\def\CONSTIF{\CONSTi_{\IF}}
\def\CONSTIFT{\CONSTi_{\IFT}}

\def\CONSTdlt{\kappa}

\def\DP{D}

\def\DV{\mathsf{D}}

\def\DVb{\breve{\DV}}

\def\DVL{\DV}

\def\DVLf{\mathbbmss{D}}

\def\DU{\mathbb{D}}

\def\deltav{\bb{\delta}}

\def\dimA{\mathtt{p}}

\def\dimtotal{\dimp^{*}}

\def\dimAs{{\dimA^{*}}}

\def\dimAs{\bar{\dimA}}

\def\dimL{\dimA}
\def\dimLe{\dimA_{\ex}}
\def\dimqLe{\dimqL_{\ex}}

\def\dimLs{\bar{\dimL}}

\def\dimH{\dimA}

\def\dimp{p}

\def\dimG{\dimA_{\GP}}

\def\dimQ{\dimA_{\QP}}

\def\dimq{q}

\def\dimqL{\mathsf{q}}

\def\dimd{d}

\def\dimLL{\dimL(\xvs)}
\def\dens{f}

\def\EfL{\E_{\lgd}}

\def\Eta{\cc{H}}

\def\err{\diamondsuit}

\def\errs{\err^{*}}

\def\errE{\err_{\Nui}}

\def\eps{\epsilon}			
\def\eps{\varepsilon}

\def\fvs{\fv}  

\def\fs{f}

\def\fba{\bar{f}}
\def\fn{g}

\def\fG{f_{\GP}}

\def\fGu{h}

\def\gaussv{\bb{\gauss}}
\def\gauss{\gamma}

\def\gp{g}

\def\gpks{\kappa}

\def\GP{G}

\def\GPY{\Gamma}
\def\GPKS{\mathcal{K}}

\def\GPa{\GP_{0}}

\def\GPT{\mathcal{G}}

\def\hKS{\hat{\KS}}

\def\IF{\mathbbmsl{F}}

\def\IFL{{\mathbbmss{F}}}
\def\IFLap{\breve{\IFL}}

\def\IFtotal{F}
\def\IFT{\mathscr{\IFtotal}}
\def\IFTb{\breve{\IFT}}

\def\IFL{\mathbbmss{F}}
\def\IFU{\breve{\IFL}}

\def\IFLba{\accentset{\circ}{\IFL}}
\def\IFCov{\Sigma}
\def\IFCovs{\IFCov^{*}}

\def\ima{z}
\def\imav{\bb{\ima}}

\def\KS{A}
\def\KSs{\KS^{*}}

\def\Kappa{\cc{K}}

\def\kullb{\cc{K}} 

\def\LT{L}
\def\LGP{\LT_{\GP}}

\def\LL{\cc{L}}

\def\lambdav{\bb{\lambda}}

\def\mm{m}

\def\muA{\mu}

\def\muH{\mu}

\def\normG{\alpha}
\def\nui{a}
\def\Nui{\mathsf{A}}
\def\nuiv{\bb{\nui}}
\def\nuiiv{\bb{\omega}}

\def\Nuiv{\bb{\Nui}}
\def\nuivs{\nuiv^{*}}
\def\Nuiv{\bb{\Nui}}
\def\nuo{\tau}
\def\nuov{\bb{\nuo}}

\def\Phiv{\bb{\Phi}}

\def\Proj{\Pi}

\def\priord{\pi}

\def\QP{Q}

\def\rhot{t}
\def\rhos{\rho^{*}}

\def\rhoIF{\rho}

\def\riskt{\cc{R}}

\def\rr{\mathtt{r}}

\def\rrbs{\bar{\rr}}

\def\rru{\rr_{\circ}}

\def\rrn{\rr}

\def\rrL{\rr}
\def\rrLs{\bar{\rr}}

\def\thetav{\bb{\theta}}

\def\Tau{T}

\def\uvc{\uv^{c}}

\def\ups{\upsilon}
\def\upsv{\bb{\ups}}

\def\upsvd{\upsv^{\circ}}
\def\upsvs{\upsv^{*}}

\def\upsvn{\upsvd}

\def\ups{\upsilon}
\def\upsv{\bb{\ups}}

\def\UV{\mathcal{U}}

\def\UVL{\UV}

\def\UVz{\UV}

\def\Ups{\varUpsilon}
\def\Upsd{\Ups^{\circ}}

\def\Upsb{\breve{\Ups}}
\def\vA{\mathtt{v}}

\def\VP{V}

\def\VV{\mathbbmsl{V}}

\def\wv{\bb{w}}

\def\xvs{\xv^{*}}
\def\xvb{\bar{\xv}}

\def\xx{\mathtt{x}}

\def\yy{\mathtt{y}}

\def\zq{z}

\def\zqe{\mathsf{z}}
\def\thetitle{Mixed Laplace approximation for marginal posterior and Bayesian inference in error-in-operator model}
\def\theruntitle {Mixed Laplace approximation}

\def\theabstract{
Laplace approximation is a very useful tool in Bayesian inference and it claims 
a nearly Gaussian behavior of the posterior.
\cite{SpLaplace2022} established some rather accurate finite sample results about
the quality of Laplace approximation in terms of the so called effective dimension \( \dimL \)
under the critical dimension constraint \( \dimL^{3} \ll n \).
However, this condition can be too restrictive for many applications like error-in-operator problem or Deep Neuronal Networks.
This paper addresses the question whether the dimensionality condition can be relaxed 
and the accuracy of approximation can be improved
if the target of estimation is low dimensional while the nuisance parameter is high or infinite dimensional. 
Under mild conditions, the marginal posterior can be approximated by a Gaussian mixture and the accuracy of the approximation only
depends on the target dimension.
Under the condition \( \dimL^{2} \ll n \) or in some special situation like semi-orthogonality, 
the Gaussian mixture can be replaced by one Gaussian distribution leading to a classical Laplace result.
The second result greatly benefits from the recent advances in Gaussian comparison from \cite{GNSUl2017}.
The results are illustrated and specified for the case of error-in-operator model.
}

\def\kwdp{62F15}
\def\kwds{62F25}

\def\thekeywords{critical dimension, mixed Gaussian approximation, Gaussian comparison, error-in-operator}

\def\thankstitle{}

\aua
{Vladimir Spokoiny}
{Spokoiny, V. }
{
    Weierstrass Institute and HU Berlin,  
    \\
    Mohrenstr. 39, 10117 Berlin, Germany
    }
{spokoiny@wias-berlin.de}
{Weierstrass-Institute and Humboldt University Berlin}
{
  Financial support by the German Research Foundation (DFG) through the Collaborative Research Center 1294 ``Data assimilation'' is gratefully acknowledged.
}

\def\HL{\mathsf{m}}
\def\biase{\mathsf{s}}
\def\biasev{\bb{s}}

\begin{document}
\thispagestyle{empty}
{
\title{\thetitle}
\theauthors

\maketitle
\begin{abstract}
{\footnotesize \theabstract}
\end{abstract}

\ifAMS
    {\par\noindent\emph{AMS 2010 Subject Classification:} Primary \kwdp. Secondary \kwds}
    {\par\noindent\emph{JEL codes}: \kwdp}

} 

\tableofcontents
\Chapter{Introduction}
\label{Sgenintr}


Laplace approximation is one of the most powerful instruments in Bayesian inference.
Let \( \lgd(\cdot) \) be a non-normalized posterior log-density, usually of the form
\begin{EQA}
	\lgd(\prmtv)
	&=&
	L(\prmtv) + \log \priord(\prmtv) 
\label{hyc7whjrfby63esswqw}
\end{EQA}
for the log-likelihood function \( L(\prmtv) \) and a prior density \( \priord(\prmtv) \).
Laplace approximation claims that, for a strongly concave function \( \lgd(x) \), 
the posterior measure \( \PfL \) with a density proportional to
\( \ex^{\lgd(\prmtv)} \) can be well approximated by a Gaussian distribution
\( \ND(\upsvs,\IFL^{-1}) \), where \( \prmtvs = \argmax_{\prmtv} \lgd(\prmtv) \) and \( \IFL = -\nabla^{2} \lgd(\prmtvs) \).
The original Laplace result was stated for the univariate case.
Recently this topic attracted a lot of attention in connection with statistical inference 
for complex high dimensional model such as nonlinear inverse problems, Deep Neuronal Networks, optimal transport, etc.
We refer to recent papers 
\cite{HeKr2022},
\cite{KGB2022},
\cite{KaRi2022},
\cite{nickl2022bayesian},
and references therein.
Classical asymptotic setup with the fixed prior and a growing sample size leads to the prominent Bernstein--von Mises Theorem
that claims prior free asymptotic normal posterior distribution 
centered at the MLE \( \tilde{\prmtv} = \argmax_{\prmtv} L(\prmtv) \).
We refer to \cite[Chapter 12]{GhVa2017} for a detailed presentation and literature overview.
However, the case of large or even huge parameter dimension requires a separate study, the Bernstein--von Mises phenomenon
does not apply any more, the prior impact could be significant. 
A particular issue is the \emph{critical dimension} meaning the dimensionality constraint under which the Laplace approximation
still applies.
\cite{SpLaplace2022} established rather precise nonasymptotic bounds on the quality of Laplace approximation
in terms of the so called \emph{Laplace effective dimension} \( \dimLs \); see Section~\ref{SeffdimLa} for a precise definition.
The main result of \cite{SpLaplace2022} requires \( \dimLs^{3} \ll n \), where \( n \) is the sample size;
see e.g. Theorem~\ref{TLaplaceTV34}. 
For high dimensional models this constraint appears to be very limiting.
Unfortunately, it seems to be impossible to relax this condition in a general situation.
This paper focuses on a slightly different problem.
Let the parameter vector \( \prmtv \) be high dimensional but we are only interested in its low dimensional component \( \targv \).
The main question under consideration is whether the accuracy of approximation of the marginal 
\( \targv \)-posterior can be controlled in terms of the \( \targv \)-dimension \( \dimLe \) 
thus relaxing the dimensionality constraint.
The whole study can be split into two big steps.
The first result of Theorem~\ref{Ppartintegrsmi} describes an approximation of the marginal distribution of the subvector \( \targv \) 
by a Gaussian mixture.
A nice feature of this result is that the accuracy of approximation only depends on the effective dimension \( \dimLe \) 
of the target variable \( \targv \) only. 
The next step of the study aims at addressing the question if the obtained mixture can be replaced by one Gaussian distribution.
The main tools include Gaussian comparison \cite{GNSUl2017} and a bias reduction in conditional optimization based on one-point orthogonalization. 
It appears that even in a rather general situation, one can relax the condition \( \dimLs^{3} \ll n \) 
to \( \dimLs^{2} \ll n \).
In some special cases like semi-orthogonality, the full parameter dimension does not show up and the marginal \( \targv \)-posterior
can be approximated by Gaussian law with accuracy \( \sqrt{\dimLe^{3}/n} \).

\subsubsection*{Motivation 1: error-in-operator problem}
\label{SCalmexampl}
Let \( \KS \) be a linear mapping (operator) of the source signal \( \targv \in \R^{\dimp} \) to the image space \( \R^{\dimq} \).
We consider the problem of inverting the relation 
\( \zv = \KS \targv \): given an image vector \( \zv \in \R^{\dimq} \), recover the corresponding 
source \( \targv \in \R^{\dimp} \).
This leads to the linear least square problem of maximizing 
the negative fidelity function \( \lgdL(\cdot) \) of the form 
\begin{EQA}
	\lgdL(\targv) 
	&=& 
	- \frac{1}{2 \sigma^{2}} \| \zv - \KS \targv \|^{2} ,
\label{o3kv6ejfgiuhbdfjkse}
\end{EQA}
where \( \sigma^{2} \) describes the image noise.
However, this approach assumes the operator \( \KS \) to be precisely known.
In many applications, this assumption is not fulfilled and the operator \( \KS \) is known up to some error.
This means that only an estimate \( \hKS \) of \( \KS \) is available.
Particular examples include all kinds of tomography, regression with random design,
instrumental regression, 
error-in-variable regression,
functional data analysis, etc.
A pragmatic plug-in approach just replaces \( \KS \) in \eqref{o3kv6ejfgiuhbdfjkse} by its estimate \( \hKS \) leading 
to the solution
\begin{EQA}
	\hat{\targv}
	&=&
	\argmin_{\targv} \| \zv - \hKS \targv \|^{2}
	=
	\bigl( \hKS^{\T} \hKS \bigr)^{-1} \hKS^{\T} \zv .
\label{wkvoigb7y634hjgtkuiif}
\end{EQA}
Unfortunately, inverting a large matrix \( \hKS^{\T} \hKS \) can be tricky, 
especially if the operator \( \KS \) is smooth.
\cite{HoRe2008} considered simultaneous wavelet estimation of the signal \( \targv \) and the operator \( \KS \).
\cite{Trabs_2018} extended these results to Bayesian inference using a parametric assumption
about the unknown operator \( \KS = \KS_{\thetav} \) and provided some examples
from imaging and linear PDE/SDEs. 
In that paper, one can also find a literature overview and more references on this topic.
This paper treats the original problem in a semiparametric setup by including the operator 
\( \KS \) in the parameter set \( \prmtv = (\targv,\KS) \). 
This leads to the new fidelity function
\begin{EQA}[c]
	- \frac{1}{2 \sigma^{2}} \| \zv - \KS \targv \|^{2} - \frac{\muA^{2}}{2 \sigma^{2}} \| \hKS - \KS \|_{\Fr}^{2} \, , 
\label{lgbhu3hdf6vjghhksdemde}
\end{EQA}
where \( \muA \) describes the level of the operator noise and \( \| \KS \|_{\Fr}^{2} = \tr(\KS^{\T} \KS) \) stands
for the squared Frobenius norm.
The use of a Gaussian prior on \( \targv \) leads to the posterior log-density
\begin{EQA}
	\lgd(\prmtv)
	&=&
	- \frac{1}{2 \sigma^{2}} \| \zv - \KS \targv \|^{2} 
	- \frac{\muA^{2}}{2 \sigma^{2}} \| \hKS - \KS \|_{\Fr}^{2} 
	- \frac{1}{2} \| \GP (\targv - \targv_{0}) \|^{2} ,
\label{dhfdfededrsdsrdsresswd}
\end{EQA} 
where the precision operator \( \GP^{2} \) is responsible for the smoothness of the source signal while 
\( \targv_{0} \) stands for a starting guess. 
The quadratic term \( - \muA^{2} (2 \sigma^{2})^{-1} \| \hKS - \KS \|_{\Fr}^{2} \)
can be viewed as a Gaussian log-density on \( \KS \).
The function \( \lgd \) is not concave in the set of variables \( \prmtv = (\targv,\KS) \).
However, it is locally concave in a vicinity of the point \( (\targv_{0},\hKS) \) under ``warm start'' condition 
which requires a reasonable quality of the guess \( \targv_{0} \) and of the pilot \( \hKS \);
see Section~\ref{SLaplErrOp} later.
The main objective is the marginal of the full dimensional posterior distribution related to the target parameter \( \targv \).
Some results about semiparametric Bernstein--von Mises Theorem are available; see e.g. 
\cite{Ca2012},
\cite{BiKl2012},
\cite{Cast2023},
and references therein.
\cite{CaRo2015} considered the BvM result for smooth functionals.
However, the results focused on asymptotic rate of contraction and asymptotic Gaussian approximation of the posterior
in some metric like \( \ell_{2} \)-Wasserstein.
This paper, similarly to \cite{Sp2022}, rather attempts to explore the quality of Laplace approximation
and to address the issue of \emph{critical dimension} of the nuisance parameters.
The results of Section~\ref{SsemiLaplace} state a surprising phase-transition result:
an increase of the dimensionality of the nuisance parameter may lead to a mixed Laplace approximation 
in place of usual Laplace approximation.
This requires to extend the results from \cite{SpLaplace2022} to the marginal distribution. 

\subsubsection*{Motivation 2. Random design regression}
\label{Srandomdesign}
In many applications, the design \( X_{1},\ldots,X_{n} \) in regression models 
\( Y_{i} = \fs(X_{i}) + \eps_{i} \) may be naturally assumed to be random, often i.i.d.
This particularly concerns econometric, biological, chemical, sociological studies etc. 
Under linear parametric assumption \( \fs(\xv) = \sum_{j} \theta_{j} \psi_{j}(\xv) = \Psi(\xv)^{\T} \thetav \) and noise homogeneity, 
the model log-likelihood reads exactly as in the case of a deterministic design
\begin{EQA}
	L(\thetav)
	&=&
	- \frac{1}{2 \sigma^{2}} \| \Yv - \Psiv(\Xv)^{\T} \thetav \|^{2} + \CONST ,
\label{4556799kbvdfee4566y} 
\end{EQA}
where \( \Psiv(\Xv) \) is the \( \dimp \times n \) matrix with entries \( \psi_{j}(X_{i}) \).
The corresponding MLE \( \tilde{\thetav} \) is again given by
\begin{EQA}
	\tilde{\thetav}
	&=&
	\bigl\{ \Psiv(\Xv) \, \Psiv(\Xv)^{\T} \bigr\}^{-1} \Psiv(\Xv) \Yv .
\label{kvfduyerftywefope0rg}
\end{EQA}
However, a study of this estimator is quite involved due to the random nature of the matrix 
\( \Psiv(\Xv) \), especially in the case of a high design dimension \( \dimd \);
see recent papers 
\cite{baLoLu2020},
\cite{ChMo2022},
and references therein 
for the case of linear regression with \( \dimd = n \).
%

Now we describe our semiparametric approach.
Suppose to be given another collection of features \( (\phi_{k}) \) for \( k=1,\ldots,\dimq \).
The use of \( (\phi_{k}) = (\psi_{j}) \) is one of possible options.
However, one can use many others including biorthogonal bases, regression trees and forests, 
last layer neurones of Deep Neuronal Networks, etc.
Define \( \Phiv(X_{i}) = (\phi_{k}(X_{i})) \in \R^{\dimq} \) and \( \Phiv(\Xv) = (\Phiv(X_{1}), \ldots,\Phiv(X_{n})) \).
Introduce two linear operators (design kernels) \( \KS \) and \( \hKS \): \( \R^{\dimp} \to \R^{\dimq} \) by 
\begin{EQA}
\label{iiu87654e3wsdrtyuiop}
	\hKS
	& \eqdef &
	\Phiv(\Xv) \, \Psiv(\Xv)^{\T}
	=
	\sumi \Phiv(X_{i}) \, \Psiv(X_{i})^{\T} ,
	\\
	\KSs
	& \eqdef &
	\E \, \Phiv(\Xv) \, \Psiv(\Xv)^{\T}
	=
	\sumi \E \, \Phiv(X_{i}) \, \Psiv(X_{i})^{\T} .
\end{EQA}
Clearly \( \E \hKS = \KSs \), so \( \hKS \) might be viewed as an empirical version of \( \KSs \).
If the \( X_{i} \)'s are i.i.d. and \( \dens_{X} \) is the design density then
\begin{EQA}
	\KSs
	&=&
	n \E \, \Phiv(X_{1}) \, \Psiv(X_{1})^{\T}
	=
	n \int \Phiv(\xv) \, \Psiv(\xv)^{\T} \, \dens_{X}(\xv) \, d\xv \, ,
\label{plkmnbvcxcvbnmnbvcvbnm}
\end{EQA}  
In general, with \( \hKS_{i} = n \Phiv(X_{i}) \, \Psiv(X_{i})^{\T} \), it holds 
\( \hKS = n^{-1} \sumi \hKS_{i} \) and also
\begin{EQA}
	\sumi \| \hKS_{i} - \KS \|_{\Fr}^{2}
	&=&
	\sumi \| \hKS_{i} - \hKS \|_{\Fr}^{2} + n \| \hKS - \KS \|_{\Fr}^{2} \, .
\label{iujwefy7wefyikjt5r433e}
\end{EQA}
As the term \( \sumi \| \hKS_{i} - \hKS \|_{\Fr}^{2} \) only depends on the design \( \Xv \),
one may equally use \( \sumi \| \hKS_{i} - \KS \|_{\Fr}^{2} \) and 
\( n \| \hKS - \KS \|_{\Fr}^{2} \)
to measure the data misfit of a possible guess \( \KS \).

Further, denote
\begin{EQA}
	Z_{k}
	& \eqdef &
	\sumi Y_{i} \, \phi_{k}(X_{i}) 
\label{dcirf7rft64rt6rfgh}
\end{EQA}
and let \( \Zv \) be the vector in \( \R^{\dimp} \) with entries \( Z_{k} \).
The LPA \( \fvs(\xv) = \Psiv(\xv)^{\T} \, \thetav \) together with \( \Yv = \fvs + \epsv \) and  \( \Uv = \Phiv(\Xv) \epsv \) implies
\begin{EQA}
	\Zv
	&=&
	\Phiv(\Xv) \fvs + 
	\Phiv(\Xv) \epsv
	=
	\Phiv(\Xv) \, \Psiv(\Xv)^{\T} \, \thetav + 
	\Phiv(\Xv) \epsv
	=
	\hKS \, \thetav + \Uv .
\label{poiuytrr43ewsdere4dr}
\end{EQA}
The vector \( \Zv \) might be viewed as new ``observations'' of the ``response'' \( \etav = \KS \thetav \)
corrupted by two sources of errors: the additive error \( \Uv \) and the error in operator \( \hKS - \KS \).
Due to definition, these errors are not independent, however, \( \E (\Uv \cond \Xv) = 0 \).
If the errors \( \eps_{i} = Y_{i} - \fs(X_{i}) \) are Gaussian independent of the design variables \( X_{i} \), then
the new ``noise'' \( \Uv \) is Gaussian conditionally on \( \Xv \).
This enables us to decouple these sources in the likelihood by considering in place of \( \| \Zv - \hKS \thetav \|^{2} \)
two fidelity terms \( \| \Zv - \etav \|^{2} \) and 
\( n \| \hKS - \KS \|_{\Fr}^{2} \) subject to the structural 
equation \( \etav = \KS \thetav \).
Using the Lagrange multiplier idea, we put all together in one expression
\begin{EQA}
	\LL(\thetav,\etav,\KS)
	& \eqdef &
	- \frac{1}{2 \sigma^{2}} \| \Zv - \etav \|^{2}
	- \frac{n}{2 \sigma_{X}^{2}} \| \hKS - \KS \|_{\Fr}^{2} 
	- \frac{\lambda}{2 \sigma^{2}} \| \etav - \KS \thetav \|^{2} \, .
\label{3ewsazsplokfr5}
\end{EQA}
The first term specifies the quality of fitting the ``data'' \( \Zv \) by the ``response'' \( \etav \).
The second term is a similar fidelity term for the operator \( \KS \) observed with the operator noise of variance 
\( \sigma_{\KS}^{2} = n^{-1} \sigma_{X}^{2} \).
Note that the operator noise is only related to the design density \( \dens_{X} \), the observation noise \( \epsv \)
does not show up here.
The last term transfers the structural equation \( \etav = \KS \thetav \) into a penalty term.
For the parameter \( \lambda \), a natural choice is \( \lambda = 1 \).
Note that this equation is a special case of \eqref{dhfdfededrsdsrdsresswd}.

\subsubsection*{This paper contribution}
Laplace approximation claiming near normality of the posterior is the basic result in Bayesian inference.
However, its validity requires some constraints on the dimensionality of the parameter space 
of the sort \( \dimLs^{3} \ll n \) for the effective dimension \( \dimLs \) of the parameter; \cite{SpLaplace2022}.
This paper studies an important special case when only a \( \dimp \)-dimensional target of the posterior is of interest.
We present two new results.
The first one only requiring \( \dimLs \ll n \) and \( \dimLe^{3} \ll n \)
ensures a full dimensional concentration of the posterior
and approximation of the marginal posterior by a \emph{Gaussian mixture} with the accuracy of order 
\( \sqrt{\dimLe^{3}/n} \).
This result reduces the original problem to some analysis for high dimensional Gaussian measures.
The second result provides some sufficient conditions for the classical \emph{Laplace approximation} of the posterior.
Recent advances in high dimensional Gaussian probability \cite{GNSUl2017}
help to relax the condition \( \dimLs^{3} \ll n \) to \( \dimLs^{2} \ll n \).
An interesting open question is a possibility of getting a multimodal non-Gaussian approximation of the marginal posterior
for the range \( n^{1/2} \ll \dimLs \ll n \). 
The general results are illustrated on the particular example of \emph{error-in-operator} model. 

Some issues important for statistical inference are not addressed in this paper. 
In particular, we do not discuss the bias induced by the prior,
the use of credible sets as frequentist confidence sets, contraction rate over smoothness classes.
However, all these issues can been addressed similarly to \cite{Sp2022,SpLaplace2022},
the calculus from Section~\ref{Spriors} can be well used for the corresponding analysis.
Necessary finite sample guarantees for penalized estimation in error-in-operator or high dimensional random design models
require a careful derivation and will be provided in the forthcoming paper.

\subsection*{Organization of the paper}
Section~\ref{SsemiLaplace} considers a general semiparametric framework and presents some
results on marginal Laplace approximation
including Theorem~\ref{Ppartintegrsmi} about mixed Laplace approximation
and Theorem~\ref{TLaplaceopg} about Laplace approximation on the class of elliptic sets.
We also provide some sufficient conditions in terms of the efficient dimension of the full parameter, under which 
the mixed Laplace approximation can be reduced back to the standard Laplace Theorem.
Section~\ref{SLaplErrOp} specifies the result to the case of an error-in-operator problem.
The appendix collects some very useful facts about Laplace approximation for the integral  
\( \int \ex^{\lgd(\prmtv)} \, d\prmtv \) for a smooth concave function \( \lgd \) as well as some useful bounds 
from high-dimensional probability.
For reference convenience we include the general results on Laplace approximation from \cite{SpLaplace2022}
in Section~\ref{STaylor} and the technical results on Gaussian quadratic forms in Section~\ref{SGaussqf}.

{
\Chapter{Laplace approximation of the marginal}
\renewcommand{\Subsection}[1]{\subsection{#1}}

\def\dimtotal{\overline{\dimp}}

\label{SsemiLaplace}
This section considers the problem of Laplace approximation when the argument \( \prmtv \) of the function \( \lgd(\prmtv) \)
is high dimensional, \( \prmtv \in \R^{\dimtotal} \), but we are interested only in the sub-vector \( \targv \) of \( \prmtv \), 
\( \targv \in \R^{\dimp} \).
We state several results about marginal posterior.
Theorem~\ref{Tsemipostconc} explains concentration properties of the marginal posterior under the bound \( \dimAs_{\GP} \ll n \)
on the effective dimension of the full parameter \( \prmtv \).
Further, Theorem~\ref{Ppartintegrsmi} presents our main result about mixed Laplace approximation
of the marginal distribution.
Later Theorem~\ref{Ppartintegrs} and Theorem~\ref{TLaplmarginQ} provide some sufficient condition which allow to replace 
the mixed Laplace approximation by a single Gaussian distribution.
The most advanced Theorem~\ref{TLaplaceop} combines the obtained bounds with the standard trick 
in semiparametric estimation based on one-point orthogonality; see e.g. \cite{bickel1993efficient}. 
The latter result allows to improve the bound on the critical dimension from \( \dimAs^{3} \ll n \)
to \( \dimAs^{2} \ll n \) for establishing a Laplace approximation for the low-dimensional marginal distribution.
However, we guess that in zone \( n^{1/2} \ll \dimAs \ll n \), a mixed Laplace approximation is possible;
see Section~\ref{ScritdimLs} for details.

We write \( \prmtv = (\targv,\nuiv) \), where \( \nuiv \) stands for the nuisance component taking its values 
in some set \( \Nui \subseteq \R^{\dimq} \), so that \( \dimtotal = \dimp + \dimq \).
Denote \( \Ups = \R^{\dimp} \times \Nui \), 
\begin{EQA}
	\prmtvs 
	&=& 
	\argmax_{\prmtv} \lgd(\prmtv) , 
	\qquad 
	\prmtvs = (\targvs,\nuivs) .
\label{iovikjfcijeuev76rjibdkjd}
\end{EQA}
The full dimensional approach of
\ifapp{Section~\ref{STaylor}}{\cite{SpLaplace2022}}
 considers the measure \( \PfL \) on \( \Ups \) 
whose density is proportional to \( \ex^{\lgd(\upsv)} \),
and approximates it by a 
Gaussian measure.
To be more specific, assume the concavity condition \nameref{LLf0ref} with the decomposition
\begin{EQA}
	\lgd(\prmtv) 
	&=& 
	\lgdL(\prmtv) - \| \GPT \prmtv \|^{2}/2
\label{0bnit74eceudjejd}
\end{EQA}
for a concave function \( \lgdL(\upsv) \) and a symmetric positive full dimensional matrix \( \GPT^{2} \).
%
Define \( \IFT = - \nabla^{2} \lgd(\prmtvs) \), \( \IFT_{0} = - \nabla^{2} \lgdL(\prmtvs) \).

\Section{Full dimensional bounds}
\label{Sfullsemi}
This section states full dimensional about the parameter \( \prmtv \).
The full Laplace effective dimension is 
\begin{EQA}
	\dimLs
	& \eqdef &
	\tr\bigl( \IFT^{-1} \IFT_{0} \bigr) .
\label{c7djneifvf7e4kghlf}
\end{EQA}
With \( \xx \) fixed, define
\begin{EQA}
	\rrLs 
	&=& 
	2 \sqrt{\dimLs} + \sqrt{2\xx} \, . 
\label{5dhnjduye3jdfksgtwhbnd}
\end{EQA}
Theorem~\ref{TLaplaceTV34} yields the following result.

\begin{proposition}
\label{PfullLaplOPconc}
Let \( \lgd(\prmtv) \) follow \eqref{0bnit74eceudjejd}
and satisfy \nameref{LLtS3ref} with \( \prmtv = \prmtvs \), 
\( \HL^{2} = n^{-1} \IFT_{0} = - n^{-1} \nabla^{2} \lgdL(\prmtvs) \), and \( \rr = \amax^{-1} \rrLs \);
see \eqref{5dhnjduye3jdfksgtwhbnd}.
With \( \amax = 2/3 \), let 
\begin{EQA}
	\dltwbss
	\eqdef
	\frac{\hmax_{3} \, \amax^{-1} \rrLs}{n^{1/2}}
	& \leq &
	\frac{3}{4} \, .
\label{mfytf623yuedjvft6op}
\end{EQA}
Then
\begin{EQA}[ccl]
	\PfL\bigl( \| \IFT_{0}^{1/2} (\prmtv - \prmtvs) \| > \amax^{-1} \rrLs \bigr)
	& \leq &
	\ex^{-\xx} .
\label{vfocb94eiu8jry67fgop}
\end{EQA}
Moreover, if 
\begin{EQA}
	\frac{\hmax_{3} \, \amax^{-1} \rrLs \,\, \dimLs}{n^{1/2}}
	& \leq &
	2 ,
\label{0hcde4dft3igthg94yr5twewse}
\end{EQA}
then 
\begin{EQA}
	\sup_{A \in \BBB(\R^{\dimp})} \bigl| \PfL(\prmtv - \prmtvs \in A) - \P(\IFT^{-1/2} \gaussv \in A) \bigr|
	& \leq &
	4(\err_{3} + \ex^{-\xx}) 
\label{cbc5dfedrdwewwgergse}
\end{EQA}
with \( \gaussv \) standard normal in \( \R^{\dimtotal} \) and
\begin{EQA}
	\err_{3}
	&=&
	\frac{\hmax_{3} (\dimLs+1)^{3/2}}{4 (1 - \dltwbss/3)^{3/2} n^{1/2}} 
	\leq 
	\frac{\hmax_{3} \, (\dimLs+1)^{3/2}}{2 n^{1/2}} \, .
\label{5qw7dyf4e4354coefw9dufihse}
\end{EQA}
\end{proposition}

Let \( \IFU^{-1} = (\IFT^{-1})_{\targv\targv} \) be the \( \targv\targv \)-block of \( \IFT^{-1} \). 
If \( \Zv \sim \PfL \) and \( (\Xv,\Nuiv) \) are the components of \( \Zv \), then
the marginal \( \Xv \) of \( \Zv \) is also nearly Gaussian \( \ND(\targvs,\IFU^{-1}) \).
However, due to \eqref{5qw7dyf4e4354coefw9dufihse}, 
the accuracy of the full dimensional approximation deteriorates polynomially with the total
effective dimension \( \dimLs = \dimLs(\prmtvs) \) and requires \( \dimLs \ll n^{1/3} \).
The value \( \dimLs \) can be large even after dimensionality reduction caused by a penalty term.
Unfortunately, one cannot drop the condition \( \dimLs \ll n^{1/3} \), it seems to be inherent for the problem at hand. 
However, if we are interested in the low dimensional component \( \targv \) only, 
one can try to integrate out the remaining components and to relax the condition on the full parameter dimension.

\Subsection{Orthogonal case}
This section illustrates the principal idea on the special \emph{orthogonal} case when the function 
\( \lgd(\prmtv) = \lgd(\targv,\nuiv) \) satisfies the condition
\begin{EQA}
	\nabla_{\nuiv} \nabla_{\targv} \lgd(\targv,\nuiv)
	& \equiv &
	0.
\label{mnldflsytc4e4wwqfg}
\end{EQA}
In this case the function \( \lgd \) can be decomposed as
\( \lgd(\targv,\nuiv) = \lgd_{1}(\targv) + \lgd_{2}(\nuiv) \) for some functions 
\( \lgd_{1}(\targv) \) and \( \lgd_{2}(\nuiv) \).
When considering the \( \targv \)-marginal distribution of \( \PfL \),
the nuisance variable \( \nuiv \) can be ignored.

\begin{lemma}
\label{LortoLapl}
Let \( \lgd \) be twice continuously differentiable and satisfy \eqref{mnldflsytc4e4wwqfg}.
Then for any \( A \in \BBB(\R^{\dimp}) \), it holds
\( \PfL(\targv \in A) = \P_{\lgd_{1}}(\targv \in A) \). 
\end{lemma}

\begin{proof}
The condition \( \nabla_{\nuiv} \nabla_{\targv} \, \lgd(\targv,\nuiv) \equiv 0 \) implies the 
decomposition \( \lgd(\targv,\nuiv) = \lgd_{1}(\targv) + \lgd_{2}(\nuiv) \).
Now the result follows by the Fubini Theorem.
\end{proof}

This result allows to apply the result about Laplace approximation of \( \PfL \) to the function \( \lgd_{1}(\targv) \)
of the target variable \( \targv \), the corresponding accuracy depends also on the target dimension only.
Unfortunately the orthogonality condition \eqref{mnldflsytc4e4wwqfg} is too restrictive.
Later we discuss what can be stated if this condition is violated.

\Subsection{Conditional Laplace approximation}
The basic idea of our approach is, for any fixed value of the nuisance variable \( \nuiv \), to consider the Laplace approximation of
\( \lgd_{\nuiv}(\targv) = \lgd(\targv,\nuiv) \) as a function of \( \targv \) only.
We follow the line of Section~\ref{STaylor}.
Define 
\begin{EQA}[rcl]
	\targv_{\nuiv}
	& \eqdef &
	\argmax_{\targv} \lgd_{\nuiv}(\targv),
\label{jfoiuy2wedfv7tr2qsdzxdf}
	\\
	\IFL_{\nuiv}
	& \eqdef &
	- \nabla^{2} \lgd_{\nuiv}(\targv_{\nuiv}) 
	=
	- \nabla_{\targv\targv}^{2} \lgd(\targv_{\nuiv},\nuiv) .
\label{kfdy8d3fgfnm0lkgrr3}
\end{EQA}
First we state a conditional result: 
for each \( \nuiv \), the measure \( \P_{\lgd,\nuiv} \) on \( \R^{\dimp} \) with the density proportional to \( \ex^{\lgd_{\nuiv}(\targv)} \)
can be well approximated by the Gaussian measure \( \ND(\targv_{\nuiv},\IFL_{\nuiv}^{-1}) \).
The accuracy of approximation corresponds to the dimension of the target variable \( \targv \) only.
Moreover, under natural conditions, this result can be stated uniformly over the set \( \nuiv \in \Nui \).
This implies that the \( \targv \)-marginal of \( \PfL \) can be well approximated by the mixture of 
\( \ND(\targv_{\nuiv},\IFL_{\nuiv}^{-1}) \). 
Under the additional conditions \( \targv_{\nuiv} \approx \targv_{\nuivs} \),  
\( \IFL_{\nuiv} \approx \IFL \eqdef \IFL_{\nuivs} \), 
this mixture can be replaced by the Gaussian distribution \( \ND(\targvs,\IFL^{-1}) \).

Now we present our conditions.
The first one replaces the strong concavity of the full dimensional function \( f \) by 
concavity of each partial function \( \lgd_{\nuiv} \), \( \nuiv \in \Nui \).

\begin{description}
    \item[\label{LLfyref} \( \bb{(\mathcal{C}_{\nuiv})} \)]
      \textit{For any \( \nuiv \in \Nui \), there exists an operator \( \DVL_{\nuiv}^{2} \leq \IFL_{\nuiv} \) in \( \R^{\dimp} \) such that the function 
\begin{EQA}
	\lgdL_{\nuiv}(\targv_{\nuiv} + \uv)
	& \eqdef &
	\lgd_{\nuiv}(\targv_{\nuiv} + \uv) + \frac{1}{2} \| \IFL_{\nuiv}^{1/2} \uv \|^{2} - \frac{1}{2} \| \DVL_{\nuiv} \uv \|^{2}
\label{fTvufupumDTpDfx}
\end{EQA}
is concave in \( \uv \in \R^{\dimp} \). 
      }
\end{description}

If \( \lgd(\prmtv) = \lgdL(\prmtv) - \| \GPT \prmtv \|^{2}/2 \) 
with \( \| \GPT \prmtv \|^{2} = \| \GP \targv \|^{2} + \| \GPY \nuiv \|^{2} \) and 
\( \lgdL(\targv,\nuiv) \) concave in \( \targv \), then one can use 
\( \DVL_{\nuiv}^{2} = - \nabla_{\targv\targv}^{2} \, \lgdL(\targv_{\nuiv},\nuiv) \) 
while \( \IFL_{\nuiv} = \DVL_{\nuiv}^{2} + \GP^{2} \).
\nameref{LLfyref} enables us to define for each \( \nuiv \in \Nui \) the effective target dimension \( \dimL_{\nuiv} \),
the corresponding radius \( \rrL_{\nuiv} \), and the local vicinity \( \UVL_{\nuiv} \):
\begin{EQ}[rcl]
	\dimL_{\nuiv} 
	&=& 
	\tr\bigl( \DVL_{\nuiv}^{2} \IFL_{\nuiv}^{-1} \bigr),
	\\
	\rrL_{\nuiv}
	& = &
	2 \sqrt{\dimL_{\nuiv}} + \sqrt{2 \xx} , 
	\\
	\UVL_{\nuiv}
	&=&
	\bigl\{ \uv \colon \| \DVL_{\nuiv} \uv \| \leq \rrL_{\nuiv} \bigr\} .
\label{1w3rT2spTs2xlifxidey}
\end{EQ}
Later we assume that each function \( \lgd_{\nuiv} \) satisfies the smoothness conditions 
\nameref{LLseS3ref}, \nameref{LLsoS3ref}, and, if necessary, \nameref{LLseS4ref} with 
\( \rr = \amax^{-1} \rr_{\nuiv} \) for \( \amax \geq 2/3 \); 
see Section~\ref{Ssemiopt}.

\Subsection{Concentration of the conditional and marginal distribution}
We start with a concentration property for the target component.
\eqref{poybf3679jd532ff2} of Theorem~\ref{TLaplaceTV} applied to each \( \lgd_{\nuiv} \) 
yields the conditional bound on the tail probability: 
for the measure \( \P_{\lgd,\nuiv} \) with the density proportional to \( \exp \lgd_{\nuiv}(\targv) \), it holds
\begin{EQA}
	\P_{\lgd,\nuiv}\bigl( \| \DVL_{\nuiv} (\Xv - \targv_{\nuiv}) \| > \amax^{-1} \rr_{\nuiv} \bigr)
	& \leq &
	\ex^{-\xx} ,
	\qquad
	\nuiv \in \Nui .
\label{poybf3679jd532ff2se}
\end{EQA}
If we also succeed to control the variability of \( \DVL_{\nuiv}^{2} \) and \( \| \DVL_{\nuiv} \, (\targv_{\nuiv} - \targvs) \| \) in \( \nuiv \),
then the conditional bound \eqref{poybf3679jd532ff2se} would imply an unconditional one.
With \( \DVL^{2} = \DVL_{\nuivs}^{2} \), global Fr\'echet smoothness \eqref{c6ceyecc5e5etctwhcyegwc} of \nameref{LLpS3ref} implies 
\( \| \DVL^{-1} \, \DVL_{\nuiv}^{2} \, \DVL^{-1} \| \leq 1 + \dltwbss \)
for \( \nuiv \) in \( \Nui \) and \( \dltwbss \) from \eqref{8uiiihkkcrdrdteggcjcjd}; see Lemma~\ref{LvarDVa}.
In the next result we do not require global Fr\'echet smoothness.
Instead, we simply bound the variability of \( \DVL_{\nuiv}^{2} \):
for some fixed \( \CONSTdv \geq 1 \)
\begin{EQA}
	\CONSTdv^{-2} \, \DVL^{2}
	& \leq &
	\DVL_{\nuiv}^{2}
	\leq 
	\CONSTdv^{2} \, \DVL^{2} ,
	\qquad
	\nuiv \in \Nui .
\label{cvhbfte4tg3e3tw42qe4wergh}
\end{EQA}

\begin{theorem}
\label{Tsemipostconc}
Suppose \nameref{LLfyref}, \nameref{LLseS3ref}, and \eqref{cvhbfte4tg3e3tw42qe4wergh} for all \( \nuiv \in \Nui \).
Let also  
\begin{EQ}[rcl]
	\sup_{\nuiv \in \Nui} \| \DVL (\targv_{\nuiv} - \targvs) \|
	& \leq &
	\biase \, ,
\label{bijberg36g5r43tg}
\end{EQ}
and for any \( \nuiv \in \Nui \)
\begin{EQA}
	\dltwb_{3,\nuiv}
	=
	\frac{\hmax_{3} \, \rr_{\nuiv}}{n^{1/2}}
	& \leq &
	1/3 .
\label{bijberg36g5r43tgw}
\end{EQA}
Then
\begin{EQA}
	\PfL\bigl( \| \DVL (\Xv - \targvs) \| > \CONSTdv^{2} \, \amax^{-1} \rr_{\nuivs} + \biase \bigr)
	& \leq &
	\ex^{-\xx} .
\label{poybf3679jd532ff2sem}
\end{EQA}
\end{theorem}

\begin{proof}
Note first that by \eqref{cvhbfte4tg3e3tw42qe4wergh} and \eqref{bijberg36g5r43tg}
\begin{EQA}
	\| \DVL (\Xv - \targvs) \|
	& \leq &
	\| \DVL \, \DVL_{\nuiv}^{-1} \, \DVL_{\nuiv} \, (\Xv - \targv_{\nuiv}) \| 
	+ \| \DVL (\targv_{\nuiv} - \targvs) \|
	\\
	& \leq &
	\CONSTdv \, \| \DVL_{\nuiv} \, (\Xv - \targv_{\nuiv}) \| + \biase .
\label{vc8ghkhg56i5rhede432q2wer}
\end{EQA}
Moreover, \eqref{cvhbfte4tg3e3tw42qe4wergh} implies \( \rr_{\nuiv} \leq \CONSTdv \, \rr_{\nuivs} \) 
and the statement follows from \eqref{poybf3679jd532ff2se}.
\end{proof}

A nice feature of the bound \eqref{poybf3679jd532ff2sem} is that the concentration radius \( \rr_{\nuivs} \) corresponds to the dimension
of the target variable \( \Xv \) only.
Also, under \nameref{LLseS3ref}, the condition \( \dimL_{\nuiv} \ll n \) implies 
\( \hmax_{3} \rr_{\nuiv} n^{-1/2} \ll 1 \).
However, the result \eqref{poybf3679jd532ff2sem} becomes meaningful only if we can bound the semiparametric bias term \( \biase \); 
see \eqref{bijberg36g5r43tg}.
The condition \( \biase \lesssim \rr_{\nuivs} \) yields the following corollary.

\begin{corollary}
\label{CTsemipostconc}
Suppose the conditions of Theorem~\ref{Tsemipostconc} and let \( \biase \leq \CONSTbias \, \rr_{\nuivs} \).
Then
\begin{EQA}
	\PfL\bigl( \| \DVL (\Xv - \targv_{\nuiv}) \| > (\CONSTdv^{2} \, \amax^{-1} + \CONSTbias) \, \rr_{\nuivs} \bigr)
	& \leq &
	\ex^{-\xx} .
\label{poybf3679jd532ff2semc}
\end{EQA}
\end{corollary}
Usually this result is applied in combination with the so called ``one-point orthogonality'' condition;
see Section~\ref{SonepointLapl} later. 

\Subsection{Approximation by a Gaussian mixture}
\label{SmixedLaplace}

Now we turn to approximation of the marginal distribution of \( \Xv \) by a Gaussian mixture.
\eqref{ufgdt6df5dtgededsxd23gjg} of Theorem~\ref{TLaplaceTV} yields a TV-bound for each condition distribution of
\( \Xv \) given \( \nuiv \): 
for the measure \( \P_{\lgd,\nuiv} \) with the density proportional to \( \exp \lgd_{\nuiv}(\targv) \), it holds
\begin{EQA}
	\TV\bigl( \P_{\lgd,\nuiv}, \ND(\targv_{\nuiv},\IFL_{\nuiv}^{-1}) \bigr)
	& \leq &
	\frac{2 (\err_{\nuiv} + \ex^{-\xx})}{1 - \err_{\nuiv} - \ex^{-\xx}}
	\leq 
	4 (\err_{\nuiv} + \ex^{-\xx}) ,
\label{jvibhntg36fgfgftmkf}
\end{EQA}
where with \( \dltwb_{\nuiv} \) from \eqref{bijberg36g5r43tgw}
\begin{EQA}
	\err_{\nuiv}
	&=&
	\frac{0.75 \, \dltwb_{\nuiv} \, \dimL_{\nuiv}}{1 - \dltwb_{\nuiv}} \, ;
\label{juytr90f2dzaryjhfyfse}
\end{EQA}
see \eqref{juytr90f2dzaryjhfyf}.
Moreover, under the self-concordance condition \nameref{LLseS3ref}, one can use
\begin{EQA}
	\err_{\nuiv}
	&=&
	\frac{\hmax_{3}}{2} \, \sqrt{\frac{(\dimL_{\nuiv}+1)^{3}}{n}} \, ;
\label{scdugfdwyd2wywy26e6dese}
\end{EQA}
cf. \eqref{5qw7dyf4e4354coefw9dufih}.
This is again a nice bound which only involves the target dimension \( \dimL_{\nuiv} \).
However, the approximating Gaussian distribution varies with \( \nuiv \).
Introduce
\begin{EQ}[rcl]
	\feta_{\nuiv}
	& \eqdef &
	\max_{\targv} \lgd_{\nuiv}(\targv) - \lgd(\prmtvs)
	=
	\lgd_{\nuiv}(\targv_{\nuiv}) - \lgd(\prmtvs) 
	=
	\lgd(\targv_{\nuiv},\nuiv) - \lgd(\prmtvs) ,
	\\
	\discr_{\nuiv} 
	& \eqdef & 
	\frac{1}{2} \log \det (\IFL_{\nuiv}^{-1} \IFL)  .
\label{ge8qwefygw3qytfyju8qfhbd}
\end{EQ}
By definition, the \emph{deficiency} \( - \feta_{\nuiv} \) between the global maximum of \( \lgd(\prmtv) \) and the partial maximum
of \( \lgd(\targv,\nuiv) \) w.r.t. \( \targv \) for \( \nuiv \) fixed is always nonnegative.
Later we show that the study can be limited to \( \nuiv \in \Nui \) with \( \feta_{\nuiv} \) not too big in absolute value.
The first result states an approximation of the \( \targv \)-marginal by a mixture of normals.

\begin{theorem}
\label{Ppartintegrsmi}
Let \( f(\prmtvs) = \sup_{\prmtv} f(\prmtv) \).
Suppose \nameref{LLfyref} and \nameref{LLseS3ref} for all \( \nuiv \in \Nui \) and assume that
\begin{EQA}
	\sup_{\nuiv \in \Nui} \frac{\hmax_{3} \, \rr_{\nuiv} \, \dimL_{\nuiv}}{n^{1/2}}
	& \leq &
	\frac{3}{4} \, .
\label{0hcde4dft3igthg94yr5twewsemi}
\end{EQA}
Define \( \err_{\nuiv} \) by \eqref{scdugfdwyd2wywy26e6dese} and let also with \( \feta_{\nuiv} \) 
and \( \discr_{\nuiv} \) from \eqref{ge8qwefygw3qytfyju8qfhbd}
\begin{EQA}
	\errE
	& \eqdef &
	\frac{\int_{\Nui} \err_{\nuiv} \, \ex^{\feta_{\nuiv} + \discr_{\nuiv}} \, d\nuiv}
		 {\int_{\Nui} \ex^{\feta_{\nuiv} + \discr_{\nuiv}} \, d\nuiv} 
	< 
	\frac{1}{2} - \ex^{- \xx} \, .
\label{hetswffdgdgkhfeddcfwkr2mi}
\end{EQA}
Then for any function \( g(\targv) \) satisfying \( |g(\targv)| \leq 1 \), 
it holds with \( \gaussv \) standard normal
\begin{EQA}
	\left| \frac{\int_{\Ups} \ex^{ \lgd(\prmtv) } \, g(\targv) \, d\prmtv}{\int_{\Ups} \ex^{ \lgd(\prmtv)} \, d\prmtv}
		- \frac{\int_{\Nui} \ex^{\feta_{\nuiv} + \discr_{\nuiv}} \, \E g(\targv_{\nuiv} + \IFL_{\nuiv}^{-1/2} \gaussv) \, d\nuiv}
		{\int_{\Nui} \ex^{\feta_{\nuiv} + \discr_{\nuiv}} \, d\nuiv}
	\right|
	& \leq &
	\frac{2 (\errE + \ex^{- \xx})}{1 - \errE - \ex^{- \xx}}
	\, . 
	\qquad 
\label{iUtEefxwgudwmnmi}
\end{EQA}
\end{theorem}

\begin{remark}
Obviously \( \errE \) from \eqref{hetswffdgdgkhfeddcfwkr2mi} satisfies
\begin{EQA}
	\errE
	& \leq &
	\errs
	\eqdef
	\sup_{\nuiv \in \Nui} \err_{\nuiv} \, .
\label{mgfhygvfyv3rfdr4rfygg8}
\end{EQA}
However, the term \( \feta_{\nuiv} \) is always negative and it may decrease almost quadratically as \( \nuiv \) 
deviates from zero. 
Thus, bound \eqref{hetswffdgdgkhfeddcfwkr2mi} on \( \errE \) is weaker than the bound on \( \errs \).
\end{remark}

\begin{proof}
For any \( \prmtv = (\targv,\nuiv) \), the definitions of \( \lgd_{\nuiv}(\targv) = \lgd(\targv,\nuiv) \) and of 
\( \feta_{\nuiv} = \lgd_{\nuiv}(\targv_{\nuiv}) - \lgd(\prmtvs) \) yield a nice decomposition
\begin{EQA}
	\lgd(\prmtv) - \lgd(\prmtvs)
	& = &
	\feta_{\nuiv} + \lgd_{\nuiv}(\targv) - \lgd_{\nuiv}(\targv_{\nuiv})
\label{fwxfxeDeu22A}
\end{EQA}
which enables us to operate with each function \( \lgd_{\nuiv}(\targv) \) for \( \nuiv \in \Nui \) independently:
\begin{EQA}
	\int_{\Ups} \ex^{ \lgd(\prmtv) - \lgd(\prmtvs) } \, g(\uv) \, d\prmtv
	&=&
	\int_{\Nui} \ex^{\feta_{\nuiv}} \left( \int \ex^{\lgd_{\nuiv}(\targv) - \lgd_{\nuiv}(\targv_{\nuiv}) } \, g(\targv) \, d\targv  \right)\, d\nuiv 
	\\
	&=&
	\int_{\Nui} \ex^{\feta_{\nuiv}} \left( \int \ex^{\lgd_{\nuiv}(\targv_{\nuiv} + \uv) - \lgd_{\nuiv}(\targv_{\nuiv}) } \, 
		g(\targv_{\nuiv} + \uv) \, d\uv  \right)\, d\nuiv 
	\\
	&=&
	\int_{\Nui} \ex^{\feta_{\nuiv}} \left( \int \ex^{\lgd_{\nuiv}(\targv_{\nuiv} ; \uv) } \, g(\targv_{\nuiv} + \uv) \, d\uv  \right)\, d\nuiv.
	\qquad
\label{iWefxwgudwdedu}
\end{EQA}
By the arguments from the proof of Theorem~\ref{TLaplaceTV}, see \eqref{1erriexmH22du22} and \eqref{1erriexmH22du22u}, 
for any \( \nuiv \in \Nui \)
\begin{EQA}
	\left| \frac{\int \ex^{\lgd_{\nuiv}(\targv_{\nuiv} ; \uv) } \, g(\targv_{\nuiv} + \uv) \, d\uv 
				 - \int \ex^{- \| \IFL_{\nuiv}^{1/2} \uv \|^{2}/2 } \, g(\targv_{\nuiv} + \uv) \, d\uv}
		  		{\int \ex^{- \| \IFL_{\nuiv}^{1/2} \uv \|^{2}/2 } \, d\uv}
	\right|
	& \leq & 
	\err_{\nuiv} + \ex^{- \xx} .
\label{iUtEefxwgudwmnb}
\end{EQA}
Integrating w.r.t. \( \nuiv \) yields by \eqref{iWefxwgudwdedu}
\begin{EQA}
	&& \nquad
	\left| \int \ex^{ \lgd(\prmtv) - \lgd(\prmtvs) } \, g(\uv) \, d\uv \, d\nuiv
				 - \iint \ex^{\feta_{\nuiv}} \, \ex^{ - \| \IFL_{\nuiv}^{1/2} \uv \|^{2}/2 } \, g(\targv_{\nuiv} + \uv) \, d\uv \, d\nuiv
	\right|
	\\
	& \leq & 
	\iint \ex^{\feta_{\nuiv}} \, \ex^{ - \| \IFL_{\nuiv}^{1/2} \uv \|^{2}/2 } \, (\err_{\nuiv} + \ex^{- \xx}) \, d\targv \, d\nuiv \, .
\label{90fufr6r4y3er5t5er56rk}
\end{EQA}
The same bound applies with \( g(\uv) \equiv 1 \).
By definition,
\begin{EQA}
	\int \ex^{- \| \IFL_{\nuiv}^{1/2} \uv \|^{2}/2 } \, d\uv
	&=&
	\ex^{\discr_{\nuiv}} \int \ex^{- \| \IFL^{1/2} \uv \|^{2}/2 } \, d\uv 
\label{cw9wef4e3784e76fet67erf76}
\end{EQA}
and \eqref{iUtEefxwgudwmnmi} follows.
\end{proof}

The full dimensional concentration result of Theorem~\ref{TLaplaceTV34} enables us to limit \( \prmtv \) to the local vicinity
of the point \( \upsvs \)
\begin{EQA}
	\Ups_{0} 
	&=&
	\bigl\{ \prmtv \colon \| \IFT_{0}^{1/2} (\prmtv - \prmtvs) \| \leq \amax^{-1} \, \rrLs \bigr\} 
\label{vcw3klorvfy7ew7yu3}
\end{EQA}
for \( \amax = 2/3 \) and some specific \( \rrLs \).
By \( \Nui_{0} \) we denote the projection of \( \Ups_{0} \) on \( \Nui \):
\begin{EQA}
	\Nui_{0}
	&=&
	\{ \nuiv \in \Nui \colon (\targv,\nuiv) \in \Ups_{0} \text{ for some } \targv \} .
\label{4jhcc6t5e35cdjnheu7y}
\end{EQA}

\begin{theorem}
\label{Ppartintegrloc}
Suppose the conditions of Theorem~\ref{Ppartintegrsmi} and of Theorem~\ref{TLaplaceTV34}.
Then mixed Laplace approximation \eqref{iUtEefxwgudwmnmi} applies 
with \( \Nui_{0} \) in place of \( \Nui \).
\end{theorem}

\Subsection{Laplace approximation}
The result of Theorem~\ref{Ppartintegrsmi} is very useful because 
it replaces the original problem by a very particular problem for Gaussian measures only.
Namely, we are now interested in approximating a Gaussian mixture by one Gaussian distribution.
The Gaussian mixture in the approximation \eqref{iUtEefxwgudwmnmi} is characterized by the collection of 
the conditional mean \( \targv_{\nuiv} \), variance \( \IFL_{\nuiv}^{-1} \), and of the deficiencies 
\( \feta_{\nuiv} \) for each \( \nuiv \in \Nui \).
Moreover, Theorem~\ref{Ppartintegrloc} enables us to consider only \( \nuiv \in \Nui_{0} \).
The next question is about possibility of using one Gaussian distribution \( \ND(\targvs,\IFL^{-1}) \) instead of the mixture of 
\( \ND(\targv_{\nuiv},\IFL_{\nuiv}^{-1}) \).
For this, we have to control variability of the parameters \( \targv_{\nuiv} \) and \( \IFL_{\nuiv}^{-1} \). 
Given a test function \( g(\uv) \), introduce 
in addition to \( \feta_{\nuiv} \) and \( \discr_{\nuiv} \) from \eqref{ge8qwefygw3qytfyju8qfhbd} 
one more function of the argument \( \nuiv \):
\begin{EQA}[rcl]
	\Delta_{g,\nuiv}
	& \eqdef &
	\E g(\targv_{\nuiv} + \IFL_{\nuiv}^{-1/2} \gaussv) - \E g(\targv_{\nuivs} + \IFL_{\nuivs}^{-1/2} \gaussv) \, .
\label{hwef78ywsfevd2fdft5rwedtcdbsww}
\end{EQA}
Clearly the value \( \Delta_{g,\nuiv} \) is uniquely determined by 
\( \targv_{\nuiv} - \targvs \) and \( \discr_{\nuiv} \) from \eqref{ge8qwefygw3qytfyju8qfhbd}.

\begin{theorem}
\label{Ppartintegrs}
Suppose the conditions of Theorem~\ref{Ppartintegrsmi}.
Fix a function \( g(\targv) \) with \( |g(\targv)| \leq 1 \).
With \( \Delta_{g,\nuiv} \) from \eqref{hwef78ywsfevd2fdft5rwedtcdbsww}, 
\( \feta_{\nuiv} \) and \( \discr_{\nuiv} \) from \eqref{ge8qwefygw3qytfyju8qfhbd},
define
\begin{EQA}
	\Delta_{g}
	& \eqdef &
	\left| 
		\frac{\int_{\Nui} \Delta_{g,\nuiv} \, \ex^{\feta_{\nuiv} + \discr_{\nuiv}} \, d\nuiv}
		 	 {\int_{\Nui} \ex^{\feta_{\nuiv} + \discr_{\nuiv}} \, d\nuiv} 
	\right| ,
\label{hetswffdgdgkhfeddcfwkr2}
\end{EQA}
and let \( \errE + \Delta_{g} < 1/2 - \ex^{-\xx} \); see \eqref{hetswffdgdgkhfeddcfwkr2mi}.
Then 
\begin{EQA}
	\left| \frac{\int_{\Ups} \ex^{ \lgd(\prmtv) } \, g(\targv) \, d\prmtv}{\int_{\Ups} \ex^{ \lgd(\prmtv)} \, d\prmtv}
		- \E g(\targvs + \IFL_{\nuivs}^{-1/2} \gaussv)
	\right|
	& \leq &
	\frac{2 (\errE + \Delta_{g} + \ex^{- \xx})}{1 - \errE - \Delta_{g} - \ex^{- \xx}} 
	\, .
\label{5de43edewef2wserdfwded}
\end{EQA}
Under the conditions of Theorem~\ref{TLaplaceTV34}, the result applies with \( \Nui_{0} \) in place of \( \Nui \).
\end{theorem}

\begin{proof}
Bound \eqref{90fufr6r4y3er5t5er56rk} and definition \eqref{hwef78ywsfevd2fdft5rwedtcdbsww} imply
\begin{EQA}
	&& \nquad
	\left|  
		\frac{1}{\int \ex^{- \| \IFL^{1/2} \uv \|^{2}/2 } \, d\uv} \, \int_{\Ups} \ex^{ \lgd(\prmtv) - \lgd(\prmtvs) } \, g(\uv) \, d\prmtv
		-
		\E g(\targvs + \IFL_{\nuivs}^{-1/2} \gaussv) 
		\int_{\Nui} \ex^{\feta_{\nuiv} + \discr_{\nuiv}} \, d\nuiv 
	\right|
	\\
	& \leq &
	\left| \int_{\Nui} ( \err_{\nuiv} + \ex^{- \xx}  + \Delta_{g,\nuiv}) \, \ex^{\feta_{\nuiv} + \discr_{\nuiv}} \, d\nuiv \right| 
	\leq 
	(\errE + \Delta_{g} + \ex^{-\xx}) \int_{\Nui} \ex^{\feta_{\nuiv} + \discr_{\nuiv}} \, d\nuiv .
\label{hfdu7fre224etyfg8frgkdfrf}
\end{EQA}
The same bound applies to \( g(\targv) \equiv 1 \) and the result follows as in the proof of Proposition~\ref{PunbintLapl}.
\end{proof}

The bound of Theorem~\ref{Ppartintegrsmi} only assumes \( |g(\targv)| \leq 1 \), otherwise \( g(\cdot) \) can be any
measurable function of the target parameter \( \targv \).
In the contrary, the error term \( \err_{g,\nuiv} \) in  \eqref{hwef78ywsfevd2fdft5rwedtcdbsww} 
strongly relies on the function \( g(\cdot) \). 
One way for stating auniform result can be based on Pinsker's inequality. 
For two measures \( \ND_{\nuivs} = \ND(\targvs,\IFL_{\nuivs}^{-1}) \) and 
\( \ND_{\nuiv} = \ND(\targv_{\nuiv},\IFL_{\nuiv}^{-1}) \) and
for any function \( g(\cdot) \) bounded by 1, it holds  
\begin{EQA}
	|\Delta_{g,\nuiv}|
	& \leq &
	\TV( \ND_{\nuivs},\ND_{\nuiv} )
	\leq 
	\sqrt{\kullb( \ND_{\nuivs},\ND_{\nuiv} )/2} \, ,
\label{fvirnhryrfd656ege4gt}
\end{EQA}
where \( \TV( \ND_{\nuivs},\ND_{\nuiv} ) \) is the total variation distance between \( \ND_{\nuivs} \) and \( \ND_{\nuiv} \),
while \( \kullb( \ND_{\nuivs},\ND_{\nuiv} ) \) is the Kullback-Leibler divergence.
For two Gaussian measures \( \ND_{\nuivs},\ND_{\nuiv} \), it holds with \( \IFL = \IFL_{\nuivs} \)
\begin{EQA}
	\kullb( \ND_{\nuivs},\ND_{\nuiv} )
	&=&
	\frac{1}{2} \bigl\{ \| \IFL^{1/2} (\targv_{\nuiv} - \targvs) \|^{2} 
	+ \tr (\IFL_{\nuiv}^{-1} \IFL - \Id_{\dimp}) + \log \det (\IFL_{\nuiv}^{-1} \IFL) \bigr\} .
\label{h0bede4bweve7yjdyeghyse}
\end{EQA}
Moreover, if the matrix \( \BB_{\nuiv} = \IFL_{\nuiv}^{-1/2} \IFL \, \IFL_{\nuiv}^{-1/2} - \Id_{\dimp} \) satisfies 
\( \| \BB_{\nuiv} \| \leq 2/3 \) then
\begin{EQA}
	\TV( \ND_{\nuivs},\ND_{\nuiv} )
	& \leq &
	\frac{1}{2} \Bigl( \| \IFL^{1/2} (\targv_{\nuiv} - \targvs) \| + \sqrt{\tr \BB_{\nuiv}^{2}} \Bigr) .
\label{dferrfwbvf6nhdnghfkeiryse}
\end{EQA}
This statement and Theorem~\ref{Ppartintegrs} imply a bound for the total variation distance 
between the \( \targv \)-marginal of \( \PfL \) and the Gaussian approximation \( \ND(\targvs,\IFL^{-1}) \).
However, Pinsker's inequality might be too rough.
Particularly, dependence on \( \| \IFL^{1/2} (\targv_{\nuiv} - \targvs) \| \) can be problematic.

\Subsection{Gaussian comparison and elliptic credible sets}
Bound \eqref{fvirnhryrfd656ege4gt} can be drastically improved by using recent results on Gaussian comparison
if we limit ourselves to some special class of functions \( g \)
of the form \( g(\uv) = \Ind(\| \QP \uv \| \leq \rr) \); see Section~\ref{SmainresGC}.
Define
\begin{EQA}
	\Delta_{\nuiv}
	&=&
	\frac{1}{\| \QP \, \IFL^{-1} \, \QP^{\T} \|_{\Fr}}
	\left( 
		\| \QP (\IFL^{-1} - \IFL_{\nuiv}^{-1}) \QP^{\T} \|_{1} + \| \QP (\targv_{\nuiv} - \targvs) \|^{2} 
	\right) ,
\label{hysfdt7dycdgds6r}
	\\
	\Delta_{\Nui}
	& \eqdef &
	\frac{\int_{\Nui} \Delta_{\nuiv} \, \ex^{\feta_{\nuiv} + \discr_{\nuiv}} \, d\nuiv}
	 	 {\int_{\Nui} \ex^{\feta_{\nuiv} + \discr_{\nuiv}} \, d\nuiv} \, .
\label{hetswffdgdgkhfeddcfwkr2D}
\end{EQA}

\begin{theorem}
\label{TLaplmarginQ}
Suppose the conditions of Theorem~\ref{Ppartintegrsmi}.
For any \( \QP \colon \R^{\dimp} \to \R^{\dimq} \), it holds with \( \gaussv \sim \ND(0,\Id) \)
\begin{EQA}
	\sup_{\rr > 0}
	\Bigl| \PfL\bigl( \| \QP (\Xv - \targvs) \| \leq \rr \bigr) - \P\bigl( \| \QP \, \IFL^{-1/2} \gaussv \| \leq \rr \bigr) \Bigr|
	& \leq &
	\frac{2 (\err_{\QP} + \ex^{- \xx})}{1 - \err_{\QP} - \ex^{- \xx}}
	\, ,
	\qquad
\label{jbfi8jfbe32hdfvdfrdfg}
\end{EQA}
where 
\begin{EQA}
	\err_{\QP} 
	&=& 
	\errE + \CONST \Delta_{\Nui} 
\label{2FGY888765FGHBJYS}
\end{EQA}
with \( \errE \) from \eqref{hetswffdgdgkhfeddcfwkr2mi}, 
\( \Delta_{\Nui} \) from \eqref{hetswffdgdgkhfeddcfwkr2D}, and some absolute constant \( \CONST \).
Under conditions of Theorem~\ref{TLaplaceTV34}, one can use \( \Nui_{0} \) from \eqref{4jhcc6t5e35cdjnheu7y} 
in place of \( \Nui \).
\end{theorem}

As a characteristic example, consider \( \QP = \DVL \).
Then Laplace approximation \eqref{jbfi8jfbe32hdfvdfrdfg} requires small values of 
\( \| \DVL (\IFL^{-1} - \IFL_{\nuiv}^{-1}) \DVL \|_{1}/ \sqrt{\dimL_{\nuivs}} \) 
and \( \| \DVL (\targv_{\nuiv} - \targvs) \| \) for all \( \nuiv \in \Nui_{0} \).

\Subsection{Variance homogenization}

To make the result of Theorem~\ref{TLaplmarginQ} meaningful, we have to get the rid of the bias 
\( \biasev_{\nuiv} = \targv_{\nuiv} - \targvs = \targv_{\nuiv} - \targv_{\nuivs} \) and of the varying precision
matrix \( \IFL_{\nuiv} \).
This section focuses on the latter.
Variability of \( \IFL_{\nuiv} \) is measured 
%
by the first term \( \| \QP (\IFL^{-1} - \IFL_{\nuiv}^{-1}) \QP^{\T} \|_{1} \) in \eqref{hysfdt7dycdgds6r}.
Now we can restate the approximation bounds \eqref{iUtEefxwgudwmnmi} and \eqref{5de43edewef2wserdfwded} using Gaussian mixture 
with a homogeneous variance \( \IFL^{-1} \).

\begin{theorem}
\label{Ppartintegrsho}
Suppose the conditions of Theorem~\ref{Ppartintegrsmi} and of Theorem~\ref{TLaplaceTV34}.
For any \( \QP \colon \R^{\dimp} \to \R^{\dimq} \), it holds with \( \gaussv \) standard normal
and \( \Nui_{0} \) from \eqref{4jhcc6t5e35cdjnheu7y}
\begin{EQA}
	&& \nquad
	\sup_{\rr > 0}
	\left| 
		\frac{\int_{\Ups} \ex^{ \lgd(\prmtv) } \, \Ind(\| \QP \targv \| \leq \rr) \, d\prmtv}
			 {\int_{\Ups} \ex^{ \lgd(\prmtv)} \, d\prmtv}
		- \frac{\int_{\Nui_{0}} \ex^{\feta_{\nuiv} + \discr_{\nuiv}} \, \P(\| \targv_{\nuiv} + \IFL^{-1/2} \gaussv \| \leq \rr) \, d\nuiv}
		{\int_{\Nui_{0}} \ex^{\feta_{\nuiv} + \discr_{\nuiv}} \, d\nuiv}
	\right|
	\\
	& \leq &
	\frac{2 (\err_{\Nui_{0}} + \CONST \Delta_{\IFL} + \ex^{- \xx})}{1 - \err_{\Nui_{0}} - \CONST \Delta_{\IFL} - \ex^{- \xx}}
	\leq 
	4 (\err_{\Nui_{0}} + \CONST \Delta_{\IFL} + \ex^{- \xx}) ,
	\qquad 
\label{iUtEefxwgudwmnmiho}
\end{EQA}
where \( \CONST \) is an absolute constant while with \( \dltwbss \) from \eqref{mfytf623yuedjvft6op}
\begin{EQA}
	\Delta_{\IFL}
	&=& 
	\frac{1}{\| \QP \, \IFL^{-1} \, \QP^{\T} \|_{\Fr}} \,\, 
	\frac{\int_{\Nui_{0}} \| \QP (\IFL^{-1} - \IFL_{\nuiv}^{-1}) \QP^{\T} \|_{1} \, \ex^{\feta_{\nuiv} + \discr_{\nuiv}} \, d\nuiv}
	 	 {\int_{\Nui_{0}} \ex^{\feta_{\nuiv} + \discr_{\nuiv}} \, d\nuiv} 
	\leq 
	\frac{\dltwbss}{1 - \dltwbss} \, \frac{\tr (\QP \, \IFL^{-1} \QP^{\T})}{\| \QP \, \IFL^{-1} \, \QP^{\T} \|_{\Fr}} 
	\, .
\label{2FGY888765FGHBJYSho}
\end{EQA}
Furthermore, define 
\begin{EQA}
	\BB_{\QP} 
	&=& 
	\frac{1}{\| \QP \, \IFL^{-1} \, \QP^{\T} \|} \,\, \QP \, \IFL^{-1} \, \QP^{\T}   ,
	\qquad 
	\dimL_{\QP} = \tr(\BB_{\QP}) , 
\label{5chcf7f6renecejeducjddh}
\end{EQA}
and suppose that 
\begin{EQA}
	\| \BB_{\QP} \|_{\Fr}^{2}
	&=&
	\tr \BB_{\QP}^{2}
	\geq 
	\CONSTi_{0}^{2} \tr \BB_{\QP}
	=
	\CONSTi_{0}^{2} \, \dimL_{\QP}
\label{cfirj4efd76ejhfcuew2djkd}
\end{EQA} 
with some fixed positive constant \( \CONSTi_{0} \).
Then the result \eqref{iUtEefxwgudwmnmiho} applies with 
\begin{EQA}
	\Delta_{\IFL}
	& \leq &
	\frac{\dltwbss \, \sqrt{\dimL_{\QP}}}{\CONSTi_{0} (1 - \dltwbss)} 
	=
	\frac{\hmax_{3} \, \amax^{-1} \, \rrLs \, \sqrt{\dimL_{\QP}}}{\CONSTi_{0} (1 - \dltwbss) \sqrt{n}} \, .
\label{4uy7jffde223e44rbjmkj}
\end{EQA}
\end{theorem}

\begin{proof}
The first statement follows from Theorem~\ref{Ppartintegrsmi} similarly to the proof of Theorem~\ref{TLaplmarginQ}
using Gaussian comparison results from Section~\ref{SmainresGC}.
It remains to evaluate \( \| \QP (\IFL^{-1} - \IFL_{\nuiv}^{-1}) \QP^{\T} \|_{1} \).
We present a simple upper bound.

\begin{lemma}
\label{LdiscrLapl}
Assume the conditions of Proposition~\ref{PfullLaplOPconc}.
Let \( \nuiv \in \Nui_{0} \). 
Then for any \( \QP \colon \R^{\dimp} \to \R^{\dimq} \), it holds with \( \IFL = \IFL_{\nuivs} \)
\begin{EQA}
	\| \QP (\IFL^{-1} - \IFL_{\nuiv}^{-1}) \QP^{\T} \|_{1}
	& \leq &
	\frac{\dltwbss}{1 - \dltwbss} \, \tr (\QP \, \IFL^{-1} \QP^{\T}) .
	\qquad
\label{8vckmkfiw3efc5tretywy1}
\end{EQA}
In particular, with \( \QP = \DVL \) from 
\begin{EQA}
	\| \DVL (\IFL^{-1} - \IFL_{\nuiv}^{-1}) \DVL \|_{1}
	& \leq &
	\frac{\dltwbss \, \dimL_{\nuivs}}{1 - \dltwbss} \,\, .
\label{8vckmkfiw3efc5tretywy2}
\end{EQA}
\end{lemma}

\noindent
\emph{Proof.}
Lemma~\ref{LfreTay} implies for \( \prmtv_{\nuiv} = (\targv_{\nuiv},\nuiv) \) and 
\( \IFT_{\nuiv} = \IFT(\prmtv_{\nuiv}) \)
\begin{EQA}
	\frac{1}{1 + \dltwbss} \, \IFT^{-1}
	& \leq &
	\IFT_{\nuiv}^{-1}
	\leq 
	\frac{1}{1 - \dltwbss} \, \IFT^{-1} \, .
\label{ucj8fcfu7wwrehdfic}
\end{EQA}
Clearly the same bound applies to the \( \targv\targv \)-blocks \( \IFL^{-1} \) and \( \IFL_{\nuiv}^{-1} \):
\begin{EQA}
	\frac{1}{1 + \dltwbss} \, \IFL^{-1}
	& \leq &
	\IFL_{\nuiv}^{-1}
	\leq 
	\frac{1}{1 - \dltwbss} \, \IFL^{-1} \, .
\label{ucj8fcfu7wwreckeedujdfiw}
\end{EQA}
%
Now 
\begin{EQA}
	\| \QP (\IFL_{\nuiv}^{-1} - \IFL^{-1}) \QP^{\T} \|_{1}
	& \leq &
	\Bigl\| \QP \Bigl( \frac{\IFL^{-1}}{1 - \dltwbss} - \IFL^{-1} \Bigr) \QP^{\T} \Bigr\|_{1}
	\leq 
	\frac{\dltwbss}{1 - \dltwbss} \, \tr (\QP \, \IFL^{-1} \QP^{\T}) ,
\label{8vckmkfiw3efc5tretywy12}
\end{EQA}
and \eqref{8vckmkfiw3efc5tretywy1} follows.
\end{proof}

\begin{remark}
Lemma~\ref{LdiscrLapl} demonstrates advantage of using the Gaussian comparison bound in terms of \eqref{hysfdt7dycdgds6r}
instead of KL-bound \eqref{fvirnhryrfd656ege4gt} and \eqref{h0bede4bweve7yjdyeghyse}.
The quantities \( \log \det (\IFL_{\nuiv}^{-1} \IFL_{\nuivs}) \) or
\( \tr(\IFL_{\nuiv}^{-2} \IFL_{\nuivs}^{2}) \) can be evaluated using \eqref{ucj8fcfu7wwreckeedujdfiw},
however the related bound would involve the target dimension \( \dimp \) in place of the effective
dimension \( \dimL_{\nuivs} \).
If the target variable is low dimensional and \( \dimp \) is small, this difference is not critical.
\end{remark}

\begin{remark}
In an important special case of a finite dimensional target \( \targv \), the value \( \dimL_{\QP} \) 
is bounded accordingly.
Therefore, the result \eqref{iUtEefxwgudwmnmiho} only requires \( \dltwbss \ll 1 \) or \( \dimLs \ll n \).
\end{remark}

\Subsection{A bound for the bias under one-point orthogonality}
\label{SonepointLapl}
The result of Theorem~\ref{Ppartintegrsho} claims an approximation of the posterior measure \( \PfL \)
by a homogeneous mixture of Gaussian distributions \( \ND(\targv_{\nuiv},\IFL^{-1}) \) with mean \( \targv_{\nuiv} \)
and a constant variance \( \IFL^{-1} \).
This section focuses on the bias component \( \biasev_{\nuiv} = \targv_{\nuiv} - \targvs \).
The main problem in establishing a sensitive bound for marginal Laplace approximation is
that the bias vector \( \biasev_{\nuiv} \) is nearly linear in \( \nuiv \);
see Lemma~\ref{Lpartmaxq} for the case of a quadratic function \( \lgd(\prmtv) \).
The use of the one-point orthogonality device described Section~\ref{SfullLaplco}, allows to kill the linear term and
yields a better bound on the norm of the bias \( \biasev_{\nuiv} \).
\emph{One-point} orthogonality condition reads
\begin{EQA}
	\nabla_{\nuiv} \nabla_{\targv} \lgd(\prmtvs)
	& = &
	0.
\label{mnldflsytc4e4wwqfgfu}
\end{EQA}
This condition can always be enforced by a linear transform of the nuisance parameter; see Section~\ref{Sonepointtrans},
Lemma~\ref{Onepointorto}.
It obviously yields that the negative Hessian \( \IFT = - \nabla^{2} \lgd(\prmtvs) \)
is block-diagonal
and the \( \targv\targv \)-block \( \IFU^{-1} \) of \( \IFT^{-1} \) coincides with \( \IFL^{-1} \).

\begin{theorem}
\label{TLaplaceop}
Let \( \lgd(\prmtv) \) follow \eqref{0bnit74eceudjejd} and condition \eqref{mnldflsytc4e4wwqfgfu} be fulfilled.
Assume  
\nameref{LLf0ref}
and \nameref{LLsS3ref} with \( \prmtv = \prmtvs \), 
\( \HL^{2} = n^{-1} \IFT_{0} \), and \( \rr = \amax^{-1} \rrLs \).
Let also  
\begin{EQA}[c]
	\dltwbss
	\eqdef
	\frac{\hmax_{3} \, \amax^{-1} \, \rrLs}{n^{1/2}} 
	\leq 1/3 .
\label{7jk98k9kjdfyenfc76e4kop}
\end{EQA}
Given \( \QP \), assume \eqref{cfirj4efd76ejhfcuew2djkd}.
Then with \( \gaussv_{\IFL} \sim \ND(0,\IFL^{-1}) \)
\begin{EQA}
	&& \nquad
	\sup_{\rr > 0}
	\Bigl| \PfL\bigl( \| \QP (\Xv - \targvs) \| \leq \rr \bigr) - \P\bigl( \| \QP \, \gaussv_{\IFL} \| \leq \rr \bigr) \Bigr|
	\\
	& \lesssim &
		\frac{\hmax_{3} \, \rrL_{\nuivs} \, \dimL_{\nuivs}}{\sqrt{n}}
		+ \frac{\hmax_{3} \, \rrLs \, \sqrt{\dimL_{\QP}}}{\sqrt{n}} 
		+ \frac{\hmax_{3}^{2} \, {\rrLs}^{4}}{n \, \sqrt{\dimL_{\QP}}}
		+ \ex^{- \xx} .
\label{jbfi8jfbe32hdfvdfrdfgop}
\end{EQA}
\end{theorem}

\begin{proof}
It holds by \eqref{8vckmkfiw3efc5tretywy1} of Lemma~\ref{LdiscrLapl} and \eqref{cfirj4efd76ejhfcuew2djkd}
\begin{EQA}
	\frac{\| \QP (\IFL^{-1} - \IFL_{\nuiv}^{-1}) \QP^{\T} \|_{1}}{\| \QP \, \IFL^{-1} \, \QP \|_{\Fr}}	
	& \leq &
	\frac{\dltwbss}{1 - \dltwbss} \, \frac{\tr(\QP \, \IFL^{-1} \, \QP)}{\| \QP \, \IFL^{-1} \, \QP \|_{\Fr}} 
	\leq 
	\frac{\dltwbss \, \sqrt{\dimL_{\QP}}}{\CONSTi_{0} (1 - \dltwbss)} 
	=
	\frac{\hmax_{3} \, \amax^{-1} \, \rrLs \, \sqrt{\dimL_{\QP}}}{\amax \,\CONSTi_{0} (1 - \dltwbss) \sqrt{n}} \, 
\label{8vckmkfiw3efc5tretywy13}
\end{EQA}
Proposition~\ref{LsemiAv} in Section~\ref{SortLaplsemi} implies with \( \rru =	\amax^{-1} \rrLs \)
\begin{EQA}
	\| \QP \, \biasev_{\nuiv} \|
	& \leq &
	\frac{\hmax_{3} \, \rru^{2}}{\sqrt{n}} \, \| \QP \, \IFL^{-1} \QP^{\T} \|^{1/2} \, ,
\label{36gfijh94ejdvtwekop}
\end{EQA}
and
\begin{EQA}
	\frac{\| \QP \, \biasev_{\nuiv} \|^{2}}{\| \QP \, \IFL^{-1} \, \QP \|_{\Fr}}
	& \leq &
	\frac{\hmax_{3}^{2} \, \rru^{2} \, \| \QP \, \IFL^{-1} \, \QP \|}{n \, \| \QP \, \IFL^{-1} \, \QP \|_{\Fr}} 
	\leq 
	\frac{\hmax_{3}^{2} \, \rru^{2}}{n \, \CONSTi_{0} \, \sqrt{\dimL_{\QP}}} \, .
\label{xcyeyuhw3njesjcuewhu}
\end{EQA}
This yields the assertion by Theorem~\ref{Ppartintegrsmi} and Theorem~\ref{Ppartintegrsho}.
\end{proof}

\Subsection{Laplace approximation: a general bound }
\label{SfullLaplco}
One-point orthogonality condition \eqref{mnldflsytc4e4wwqfgfu} plays an important role in  the results of Section~\ref{SonepointLapl}.
Unfortunately, this condition is rarely fulfilled for the original model function \( \lgd(\prmtv) \).
However, a simple trick based on a linear transform of the nuisance variable allows to ensure this condition 
in a rather general situation.
We again consider the setup \eqref{0bnit74eceudjejd}.
Instead of \eqref{mnldflsytc4e4wwqfgfu}, we only impose a mild \emph{separability condition}: 
the full dimensional information matrix
\( \IFT = \IFT(\prmtvs) = - \nabla^{2} \lgd(\prmtvs) \) can be bounded from below by the block-diagonal matrix 
with the blocks 
\( - \nabla_{\targv\targv}^{2} \lgd(\prmtvs) \) and \( - \nabla_{\nuiv\nuiv}^{2} \lgd(\prmtvs) \).
For a precise formulation, consider the block representation of \( \IFT \):
\begin{EQA}
	\IFT
	& = &
	\begin{pmatrix}
	\IFT_{\targv\targv} & \IFT_{\targv\nuiv}
	\\
	\IFT_{\nuiv\targv} & \IFT_{\nuiv\nuiv}
	\end{pmatrix} ,
	\qquad
	\IFT_{\nuiv\targv} = \IFT_{\targv\nuiv}^{\T} .
\label{hwe78yf2diwe76tfw67etfwtbop}
\end{EQA}

\begin{description}
	\item[\label{IFSref} \( \bb{(\IFT)} \)]
	\emph{It holds \( \IFT_{\targv\targv} > 0 \), \( \IFT_{\nuiv\nuiv} > 0 \), and for some \( \rhoIF = \rhoIF(\IFT) < 1 \)
	} 
\begin{EQA}
	\| 
	\IFT_{\targv\targv}^{-1/2} \IFT_{\targv\nuiv} \, \IFT_{\nuiv\nuiv}^{-1} \, \IFT_{\nuiv\targv} \, \IFT_{\targv\targv}^{-1/2} 
	\|
	& \leq &
	\rhoIF .
\label{vcjcfvedt7wdwhesqgghwqLsL}
\end{EQA}
\end{description}

\noindent
Define the efficient semiparametric Fisher matrix
\begin{EQA}[c]
	\IFU
	\eqdef 
	\IFT_{\targv\targv} - \IFT_{\targv\nuiv} \, \IFT_{\nuiv\nuiv}^{-1} \, \IFT_{\nuiv\targv} \, ;
\label{5hfjdieyedyiwhihidh3i3ce}
\end{EQA}
see \eqref{n6d2fgujrt6fgfrjrd}.
Lemma~\ref{Lidentsemi} implies
\begin{EQA}
	(1 - \rhoIF) \IFT_{\targv\targv} 
	\leq 
	\IFU
	& \leq &
	\IFT_{\targv\targv} \, .
\label{ijffd3e6dijuwsejk}
\end{EQA}
Moreover, \nameref{IFSref} can be written as 
\( \IFT_{\targv\nuiv} \, \IFT_{\nuiv\nuiv}^{-1} \, \IFT_{\nuiv\targv} \leq \rhoIF \IFT_{\targv\targv} \), that is,
\eqref{ijffd3e6dijuwsejk} is equivalent to \nameref{IFSref}.
Define also \( \CFT = \IFT_{\nuiv\nuiv}^{-1} \, \IFT_{\nuiv\targv} \).
Lemma~\ref{Onepointorto} explains how the linear transform 
\begin{EQ}[rcl]
	\nuov 
	&=& 
	\nuiv + \IFT_{\nuiv\nuiv}^{-1} \, \IFT_{\nuiv\targv} \, (\targv - \targvs) 
	=
	\nuiv + \CFT \, (\targv - \targvs) ,
	\\
	\lgdb(\targv,\nuov)
	&=&
	\lgd(\targv,\nuiv) 
	=
	\lgd(\targv,\nuov - \CFT \, (\targv - \targvs)) .
\label{0vcucvf653nhftdreyenop}
\end{EQ}
ensures one-point orthogonality of the transformed function \( \lgdb(\targv,\nuov) \).
Moreover, all the important concavity and smoothness conditions on the function \( \lgd(\prmtv) \)
are transferred to the function \( \lgdb(\targv,\nuov) \).
In particular, the self-concordance condition \nameref{LLsS3ref} is not dramatically affected 
by the linear transform \eqref{0vcucvf653nhftdreyenop},
perhaps, the constant \( \hmax_{3} \) has to be slightly updated due to a change of the shape of the
local elliptic vicinity of \( \prmtvs \).
The same holds for the Laplace effective dimensions \( \dimLs \) and \( \dimL_{\nuiv} \).
Now everything is prepared for applying the result of Theorem~\ref{TLaplaceop} 
to the function \( \lgdb(\targv,\nuov) \).
The major change is in using the Gaussian approximation \( \ND(\targvs,\IFU^{-1}) \)
in place of \( \ND(\targvs,\IFL^{-1}) \).

\begin{theorem}
\label{TLaplaceopg}
Let \( \lgd(\prmtv) \) follow \eqref{0bnit74eceudjejd}.
Assume  
\nameref{LLf0ref},
\nameref{IFSref},
and \nameref{LLsS3ref} with \( \prmtv = \prmtvs \), 
\( \HL^{2} = n^{-1} \IFT_{0} \), and \( \rr = \amax^{-1} \rrLs \) following \eqref{7jk98k9kjdfyenfc76e4kop}.
For \( \QP \) satisfying \eqref{cfirj4efd76ejhfcuew2djkd}, it holds
with \( \gaussv_{\IFU} \sim \ND(0,\IFU^{-1}) \)
\begin{EQA}
	&& \nquad
	\sup_{\rr > 0}
	\Bigl| \PfL\bigl( \| \QP (\Xv - \targvs) \| \leq \rr \bigr) - \P\bigl( \| \QP \, \gaussv_{\IFU} \| \leq \rr \bigr) \Bigr|
	\\
	& \lesssim &
		\frac{\hmax_{3} \, \rrL_{\nuivs} \, \dimL_{\nuivs}}{\sqrt{n}}
		+ \frac{\hmax_{3} \, \rrLs \, \sqrt{\dimL_{\QP}}}{\sqrt{n}} 
		+ \frac{\hmax_{3}^{2} \, {\rrLs}^{4}}{n \, \sqrt{\dimL_{\QP}}}
		+ \ex^{- \xx} .
\label{jbfi8jfbe32hdfvdfrdfgopU}
\end{EQA}
\end{theorem}

\Subsection{Critical dimension in marginal Laplace approximation}
\label{ScritdimLs}
The result of Theorem~\ref{TLaplaceop} yields an important conclusion.
With \( \rrLs \asymp \dimLs^{1/2} \), the condition \( \Delta_{\nuiv} \ll 1 \) from \eqref{hysfdt7dycdgds6r} 
means that
\begin{EQA}[c]
	\hmax_{3}^{2} \,\, \dimL_{\QP} \,\, \dimLs \ll n, 
	\qquad 
	\hmax_{3}^{2} \,\, {\dimLs}^{2} \ll n \, \sqrt{\dimL_{\QP}} \, .
\label{c0cle3fc98vkru8ecm}
\end{EQA}
The marginal Laplace approximation requires in addition 
\( \hmax_{3} \, \rr_{\nuivs} \, \dimL_{\nuivs} \asymp 
\hmax_{3} \, \dimL_{\nuivs}^{3/2} \ll n^{1/2} \).
In particular, if the target component is low dimensional and the factors \( \dimL_{\nuivs}, \dimL_{\QP} \) can be ignored, 
condition \eqref{c0cle3fc98vkru8ecm} reads as
\( \hmax_{3}^{2} \,\, \dimLs^{2} \ll n \) which improves the full dimensional condition \( \hmax_{3}^{2} \,\, \dimLs^{3} \ll n \).
This improvement is obtained by a non-trivial combination of the Gaussian mixture approximation of Theorem~\ref{Ppartintegrsmi}
and the Gaussian comparison bound used in Theorem~\ref{TLaplmarginQ}.

\section{Error-in-operator model}
\label{SLaplErrOp}

This section aims at applying the general results of Section~\ref{SsemiLaplace} to the special case
of error-in-operator model.
The main focus is on the properties of the marginal posterior such as concentration and Laplace approximation
under possibly weak conditions on the full parameter dimension.
In the contrary to \cite{Trabs_2018}, the issues like contraction rate under smoothness conditions on the source signal
\( \targv \) and the operator \( \KS \) are not discussed. 
However, the main result from Theorem~\ref{TLaplaceEO} claims that the marginal posterior behaves as 
in the classical linear inverse problem with the true known operator.
This enables us to reduce the remaining questions to the well studied case of a linear inverse problem; 
see e.g. \cite{KnVaZa2011}, \cite{KnSzVa2016}.

Given a vector \( \zv \in \R^{n} \), the task is to invert the relation \( \zv = \KS \targv + \epsv \) when 
the linear operation \( \KS \) is not known and only a pilot \( \hKS \) is available.
Following the suggestion in Section~\ref{SCalmexampl}, we treat this task as a semiparametric problem 
of recovering \( \targv \) while \( \KS \) serves as nuisance parameter.
Given an image vector \( \zv \in \R^{\dimq} \) 
and a pilot \( \hKS \colon \R^{\dimp} \to \R^{\dimq} \), consider the function
\begin{EQA}
	\lgd(\targv,\KS)
	&=&
	- \frac{1}{2} \| \zv - \KS \targv \|^{2} 
	- \frac{\muA^{2}}{2} \| \hKS - \KS \|_{\Fr}^{2} 
	- \frac{1}{2} \| \GP \targv \|^{2} 
	- \frac{1}{2} \| \KS \|_{\GPKS}^{2} \, .
\label{bh2gfeft65t6f5rh3wqe3r}
\end{EQA}
Here \( \targv \in \R^{\dimp} \) is the variable of interest 
while \( \KS \in \R^{\dimq \times \dimp} \) 
is an operator from \( \R^{\dimp} \) to \( \R^{\dimq} \) serving as a nuisance variable.
The factor \( \muA \) in the fidelity term \( \muA^{2} \| \hKS - \KS \|_{\Fr}^{2} \) 
scales the operator noise. 
As usual, smoothness of the source signal \( \targv \) is controlled by the penalty term 
\( \| \GP \targv \|^{2} \).
One can easily extend the derivation for a penalty \( \| \GP (\targv - \targv_{0}) \|^{2}/2 \) 
in which \( \targv_{0} \) is an initial guess for \( \targv \).
Equivalently, one can treat this penalty as log-density of the Gaussian prior 
\( \ND(\targv_{0},\GP^{-2}) \) on \( \targv \).
The norm \( \| \cdot \|_{\GPKS} \) reflects our prior information about operator smoothness.
Note that the fidelity \( \muA^{2} \| \hKS - \KS \|_{\Fr}^{2} \) already collects 
some prior information about \( \KS \), so an additional penalization
is not ultimate. 
However, this term does not mimic any smoothness of the operator \( \KS \).
We apply a construction reflecting the standard ``approximation spaces'' approach;
see \cite{HoRe2008},
\cite{Trabs_2018}.
Consider \( \GPKS = \blk\{ \GPKS_{1}, \ldots,\GPKS_{\dimq} \} \),
where \( \GPKS_{m} \) is a positive symmetric operator in \( \R^{\dimp} \).
For \( \KS = (\KS_{1},\ldots,\KS_{\dimq})^{\T} \), this yields
\begin{EQA}
	\| \KS \|_{\GPKS}^{2}
	&=&
	\sum_{m=1}^{\dimq} \| \GPKS_{m} \KS_{m} \|^{2} 
	\, .
\label{4f8g8trthwcibgitiehr}
\end{EQA}
The ridge regression case corresponds to
\( \GP^{2} = \gp^{2} \Id_{\dimp} \) and \( \GPKS_{m}^{2} \equiv 0 \).

The main question under study is the set of sufficient conditions ensuring Laplace approximation 
for the marginal distribution for the target parameter \( \targv \) when considering the unknown operator \( \KS \)
as nuisance. 
Define
\begin{EQA}
	(\targvs,\KSs)
	&=&
	\argmax_{(\targv,\KS)} \lgd(\targv,\KS) .
\label{b983mfguy7gvu6tyfdtde}
\end{EQA}
Because of the product term \( \KS \targv \) in \eqref{bh2gfeft65t6f5rh3wqe3r}, 
this problem is not quadratic in the scope of parameters 
\( \prmtv = (\targv,\KS) \) and a closed form solution is not available. 
In particular, \( \targvs \) does not coincide with the plug-in solution \( \targv_{\hKS} \)
\begin{EQA}
	\targv_{\hKS}
	&=&
	\bigl( \hKS^{\T} \hKS + \GP^{2} \bigr)^{-1} \hKS^{\T} \zv  .
\label{jjgxd556dgh3wecuweuw}
\end{EQA}

\Subsection{Identifiability, warm start, efficient information matrix}
The original relation \( \zv = \KS \targv \) does not allow to recover \( \targv \) when \( \KS \) is unknown.
The additional information about \( \KS \) in form of a pilot \( \hKS \) or, alternatively, a prior on \( \KS \),
improves this lack of identifiability.
However, the quality of the pilot measured by the multiplicative factor \( \muA \) in \eqref{bh2gfeft65t6f5rh3wqe3r} 
is important.
The smaller is the operator noise, or, equivalently, the better is the accuracy 
of the pilot \( \hKS \), the larger \( \muA \) should be.
Still we have an issue with identifiability for \( \targv \) very large due to multiplicative
term \( \KS \targv \) in \( \| \zv - \KS \targv \|^{2}/2 \): 
even a small error in the operator \( \KS \)
may result in a large shift of the solution \( \targv \).
Proposition~\ref{PIFSrefEiO} later in this section helps to address the identifiability question using the so called 
``warm start'' condition.
We also compute the efficient semiparametric matrix \( \IFU \) and check condition \nameref{IFSref}.

\begin{proposition}
\label{PIFSrefEiO}
Let \( \GP^{2} \geq \GPa^{2} \) for some \( \GPa^{2} \) positive semi-definite.
Define for \( \rhoIF < 1 \)
\begin{EQA}
	\Upsd
	& \eqdef &
	\Bigl\{ (\targv,\KS) \colon  
	4 \| \targv \|^{2} \leq \rhoIF \, \muA^{2} ,
	\, \,
	4 \| \KS \targv - \imav \|^{2} \Id_{\dimp} \leq \rhoIF \, \muA^{2} \, (\KS^{\T} \KS + 2 \GPa^{2})
	\Bigr\}.
	\qquad
\label{23fcgftfgwbsnhewretge}
\end{EQA}
Let the point \( (\targvs,\KSs) \) belong to a set \( \Upsd \).
Then \nameref{IFSref} is fulfilled for \( \IFT = \IFT(\prmtv) \) and all \( \prmtv \in \Upsd \) with the same \( \rhoIF \), 
and \( \lgd(\prmtv) \) is strongly concave on \( \Upsd \).
\end{proposition}

\begin{proof}
Let \( \prmtv = (\targv,\KS) \) and \( \KS_{m}^{\T} \) be the rows of \( \KS \).
The function \( \lgd(\targv,\KS) \) from \eqref{bh2gfeft65t6f5rh3wqe3r} can be represented as 
\begin{EQA}
	\lgd(\targv,\KS)
	&=&
	- \frac{1}{2} \sum_{m=1}^{\dimq} \bigl( 
		| z_{m} - \KS_{m}^{\T} \targv |^{2} 
		+ \| \hat{\KS}_{m} - \KS_{m} \|^{2} 
		+ \| \GPKS_{m} \KS_{m} \|^{2} 
	\bigr)
	- \frac{1}{2} \| \GP \targv \|^{2} \, .
\label{bh2gfeft65t6f5rh3wqe3rm1M}
\end{EQA}
It holds
\begin{EQA}
	\IFT_{\targv\targv}(\prmtv)
	=
	- \nabla_{\targv\targv}^{2} \lgd(\targv,\KS)
	&=&
	\KS^{\T} \KS + \GP^{2}
	=
	\sum_{m=1}^{\dimq} \KS_{m} \, \KS_{m}^{\T} + \GP^{2}
\label{un4v7frue43kgv8uhkj}
\end{EQA}
and for each \( m = 1, \ldots,\dimq \) with 
\( \HFblk_{m}^{2} = \muA^{2} \Id_{\dimp} + \GPKS_{m}^{2} \)
\begin{EQ}[rcccl]
	- \nabla_{\KS_{m} \KS_{m}}^{2} \lgd(\targv,\KS)
	&=&
	\targv \targv^{\T} + \muA^{2} \Id_{\dimp} + \GPKS_{m}^{2} 
	& \eqdef &
	\targv \targv^{\T} + \HFblk_{m}^{2} \, ,
	\\
	- \nabla_{\targv} \nabla_{\KS_{m}} \lgd(\targv,\KS)
	&=&
	(\KS_{m}^{\T} \targv - z_{m}) \Id_{\dimp} + \KS_{m}\, \targv^{\T} 
	& \eqdef &
	\AFblk_{m}(\targv,\KS_{m}) \, .
\label{8dfnjey63f6ytjutyswwwq}
\end{EQ}
Block-diagonal structure of \( \IFT_{\KS\KS}(\prmtv) \) yields
\begin{EQA}
	\IFT_{\targv\KS}(\prmtv) \, \IFT_{\KS\KS}^{-1}(\prmtv) \, \IFT_{\KS\targv}(\prmtv) 
	&=& 
	\sum_{m=1}^{\dimq} \AFblk_{m}(\targv,\KS_{m}) \, \IFT_{\KS_{m}\KS_{m}}^{-1}(\prmtv) \, \AFblk_{m}^{\T}(\targv,\KS_{m}) \, . 
\label{5hfjdieyhidh3i3ceaL}
\end{EQA}
Moreover, \( \IFT_{\KS_{m}\KS_{m}}^{-1}(\prmtv) = (\targv \targv^{\T} + \muA^{2} \Id_{\dimp} + \GPKS_{m}^{2})^{-1} 
\leq (\muA^{2} \Id_{\dimp} + \GPKS_{m}^{2})^{-1} \leq \muA^{-2} \Id_{\dimp} \)
and 
\begin{EQA}
	&& \nquad
	\AFblk_{m}(\targv,\KS_{m}) \, \IFT_{\KS_{m}\KS_{m}}^{-1}(\prmtv) \, \AFblk_{m}^{\T}(\targv,\KS_{m})
	\leq 
	2 (\KS_{m}^{\T} \targv - \ima_{m})^{2} (\muA^{2} \Id_{\dimp} + \GPKS_{m}^{2})^{-1} 
	\\
	&&
	+ \, 2 \, \targv^{\T} (\targv \targv^{\T} + \muA^{2} \Id_{\dimp} + \GPKS_{m}^{2})^{-1} \targv \, \KS_{m} \, \KS_{m}^{\T}
	\\
	& \leq &
	\frac{2 (\KS_{m}^{\T} \targv - \ima_{m})^{2}}{\muA^{2}} \, \Id_{\dimp} 
	+ \frac{2 \| \targv \|^{2}}{\| \targv \|^{2} + \muA^{2}} \, \KS_{m} \, \KS_{m}^{\T} \, 
\label{4tge8vuyew3geywsxcgdmL}
\end{EQA}
yielding by \eqref{23fcgftfgwbsnhewretge}
\begin{EQA}
	\IFT_{\targv\KS}(\prmtv) \, \IFT_{\KS\KS}^{-1}(\prmtv) \, \IFT_{\KS\targv}(\prmtv)
	& \leq &
	\frac{2}{\muA^{2}} \| \imav - \KS \targv \|^{2} \Id_{\dimp} 
	+ \frac{2 \| \targv \|^{2}}{\| \targv \|^{2} + \muA^{2}} \, \KS^{\T} \KS 
	\\
	& \leq &
	\frac{\rhoIF}{2} (\KS^{\T} \KS + 2 \GPa^{2}) + \frac{\rhoIF}{2 (\rhoIF/4 + 1)} \, \KS^{\T} \KS .
\label{crduybv7e44er7ftyuL}
\end{EQA}
and for any \( \prmtv \in \Upsd \)
\begin{EQA}
	\IFU(\prmtv)
	& \geq &
	(1 - \rhoIF) \IFT_{\targv\targv}(\prmtv) .
\label{5xcytds65def5c5dwetw}
\end{EQA}
By \eqref{ijffd3e6dijuwsejk}, this yields \nameref{IFSref}.
Moreover, \( \IFT_{\targv\targv}(\prmtv) \geq \GP^{2} \) by \eqref{un4v7frue43kgv8uhkj}, and also
by \eqref{8dfnjey63f6ytjutyswwwq}, it holds
\( \IFT_{\KS\KS} \geq \blk\{ \HFblk_{1}^{2},\ldots,\HFblk_{\dimq}^{2} \} \)
with \( \HFblk_{m}^{2} = \muA^{2} \Id_{\dimp} + \GPKS_{m}^{2} \geq \muA^{2} \Id_{\dimp} \).
This and \nameref{IFSref} yields strong concavity of \( \lgd(\prmtv) \) on \( \Upsd \); see Lemma~\ref{Lidentsemi}.
\end{proof}

We already mentioned that the fidelity term \( \| \zv - \KS \targv \|^{2} \) is not convex 
in the scope of variables \( \targv,\KS \).
Proposition~\ref{PIFSrefEiO} shows that
adding the quadratic penalties helps to improve the situation, at least, locally.
Introduce the function
\begin{EQA}
	\lgdL(\prmtv)
	=
	\lgdL(\targv,\KS)
	&=&
	- \frac{1}{2} \| \zv - \KS \targv \|^{2} 
	- \frac{\muA^{2}}{2} \| \hKS - \KS \|_{\Fr}^{2} 
	- \frac{1}{2} \| \GPa \targv \|^{2} 
\label{bh2gfeft65t6f5rhKS0r}
\end{EQA}
with \( \GPa^{2} \leq \GP^{2} \).
In fact, we proved in Proposition~\ref{PIFSrefEiO} that it is strongly concave in \( (\targv,\KS) \) 
on the set \( \Upsd \) from \eqref{23fcgftfgwbsnhewretge}.
Later we denote
\begin{EQA}
	\IFT_{0}(\prmtv)
	&=&
	- \nabla^{2} \lgdL(\prmtv) .
\label{54hjvbn8uyi95eywe7edu}
\end{EQA}

\Subsection{Full and target efficient dimension}
The results about posterior concentration heavily rely on 
the \emph{full effective dimension} \( \dimLs(\prmtv) \)
and the \emph{target effective dimension} \( \dimLe(\prmtv) \) defined as
\begin{EQA}
\label{6hv643e4jhf98rwqikofu}
	\dimLs(\prmtv) 
	&=& 
	\tr\bigl\{ \IFT_{0}(\prmtv) \, \IFT^{-1}(\prmtv) \bigr\} ,
	\\
	\dimLe(\prmtv) 
	&=& 
	\tr\bigl\{ \IFT_{0,\targv\targv}(\prmtv) \, \IFT_{\targv\targv}^{-1}(\prmtv) \bigr\} 
	=
	\tr\bigl\{ (\KS^{\T} \KS + \GPa^{2}) (\KS^{\T} \KS + \GP^{2})^{-1} \bigr\}.
\label{6hv643e4jhf98rwqiko}
\end{EQA}
As \( \dimLe(\prmtv) \) only depends on \( \KS \), we write \( \dimLe(\KS) \) instead of \( \dimLe(\prmtv) \).
It appears that the full effective dimension \( \dimLs(\prmtv) \) can be bounded from above by 
the sum of the target effective dimension and the \emph{nuisance effective dimension} \( \dimqLe \) defined as
\begin{EQA}
	\dimqLe
	& \eqdef &
	\sum_{m=1}^{\dimq} \muA^{2} \, \tr( \muA^{2} \Id_{\dimp} + \GPKS_{m}^{2})^{-1} .
\label{6uviru8ie3ikfgvor56wujh}
\end{EQA}
The next lemma quantifies this statement.

\begin{lemma}
\label{LfulldimEoP}
Let 
\( \Upsd \) be given by \eqref{23fcgftfgwbsnhewretge}.
It holds for \( \dimLs(\prmtv) \) from \eqref{6hv643e4jhf98rwqikofu}
and any \( \prmtv = (\targv,\KS) \in \Upsd \) with \( \dimqLe \) from \eqref{6uviru8ie3ikfgvor56wujh} 
\begin{EQA}
	\dimLs(\prmtv)
	& \leq &
	\frac{\dimLe(\KS)}{1 - \rhoIF}
	+ \frac{1}{1 - \rhoIF} \sum_{m=1}^{\dimq} 
		\tr\bigl\{ (\targv \targv^{\T} + \muA^{2} \Id_{\dimp}) 
		(\targv \targv^{\T} + \muA^{2} \Id_{\dimp} + \GPKS_{m}^{2})^{-1} \bigr\} 
	\\
	& \leq &
	\frac{\dimLe(\KS)}{1 - \rhoIF}
	+ \frac{(1 + \rhoIF/4) \dimqLe}{1 - \rhoIF} \, .
\label{kifdgegviued0uvueud}
\end{EQA}
\end{lemma}

\begin{proof}
Proposition~\ref{PIFSrefEiO} and Lemma~\ref{Lidentsemi} yield
\begin{EQA}
	\dimLs(\prmtv)
	& \leq &
	\frac{1}{1 - \rhoIF} \, 
	\tr\bigl\{ \IFT_{0}(\prmtv) \, \blk\{ \IFT_{\targv\targv}^{-1}(\prmtv), \IFT_{\KS\KS}^{-1}(\prmtv)\} \bigr\}	
	\\
	& \leq &
	\frac{1}{1 - \rhoIF} \, \tr\bigl\{ \IFT_{0,\targv\targv}(\prmtv) \, \IFT_{\targv\targv}^{-1}(\prmtv) \bigr\}
	+ \frac{1}{1 - \rhoIF} \, \tr\bigl\{ \IFT_{0,\KS\KS}(\prmtv) \, \IFT_{\KS\KS}^{-1}(\prmtv) \bigr\} .
\label{89ciodwe3uedfcgtwhcjh}
\end{EQA}
Block structure of \( \IFT_{\KS\KS}(\prmtv) \) and \( \IFT_{0,\KS\KS}(\prmtv) \) due to 
\eqref{8dfnjey63f6ytjutyswwwq} yields
\begin{EQA}
	\tr\bigl\{ \IFT_{0,\KS\KS}(\prmtv) \, \IFT_{\KS\KS}^{-1}(\prmtv) \bigr\}
	& \leq &
	\sum_{m=1}^{\dimq} 
		\tr\bigl\{ (\targv \targv^{\T} + \muA^{2} \Id_{\dimp}) 
		(\targv \targv^{\T} + \muA^{2} \Id_{\dimp} + \GPKS_{m}^{2})^{-1} \bigr\}
\label{c7euer8bn8gyt895tr8i}
\end{EQA}
and the assertions follow in view of the bound \( \| \targv \|^{2} \leq \rhoIF \, \muA^{2} / 4 \).
\end{proof}

By Lemma~\ref{LfulldimEoP},
the full effective dimension exceeds the target effective dimension by the value of order
\( \dimqLe \).
An important message from this result is as follows: 
if the value \( \dimqLe \) from \eqref{6uviru8ie3ikfgvor56wujh} and thus, 
the full effective dimension \( \dimLs(\prmtv) \) is 
of the same order as the target effective dimension \( \dimLe(\prmtv) \), then
even full dimensional concentration and Laplace approximation deliver the desired quality
corresponding to the target effective dimension.
In particular, the plug-in procedure with \( \hKS \) in place of \( \KS \) is nearly efficient.

Without penalization of the operator \( \KS \) or for a small penalization with
\( \max_{m} \| \GPKS_{m}^{2} \| \ll \muA^{2} \), one can bound
\begin{EQA}
	\dimqLe
	& \approx &
	\dimp \, \dimq .
\label{od63d6e36dhd3oeffr3fgg}
\end{EQA}
Although \( \muA^{2} \) can be large, as in the random regression case with \( \muA^{2} \asymp n \),
the bound \eqref{od63d6e36dhd3oeffr3fgg} is not satisfactory because
it is not dimension free and involve the ambient dimension \( \dimp,\dimq \). 
The value \( \dimqLe \) and hence, \( \dimLs \) can be drastically reduced by using 
smoothness properties of the signal \( \targv \) and of the operator \( \KS \)
given in terms of the penalty \( \| \KS \|_{\GPKS}^{2}/2 \) for
\( \GPKS = \blk\{ \GPKS_{1}, \ldots,\GPKS_{\dimq} \} \).
As smoothness of \( \targv \) is anyway described by the penalty term \( \| \GP \targv \|^{2}/2 \),
it is natural to take \( \GPKS_{m}^{2} = \gpks_{m}^{2} \GP^{2} \).
The growth of the factors \( \gpks_{m}^{2} \) describes smoothness properties of 
the image of \( \KS \).

As usual, suppose that \( \gpks_{m}^{2} \) grow sufficiently fast, e.g. polynomially 
or exponentially, and similarly for the ordered eigenvalues \( \gp_{j}^{2} \) of \( \GP^{2} \).
Then for each \( m \) we can define \( j_{m} \) as the largest index \( j \) with 
\( \gp_{j}^{2} \leq \muA^{2} \, \gpks_{m}^{-2} \). 
For a polynomial growth of the \( g_{j}^{2} \), simple calculus show that 
\begin{EQA}
	\tr( \muA^{2} \Id_{\dimp} + \GPKS_{m}^{2})^{-1}
	=
	\muA^{-2} \, \muA^{2} \, \gpks_{m}^{-2} \tr( \muA^{2} \gpks_{m}^{-2} \Id_{\dimp} + \GP^{2})^{-1}
	& \asymp &
	\muA^{-2} \, j_{m} \, ;
\label{hdcgcabt78ewf78f786r738}
\end{EQA}
see Section~\ref{Seffdima}.
Therefore,
\begin{EQA}
	\dimqLe
	& \asymp &
	\sum_{m} j_{m} \, .
\label{ftrwdfygdusuwqys2h2}
\end{EQA}

\Subsection{Full dimensional concentration of the posterior}
We are now about to state the main results for the error-in-operator model. 
Denote \( \IFT_{0} = \IFT_{0}(\upsvs) \), \( \dimLs = \dimLs(\prmtvs) \) 
and define for \( \amax  = 2/3 \) the concentration set
\begin{EQ}[rcl]
	\CA
	&=&
	\bigl\{ \prmtv \colon \| \IFT_{0}^{1/2} (\prmtv - \prmtvs) \| \leq \amax^{-1} \rrbs \bigr\},
	\\
	\rrbs
	&=&
	2 \, \dimLs^{1/2} + (2\xx)^{1/2} . 
\label{4etgdf5e5w3e4ed4w3e453}
\end{EQ}

\begin{theorem}
\label{TEiOconc}
Let the set \( \CA \) 
from \eqref{4etgdf5e5w3e4ed4w3e453} satisfy \( \CA \subset \Upsd \) 
for \( \Upsd \) from \eqref{23fcgftfgwbsnhewretge}.
If \( \hmax_{3} \, \rrbs \, n^{-1/2} \leq 1/3 \) with \( \hmax_{3} = 6 \muA^{-1} \), then 
\begin{EQA}
	\PfL\bigl( \| \IFT_{0}^{1/2} (\prmtv - \prmtvs) \| > \amax^{-1} \rrLs \bigr)
	& \leq & 
	\ex^{-\xx}.
\label{32ghgh9nghbuujii9o}
\end{EQA}
Moreover, 
\begin{EQA}
	\PfL\bigl( (1 - \rhoIF) \muA \| \KS - \KSs \|_{\Fr} > \amax^{-1} \rrLs \bigr)
	& \leq & 
	\ex^{-\xx}.
\label{vsgx4s4rw2chswtghxg}
\end{EQA}
\end{theorem}

\begin{proof}
Later we show that \( \lgdL(\prmtv) \) fulfills the full dimensional smoothness condition \nameref{LLtS3ref}
on the set \( \Upsd \) with \( \hmax_{3} = 6 \muA^{-1} \),
and the result follows from \eqref{poybf3679jd532ff2} of Theorem~\ref{TLaplaceTV}. 
\end{proof}

\Subsection{Marginal posterior: concentration and Laplace approximation}
\label{SLaplaceEiO}
This section presents our main results about Laplace approximation of the marginal posterior.
We apply the general results of Section~\ref{SfullLaplco} to the considered setup after restricting 
the parameter space to \( \Upsd \).
Remind the notation \( \IFL(\prmtv) = \IFT_{\targv\targv}(\prmtv) = \KS^{\T} \KS + \GP^{2} \);
see \eqref{5xcytds65def5c5dwetw}.
As in Section~\ref{SfullLaplco}, for \( \dimLe(\KS) = \dimLe(\prmtv) \) from \eqref{6hv643e4jhf98rwqiko}, define 
\begin{EQA}
	\rrL(\KS)
	&=&
	2 \sqrt{\dimLe(\KS)} + \sqrt{2 \xx} \, .
\label{7yuhjcuye64e376fbgtrer}
\end{EQA}
Notation \( \IFL = \IFT_{\targv\targv}(\prmtvs) \) and 
\( \IFU = \IFU(\prmtvs) \) will be used.
One good news is that the efficient semiparametric matrix \( \IFU = \IFU(\prmtvs) \) from \eqref{5hfjdieyedyiwhihidh3i3ce}
satisfies \( (1 - \rhoIF) \IFL \leq \IFU \leq \IFL \). 

With \( \GPa \) shown in \eqref{bh2gfeft65t6f5rhKS0r}, define the \emph{effective sample size} \( n \) by the equation
\begin{EQA}
	n^{-1}
	&=&
	\| ({\KSs}^{\T} \KSs + \GPa^{2})^{-1} \| .
\label{jdndhnnfer445w6dhdy}
\end{EQA}
Note that this value does not depend on \( \muA \).
For a linear mapping \( \QP \) in \( \R^{\dimp} \), define
\begin{EQA}
	\BB_{\QP} 
	&=& 
	\frac{1}{\| \QP \, \IFL^{-1} \, \QP^{\T} \|} \,\, \QP \, \IFL^{-1} \, \QP^{\T}   ,
	\qquad 
	\dimQ = \tr(\BB_{\QP}) .
\label{5chcf7f6renecejeducjddhEO}
\end{EQA}
Theorem~\ref{TLaplaceopg} implies the following bound.

\begin{theorem}
\label{TLaplaceEO}
Assume the conditions of Theorem~\ref{TEiOconc}.
Let \( \QP \) be such that  
\begin{EQA}
	\| \BB_{\QP} \|_{\Fr}^{2}
	&=&
	\tr \BB_{\QP}^{2}
	\geq 
	\CONSTi_{0}^{2} \tr \BB_{\QP}
	=
	\CONSTi_{0}^{2} \, \dimQ
\label{cfirj4efd76ejhfcuew2djkdEO}
\end{EQA} 
with some fixed positive constant \( \CONSTi_{0} \).
For \( \IFU = \IFU(\prmtvs) \), it holds with \( \gaussv_{\IFU} \sim \ND(0,\IFU^{-1}) \)
\begin{EQA}
	&& \nquad
	\sup_{\rr > 0}
	\Bigl| \PfL\bigl( \| \QP (\Xv - \targvs) \| \leq \rr \bigr) - \P\bigl( \| \QP \, \gaussv_{\IFU} \| \leq \rr \bigr) \Bigr|
	\\
	& \lesssim &
		\frac{\muA^{-1} \, \rrL(\KSs) \, \dimLe(\KSs)}{\sqrt{n}}
		+ \frac{\muA^{-1} \, \rrLs \, \sqrt{\dimQ}}{\sqrt{n}} 
		+ \frac{\muA^{-2} \, {\rrLs}^{4}}{n \, \sqrt{\dimQ}}
		+ \ex^{- \xx} .
\label{jbfi8jfbe32hdfvdfrdfgopEO}
\end{EQA}
\end{theorem}

\begin{proof}
Consider \( \lgdL(\prmtv) \) from \eqref{bh2gfeft65t6f5rhKS0r}.
It obviously holds
\begin{EQA}
	- \nabla^{2}_{\targv\targv} \lgdL(\targv,\KS)
	&=&
	\DVL^{2}(\KS)
	=
	\KS^{\T} \KS + \GPa^{2} \, ,
\label{89bkekerivbmdsfjbtlDV}
	\\
	- \nabla^{2}_{\KS\KS} \lgdL(\targv,\KS) 
	&=&
	\blk\{ \targv \targv^{\T} + \muA^{2} \Id_{\dimp}, \ldots ,
		\targv \targv^{\T} + \muA^{2} \Id_{\dimp} \};
\label{w7ucf7fryeruhe4f6wedef}
\end{EQA}
see \eqref{8dfnjey63f6ytjutyswwwq}.
It is natural to describe the local geometry at \( \prmtv = (\targv,\KS) \) by 
\begin{EQA}
	\HL^{2}(\prmtv) 
	&=& 
	n^{-1} \blk\{ \DVL^{2}(\KS), \muA^{2} \Id \} 
\label{iiii774323dffdwqqdh}
\end{EQA}
with \( n \) from \eqref{jdndhnnfer445w6dhdy}.
The next step is in checking \nameref{LLtS3ref} and \nameref{LLtS4ref}.

\begin{lemma}
\label{LS3checkEO}
Let \( \muA \leq \muA \).
For the functions \( \lgd(\prmtv) \) and \( \lgdL(\prmtv) \), conditions 
\nameref{LLtS3ref} and \nameref{LLtS4ref} hold at any \( \prmtv \in \Upsd \) with 
\( \HL^{2}(\prmtv) \) from \eqref{iiii774323dffdwqqdh} and
\begin{EQA}
	\hmax_{3}
	&=&
	6 \muA^{-1} ,
	\qquad
	\hmax_{4}
	=
	3 \muA^{-2} .
\label{gdjddheebfyygredtgebegxg}
\end{EQA}
\end{lemma}

\begin{proof}
Fix \( \uv = (\xiv,\Uv) \) with 
\( \Uv^{\T} = (\nuiiv_{1},\ldots,\nuiiv_{\dimq}) \) such that
\begin{EQA}
	n \, \uv^{\T} \HL^{2}(\prmtv) \, \uv
	=
	\xiv^{\T} (\KS^{\T} \KS + \GPa^{2}) \xiv
	+ \sum_{m=1}^{\dimq} 
	\bigl( |\targv^{\T} \nuiiv_{m}|^{2} + \muA^{2} \| \nuiiv_{m} \|^{2} \bigr)
	& = &
	\rr^{2} .
	\qquad
\label{f564hjfgujg8747hfy76ehe}
\end{EQA}
Consider \( \lgdL(\prmtv + t \uv) = \lgdL(\targv + t \xiv,\KS + t \Uv) \).
The fourth derivative \( \nabla^{4} \lgdL(\prmtv + t \uv) \) in \( t \) does not depend on \( \prmtv \) and
it holds for any \( t \)
\begin{EQA}
	- \frac{d^{4}}{dt^{4}} \lgdL(\prmtv + t \uv) 
	& = &
	\langle \nabla^{4} \lgdL(\prmtv), \uv^{\otimes 4} \rangle
	=
	12  \sum_{m=1}^{\dimq} (\xiv^{\T} \nuiiv_{m})^{2} .
\label{ywkdofcfyweyenvcuduye}
\end{EQA}
By \eqref{f564hjfgujg8747hfy76ehe}, \eqref{jdndhnnfer445w6dhdy},
and the inequality \( 4ab \leq (a + b)^{2} \), it holds 
\begin{EQA}
	&& \nquad
	4 \muA^{2} \sum_{m=1}^{\dimq} (\xiv^{\T} \nuiiv_{m})^{2}
	\leq 
	4 \| \xiv \|^{2} \sum_{m=1}^{\dimq} \muA^{2} \| \nuiiv_{m} \|^{2}
	\\
	& \leq &
	4 \,\| (\KS^{\T} \KS + \GPa^{2})^{-1} \| \,
	\, \xiv^{\T} (\KS^{\T} \KS + \GPa^{2}) \xiv \, \,
	\sum_{m=1}^{\dimq} \muA^{2} \| \nuiiv_{m} \|^{2}
	\\
	& \leq &
	\| (\KS^{\T} \KS + \GPa^{2})^{-1} \| \, \rr^{4}
	=
	\rr^{4}/n \, .
\label{ugfr5544eehhjkudjdebh}
\end{EQA}
This easily implies \nameref{LLtS4ref} with \( \hmax_{4} \) from \eqref{gdjddheebfyygredtgebegxg}.

Now we proceed with \nameref{LLtS3ref}.
Definition yields
\begin{EQA}
	- \frac{d^{3}}{dt^{3}} \lgdL(\targv + t \xiv,\KS + t \Uv) \bigg|_{t=0}
	&=&
	6 \sum_{m=1}^{\dimq} (\targv^{\T} \nuiiv_{m} + \KS_{m}^{\T} \xiv ) \xiv^{\T} \nuiiv_{m} \, .
\label{ywkdofcfyweyenvcuduye}
\end{EQA}
The use of \eqref{f564hjfgujg8747hfy76ehe} yields 
similarly to \eqref{ugfr5544eehhjkudjdebh} with \( s = \muA^{-1} \rr \, n^{-1/2}/2 \)
\begin{EQA}
	&& \nquad
	2 \sum_{m=1}^{\dimq} \targv^{\T} \nuiiv_{m} \,\, \xiv^{\T} \nuiiv_{m}
	\leq 
	s \, \sum_{m=1}^{\dimq} (\targv^{\T} \nuiiv_{m})^{2} 
	+ s^{-1} \sum_{m=1}^{\dimq} (\xiv^{\T} \nuiiv_{m})^{2} 
	\\
	& \leq &
	s \,\rr^{2} + s^{-1} \muA^{-2} \rr^{4} /(4n)
	\leq 
	\muA^{-1} \, \rr^{3} \, n^{-1/2} .
\label{poiuytr5tgbcvfgyuhwe7}
\end{EQA}
Similarly, 
\begin{EQA}
	&& \nquad
	2 \sum_{m=1}^{\dimq} \KS_{m}^{\T} \xiv \,\, \xiv^{\T} \nuiiv_{m}
	\leq
	s \sum_{m=1}^{\dimq} (\KS_{m}^{\T} \xiv)^{2} 
	+ s^{-1} \sum_{m=1}^{\dimq} (\xiv^{\T} \nuiiv_{m})^{2}
	\\
	& \leq &
	s \, \, \xiv^{\T} \KS^{\T} \KS \, \xiv
	+ s^{-1} \, \muA^{-2} \rr^{4} /(4n)
	\\
	& \leq &
	s \, \rr^{2} + s^{-1} \muA^{-2} \rr^{4} /(4n)
	\leq 
	\muA^{-1} \, \rr^{3} \, n^{-1/2} .
\label{njio4h0fke7vuf4jel}
\end{EQA}
Summing up the latter bounds yields 
\begin{EQA}
	\bigl| \langle \nabla^{3} \lgdL(\prmtv), \uv^{\otimes 3} \rangle \bigr|
	& \leq &
	6 \muA^{-1} \, \rr^{3} \, n^{-1/2} 
\label{7rr3ecrxewwfwycel25g}
\end{EQA}
and \nameref{LLtS3ref} follows with \( \hmax_{3} = 6 \muA^{-1} \).
\end{proof}

Lemma~\ref{LS3checkEO} ensures \nameref{LLtS3ref} with \( \hmax_{3} \asymp \muA^{-1} \)
and the accuracy bound \eqref{jbfi8jfbe32hdfvdfrdfgopEO} of the theorem follows from Theorem~\ref{TLaplaceopg}.
\end{proof}

\Subsection{Critical dimension}
Now we discuss the issue of critical dimension.
With \( n \) from \eqref{jdndhnnfer445w6dhdy} and \( \dimLs \) from Lemma~\ref{LfulldimEoP},
Theorem~\ref{TEiOconc} requires \( \muA^{-1} \rrLs \ll n^{1/2} \) or 
\begin{EQA}
	\muA^{-2} \, {\rrLs}^{2}
	\approx
	\muA^{-2} \, \dimLs 
	\asymp
	\muA^{-2} (\dimLe + \dimqLe) 
	& \ll &
	n .
\label{7sswxtwy2hdywyhjjjtrr}
\end{EQA}
For the unpenalized case with \( \GPKS = 0 \), it leads to
\( \muA^{-2} \dimqLe = \muA^{-2} \dimp \, \dimq \ll n \); see \eqref{od63d6e36dhd3oeffr3fgg}.
A particular case with \( \muA^{-2} = \dimp = \dimq = n \) as in random design regression
is nearly included up to a small factor. 
The use of penalization by \( \| \KS \|_{\GPKS}^{2} \) results in the updated condition
\begin{EQA}
	\sum_{m=1}^{\dimq} \tr( \muA^{2} \Id_{\dimp} + \GPKS_{m}^{2})^{-1}
	& \ll &
	n .
\label{hd6dy3y366r646gsjjdueu}
\end{EQA}
This penalization allows to incorporate the cases with \( \dimp = \dimq \gg n \)
and even \( \dimp = \dimq = \infty \).

The condition on critical dimension in Theorem~\ref{TLaplaceEO} is even more striking.
Suppose that \( \QP = \DVL \) and \( \dimQ \approx \dimLe \).
Theorem~\ref{TLaplaceEO} applies under the conditions
\begin{EQA}
	\muA^{-2} \, \rrLs^{2} \, \dimLe
	& \ll &
	n,
	\qquad
	\muA^{-2} \, {\rrLs}^{4}
	\ll
	n \sqrt{\dimLe} \, .
\label{78jxdtdrtsg5rgghhye}
\end{EQA}
The use of \( \muA^{-2} \, {\rrLs}^{2} \approx \muA^{-2} \, \dimLs \asymp \muA^{-2} (\dimLe + \dimqLe) \)
and \eqref{6uviru8ie3ikfgvor56wujh} leads to
\begin{EQA}
	\dimLe \sum_{m=1}^{\dimq} \tr( \muA^{2} \Id_{\dimp} + \GPKS_{m}^{2})^{-1}
	& \ll &
	n ,
	\\
	\frac{\muA}{\dimLe^{1/4}} \sum_{m=1}^{\dimq} \tr( \muA^{2} \Id_{\dimp} + \GPKS_{m}^{2})^{-1}
	& \ll &
	\sqrt{n} \, .
\label{dd523tdhfcjcufurydsrw5w}
\end{EQA}
Without operator penalty, that is, for \( \GPKS_{m} \equiv 0 \) this becomes 
\begin{EQA}
	\dimLe \, \muA^{-2} \,\, \dimp \, \dimq
	& \ll &
	n,
	\\
	\dimLe^{1/4} \, \muA^{-1} \, \dimp \, \dimq
	& \ll &
	\sqrt{n} \, .
\label{tdrwgwdwey3wydhscgfghm}
\end{EQA}
While the first condition can be nearly fulfilled for \( \dimLe \) bounded 
and \( \dimp = \dimq \asymp n \),
the second condition fails completely when \( n \) is large.
However, the use of a proper operator \( \GPKS^{2} \) ensures the required ``critical dimension'' condition;
see Section~\ref{Spriors} for more details.

}

\newpage
\appendix


\Chapter{Local smoothness conditions}
\label{Slocalsmooth}
This section discusses different local smoothness characteristics of a 
multivariate function \( f(\upsv) = \E L(\upsv) \), \( \upsv \in \R^{\dimp} \).

\Section{Smoothness and self-concordance in Gateaux sense}
Below we assume 
the function \( f(\upsv) \) to be strongly concave with the negative Hessian 
\( \IFL(\upsv) \eqdef - \nabla^{2} f(\upsv) \in \Matr_{\dimp} \) positive definite. 
Also assume \( f(\upsv) \) three or sometimes even four times Gateaux differentiable in \( \upsv \in \Ups \).
For any particular direction \( \uv \in \R^{\dimp} \), we consider the univariate function 
\( f(\upsv + t \uv) \) and measure its smoothness in \( t \).
Local smoothness of \( f \) will be described by the relative error of the Taylor expansion 
of the third or four order.
Namely, define
\begin{EQ}[rcl]
	\dltw_{3}(\upsv,\uv) 
	&=& 
	f(\upsv + \uv) - f(\upsv) - \langle \nabla f(\upsv), \uv \rangle 
	- \frac{1}{2} \langle \nabla^{2} f(\upsv), \uv^{\otimes 2} \rangle , 
	\\
	\dltwd_{3}(\upsv,\uv) 
	&=&
	\langle \nabla f(\upsv + \uv), \uv \rangle - \langle \nabla f(\upsv), \uv \rangle 
	- \langle \nabla^{2} f(\upsv), \uv^{\otimes 2} \rangle \, ,
\label{dltw3vufuv12f2ga}
\end{EQ}
and
\begin{EQA}
	\dltw_{4}(\upsv,\uv)
	& \eqdef &
	f(\upsv + \uv) - f(\upsv) - \langle \nabla f(\upsv), \uv \rangle 
	- \frac{1}{2} \langle \nabla^{2} f(\upsv), \uv^{\otimes 2} \rangle
	- \frac{1}{6} \langle \nabla^{3} f(\upsv), \uv^{\otimes 3} \rangle .
\label{hvcduywgedfuyg2y1y35e3wweg}
\end{EQA}
Now, for each \( \upsv \), suppose to be given a positive symmetric operator 
\( \DVL(\upsv) \in \Matr_{\dimp} \) with \( \DVL^{2}(\upsv) \leq  \IFL(\upsv) = - \nabla^{2} f(\upsv) \)
defining a local metric and a local vicinity around \( \upsv \):
\begin{EQA}
	\UVz(\upsv)
	&=&
	\bigl\{ \uv \in \R^{\dimp} \colon \| \DVL(\upsv) \uv \| \leq \rr \bigr\}
\label{ed7sycf7wedwgedq2ftwdfgtv}
\end{EQA}
for some radius \( \rr \).

Local smoothness properties of \( f \) are given via the quantities
\begin{EQA}[rcccl]
    \dltwb(\upsv)
    & \eqdef &
    \sup_{\uv \colon \| \DVL(\upsv) \uv \| \leq \rr} \,
    \frac{2 |\dltw_{3}(\upsv,\uv)|}{\| \DVL(\upsv) \uv \|^{2}} 
    \,\, ,
    \qquad
    \dltwbd(\upsv)
    & \eqdef &
    \sup_{\uv \colon \| \DVL(\upsv) \uv \| \leq \rr} \, \frac{2 |\dltwd_{3}(\upsv,\uv)|}{\| \DVL(\upsv) \uv \|^{2}} \,\, . 
    \qquad
\label{dtb3u1DG2d3GPg}
\end{EQA}
The Taylor expansion yields for any \( \uv \) with \( \| \DVL(\upsv) \uv \| \leq \rr \)
\begin{EQ}[rcccl]
	\bigl| \dltw_{3}(\upsv,\uv) \rangle \bigr|
	& \leq &
	\frac{\dltwb(\upsv)}{2} \| \DVL(\upsv) \uv \|^{2} 
	\, ,
	\qquad
	\bigl| \dltwd_{3}(\upsv,\uv) \bigr|
	& \leq &
	\frac{\dltwbd(\upsv)}{2} \| \DVL(\upsv) \uv \|^{2}
	\, .
	\qquad
\label{dta3u1DG2d3GPa1g}
\end{EQ}
%
The introduced quantities \( \dltwb(\upsv) \), \( \dltwbd(\upsv) \) 
strongly depend on the radius \( \rr \) of the local vicinity \( \UVz(\upsv) \).
The results about Laplace approximation can be improved 
provided a homogeneous upper bound on the error of Taylor expansion. 
Assume a subset \( \Upsd \) of \( \Ups \) to be fixed.

\begin{description}
    \item[\label{LL3tref} \( \bb{(\mathcal{T}_{3})} \)]
      \textit{There exists \( \dltwu_{3} \) such that for all \( \upsv \in \Upsd \)}
\begin{EQA}
	\bigl| \dltw_{3}(\upsv,\uv) \bigr|
	& \leq &
	\frac{\dltwu_{3}}{6} \| \DVL(\upsv) \, \uv \|^{3} \, ,
	\quad
	\bigl| \dltwd_{3}(\upsv,\uv) \bigr|
	\leq 
	\frac{\dltwu_{3}}{2} \| \DVL(\upsv) \, \uv \|^{3} \, ,
	\quad
	\uv \in \UVz(\upsv).
\label{bd3xu16f3uo3st}
\end{EQA}
\end{description}
 
\begin{description}
    \item[\label{LL4tref} \( \bb{(\mathcal{T}_{4})} \)]
      \textit{ There exists \( \dltwu_{4} \) such that for all \( \upsv \in \Upsd \)}
\begin{EQA}
	\bigl| \dltw_{4}(\upsv,\uv) \bigr|
	& \leq &
	\frac{\dltwu_{4}}{24} \| \DVL(\upsv) \, \uv \|^{4} \, ,
	\qquad
	\uv \in \UVz(\upsv).
\label{1mffmxum5st}
\end{EQA}
\end{description}

\begin{lemma}
\label{LdltwLa3t}
Under \nameref{LL3tref},
the values \( \dltwb(\upsv) \) and \( \dltwbd(\upsv) \) from \eqref{dtb3u1DG2d3GPg} satisfy
\begin{EQA}[rcccl]
\label{gtcdsftdffrvsewsea}
	\dltwb(\upsv)
	& \leq &
	\frac{\dltwu_{3} \, \rr}{3 } \, ,
	\qquad
	\dltwbd(\upsv)
	& \leq &
	{\dltwu_{3} \, \rr} \, ,
	\qquad
	\upsv \in \Upsd .
\label{gtcdsftdfvtwdsefhfdvfrvsewseG}
\end{EQA}
\end{lemma}

\begin{proof}
For any \( \uv \in \UVz(\upsv) \) with \( \| \DVL(\upsv) \uv \| \leq \rr \)
\begin{EQA}
	\bigl| \dltw_{3}(\upsv,\uv) \bigr|
	& \leq &
	\frac{\dltwu_{3}}{6} \, \| \DVL(\upsv) \uv \|^{3} 
	\leq 
	\frac{\dltwu_{3} \, \rr}{6} \, \| \DVL(\upsv) \uv \|^{2},
\label{jrgeteteer2234587654}
\end{EQA}
and the bound for \( \dltwb(\upsv) \) follows.
The proof for \( \dltwbd(\upsv) \) is similar.
\end{proof}

The values \( \dltwu_{3} \) and \( \dltwu_{4} \) are usually very small.
Some quantitative bounds are given later in this section
under the assumption that the function \( f(\upsv) = \E \LGP(\upsv) \) can be written in the form \( - f(\upsv) = n \hL(\upsv) \) 
for a fixed smooth function \( h(\upsv) \) with the Hessian \( \nabla^{2} \hL(\upsv) \). 
The factor \( n \) has meaning of the sample size%
\ifapp{; see \Chname \ref{ScritdimMLE} or \Chname \ref{SGBvM}.}{.}

\begin{description}
    \item[\label{LLtS3ref} \( \bb{(\mathcal{S}_{3})} \)]
      \emph{ \( - f(\upsv) = n \hL(\upsv) \) for \( \hL(\upsv) \) convex with 
      \( \nabla^{2} \hL(\upsv) \geq \HL^{2}(\upsv) = \DVL^{2}(\upsv)/n \) 
      and
\begin{EQA}
	\sup_{\uv \colon \| \HL(\upsv) \uv \| \leq \rr/\sqrt{n}} 
	\frac{\bigl| \langle \nabla^{3} \hL(\upsv + \uv), \uv^{\otimes 3} \rangle \bigr|}{\| \HL(\upsv) \uv \|^{3}}
	& \leq &
	\hmax_{3} \, .
\end{EQA}
}
    \item[\label{LLtS4ref} \( \bb{(\mathcal{S}_{4})} \)]
      \emph{ the function \( \hL(\cdot) \) satisfies \nameref{LLtS3ref} and  
\begin{EQA}
	\sup_{\uv \colon \| \HL(\upsv) \uv \| \leq \rr/\sqrt{n}}
	\frac{\bigl| \langle \nabla^{4} \hL(\upsv + \uv), \uv^{\otimes 4} \rangle \bigr|}{\| \HL(\upsv) \uv \|^{4}}
	& \leq &
	\hmax_{4} \, .
\end{EQA}
}
\end{description}

\noindent
\nameref{LLtS3ref} and \nameref{LLtS4ref}
are local versions of the so called self-concordance condition; see \cite{Ne1988}.
In fact, they require that each univariate function \( \hL(\upsv + t \uv) \) of \( t \in \R \)
is self-concordant with some universal constants \( \hmax_{3} \) and \( \hmax_{4} \).
Under \nameref{LLtS3ref} and \nameref{LLtS4ref}, we can use \( \DVL^{2}(\upsv) = n \, \HL^{2}(\upsv) \) 
and easily bound the values 
\( \dltw_{3}(\upsv,\uv) \), \( \dltw_{4}(\upsv,\uv) \), and \( \dltwb(\upsv) \), \( \dltwbd(\upsv) \).

\begin{lemma}
\label{LdltwLaGP}
Suppose \nameref{LLtS3ref}.
Then 
\nameref{LL3tref} follows with \( \dltwu_{3} = \hmax_{3} n^{-1/2} \).
Moreover, for \( \dltwb(\upsv) \) and \( \dltwbd(\upsv) \) from \eqref{dtb3u1DG2d3GPg}, it holds
\begin{EQA}[rcccl]
	\dltwb(\upsv)
	& \leq &
	\frac{\hmax_{3} \, \rr}{3 n^{1/2}} \, ,
	\qquad
	\dltwbd(\upsv)
	& \leq &
	\frac{\hmax_{3} \, \rr}{n^{1/2}} \, .
\label{gtcdsftdfvtwdsefhfdvfrvsewseGP}
\end{EQA}
Also \nameref{LL4tref} follows from \nameref{LLtS4ref} with \( \dltwu_{4} = \hmax_{4} n^{-1} \).
\end{lemma}

\begin{proof}
For any \( \uv \in \UVz(\upsv) \) and \( t \in [0,1] \), by the Taylor expansion of the third order
\begin{EQA}
	|\dltw(\upsv,\uv)|
	& \leq &
	\frac{1}{6} \bigl| \langle \nabla^{3} f(\upsv + t \uv), \uv^{\otimes 3} \rangle \bigr|
	=
	\frac{n}{6} \, \bigl| \langle \nabla^{3} \hL(\upsv + t \uv), \uv^{\otimes 3} \rangle \bigr|
	\leq 
	\frac{n \, \hmax_{3}}{6} \, \| \HL(\upsv) \uv \|^{3} 
	\\
	&=&
	\frac{n^{-1/2} \, \hmax_{3}}{6} \, \| \DVL(\upsv) \uv \|^{3}
	\leq 
	\frac{n^{-1/2} \, \hmax_{3} \, \rr}{6} \, \| \DVL(\upsv) \uv \|^{2} \, .
\label{jrgeteteer2234587654}
\end{EQA}
This implies \nameref{LL3tref} as well as \eqref{gtcdsftdfvtwdsefhfdvfrvsewseGP}; see \eqref{dta3u1DG2d3GPa1g}.
The statement about \nameref{LL4tref} is similar.
\end{proof}

\Section{Fr\'echet derivatives and smoothness of the Hessian}
\label{SlocsmooFr}
For evaluation of the bias, 
we also need stronger smoothness conditions in the Fr\'echet sense.
Let \( f \) be a strongly concave function. 
Essentially we need some continuity of the negative Hessian \( \IFL(\upsv) = - \nabla^{2} f(\upsv) \).
For \( \upsv \in \Ups \), 
define \( \DVL(\upsv) = \IFL^{1/2}(\upsv) \) and
\begin{EQA}
	\dltwbss(\upsv)
    & \eqdef & 
    \sup_{\uv \colon \| \DVL(\upsv) \uv \| \leq \rr} \,\, \sup_{\gammav \in \R^{\dimp}} \,\, 
    \frac{|\langle \IFL(\upsv + \uv) - \IFL(\upsv), \gammav^{\otimes 2} \rangle|}{\| \DVL(\upsv) \gammav \|^{2} } \, .
\label{jcxhydtferyufgy7tfdsy7ft}
\end{EQA}
This definition of \( \dltwbss(\upsv) \) is, of course, stronger than
the one-directional definition of \( \dltwb(\upsv) \) in \eqref{dtb3u1DG2d3GPg}.
However, in typical examples these quantities \( \dltwb(\upsv) \) and \( \dltwbss(\upsv) \) are similar%
\ifapp{; see the examples from \Chname \ref{Saniligot} and \Chname \ref{SGBvM}.}{.}

We also present a Fr\'echet version of \nameref{LLtS3ref}.
\begin{description}

    \item[\label{LLsS3ref} \( \bb{(\mathcal{S}_{3}^{+})} \)]
      \emph{\( - f(\upsv) = n \hL(\upsv) \) with \( \hL(\cdot) \) strongly concave, 
      \( \HL^{2}(\upsv) \leq \nabla^{2} \hL(\upsv) \), 
      and
\begin{EQA}
	\sup_{\| \HL(\upsv) \uv \| \leq \rr/\sqrt{n}}  \,\, \sup_{\gammav \in \R^{\dimp}} \,\,
	\frac{\bigl| \langle \nabla^{3} \hL(\upsv + \uv), \uv \otimes \gammav^{\otimes 2} \rangle \bigr|}
		 {\| \HL(\upsv) \uv \|\, \| \HL(\upsv) \gammav \|^{2}}
	& \leq &
	\hmax_{3} \, .
\end{EQA}
}
\end{description}

\begin{lemma}
\label{LfreTay}
For \( \dltwbss(\upsv) \) from \eqref{jcxhydtferyufgy7tfdsy7ft} and any \( \uv \) with 
\( \| \IFL^{1/2}(\upsv) \uv \| \leq \rr \), it holds
\begin{EQA}
	\| \IFL^{-1/2}(\upsv) \, \IFL(\upsv + \uv) \, \IFL^{-1/2}(\upsv) - \Id_{\dimp} \|
	& \leq &
	\dltwbss(\upsv) .
\label{ghdrd324ee4ew222w3ew}
\end{EQA}
Moreover, \nameref{LLsS3ref} yields \( \dltwbss(\upsv) \leq \hmax_{3} \, \rr \, n^{-1/2} \) and
for any \( \uv \) with \( \| \IFL^{1/2}(\upsv) \uv \| \leq \rr \)
\begin{EQA}
	\| \IFL^{-1/2}(\upsv) \, \IFL(\upsv + \uv) \, \IFL^{-1/2}(\upsv) - \Id_{\dimp} \|
	& \leq &
	\frac{\hmax_{3}}{n^{1/2}} \, \| \IFL^{1/2}(\upsv) \uv \|
	\leq 
	\frac{\hmax_{3} \, \rr}{n^{1/2}} \, .
\label{ghdrd324ee4ew222w3ewss}
\end{EQA}
\end{lemma}

\begin{proof}
Denote 
\( \Delta(\uv) = \IFL(\upsv + \uv) - \IFL(\upsv) \).
Then by \eqref{jcxhydtferyufgy7tfdsy7ft} for any \( \gammav \in \R^{\dimp} \) with \( \deltav = \IFL^{-1/2}(\upsv) \gammav \)
\begin{EQA}
	&& \nquad
	\bigl| \bigl\langle \IFL^{-1/2}(\upsv) \, \IFL(\upsv + \uv) \, \IFL^{-1/2}(\upsv) - \Id_{\dimp}, \gammav^{\otimes 2} \bigr\rangle \bigr|
	=
	\bigl| \langle \Delta(\uv), \deltav^{\otimes 2} \rangle \bigr|
	\\
	& \leq &
	\dltwbss(\upsv) \, \| \IFL^{1/2}(\upsv) \, \deltav \|^{2} 
	=
	\dltwbss(\upsv) \, \| \gammav \|^{2} .
\label{9487654r5tghasdfg}
\end{EQA}
This yields \eqref{ghdrd324ee4ew222w3ew}.
The arguments from Lemma~\ref{LdltwLaGP} yield bound \eqref{ghdrd324ee4ew222w3ewss}.
\end{proof}

Define now for any \( \uv \)
\begin{EQA}
	\IFLba(\upsv;\uv)
	& \eqdef &
	\int_{0}^{1} \IFL(\upsv + t \uv) \, dt \, .
\label{vhudfg7sdfyt7s8hgbfuhbugrh9}
\end{EQA}

\begin{lemma}
\label{LfreTayb}
For any \( \uv \in \R^{\dimp} \) with 
\( \| \IFL^{1/2}(\upsv) \uv \| \leq \rr \)
\begin{EQ}[rcl]
	\bigl\| \IFL^{-1/2}(\upsv) \, \IFLba(\upsv ; \uv) \, \IFL^{-1/2}(\upsv) - \Id_{\dimp} \bigr\|
	& \leq &
	\dltwbss(\upsv) 
	\, .
\label{dsvfhgevrfgtwfyikoub}
\end{EQ}
Moreover,
\begin{EQA}[rcccl]
	\{ 1 - \dltwbss(\upsv) \} \IFL(\upsv)
	& \leq &
	\IFLba(\upsv ; \uv)
	& \leq &
	\{ 1 + \dltwbss(\upsv) \} \IFL(\upsv) \, .
\label{5chfdc7e3yvc5ededww}
\end{EQA}
\end{lemma}

\begin{proof}
For any \( \gammav \in \R^{\dimp} \) and \( t \in [0,1] \), the definition \eqref{jcxhydtferyufgy7tfdsy7ft} implies
\begin{EQA}
	\bigl| \bigl\langle \IFL(\upsv + t \uv) - \IFL(\upsv), \gammav^{\otimes 2} \bigr\rangle \bigr|
	& \leq &
	\dltwbss(\upsv) \, \| \IFL^{1/2}(\upsv) \, \gammav \|^{2} .
\label{09vbfkmdior329iufe9o}
\end{EQA}
This obviously yields under \nameref{LLsS3ref}
\begin{EQA}
	\bigl| \bigl\langle \IFLba(\upsv; \uv) - \IFL(\upsv), \gammav^{\otimes 2} \bigr\rangle \bigr|
	& \leq &
	\frac{\hmax_{3} \, \rr}{n^{1/2}} \, \| \IFL^{1/2}(\upsv) \, \gammav \|^{2} \int_{0}^{1} t \, dt
	=
	\frac{\hmax_{3} \, \rr}{2 n^{1/2}} \, \| \IFL^{1/2}(\upsv) \, \gammav \|^{2} 
\label{09vbfkmdior329iufe9ob}
\end{EQA}
and \eqref{dsvfhgevrfgtwfyikoub} follows as in Lemma~\ref{LfreTay}.
\end{proof}

If \( f(\upsv) \) is quadratic, then \( \nabla f(\upsv) \) is linear while 
\( \IFL \equiv - \nabla^{2} f(\upsv) \) is constant and we obtain the identity
\begin{EQA}
	\IFL^{-1} \bigl\{ \nabla f(\upsv + \uv) - \nabla f(\upsv) \bigr\}
	&=&
	- \uv .
\label{7dhc5fye4378fjfrurivj}
\end{EQA}
Now we consider the case of a smooth \( f(\cdot) \) satisfying \eqref{jcxhydtferyufgy7tfdsy7ft} 
or \nameref{LLsS3ref}.

\begin{lemma}
\label{LfreTayc}
Suppose the conditions of Lemma~\ref{LfreTay}. 
For any \( \uv \) with \( \| \IFL^{1/2}(\upsv) \uv \| \leq \rr \)
\begin{EQA}
	\nabla f(\upsv + \uv) - \nabla f(\upsv) 
	&=&
	- \IFLba(\upsv ; \uv) \, \uv .
\label{7jf98vgiey47fvrtirdvbk}
\end{EQA}
Moreover, for any positive definite symmetric operator  
\( \QP \colon \R^{\dimp} \to \R^{\dimp} \), it holds
\begin{EQA}[rcl]
	\bigl\| \QP \, \bigl\{ \nabla f(\upsv + \uv) - \nabla f(\upsv) \bigr\} \bigr\|
	& \leq &
	\bigl\{ 1 + \dltwbss(\upsv) \bigr\} \, \| \QP \, \IFL(\upsv) \, \uv \| ,
	\\
	\| \QP \, \IFL(\upsv) \, \uv \|
	& \leq &
	\frac{1}{1 - \dltwbss(\upsv)} \, 
	\bigl\| \QP \, \bigl\{ \nabla f(\upsv + \uv) - \nabla f(\upsv) \bigr\} \bigr\| ,
\label{5chfdc7e3yvc5ededwc}
\end{EQA}
and
\begin{EQA}
	\bigl\| \QP \, \bigl\{ \nabla f(\upsv + \uv) - \nabla f(\upsv) + \IFL(\upsv) \, \uv \bigr\} \bigr\|
	& \leq &
	\dltwbss(\upsv) \, \| \QP \, \IFL(\upsv) \, \uv \| .
\label{5chfdc7e3yvc5ededwcuv}
\end{EQA}
\end{lemma}

\begin{proof}
First we check \eqref{7jf98vgiey47fvrtirdvbk}.
Indeed, for any \( \gammav \in \R^{\dimp} \), consider the univariate function 
\( h(t) = \langle \nabla f(\upsv + t \uv) - \nabla f(\upsv), \gammav \rangle \).
The the statement follows from definition \eqref{vhudfg7sdfyt7s8hgbfuhbugrh9} 
and the identity
\( h(1) - h(0) = \int_{0}^{1} h'(t) \, dt \).
Further, it holds by \eqref{7jf98vgiey47fvrtirdvbk} and \eqref{5chfdc7e3yvc5ededww}
\begin{EQA}
	&& \nquad
	\bigl\| \QP \, \bigl\{ \nabla f(\upsv + \uv) - \nabla f(\upsv) \bigr\} \bigr\|
	=
	\| \QP \, \IFLba(\upsv ; \uv) \, \uv \|
	=
	\bigl\| \QP \, \IFLba(\upsv ; \uv) \IFL^{-1}(\upsv) \, \QP^{-1} \, \QP \, \IFL(\upsv) \, \uv \bigr\|
	\\
	& \leq &
	\bigl\| \QP \, \IFLba(\upsv ; \uv) \IFL^{-1}(\upsv) \, \QP^{-1} \bigr\| \, \,
	\| \QP \, \IFL(\upsv) \, \uv \|
	\leq 
	\bigl\{ 1 + \dltwbss(\upsv) \bigr\} \, \| \QP \, \IFL(\upsv) \, \uv \|
\label{9ckjf7rf74er4re7yfwweee}
\end{EQA}
and
\begin{EQA}
	\| \QP \uv \|
	&=&
	\bigl\| \QP \, \IFLba^{-1}(\upsv ; \uv) \, \bigl\{ \nabla f(\upsv + \uv) - \nabla f(\upsv) \bigr\} \bigr\|
	\\
	& = &
	\bigl\| \QP \, \IFLba^{-1}(\upsv ; \uv) \, \IFL(\upsv) \, \QP^{-1} \, 
		\QP \, \IFL^{-1}(\upsv) \, \bigl\{ \nabla f(\upsv + \uv) - \nabla f(\upsv) \bigr\} 
	\bigr\|
	\\
	& \leq &
	\bigl\| \QP \, \IFLba^{-1}(\upsv ; \uv) \, \IFL(\upsv) \, \QP^{-1} \bigr\|
	\,\,
	\bigl\| \QP \, \IFL^{-1}(\upsv) \, \bigl\{ \nabla f(\upsv + \uv) - \nabla f(\upsv) \bigr\} \bigr\|
	\\
	& \leq &
	\frac{1}{1 - \dltwbss(\upsv)} \, 
	\bigl\| \QP \, \IFL^{-1}(\upsv) \, \bigl\{ \nabla f(\upsv + \uv) - \nabla f(\upsv) \bigr\} \bigr\|
\label{oiuygf5tgh98y3eyhhui89}
\end{EQA}
as required.
Similarly 
\begin{EQA}
	&& \nquad
	\bigl\| \QP \bigl\{ \nabla f(\upsv + \uv) - \nabla f(\upsv) + \IFL(\upsv) \, \uv \bigr\} \bigr\|
	=
	\bigl\| \QP \bigl\{ \IFLba(\upsv ; \uv) - \IFL(\upsv) \bigr\} \uv \bigr\|
	\leq 
	\dltwbss(\upsv) \, \| \QP \, \IFL(\upsv) \, \uv \|
\label{9dhdfre64fvieefhvfue}
\end{EQA}
and \eqref{5chfdc7e3yvc5ededwcuv} follows from \eqref{5chfdc7e3yvc5ededww}.
\end{proof}


\Section{Optimization after linear or quadratic perturbation}
\label{Squadnquad}

Let \( \fs(\upsv) \) be a smooth concave function, 
\begin{EQA}
	\upsvs
	&=&
	\argmax_{\upsv} \fs(\upsv),
\label{fg5hg3gf98tkj3dciryt}
\end{EQA}
and \( \IFL = \DVL^{2} = - \nabla^{2} \fs(\upsvs) \).
Later we study the question how the point of maximum and the value of maximum of \( \fs \) change if we add a linear or quadratic 
component to \( \fs \).

\Subsection{A linear perturbation}
This section studies the case of a linear change of \( \fs \).
More precisely, let another function \( \fn(\upsv) \) satisfy for some vector \( \Av \)
\begin{EQA}
	\fn(\upsv) - \fn(\upsvs) 
	&=&
	\bigl\langle \upsv - \upsvs, \Av \bigr\rangle + \fs(\upsv) - \fs(\upsvs) .
\label{4hbh8njoelvt6jwgf09}
\end{EQA}
A typical example corresponds to \( \fs(\upsv) = \E L(\upsv) \) and \( \fn(\upsv) = L(\upsv) \) 
for a random function \( L(\upsv) \) with a linear stochastic component \( \zeta(\upsv) = L(\upsv) - \E L(\upsv) \)%
\ifapp{; see \nameref{Eref}.}{.}
Then \eqref{4hbh8njoelvt6jwgf09} is satisfied with 
\begin{EQA}
	\Av 
	&=&
	\nabla \zeta .
\label{mdfwdgfwes3e4gfdyc7f}
\end{EQA}
The aim of the analysis is evaluate the values 
\begin{EQA}
	\upsvn
	& \eqdef &
	\argmax_{\upsv} \fn(\upsv),
	\qquad
	\fn(\upsvn)
	=
	\max_{\upsv} \fn(\upsv) .
\label{6yc63yhudf7fdy6edgehy} 
\end{EQA}
The results will be stated mainly in terms of the quantity \( \| \IFL^{-1/2} \Av \| \).
First we consider the case of a quadratic function \( \fs \).

\begin{lemma}
\label{Pquadquad}
Let \( \fs(\upsv) \) be quadratic with \( \nabla^{2} \fs(\upsv) \equiv - \IFL \).
If \( \fn(\upsv) \) satisfy \eqref{4hbh8njoelvt6jwgf09}, then 
\begin{EQA}
	\upsvn - \upsvs
	&=&
	\IFL^{-1} \Av,
	\qquad
	\fn(\upsvn) - \fn(\upsvs)
	=
	\frac{1}{2} \| \IFL^{-1/2} \Av \|^{2} .
\label{kjcjhchdgehydgtdtte35}
\end{EQA}
\end{lemma}

\begin{proof}
If \( \fs(\upsv) \) is quadratic, then, of course, under \eqref{4hbh8njoelvt6jwgf09}, \( \fn(\upsv) \) is quadratic as well
with \( - \nabla^{2} \fn(\upsv) \equiv - \IFL \).
This implies
\begin{EQA}
	\nabla \fn(\upsvs) - \nabla \fn(\upsvn)
	&=&
	\IFL (\upsvn - \upsvs) .
\label{dcudydye67e6dy3wujhds7}
\end{EQA}
Further, \eqref{4hbh8njoelvt6jwgf09} and \( \nabla \fs(\upsvs) = 0 \) yield \( \nabla \fn(\upsvs) = \Av \).
Together with \( \nabla \fn(\upsvn) = 0 \), this implies
\( \upsvn - \upsvs = \IFL^{-1} \Av \).
The Taylor expansion of \( \fn \) at \( \upsvn \) yields by \( \nabla \fn(\upsvn) = 0 \)
\begin{EQA}
	\fn(\upsvs) - \fn(\upsvn)
	&=&
	- \frac{1}{2} \| \IFL^{1/2} (\upsvn - \upsvs) \|^{2}
	=
	- \frac{1}{2} \| \IFL^{-1/2} \Av \|^{2} 
\label{8chuctc44wckvcuedje}
\end{EQA}
and the assertion follows.
\end{proof}
  
The next result describes the concentration properties of \( \upsvn \) from \eqref{6yc63yhudf7fdy6edgehy} in a local elliptic set
of the form
\begin{EQA}
	\CA(\rr)
	& \eqdef &
	\{ \upsv \colon \| \IFL^{1/2} (\upsv - \upsvs) \| \leq \rr \} ,
\label{0cudc7e3jfuyvct6eyhgwe}
\end{EQA}
where \( \rr \) is slightly larger than \( \| \IFL^{-1/2} \Av \| \).

\begin{proposition}
\label{Pconcgeneric}
Let \( \fs(\upsv) \) be a strongly concave function with \( \fs(\upsvs) = \max_{\upsv} \fs(\upsv) \)  
and \( \IFL = - \nabla^{2} \fs(\upsvs) \).
Let further \( \fn(\upsv) \) and \( \fs(\upsv) \) be related by \eqref{4hbh8njoelvt6jwgf09} with some vector \( \Av \).
Fix \( \amax \leq 2/3 \) and \( \rrn \) such that \( \| \IFL^{-1/2} \Av \| \leq \amax \, \rrn \).
Suppose now that \( \fs(\upsv) \) satisfy \eqref{dtb3u1DG2d3GPg} for \( \upsv = \upsvs \), 
\( \DVL(\upsvs) = \IFL^{1/2} = \DVL \), 
and \( \dltwbd \) such that 
\begin{EQA}
	1 - \amax - \dltwbd
	& \geq &
	0 .
\label{rrm23r0ut3ua}
\end{EQA}
Then for \( \upsvn \) from \eqref{6yc63yhudf7fdy6edgehy}, it holds 
\begin{EQA}
	\| \IFL^{1/2} (\upsvn - \upsvs) \|  
	& \leq &
	\rrn \, . 
\label{rhDGtuGmusGU0a}
\end{EQA}
\end{proposition}

\begin{proof}
With \( \DVL = \IFL^{1/2} \),
the bound \( \| \DVL^{-1} \Av \| \leq \amax \, \rrn \) implies for any \( \uv \)
\begin{EQA}
	\bigl| \langle \Av, \uv \rangle \bigr|
	& = &
	\bigl| \langle \DVL^{-1} \Av, \DVL \uv \rangle \bigr|
	\leq 
	\amax \, \rrn \| \DVL \uv \| \, .
\label{LLoDGm1nzua}
\end{EQA}
If \( \| \DVL \uv \| > \rrn \), then \( \rrn \| \DVL \uv \| \leq \| \DVL \uv \|^{2} \). 
Therefore,
\begin{EQA}
	\bigl| \langle \Av, \uv \rangle \bigr|
	& \leq &
	\amax \| \DVL \uv \|^{2} \, ,
	\qquad
	\| \DVL \uv \| > \rrn \, .
\label{LLoDGm1nzun}
\end{EQA}
Let \( \upsv \) be a point on the boundary of the set \( \CA(\rrn) \) from \eqref{0cudc7e3jfuyvct6eyhgwe}.
We also write \( \uv = \upsv - \upsvs \).
The idea is to show that the derivative  \( \frac{d}{dt} \fn(\upsvs + t \uv) < 0 \) 
is negative for \( t > 1 \).
Then all the extreme points of \( \fn(\upsv) \) are within \( \CA(\rrn) \).
We use the decomposition
\begin{EQA}
	\fn(\upsvs + \rhot \uv) - \fn(\upsvs)
	&=&
	\bigl\langle \Av, \uv \bigr\rangle \, \rhot 
	+ \fs(\upsvs + \rhot \uv) - \fs(\upsvs) .
\label{LGtsGtuLGtsa}
\end{EQA}
With \( \fGu(t) = \fs(\upsvs + \rhot \uv) \), it holds
\begin{EQA}
	\frac{d}{d \rhot} \fGu(\upsvs + \rhot \uv)
	&=&
	\bigl\langle \Av, \uv \bigr\rangle + \fGu'(\rhot) .
\label{frddtLtGstua}
\end{EQA}
By definition of \( \upsvs \), it also holds \( \fGu'(0) = 0 \).
The identity \( \nabla^{2} \fs(\upsvs) = - \DVL^{2} \) yields \( \fGu''(0) = - \| \DVL \uv \|^{2} \).
Bound \eqref{dta3u1DG2d3GPa1g} implies for \( | \rhot | \leq 1 \)
\begin{EQA}
	\bigl| \fGu'(\rhot) - \rhot \fGu''(0) \bigr|
	&=&
	\bigl| \fGu'(\rhot) - \fGu'(0) - \rhot \fGu''(0) \bigr|
	\leq 
	\rhot^{2} \, \bigl| \fGu''(0) \bigr| \, \dltwbd \, .
\label{fptfp0fpttfpp13a}
\end{EQA}
For \( \rhot = 1 \), we obtain by \eqref{rrm23r0ut3ua}
\begin{EQA}
	\fGu'(1) 
	& \leq &
	\fGu''(0) - \fGu''(0) \, \dltwbd
	=
	- \bigl| \fGu''(0) \bigr| (1 -  \dltwbd)
	< 0 .
\label{fp1fpp13d3rGa}
\end{EQA}
Moreover, concavity of \( \fGu(\rhot) \) and \( \fGu'(0) = 0 \) imply that \( \fGu'(\rhot) \) decreases in 
\( \rhot \) for \( \rhot > 1 \).
Further, summing up the above derivation yields 
\begin{EQA}
	\frac{d}{dt} \fGu(\upsvs + \rhot \uv) \Big|_{\rhot=1}
	& \leq &
	- \| \DVL \uv \|^{2} (1 - \amax - \dltwbd)
	< 0 .
\label{ddtLGtstu33a}
\end{EQA}
As \( \frac{d}{d \rhot} \fGu(\upsvs + \rhot \uv) \) decreases with \( \rhot \geq 1 \) together with 
\( \fGu'(\rhot) \) due to \eqref{frddtLtGstua}, the same applies to all such \( \rhot \).
This implies the assertion.
\end{proof}

The result of Proposition~\ref{Pconcgeneric} allows to localize the point \( \upsvn = \argmax_{\upsv} \fn(\upsv) \)
in the local vicinity \( \CA(\rrn) \) of \( \upsvs \).
The use of smoothness properties of \( \fn \) or, equivalently, of \( \fs \), in this vicinity helps to obtain
rather sharp expansions for \( \upsvn - \upsvs \) and for \( \fn(\upsvn) - \fn(\upsvs) \); cf. \eqref{kjcjhchdgehydgtdtte35}.

\begin{proposition}
\label{PFiWigeneric}
Under the conditions of Proposition~\ref{Pconcgeneric}
\begin{EQ}[rcccl]
    - \frac{\dltwb}{1 + \dltwb} \| \DVL^{-1} \Av \|^{2}
    & \leq &
    2 \fn(\upsvn) - 2 \fn(\upsvs) 
    - \| \DVL^{-1} \Av \|^{2}
    & \leq &
    \frac{\dltwb}{1 - \dltwb} \| \DVL^{-1} \Av \|^{2} \, .
    \qquad
\label{3d3Af12DGttGa}
\end{EQ}
Also
\begin{EQ}[rcl]
    \| \DVL (\upsvn - \upsvs) - \DVL^{-1} \Av \|^{2}
    & \leq &
    \frac{3 \dltwb}{(1 - \dltwb)^{2}} \, \| \DVL^{-1} \Av \|^{2} \, ,
    \\
    \| \DVL (\upsvn - \upsvs) \|
    & \leq &
    \frac{1 + \sqrt{2 \dltwb}}{1 - \dltwb} \, \| \DVL^{-1} \Av \| \, .
\label{DGttGtsGDGm13rGa}
\end{EQ}
\end{proposition}

\begin{proof}
By \eqref{dtb3u1DG2d3GPg}, for any \( \upsv \in \CA(\rrn) \)
\begin{EQA}
	\Bigl| 
		\fs(\upsvs) - \fs(\upsv) - \frac{1}{2} \| \DVL (\upsv - \upsvs) \|^{2} 
	\Bigr|
	& \leq &
	\frac{\dltwb}{2} \| \DVL (\upsv - \upsvs) \|^{2} .
\label{d3GrGELGtsG12}
\end{EQA}
Further, 
\begin{EQA}[rcl]
	&& \nquad
	\fn(\upsv) - \fn(\upsvs) - \frac{1}{2} \| \DVL^{-1} \Av \|^{2}
	\\
	&=&
	\bigl\langle \upsv - \upsvs, \Av \bigr\rangle
	+ \fs(\upsv) - \fs(\upsvs) - \frac{1}{2} \| \DVL^{-1} \Av \|^{2} 
	\\
	&=&
	- \frac{1}{2} \bigl\| \DVL (\upsv - \upsvs) - \DVL^{-1} \Av \bigr\|^{2}
	+ \fs(\upsv) - \fs(\upsvs) + \frac{1}{2} \| \DVL (\upsv - \upsvs) \|^{2} .
	\qquad 
\label{12ELGuELusG}
\end{EQA}
As \( \upsvn \in \CA(\rrn) \) and it maximizes \( \fn(\upsv) \), we derive by \eqref{d3GrGELGtsG12}
\begin{EQA}
	&& \nquad
	\fn(\upsvn) - \fn(\upsvs) - \frac{1}{2} \| \DVL^{-1} \Av \|^{2}
	=
	\max_{\upsv \in \CA(\rrn)} 
	\Bigl\{ 
		\fn(\upsv) - \fn(\upsvs) - \frac{1}{2} \| \DVL^{-1} \Av \|^{2} 
	\Bigr\}
	\\
	& \leq &
	\max_{\upsv \in \CA(\rrn)} 
	\Bigl\{ 
		- \frac{1}{2} \bigl\| \DVL (\upsv - \upsvs) - \DVL^{-1} \Av \bigr\|^{2} 
		+ \frac{\dltwb}{2} \| \DVL (\upsv - \upsvs) \|^{2}
	\Bigr\} .
\label{d3G1212222B} 
\end{EQA}
Further, \( \max_{\uv} \bigl\{ \dltwb \| \uv \|^{2} - \| \uv - \xiv \|^{2} \bigr\} = \frac{\dltwb}{1 - \dltwb} \| \xiv \|^{2} \)
for \( \dltwb \in [0,1) \) and \( \xiv \in \R^{\dimp} \), 
yielding
\begin{EQA}
	\fn(\upsvn) - \fn(\upsvs) - \frac{1}{2} \| \DVL^{-1} \Av \|^{2}
	& \leq &
	\frac{\dltwb}{2(1 - \dltwb)} \| \DVL^{-1} \Av \|^{2} . 
\label{fd3G122B2} 
\end{EQA}
Similarly 
\begin{EQA}
	&& \nquad \quad
	\fn(\upsvn) - \fn(\upsvs) - \frac{1}{2} \| \DVL^{-1} \Av \|^{2}
	\geq 
	\max_{\upsv \in \CA(\rrn)} 
	\Bigl\{ 
		- \frac{1}{2} \bigl\| \DVL (\upsv - \upsvs) - \DVL^{-1} \Av \bigr\|^{2} 
		- \frac{\dltwb}{2} \| \DVL (\upsv - \upsvs) \|^{2}
	\Bigr\}
	\\
	&=&
	- \frac{\dltwb }{2(1 + \dltwb)} \, \| \DVL^{-1} \Av \|^{2} . 
	\qquad \quad
\label{fd3G122B2m} 
\end{EQA}
These bounds imply 
\eqref{3d3Af12DGttGa}.

Now we derive similarly to \eqref{12ELGuELusG} that for \( \upsv \in \CA(\rrn) \)
\begin{EQA}
	\fn(\upsv) - \fn(\upsvs) 
	& \leq &
	\bigl\langle \upsv - \upsvs, \Av \bigr\rangle
	- \frac{1 - \dltwb}{2} \| \DVL (\upsv - \upsvs) \|^{2} .
\label{LGvLGvsGf1d3G2}
\end{EQA}
A particular choice \( \upsv = \upsvn \) yields
\begin{EQA}
	\fn(\upsvn) - \fn(\upsvs) 
	& \leq &
	\bigl\langle \upsvn - \upsvs, \Av \bigr\rangle
	- \frac{1 - \dltwb}{2} \| \DVL (\upsvn - \upsvs) \|^{2} .
\label{21GsvtvGDG3G2}
\end{EQA}
Combining with \eqref{fd3G122B2m} allows to bound
\begin{EQA}
	&& \nquad
	\bigl\langle \upsvn - \upsvs, \Av \bigr\rangle
	- \frac{1 - \dltwb}{2} \| \DVL (\upsvn - \upsvs) \|^{2} 
	- \frac{1}{2} \| \DVL^{-1} \Av \|^{2}
	\geq 
	- \frac{\dltwb}{2(1 + \dltwb)} \| \DVL^{-1} \Av \|^{2} .
\label{2m1DGd3G123G}
\end{EQA}
Further, for \( \xiv = \DVL^{-1} \Av \), 
\( \uv = \DVL (\upsvn - \upsvs) \), and 
\( \dltwb \in [0,1/3] \), the inequality 
\begin{EQA}
	2 \bigl\langle \uv, \xiv \bigr\rangle - (1 - \dltwb) \| \uv \|^{2} - \| \xiv \|^{2}
	& \geq &
	- \frac{\dltwb}{1 + \dltwb} \| \xiv \|^{2} 
\label{dtxi2fd1d22}
\end{EQA}
implies 
\begin{EQA}
	\bigl\| \uv - \frac{1}{1-\dltwb} \xiv \bigr\|^{2}
	& \leq &
	\frac{2 \dltwb}{(1 + \dltwb) (1 - \dltwb)^{2}} \| \xiv \|^{2}
\label{uv11wxi22w1w}
\end{EQA}
yielding for \( \dltwb \leq 1/3 \)
\begin{EQA}
	\| \uv - \xiv \|
	& \leq &
	\biggl( \dltwb + \sqrt{\frac{2 \dltwb}{1 + \dltwb}} \biggr)
	\frac{\| \xiv \|}{1 - \dltwb}
	\leq 
	\frac{\sqrt{3 \dltwb} \| \xiv \|}{1 - \dltwb} \, ,
	\\
	\| \uv \|
	& \leq &
	\biggl( 1 + \sqrt{\frac{2 \dltwb}{1 + \dltwb}} \biggr)
	\frac{\| \xiv \|}{1 - \dltwb}
	\leq 
	\frac{1 + \sqrt{2 \dltwb} \| \xiv \|}{1 - \dltwb} \, ,
\label{uxiBws2w1w31w}
\end{EQA}
and \eqref{DGttGtsGDGm13rGa} follows.
\end{proof}

\Subsection{Quadratic penalization}
Here we discuss the case when \( \fn(\upsv) - \fs(\upsv) \) is quadratic.
The general case can be reduced to the situation with \( \fn(\upsv) = \fs(\upsv) - \| \GP \upsv \|^{2}/2 \).
To make the dependence of \( \GP \) more explicit, denote 
\begin{EQA}
	\fG(\upsv) 
	&=& 
	\fs(\upsv) - \| \GP \upsv \|^{2}/2 .
\label{8cvkfc9fujf6jnmcer4cd}
\end{EQA}
With \( \upsvs = \argmax_{\upsv} \fs(\upsv) \) and \( \upsvs_{\GP} = \argmax_{\upsv} \fG(\upsv) \),
we study the bias \( \upsvs_{\GP} - \upsvs \) induced by this penalization
under Fr\'echet-type smoothness conditions. 
To get some intuition, consider first the case of a quadratic function \( \fs(\upsv) \).

\begin{lemma}
\label{Lbiasquadgen}
Let \( \fs(\upsv) \) be quadratic with \( \IFL \equiv - \nabla^{2} \fs(\upsv) \) and \( \Av_{\GP} \equiv - \GP^{2} \upsvs \).
Then it holds with \( \IFL_{\GP} = \IFL + \GP^{2} \)
\begin{EQA}[rcccl]
	\upsvs_{\GP} - \upsvs
	&=&
	\IFL_{\GP}^{-1} \Av_{\GP} 
	&=&
	- \IFL_{\GP}^{-1} \GP^{2} \upsvs,
	\\
	\fG(\upsvs_{\GP}) - \fG(\upsvs)
	&=&
	\frac{1}{2} \| \IFL_{\GP}^{-1/2} \Av_{\GP} \|^{2}
	&=&
	\frac{1}{2} \| \IFL_{\GP}^{-1/2} \GP^{2} \upsvs \|^{2} \, .
\label{kjfydf554dfwertdgwdf}
\end{EQA}
\end{lemma} 

\begin{proof}
Quadraticity of \( \fs(\upsv) \) implies quadraticity of \( \fG(\upsv) \) with \( \nabla^{2} \fG(\upsv) \equiv - \IFL_{\GP} \).
This implies
\begin{EQA}
	\nabla \fG(\upsvs) - \nabla \fG(\upsvs_{\GP})
	&=&
	\IFL_{\GP} (\upsvs_{\GP} - \upsvs) .
\label{dcudydye67e6dy3wujhdqu}
\end{EQA}
Further, \( \nabla \fs(\upsvs) = 0 \) yielding \( \nabla \fG(\upsvs) = \Av_{\GP} = - \GP^{2} \upsvs \).
Together with \( \nabla \fG(\upsvs_{\GP}) = 0 \), this implies
\( \upsvs_{\GP} - \upsvs = \IFL_{\GP}^{-1} \Av_{\GP} \).
The Taylor expansion of \( \fG \) at \( \upsvs_{\GP} \) yields with \( \DVL_{\GP} = \IFL_{\GP}^{1/2} \) 
\begin{EQA}
	\fG(\upsvs) - \fG(\upsvs_{\GP})
	&=&
	- \frac{1}{2} \| \DVL_{\GP} (\upsvs_{\GP} - \upsvs) \|^{2}
	=
	- \frac{1}{2} \| \DVL_{\GP}^{-1} \Av_{\GP} \|^{2} 
\label{8chuctc44wckvcuedjequ}
\end{EQA}
and the assertion follows.
\end{proof}

Now we turn to the general case with \( \fs \) smooth in Fr\'echet sense.
Let \( \IFL(\upsv) = - \nabla^{2} \fs(\upsv) \).
For \( \upsv \in \Ups \) and \( \uv \in \R^{\dimp} \), define
\begin{EQA}
	\IFLba(\upsv;\uv)
	& \eqdef &
	\int_{0}^{1} \IFL(\upsv + t \uv) \, dt \, ;
\label{vhudfg7hgbfuhbugr}
\end{EQA}
cf. \eqref{vhudfg7sdfyt7s8hgbfuhbugrh9}.
Similarly define 
\( \IFL_{\GP}(\upsv) = - \nabla^{2} \fG(\upsv) = \IFL(\upsv) + \GP^{2} \) and
\begin{EQA}
	\IFLba_{\GP}(\upsv;\uv)
	& \eqdef &
	\int_{0}^{1} \IFL_{\GP}(\upsv + t \uv) \, dt
	=
	\IFLba(\upsv;\uv) + \GP^{2} .
\label{virjtfgy6f53e5fuhuyyuu}
\end{EQA}

\begin{lemma}
\label{LfreTgrad}
For the vector \( \biasv_{\GP} \eqdef \upsvs_{\GP} - \upsvs \), define 
\( \IFLba_{\GP} = \IFLba_{\GP}(\upsvs;\biasv_{\GP}) \).
Then
\begin{EQA}
	\upsvs_{\GP} - \upsvs
	&=&
	\IFLba_{\GP}^{-1} \, \Av_{\GP} \, .
\label{wivufhu338fvjed3jfufb}
\end{EQA}
\end{lemma} 

\begin{proof}
First we show for any \( \upsv \in \Ups \) and \( \uv \in \R^{\dimp} \) that
\begin{EQA}
	\nabla \fG(\upsv + \uv) - \nabla \fG(\upsv) 
	& = &
	- \IFLba_{\GP}(\upsv;\uv) \, \uv .
\label{dfhhjvqwbhewjmaqjmqjew}
\end{EQA}
Indeed, for any \( \gammav \in \R^{\dimp} \), consider the univariate function 
\( h(t) = \langle \nabla \fG(\upsv + t \uv) - \nabla \fG(\upsv), \gammav \rangle \).
Statement \eqref{dfhhjvqwbhewjmaqjmqjew} follows from definitions \eqref{vhudfg7hgbfuhbugr}, \eqref{virjtfgy6f53e5fuhuyyuu}, 
and the identity
\( h(1) - h(0) = \int_{0}^{1} h'(t) \, dt \).

Further, the definition \( \upsvs = \argmax_{\upsv} \fs(\upsv) \) yields \( \nabla \fs(\upsvs) = 0 \) and
\begin{EQA}
	\nabla \fG(\upsvs)
	&=& 
	\nabla f(\upsvs) - \GP^{2} \upsvs
	=
	\Av_{\GP} \, .
\label{dfse423q2esftcsdfgsygyb}
\end{EQA}
In view of \( \nabla \fG(\upsvs_{\GP}) = 0 \), we derive
\begin{EQA}
	\nabla \fG(\upsvs_{\GP}) - \nabla \fG(\upsvs)
	&=&
	- \Av_{\GP} \, .
\label{4hv78gfu7jh3ewo0ggt2lobhb}
\end{EQA}
Representation \eqref{dfhhjvqwbhewjmaqjmqjew} with \( \upsv = \upsvs \) and \( \uv = \biasv_{\GP} \) yields \eqref{wivufhu338fvjed3jfufb}.
\end{proof}

Representation \eqref{wivufhu338fvjed3jfufb} is very useful to bound from above the bias 
\( \upsvs_{\GP} - \upsvs \).
Indeed, assuming that \( \upsvs_{\GP} \) is in a local vicinity of \( \upsvs \) we may use 
Fr\'echet smoothness of \( \fs \) in terms of the value \( \dltwbss(\upsvs) \) from \eqref{jcxhydtferyufgy7tfdsy7ft} to approximate 
\( \IFLba_{\GP} \approx \IFL_{\GP}(\upsvs) \) and
\( \upsvs_{\GP} - \upsvs \approx - \IFL_{\GP}^{-1}(\upsvs) \, \GP^{2} \upsvs \).

\begin{proposition}
\label{Pbiasgeneric} 
Define \( \DVLf_{\GP} \) by
\begin{EQA}
	\DVLf_{\GP}^{2}
	& = &
	\IFL_{\GP}(\upsvs) ,
\label{95jioo00853rdvfyu7}
\end{EQA} 
and let \( \QP \) be a positive definite symmetric operator in \( \R^{\dimp} \) with \( \QP \leq \DVLf_{\GP} \).
Fix 
\begin{EQA}
	\biasQ_{\GP}
	&=& 
	\| \QP \, \DVLf_{\GP}^{-2} \, \GP^{2} \upsvs \| 
	=
	\| \QP \, \DVLf_{\GP}^{-2} \, \Av_{\GP} \| \, .
\label{fd9dfhy4ye6fuydfrerf}
\end{EQA} 
With \( \amax = 2/3 \), assume
\begin{EQA}
	\dltwbs_{\GP}
	\eqdef
	\sup_{\uv \colon \| \QP \uv \| \leq \amax^{-1} \biasQ_{\GP}} \,\, 
	\| \DVLf_{\GP}^{-1} \, \IFL_{\GP}(\upsvs + \uv) \, \DVLf_{\GP}^{-1} - \Id_{\dimp} \|
	& \leq &
	\frac{1}{3} \, ;
\label{ghdrd324ee4ew222gen}
\end{EQA}
cf. \eqref{ghdrd324ee4ew222w3ew}. 
Then \( \| \QP (\upsvs_{\GP} - \upsvs) \| \leq \amax^{-1} \biasQ_{\GP} \) or, equivalently,
\begin{EQA}
	\upsvs_{\GP}
	 & \in &
	 \CA_{\GP}
	 \eqdef
	 \{ \upsv \colon \| \QP (\upsv - \upsvs) \| \leq \amax^{-1} \biasQ_{\GP} \} .
\label{odf6fdyr6e4deuewjug}
\end{EQA}
Moreover, 
\begin{EQA}
	\| \QP (\upsvs_{\GP} - \upsvs - \DVLf_{\GP}^{-2} \Av_{\GP}) \|
	& \leq &
	\dltwbs_{\GP} \, \amax^{-1} \biasQ_{\GP} \, .
\label{11ma3eaelDebbgen}
\end{EQA}
\end{proposition}

\begin{proof}
First we check that \( \upsvs_{\GP} \) concentrates in the local vicinity 
\( \CA_{\GP} \) from \eqref{odf6fdyr6e4deuewjug}.
Strong concavity of \( \fG \) implies that the solution \( \upsvs_{\GP} \) exists and unique.
Let us fix any vector \( \uv \in \R^{\dimp} \) with \( \| \QP \, \uv \| = \amax^{-1} \biasQ_{\GP} \).
Due to \eqref{wivufhu338fvjed3jfufb}, 
we are looking at the solution \( \IFLba_{\GP}(\upsvs; s \uv) \, s \uv = \Av_{\GP} \) in \( s \).
It suffices to ensure that 
\begin{EQA}
	s \, \QP \, \uv 
	& = &
	\QP \, \IFLba_{\GP}(\upsvs; s \uv)^{-1} \Av_{\GP} 
\label{ihy666677erdft4e3te4yt}
\end{EQA}
is impossible for \( s \geq 1 \).
For \( s = 1 \), we can use that \( \IFLba_{\GP}(\upsvs; \uv) \geq (1 - \dltwbs_{\GP}) \DVLf_{\GP}^{2} \);
see \eqref{5chfdc7e3yvc5ededww} of Lemma~\ref{LfreTayb}.
Therefore,
\begin{EQA}
	\| \QP \, \IFLba_{\GP}^{-1}(\upsvs; \uv) \Av_{\GP} \|
	& \leq &
	\| \QP \, \IFLba_{\GP}^{-1}(\upsvs; \uv) \DVLf_{\GP}^{2} \QP^{-1} \| \, \| \QP \, \DVLf_{\GP}^{-2} \, \Av_{\GP} \|
	\\
	& \leq &
	(1 - \dltwbs_{\GP})^{-1} \biasQ_{\GP}
	<
	\amax^{-1} \biasQ_{\GP}
\label{v565tf5tfd56fd65tgtsydx}
\end{EQA}
and \eqref{ihy666677erdft4e3te4yt} with \( s=1 \) is impossible because \( \| \QP \uv \| = \amax^{-1} \biasQ_{\GP} \).
It remains to note that the matrix \( s \, \IFLba_{\GP}(\upsvs; s \uv) \) grows with \( s \) as
\begin{EQA}
	s \, \IFLba_{\GP}(\upsvs; s \uv)
	&=&
	s \int_{0}^{1} \IFL_{\GP}(\upsvs + t \, s \, \uv) \, dt
	=
	\int_{0}^{s} \IFL_{\GP}(\upsvs + t \, \uv) \, dt .
\label{r0ofigvjbh3wd9eorikcjv}
\end{EQA}
Now we bound \( \| \QP \, \biasv_{\GP} \| \) assuming that \( \| \QP (\upsvs_{\GP} - \upsvs) \| \leq \amax^{-1} \biasQ_{\GP} \)
and \eqref{ghdrd324ee4ew222gen} applies.  
Statement \eqref{5chfdc7e3yvc5ededww} of Lemma~\ref{LfreTayb} implies 
\( \IFLba_{\GP}^{-1} \leq (1 - \dltwbs_{\GP})^{-1} \DVLf_{\GP}^{2} \) and
\begin{EQA}
	\| \QP \, \biasv_{\GP} \|
	&=&
	\| \QP \, \IFLba_{\GP}^{-1} \, \DVLf_{\GP}^{2} \, \QP^{-1} \, \QP \, \DVLf_{\GP}^{-2} \, \Av_{\GP} \|
	\leq 
	\| \QP \, \IFLba_{\GP}^{-1} \, \DVLf_{\GP}^{2} \, \QP^{-1} \| \, \| \QP \, \DVLf_{\GP}^{-2} \Av_{\GP} \|
	\\
	& \leq &
	\frac{1}{1 - \dltwbs_{\GP}} \, \| \QP \, \DVLf_{\GP}^{-2} \Av_{\GP} \| 
	=
	\frac{\biasQ_{\GP}}{1 - \dltwbs_{\GP}} \, .
\label{jhfdt3dtyiu89jgrert6b}
\end{EQA}
In the same way we derive
\begin{EQA}
	\| \QP \, (\biasv_{\GP} - \DVLf_{\GP}^{-2} \Av_{\GP}) \|
	&=&
	\bigl\| \QP \, (\IFLba_{\GP}^{-1} - \DVLf_{\GP}^{-2}) \Av_{\GP} \bigr\|
	\leq 
	\frac{\dltwbs_{\GP}}{1 - \dltwbs_{\GP}} \, \| \QP \, \DVLf_{\GP}^{-2} \Av_{\GP} \| \, ,
\label{ae3DeuseDeutb}
\end{EQA}
and \eqref{11ma3eaelDebbgen} follows as well. 
\end{proof}

\begin{remark}
\label{Rbiasgeneric}
Inspection of the proofs of Proposition~\ref{Pbiasgeneric} indicates that 
the results \eqref{odf6fdyr6e4deuewjug} through \eqref{11ma3eaelDebbgen} can be restated 
with \( \DVL_{\GP}^{2} = \IFL_{\GP}(\upsvs_{\GP}) \) in place of \( \DVLf_{\GP}^{2} = \IFL_{\GP}(\upsvs) \).
\end{remark}


\newcommand{\binomb}[2]{\genfrac{}{}{0pt}{}{#1}{#2}}
\def\Bv{\bb{B}}

\Section{Conditional and marginal optimization}
\label{Ssemiopt}
This section describes the problem of sequential and marginal optimization. 
Consider a function \( \lgd(\prmtv) \) of a parameter \( \prmtv \in \R^{\dimtotal} \)
which can be represented as \( \prmtv = (\targv,\nuiv) \), where \( \targv \in \R^{\dimp} \)
is the target subvector while \( \nuiv \in \R^{\dimq} \) is a nuisance variable.

\Subsection{Partial optimization and local partial smoothness}
For any fixed value of the nuisance variable \( \nuiv \in \Nui \), consider 
\( \lgd_{\nuiv}(\targv) = \lgd(\targv,\nuiv) \) as a function of \( \targv \) only.
Below we assume that \( \lgd_{\nuiv}(\targv) \) is concave in \( \targv \) for any \( \nuiv \in \Nui \).
Define 
\begin{EQA}[rcl]
	\targv_{\nuiv}
	& \eqdef &
	\argmax_{\targv} \lgd_{\nuiv}(\targv),
\label{jfoiuy2wedfv7tr2qsdzxdfso}
	\\
	\feta_{\nuiv}
	& \eqdef &
	\max_{\targv} \lgd_{\nuiv}(\targv) - \lgd(\prmtvs)
	=
	\lgd_{\nuiv}(\targv_{\nuiv}) - \lgd(\prmtvs) .
\label{kfdy8d3fgfnm0lkgrr3so}
\end{EQA}
Also define
\begin{EQA}
	\IFL_{\nuiv}
	& \eqdef &
	- \nabla^{2} \lgd_{\nuiv}(\targv_{\nuiv}) 
	=
	- \nabla_{\targv\targv}^{2} \lgd(\targv_{\nuiv},\nuiv) .
\label{ge8qwefygw3qytfyju8qfhbdso}
\end{EQA}

Local smoothness of each function \( \lgd_{\nuiv}(\cdot) \) around \( \targv_{\nuiv} \) 
can be well described under the self-concordance property.
Let some radius \( \rr \) be fixed. 
In general it may depend on \( \nuiv \) or on the effective dimension \( \dimA_{\nuiv} \)
for \( \nuiv \in \Nui \).
We also assume that for each \( \nuiv \in \Nui \), a local metric on \( \R^{\dimp} \)
for the target variable \( \targv \) is defined by a matrix \( \HL_{\nuiv} \in \Matr_{\dimp} \).

\begin{description}
    \item[\label{LLseS3ref} \( \bb{(\mathcal{S}_{3 | \nuiv})} \)]
      \emph{For any \( \nuiv \in \Nui \), it holds 
      \( \lgd_{\nuiv}(\targv) = - n \hL_{\nuiv}(\targv) \) with a strongly convex function 
      \( \hL_{\nuiv}(\targv) \) such that 
      \( \HL_{\nuiv}^{2} \leq \nabla^{2} \hL_{\nuiv}(\targv_{\nuiv}) \), and
\begin{EQA}
	\sup_{\| \HL_{\nuiv} \uv \| \leq \rr/\sqrt{n}} \, \sup_{t \in [0,1]} 
	\frac{\bigl| \langle \nabla^{3} \hL_{\nuiv}(\targv_{\nuiv} + t \uv), \uv^{\otimes 3} \rangle \bigr|}
		 {\| \HL_{\nuiv} \uv \|^{3}}
	& \leq &
	\hmax_{3} \, .
\label{hvgtgdfcyt4e3f6tfterso}
\end{EQA}
}
    \item[\label{LLseS4ref} \( \bb{(\mathcal{S}_{4 |\nuiv})} \)]
      \emph{For any \( \nuiv \in \Nui \), the function \( h_{\nuiv}(\cdot) \) satisfies \nameref{LLseS3ref} and 
\begin{EQA}
	\sup_{\| \HL_{\nuiv} \uv \| \leq \rr/\sqrt{n}} \, \sup_{t \in [0,1]} 
	\frac{\bigl| \langle \nabla^{4} \hL_{\nuiv}(\targv_{\nuiv} + t \uv), \uv^{\otimes 4} \rangle \bigr|}
		 {\| \HL_{\nuiv} \uv \|^{4}}
	& \leq &
	\hmax_{4} \, .
\end{EQA}
}
\end{description}

Under \nameref{LLseS3ref} it holds for all \( \uv \) with \( \| \HL_{\nuiv}(\targv) \uv \| \leq \rr/\sqrt{n} \)
for the corresponding errors of the Taylor expansion
\begin{EQA}
	\dltw_{3,\nuiv}(\uv) 
	& \eqdef & 
	\lgd_{\nuiv}(\targv_{\nuiv} + \uv) - \lgd_{\nuiv}(\targv_{\nuiv}) - \bigl\langle \lgd_{\nuiv}(\targv_{\nuiv}), \uv \bigr\rangle + \| \IFL_{\nuiv}^{1/2} \uv \|^{2}/2 
	\\ 
	& \leq &
	\frac{1}{6} \, \hmax_{3} \, n \| \HL_{\nuiv} \uv \|^{3}
	\leq 
	\frac{\hmax_{3} \, \rr}{6 n^{1/2}} \| \IFL_{\nuiv}^{1/2} \uv \|^{2} ,
\label{6ghd56eg3nkcxtwhwkthhj}
\end{EQA}
and
\begin{EQA}
	\dltw_{4,\nuiv}(\uv) 
	& \eqdef & 
	\lgd_{\nuiv}(\targv) - \lgd_{\nuiv}(\targv_{\nuiv}) - \bigl\langle \lgd_{\nuiv}(\targv_{\nuiv}), \uv \bigr\rangle + \| \IFL_{\nuiv}^{1/2} \uv \|^{2}/2 
	- \bigl\langle \nabla^{3} \lgd_{\nuiv}(\targv_{\nuiv}), \uv^{\otimes 3} \bigr\rangle/6	
	\\
	& \leq & 
	\frac{1}{24} \, \hmax_{4} \, n \| \HL_{\nuiv} \uv \|^{4}
	\leq 
	\frac{\hmax_{4} \, \rr^{2}}{24 n} \| \IFL_{\nuiv}^{1/2} \uv \|^{2} ;
\label{4cfg9hh0rth47re7fd6ehj}
\end{EQA}
see Lemma~\ref{LdltwLaGP}.

Condition \nameref{LLsS3ref} can also be extended to the \( \nuiv \)-conditional framework.

\begin{description}

    \item[\label{LLsoS3ref} \( \bb{(\mathcal{S}_{3|\nuiv}^{+})} \)]
      \emph{\( f_{\nuiv}(\targv) = - n \hL_{\nuiv}(\targv) \) with \( \hL_{\nuiv}(\targv) \) strongly concave. 
      \( \HL_{\nuiv}^{2} \leq \nabla^{2} \hL_{\nuiv}(\targv_{\nuiv}) \), 
      and
\begin{EQA}
	\sup_{\| \HL_{\nuiv} \uv \| \leq \rr/\sqrt{n}}  \,\, \sup_{\gammav \in \R^{\dimp}} \,\,
	\frac{\bigl| \langle \nabla^{3} \hL_{\nuiv}(\targv_{\nuiv} + \uv), \uv \otimes \gammav^{\otimes 2} \rangle \bigr|}
		 {\| \HL_{\nuiv} \uv \|\, \| \HL_{\nuiv} \gammav \|^{2}}
	& \leq &
	\hmax_{3} \, .
\label{d6f53ye5vry4fddfgeyd}
\end{EQA}
}
\end{description}

To pull together the results of partial optimization w.r.t. \( \targv \) conditioned on \( \nuiv \),
we also need a condition on the cross derivative of \( \lgd(\targv,\nuiv) \).

\begin{description}
    \item[\label{LLpS3ref} \( \bb{(\mathcal{S}_{3,\nuiv}^{+})} \)]
      \emph{It holds \( \lgd(\prmtv) = - n \hL(\prmtv) \) with a strongly convex function 
      \( \hL(\prmtv) = \hL(\targv,\nuiv) \) satisfying \nameref{LLseS3ref} for each \( \nuiv \in \Nui \).
      Moreover,  it holds 
}
\begin{EQA}
	\sup_{\gammav \in \R^{\dimp}} \,\, 
	\sup_{\| \HL \, \wv \| \leq \rru/\sqrt{n}} \,\, 
	\sup_{t \in [0,1]} 
	\frac{| \langle \nabla_{\targv} \nabla^{2} \hL(\prmtvs + t \wv), \gammav \otimes \wv^{\otimes 2} \rangle |}
		 {\| \HL_{\nuivs} \gammav \| \, \| \HL \, \wv \|^{2}}
	& \leq &
	\hmax_{3} \, ,
	\\
	\sup_{\gammav \in \R^{\dimp}} \,\, 
	\sup_{\| \HL \, \wv \| \leq \rru/\sqrt{n}} \,\, 
	\sup_{t \in [0,1]} 
	\frac{| \langle \nabla_{\targv\targv}^{2} \nabla \hL(\prmtvs + t \wv), \gammav^{\otimes 2} \otimes \wv \rangle |}
		 {\| \HL_{\nuivs} \gammav \|^{2} \, \| \HL \, \wv \|}
	& \leq &
	\hmax_{3} \, .
\label{c6ceyecc5e5etctwhcyegwc}
\end{EQA}
\end{description}

This condition bounds the third order cross derivative in \( \targv \) and \( \nuiv \).
Of course, it suffices to bound the third full derivative of \( \hL(\prmtv) \); see \nameref{LLsS3ref}.

\Subsection{Conditional optimization and a bound on the bias}
\label{SortLaplsemi}
Here we study variability of the value \( \targv_{\nuiv} = \argmax_{\targv} \lgd(\targv,\nuiv) \)
w.r.t. the nuisance para\-meter \( \nuiv \).
It appears that local quadratic approximation of the function \( \lgd \) in a vicinity of the extreme point
\( \prmtvs \) yields a nearly linear dependence of \( \targv_{\nuiv} \) on \( \nuiv \).
We illustrate this fact on the case of a quadratic function \( \lgd(\cdot) \).
Consider the negative Hessian \( \IFT = - \nabla^{2} \lgd(\prmtvs) \) in the block form:
\begin{EQA}
	\IFT
	& \eqdef &
	- \nabla^{2} \lgd(\prmtvs)
	=
	\begin{pmatrix}
	\IFT_{\targv\targv} & \IFT_{\targv\nuiv}
	\\
	\IFT_{\nuiv\targv} & \IFT_{\nuiv\nuiv}
	\end{pmatrix} 
\label{hwe78yf2diwe76tfw67etfwtbso}
\end{EQA}
with \( \IFT_{\nuiv\targv} = \IFT_{\targv\nuiv}^{\T} \).
If \( \lgd(\prmtv) \) is quadratic then \( \IFT \) and its blocks do not depend on \( \prmtv \).

\begin{lemma}
\label{Lpartmaxq}
Let \( \lgd(\prmtv) \) be {quadratic}, {strongly concave}, and \( \nabla \lgd(\prmtvs) = 0 \).
Then
\begin{EQA}
	\targv_{\nuiv} - \targvs 
	&=& 
	- \IFT_{\targv\targv}^{-1} \, \IFT_{\targv\nuiv} \, (\nuiv - \nuivs) .
\label{0fje7fhihy84efiewkw}
\end{EQA}
\end{lemma}

\begin{proof}
The condition \( \nabla \lgd(\prmtvs) = 0 \) implies 
\( \lgd(\prmtv) = \lgd(\prmtvs) - (\prmtv - \prmtvs)^{\T} \IFT \, (\prmtv - \prmtvs)/2 \) 
with \( \IFT = - \nabla^{2} \lgd(\prmtvs) \).
For \( \nuiv \) fixed, the point \( \targv_{\nuiv} \) maximizes 
\( - (\targv - \targvs)^{\T} \IFT_{\targv\targv} \, (\targv - \targvs) /2 - (\targv - \targvs)^{\T} \IFT_{\targv\nuiv} \, (\nuiv - \nuivs) \)
and thus, \( \targv_{\nuiv} - \targvs = - \IFT_{\targv\targv}^{-1} \, \IFT_{\targv\nuiv} \, (\nuiv - \nuivs) \).
\end{proof} 

This observation \eqref{0fje7fhihy84efiewkw} is in fact discouraging 
because the bias \( \targv_{\nuiv} - \targvs \) has the same magnitude as the nuisance parameter \( \nuiv - \nuivs \).
However, the condition \( \IFT_{\targv\nuiv} = 0 \) yields \( \targv_{\nuiv} \equiv \targvs \) 
and the bias vanishes.
If \( \lgd(\prmtv) \) is not quadratic, the \emph{orthogonality} condition \( \nabla_{\nuiv} \nabla_{\targv} \, \lgd(\targv,\nuiv) \equiv 0 \)
for all \( (\targv,\nuiv) \in \WV \) still ensures a vanishing bias.

\begin{lemma}
\label{Lpartmax}
Let \( \lgd(\targv,\nuiv) \) be a continuously differentiable and 
\( \nabla_{\nuiv} \nabla_{\targv} \, \lgd(\targv,\nuiv) \equiv 0 \).
Then the point 
\( \targv_{\nuiv} = \argmax_{\targv} f(\targv,\nuiv) \) does not depend on \( \nuiv \).
\end{lemma}

\begin{proof}
The condition \( \nabla_{\nuiv} \nabla_{\targv} \, \lgd(\targv,\nuiv) \equiv 0 \) implies the decomposition 
\( \lgd(\targv,\nuiv) = \lgd_{1}(\targv) + \lgd_{2}(\nuiv) \) for some functions 
\( \lgd_{1} \) and \( \lgd_{2} \).
This in turn yields \( \targv_{\nuiv} \equiv \targvs \).
\end{proof}

As a corollary, the maximizer \( \targv_{\nuiv} \) and the corresponding negative Hessian \( \IFL_{\nuiv} = \DV_{\nuiv}^{2} \)
do not depend on \( \nuiv \).
Unfortunately, the orthogonality condition \( \nabla_{\nuiv} \nabla_{\targv} \, \lgd(\targv,\nuiv) \equiv 0 \) is too restrictive
and fulfilled only in the special additive case \( \lgd(\targv,\nuiv) = \lgd_{1}(\targv) + \lgd_{2}(\nuiv) \).

In some cases, one can check \emph{semi-orthogonality} condition 
\begin{EQA}
	\nabla_{\nuiv} \nabla_{\targv} \, \lgd(\targvs,\nuiv)  
	&=&
	0,
	\qquad
	\forall \nuiv \in \Nui \, .
\label{jdgyefe74erjscgygfydhwse}
\end{EQA}
\ifapp{A typical example is given by nonlinear regression; see Section~\ref{SnonlinLA}.}{}

\begin{lemma}
\label{Lsemiorto}
Assume \eqref{jdgyefe74erjscgygfydhwse}.
Then 
\begin{EQA}[rcccl]
	\nabla_{\targv} \, \lgd(\targvs,\nuiv)
	& \equiv &
	0,
	\qquad
	\nabla_{\targv\targv}^{2} \lgd(\targvs,\nuiv) 
	& \equiv &
	\nabla_{\targv\targv}^{2} \lgd(\targvs,\nuivs),
	\qquad
	\nuiv \in \Nui .
\label{0xyc5ftr4j43vefvuruendt}
\end{EQA}
Moreover, if \( \lgd(\targv,\nuiv) \) is concave in \( \targv \) given \( \nuiv \) then
\begin{EQA}
	\targv_{\nuiv} 
	& \eqdef & 
	\argmax_{\targv} \lgd(\targv,\nuiv) \equiv \targvs ,
	\qquad
	\forall \nuiv \in \Nui \,.
\label{irdtyrnjjutu45r7gjhr}
\end{EQA}
\end{lemma}

\begin{proof}
Consider the vector
\begin{EQA}
	\Av_{\nuiv} 
	& \eqdef & 
	\nabla_{\targv} \, \lgd(\targvs,\nuiv) \, .
\label{02h874shki7rwhuyfv}
\end{EQA}
Obviously \( \Av_{\nuivs} = 0 \).
Moreover, \eqref{jdgyefe74erjscgygfydhwse} implies that \( \Av_{\nuiv} \) 
does not depend on \( \nuiv \) and thus, vanishes everywhere.
As \( \lgd \) is concave in \( \targv \), 
this implies \( \lgd(\targvs,\nuiv) = \max_{\targv} \lgd(\targv,\nuiv) \)
and \( \targv_{\nuiv} = \argmax_{\targv} \lgd(\targv,\nuiv) \equiv \targvs \).
Similarly, \eqref{jdgyefe74erjscgygfydhwse} implies 
\( \nabla_{\nuiv} \nabla_{\targv\targv} \, \lgd(\targvs,\nuiv) \equiv 0 \) 
and \eqref{0xyc5ftr4j43vefvuruendt} follows.
Concavity of \( \lgd(\targv,\nuiv) \) in \( \targv \) for \( \nuiv \) fixed and \eqref{0xyc5ftr4j43vefvuruendt} 
imply \eqref{irdtyrnjjutu45r7gjhr}.
\end{proof}


Unfortunately, semi-orthogonality \eqref{jdgyefe74erjscgygfydhwse} condition is also rather restrictive and 
fulfilled only in special situations.
A weaker condition of 
\emph{one-point orthogonality} means \( \nabla_{\nuiv} \nabla_{\targv} \, \lgd(\targvs,\nuivs) = 0 \).
This condition is not restrictive and can always be enforced by a linear transform of the nuisance variable \( \nuiv \);
see Section~\ref{Sonepointtrans}.
Here we evaluate variability of \( \IFL_{\nuiv} = - \nabla_{\targv\targv}^{2} \lgd(\targvs,\nuiv) \)
and of the norm of \( \biasv_{\nuiv} = \targv_{\nuiv} - \targvs \) 
under this condition.
Similarly to Lemma~\ref{LfreTay}, we derive the following bound.

\begin{lemma}
\label{LvarDVa}
Assume \eqref{c6ceyecc5e5etctwhcyegwc} from \nameref{LLpS3ref} and let 
\begin{EQA}
	\dltwbss
	&=&
	\hmax_{3} \, \rru n^{-1/2} 
	\leq 
	1/3 .
\label{8uiiihkkcrdrdteggcjcjd}
\end{EQA}
Then for any \( \prmtv = (\targvs,\nuiv) \) with \( \| \HL (\prmtv - \prmtvs) \| \leq \rru/\sqrt{n} \), it holds
\begin{EQA}
	(1 - \dltwbss) \IFL_{\nuivs}
	& \leq &
	\IFL_{\nuiv}
	\leq 
	(1 + \dltwbss) \IFL_{\nuivs} .
\label{8c7c67cc6c63kdldlvvudw}
\end{EQA}
\end{lemma}

\begin{proposition}
\label{LsemiAv}
Assume one-point orthogonality \( \nabla_{\nuiv} \nabla_{\targv} \, \lgd(\targvs,\nuivs) = 0 \). 
Let also conditions \nameref{LLsoS3ref}, \nameref{LLpS3ref}
hold and the involved radii \( \rr \) and \( \rru \) satisfy
\begin{EQA}[c]
	\hmax_{3} \, \rru n^{-1/2} 
	\leq 1/3, 
	\qquad 
	\rr \geq \hmax_{3} \, n^{-1/2} \, \rru^{2} \, .
\label{7jk98k9kjdfyenfc76e4k}
\end{EQA}
Then for any \( \prmtv = (\targvs,\nuiv) \) with \( \| \HL (\prmtv - \prmtvs) \| \leq \rru/\sqrt{n} \) and any linear mapping 
\( \QP \colon \R^{\dimp} \to \R^{\dimq} \), the bias \( \biasv_{\nuiv} = \targv_{\nuiv} - \targvs \) satisfies
with \( \IFL = \IFL_{\nuivs} \)
\begin{EQA}
	\| \QP \, \biasv_{\nuiv} \|
	& \leq &
	\frac{\hmax_{3} \, \rru^{2}}{\sqrt{n}} \, \| \QP \, \IFL^{-1} \QP^{\T} \|^{1/2} \, .
\label{36gfijh94ejdvtwekoi}
\end{EQA}
In particular, with \( \QP = \IFL^{1/2} \)
\begin{EQA}
	\| \IFL^{1/2} \, \biasv_{\nuiv} \|
	& \leq &
	\frac{\hmax_{3} \, \rru^{2}}{\sqrt{n}} \, .
\label{36gfijh94ejdvtwekoD}
\end{EQA}
\end{proposition}

\begin{remark}
Condition \nameref{LLsS3ref} with \( \rr = \rru \) satisfying \( \hmax_{3} \, \rru n^{-1/2} \leq 1/3 \)
obviously implies \nameref{LLsoS3ref}, \nameref{LLpS3ref} with the same \( \rr = \rru \)
and \eqref{7jk98k9kjdfyenfc76e4k} is fulfilled as well.
\end{remark}

\begin{proof}
Fix \( \prmtv = (\targvs,\nuiv) \) with \( \| \HL (\prmtv - \prmtvs) \| \leq \rru/\sqrt{n} \).
Define the vector \( \Av_{\nuiv} \in \R^{\dimp} \) by
\begin{EQA}
	\Av_{\nuiv}
	& \eqdef &
	\nabla \lgd_{\nuiv}(\targvs) 
	=
	\nabla_{\targv} \, \lgd(\targvs,\nuiv) .
\label{8g5fgfg2frdg3ef6fh}
\end{EQA}
%
Define \( \wv = \prmtv - \prmtvs \), \( \etav = \nuiv - \nuivs \).
For bounding the vector \( \Av_{\nuiv} \), use
\( \Av_{\nuivs} = 0 \) and one-point orthogonality \( \nabla_{\nuiv} \Av_{\nuivs} = \nabla_{\nuiv} \nabla_{\targv} \, \lgd(\prmtvs) = 0 \).
This implies for \( \Bv_{\nuiv}(t) = \Av_{\nuivs + t \etav} \)
\begin{EQA}
	\Av_{\nuiv} 
	&=&
	\Bv_{\nuiv}(1) - \Bv_{\nuiv}(0) - \Bv'_{\nuiv}(0)
	=
	\int_{0}^{1} (\Bv'_{\nuiv}(t) - \Bv'_{\nuiv}(0)) \, dt
	=
	\int_{0}^{1} (1 - s) \, \Bv''_{\nuiv}(s) \, ds \, ,
\label{g8eh35fgtg7t76868jdy}
\end{EQA}
where \( \Bv'_{\nuiv}(t) = \frac{d}{dt} \Bv_{\nuiv}(t) \),
\( \Bv''_{\nuiv}(t) = \frac{d^{2}}{dt^{2}} \Bv_{\nuiv}(t) \).
Condition \nameref{LLpS3ref} ensures
\begin{EQA}
	\bigl| \langle \Bv''_{\nuiv}(s), \gammav \rangle \bigr|
	&=&
	\bigl| \bigl\langle \nabla_{\targv\nuiv\nuiv}^{3}\hL(\targvs,\nuivs + t \etav), \gammav \otimes \etav^{\otimes 2} \bigr\rangle \bigr|
	\\
	& \leq &
	\hmax_{3} \, n \, \| \HL_{\nuivs} \gammav \| \, \| \HL \, \wv \|^{2}
	\leq 
	\hmax_{3} \, n^{1/2} \, \| \IFL^{1/2} \gammav \| \, \| \HL \, \wv  \|^{2} .
\label{5cbgcfyt6webhwefuy6}
\end{EQA}
Now we use that 
\begin{EQA}
	\| \IFL^{-1/2} \Bv''_{\nuiv}(s) \|
	&=&
	\sup_{\gammav \colon \| \IFL^{1/2} \gammav \| \leq 1} \bigl| \langle \IFL^{-1/2} \Bv''_{\nuiv}(s), \IFL^{1/2} \gammav \rangle \bigr|
	=
	\sup_{\gammav \colon \| \IFL^{1/2} \gammav \| \leq 1}  \bigl| \langle \Bv''_{\nuiv}(s), \gammav \rangle \bigr|
	\\
	& \leq &
	\hmax_{3} \, n^{1/2} \, \sup_{\gammav \colon \| \IFL^{1/2} \gammav \| \leq 1}
	\| \IFL^{1/2} \gammav \| \, \| \HL \, \wv  \|^{2}
	\leq 
	\hmax_{3} \, n^{1/2} \, \| \HL \, \wv \|^{2} \, 
\label{hgvytew3hw36fdfnhrg}
\end{EQA}
yielding
\begin{EQA}
	\| \IFL^{-1/2} \Av_{\nuiv} \|
	& \leq &
	\hmax_{3} \, n^{1/2} \, \| \HL \, \wv  \|^{2} \int_{0}^{1} (1 - s) \, ds
	\leq 
	\frac{\hmax_{3} \, \rru^{2}}{2 \sqrt{n}} \, .
\label{4cbijjm09j9hhrjfifkr}
\end{EQA}
Now we bound the norm of \( \biasv_{\nuiv} \).
We use that \( \hmax_{3} \, n^{-1/2} \, \rru^{2} \leq \rr \) and \eqref{d6f53ye5vry4fddfgeyd}
of \nameref{LLsoS3ref} holds for this \( \rr \).
As \( \nabla \lgd_{\nuiv}(\targv_{\nuiv}) = 0 \), we derive with 
\( \biasv_{\nuiv} = \targv_{\nuiv} - \targvs \)
\begin{EQA}
	\Av_{\nuiv}
	&=&
	\nabla \lgd_{\nuiv}(\targvs) - \nabla \lgd_{\nuiv}(\targv_{\nuiv})
	=
	- \Bigl( \int_{0}^{1} \nabla^{2} \lgd_{\nuiv}(\targvs + t \biasv_{\nuiv}) \, dt \Bigr) \, \biasv_{\nuiv} 
	=
	\IFLba_{\nuiv} \, \biasv_{\nuiv} \, ;
	\qquad
\label{ub7eywhdtwjeyfdqxq}
\end{EQA}
see Lemma~\ref{LfreTayc}.
This yields \( \biasv_{\nuiv} = \IFLba_{\nuiv}^{-1} \Av_{\nuiv} \).
Moreover, for \( \rr_{\nuiv} = \hmax_{3} \, n^{-1/2} \, \rru^{2} \), it holds 
\( \dltwbss_{\nuiv} \eqdef \hmax_{3} \, n^{-1/2} \, \rr_{\nuiv} \leq {\dltwbss_{3}}^{2} \) for 
\( \dltwbss \eqdef \hmax_{3} \, \rru n^{-1/2} \)
implying \( (1 - \dltwbss_{\nuiv}) (1 - \dltwbss) \geq 1/2 \) and
by \eqref{5chfdc7e3yvc5ededww} of Lemma~\ref{LfreTayb} and \eqref{8c7c67cc6c63kdldlvvudw}
\begin{EQA}
	\| \QP \, \biasv_{\nuiv} \|
	& = &
	\| \QP \, \IFLba_{\nuiv}^{-1} \Av_{\nuiv} \|
	\leq 
	\frac{1}{1 - \dltwbss_{\nuiv}} \, \| \QP \, \IFL^{-1}(\nuiv) \Av_{\nuiv} \|
	\leq 
	\frac{1}{(1 - \dltwbss_{\nuiv}) (1 - \dltwbss)} \, \| \QP \, \IFL^{-1} \Av_{\nuiv} \|
	\\
	& \leq &
	\frac{1}{(1 - \dltwbss_{\nuiv}) (1 - \dltwbss)} \, 
	\| \QP \, \IFL^{-1} \QP^{\T} \|^{1/2} \, \| \IFL^{-1/2} \Av_{\nuiv} \|
	\leq 
	\frac{\hmax_{3} \, \rru^{2}}{\sqrt{n}} \, \| \QP \, \IFL^{-1} \QP^{\T} \|^{1/2} \, .
\label{36gfijh94ejdvtweko12}
\end{EQA}
This yields the assertion.
\end{proof}

\Subsection{One-point orthogonality by a linear transform}
\label{Sonepointtrans}
\emph{One-point orthogonality} condition means \( \nabla_{\nuiv} \nabla_{\targv} \, \lgd(\prmtvs) = 0 \).
This section explains how
this condition can be enforced for a general function \( \lgd \) by 
a linear transform of the nuisance parameter.
We only need a mild \emph{separability condition} which 
assumes that the full dimensional information matrix
\( \IFT(\prmtv) = - \nabla^{2} \lgd(\prmtv) \) is well posed
for \( \prmtv = \prmtvs \).

\begin{description}
	\item[\label{IFLref} \( \bb{(\IFT)} \)]
	\emph{
	The matrix \( \IFT = \IFT(\upsvs) \) is positive definite 
	and for some constant \( \rhoIF = \rhoIF(\IFT) < 1 \)} 
\begin{EQA}
	\| 
	\IFT_{\targv\targv}^{-1/2} \IFT_{\targv\nuiv} \, \IFT_{\nuiv\nuiv}^{-1} \, \IFT_{\nuiv\targv} \, \IFT_{\targv\targv}^{-1/2} 
	\|
	& \leq &
	\rhoIF
	< 
	1 .
\label{vcjcfvedt7wdwhesqgghwqLco}
\end{EQA}
\end{description}

\begin{lemma}
\label{Lidentsemi}
Suppose \nameref{IFLref}. 
Then 
\begin{EQA}
	(1 - \rhoIF) \blk\{ \IFT_{\targv\targv},\IFT_{\nuiv\nuiv} \}
	\leq 
	\IFT
	& \leq &
	(1 + \rhoIF) \blk\{ \IFT_{\targv\targv},\IFT_{\nuiv\nuiv} \}
\label{f6eh3ewfd6ehev65r43tre}
\end{EQA}
and the matrices
\begin{EQA}
	\IFTb_{\targv\targv} 
	& \eqdef &
	\IFT_{\targv\targv} - \IFT_{\targv\nuiv} \, \IFT_{\nuiv\nuiv}^{-1} \, \IFT_{\nuiv\targv} \, ,
	\qquad
	\IFTb_{\nuiv\nuiv}
	\eqdef 
	\IFT_{\nuiv\nuiv} - \IFT_{\nuiv\targv} \, \IFT_{\targv\targv}^{-1} \, \IFT_{\targv\nuiv} \, ,
\label{ljhy6furf4rfeder8jfu8rd}
\end{EQA}
satisfy
\begin{EQA}
	(1 - \rhoIF) \, \IFT_{\targv\targv}
	\leq 
	\IFTb_{\targv\targv}
	& \leq &
	\IFT_{\targv\targv} \, ,
	\qquad
	(1 - \rhoIF) \, \IFT_{\nuiv\nuiv}
	\leq 
	\IFTb_{\nuiv\nuiv}
	\leq 
	\IFT_{\nuiv\nuiv} \, .
\label{yg3w5dffctvyry4r7er7dgfe}
\end{EQA}
\end{lemma}

\begin{proof}
Define \( \DV^{2} = \IFT_{\targv\targv} \), \( \HV^{2} = \IFT_{\nuiv\nuiv} \),
\( \IFL = \blk\{ \DV^{2},\HV^{2} \} \), \( \AFblkn = \DV^{-1} \IFT_{\targv\nuiv} \, \HV^{-1} \), and consider the matrix 
\begin{EQA}
	\IFL^{-1/2} \IFT \, \IFL^{-1/2}
	&=&
	\begin{pmatrix}
		\Id_{\dimp} & \DV^{-1} \IFT_{\targv\nuiv} \, \HV^{-1}
		\\
		\HV^{-1} \IFT_{\nuiv\targv} \, \DV^{-1} & \Id_{\dimq}
	\end{pmatrix}
	=
	\begin{pmatrix}
		\Id_{\dimp} & \AFblkn
		\\
		\AFblkn^{\T} & \Id_{\dimq}
	\end{pmatrix} .
\label{gswrew35e35e35e5txzgtqw}
\end{EQA}
Condition \eqref{vcjcfvedt7wdwhesqgghwqLco} implies \( \| \AFblkn \AFblkn^{\T} \| \leq \rhoIF \) and hence,
\begin{EQA}
	1 - \rhoIF
	& \leq &
	\| \IFL^{-1/2} \IFT \, \IFL^{-1/2} \| 
	\geq 
	1 + \rhoIF .
\label{ye376evhgder52wedsytg}
\end{EQA}
Moreover,
\begin{EQA}
	\IFTb_{\targv\targv} 
	& = &
	\IFT_{\targv\targv} - \IFT_{\targv\nuiv} \, \IFT_{\nuiv\nuiv}^{-1} \, \IFT_{\nuiv\targv} 
	=
	\DV (\Id_{\dimp} - \AFblkn \AFblkn^{\T}) \DV 
	\geq 
	(1 - \rhoIF) \DV^{2} \, ,
\label{ljhy6furf4jfu8rdfcweerdd}
\end{EQA}
and similarly for \( \IFTb_{\nuiv\nuiv} \).
\end{proof}
\noindent
Note that \( \DV^{2} = \IFT_{\targv\targv} \) is the \( \targv\targv \)-block of \( \IFT = - \nabla^{2} \lgd(\prmtvs) \), while 
\( \DVb^{-2} = \IFTb_{\targv\targv}^{-1} \) is the \( \targv\targv \)-block of \( \IFT^{-1} \)
by the formula of block-inversion.
The matrices \( \DV^{2} \) and \( \DVb^{2} \) only coincide if the matrix \( \IFT \) is of block-diagonal structure.

\begin{lemma}
\label{Onepointorto}
With block representation \eqref{hwe78yf2diwe76tfw67etfwtbso},
define \( \CFT = \IFT_{\nuiv\nuiv}^{-1} \, \IFT_{\nuiv\targv} \) and
\begin{EQ}[rcl]
	\nuov 
	&=& 
	\nuiv + \IFT_{\nuiv\nuiv}^{-1} \, \IFT_{\nuiv\targv} \, (\targv - \targvs) 
	=
	\nuiv + \CFT \, (\targv - \targvs) ,
	\\
	\lgdb(\targv,\nuov)
	&=&
	\lgd(\targv,\nuiv) 
	=
	\lgd(\targv,\nuov - \CFT \, (\targv - \targvs)) .
\label{0vcucvf653nhftdreyenwe}
\end{EQ}
Then 
\begin{EQA}
	\nabla_{\nuov} \nabla_{\targv} \, \lgdb(\targv,\nuov) \Big|_{\binomb{\targv=\targvs}{\nuov=\nuivs}}
	&=&
	0
	\qquad
\label{jdgyefe74erjscgygfydhw}
	\\
	\nabla_{\targv\targv}^{2} \, \lgdb(\targv,\nuov) \Big|_{\binomb{\targv=\targvs}{\nuov=\nuivs}} 
	&=&
	\DVb^{2} 
	 \, .
\label{n6d2fgujrt6fgfrjrd}
\end{EQA}
\end{lemma}

\begin{proof}
By definition \( \lgd(\targv,\nuiv) = \lgd(\targv,\nuov - \CFT \, \targv) \), and it is straightforward to check that
\begin{EQA}
	\nabla_{\nuov} \nabla_{\targv} \, \lgdb(\targv,\nuov) \Big|_{\binomb{\targv=\targvs}{\nuov=\nuivs}}
	=
	\nabla_{\nuov} \nabla_{\targv} \, \lgd(\targv,\nuov 
	- \CFT \, (\targv - \targvs)) \Big|_{\binomb{\targv=\targvs}{\nuov=\nuivs}} 
	&=&
	0 .
	\qquad
\label{jdgyefe74erjscgygfydhwp}
\end{EQA}
Similarly
\begin{EQA}
	\nabla_{\targv\targv}^{2} \, \lgdb(\prmtvs) 
	=
	\nabla_{\targv\targv}^{2} \, \lgd(\targv,\nuov - \CFT \, (\targv - \targvs)) \Big|_{\binomb{\targv=\targvs}{\nuov=\nuivs}}
	&=&
	\IFT_{\targv\targv} - \IFT_{\targv\nuiv} \, \IFT_{\nuiv\nuiv}^{-1} \, \IFT_{\nuiv\targv} 
	=
	\DVb^{2} \, 
\label{d3gr7hyrfygieshfdhfjk}
\end{EQA}
as required.
\end{proof}

We can summarize that the linear transform \eqref{0vcucvf653nhftdreyenwe} ensures the one-point orthogonality
condition \( \nabla_{\nuov} \nabla_{\targv} \, \lgdb(\prmtvs) = 0 \) for the function 
\( \lgdb(\targv,\nuov) = \lgd(\targv,\nuov - \cc{C} \, (\targv - \targvs)) \) with 
\( - \nabla_{\targv\targv}^{2} \lgdb(\prmtvs) = \DVb^{2} \).
For this function, one can redefine all the characteristics including
\begin{EQA}
	\breve{\feta}_{\nuov}
	& \eqdef &
	\max_{\targv} \lgdb(\targv,\nuov),
	\\
	\targvb_{\nuov}
	& \eqdef &
	\argmax_{\targv} \lgdb(\targv,\nuov),
	\\
	\DVb_{\nuov}^{2}
	& \eqdef &
	- \nabla_{\targv\targv}^{2} \lgdb(\targvb_{\nuov},\nuov). 
\label{nfyehendcghegwdgqgqngd}
\end{EQA}
If \( \lgd(\targv,\nuiv) \) is quadratic, then orthogonality can be achieved by a linear transform of the nuisance parameter \( \nuiv \).
In particular, for \( \lgd(\prmtv) = \lgd(\prmtvs) - (\prmtv - \prmtvs)^{\T} \IFT \, (\prmtv - \prmtvs)/2 \) quadratic, it holds with 
\( \HV^{2} = \IFT_{\nuiv\nuiv} \)
\begin{EQA}
	\lgdb(\targv,\nuov)
	&=&
	\lgd(\prmtvs) - \frac{\| \DVb (\targv - \targvs) \|^{2}}{2}  - \frac{\| \HV (\nuov - \nuivs) \|^{2}}{2}  \, ,
\label{ovjv8734mfdgdytwndbgg}
\end{EQA}
and hence
\begin{EQA}
	\targvb_{\nuov}
	& \equiv &
	\targvs,
	\qquad
	\DVb_{\nuov}^{2}
	\equiv 
	\DVb^{2} ,
	\qquad
	\breve{\feta}_{\nuov}
	=
	- \frac{\| \HV (\nuov - \nuivs) \|^{2}}{2} \, .
\label{lfibubgue5t335r6tt7ehw}
\end{EQA}
Next we discuss another special case when the mixed derivative of 
\( \lgd(\targv,\nuiv) \) only depends on \( \targv \):
\begin{EQA}
	\nabla_{\nuiv} \nabla_{\targv} \lgd(\targv,\nuiv)
	& \equiv &
	- \dCross(\targv) ,
	\qquad
	\nabla_{\nuiv}^{2} \lgd(\targv,\nuiv)
	\equiv 
	- \dPartial^{2}(\targv),
\label{0hy533rycirudle}
\end{EQA}
for a \( \dimq \times \dimp \)-matrix function \( \dCross(\targv) \) and a positive \( \dimq \times \dimq \)-matrix function 
\( \dPartial^{2}(\targv) \).
We also write \( \dCross \) and \( \dPartial^{2} \) in place of \( \dCross(\targvs) \) and \( \dPartial^{2}(\targvs) \).
Due to the next result, this condition yields a semi-orthogonality of the transformed function \( \lgdb(\cdot) \).

\begin{lemma}
\label{LCrossdervs}
Assume \eqref{0hy533rycirudle}.
Then the function \( \lgdb(\targv,\nuov) \) from \eqref{0vcucvf653nhftdreyenwe} satisfies the semi-orthogonality condition 
\begin{EQA}
	\nabla_{\nuov} \nabla_{\targv} \lgdb(\targvs,\nuov)
	& \equiv &
	0.
\label{juifue36fhg65e3yeh}
\end{EQA}
Moreover, 
for any \( \nuov \) with \( (\targvs,\nuov) \in \Upsb \), it holds
\begin{EQA}[rcccl]
	\nabla_{\targv} \lgdb(\targvs,\nuov)
	& \equiv &
	0,
	\qquad
	\nabla_{\targv\targv}^{2} \lgdb(\targvs,\nuov)
	& \equiv &
	\nabla_{\targv\targv}^{2} \lgdb(\targvs,\nuivs) .
\label{klfu87fjhr57765645r}
\end{EQA}
If \( \lgd(\prmtv) \) is concave in \( \prmtv \) then 
\begin{EQA}
	\targvb_{\nuov}
	& \eqdef &
	\argmax_{\targv} \lgdb(\targv,\nuov)
	\equiv 
	\targvs \, .	
\label{juifue36fhg65e3yehop}
\end{EQA}
\end{lemma}

\begin{proof}
Let \( \lgdb(\targv,\nuov) \) be defined by \eqref{0vcucvf653nhftdreyenwe}.
Lemma~\ref{Onepointorto} ensures one-point orthogonality 
\( \nabla_{\nuov} \nabla_{\targv} \lgdb(\prmtvs) = \nabla_{\nuov} \nabla_{\targv} \lgdb(\targvs,\nuivs) = 0 \).
Moreover, \eqref{0hy533rycirudle} yields \( \CFT = \dPartial^{-2} \, \dCross \).
Now consider \( \nabla_{\nuov} \nabla_{\targv} \lgdb(\targvs,\nuov) \).
It holds by \eqref{0hy533rycirudle} and definition \eqref{0vcucvf653nhftdreyenwe}
\begin{EQA}
	\nabla_{\nuov} \nabla_{\targv} \lgdb(\targvs,\nuov) 
	&=&
	\nabla_{\nuov} \nabla_{\targv} \lgd(\targv,\nuov - \CFT \, (\targv - \targvs)) \bigg|_{\targv = \targvs}
	\\
	&=&
	\nabla_{\nuov} \nabla_{\targv} \lgd(\targvs,\nuov) - \nabla_{\nuov}^{2} \lgd(\targvs,\nuov) \CFT
	=
	- \dCross + \dPartial^{2} \CFT
	=
	0
\label{ucieruier3dchedndm}
\end{EQA}
yielding semi-orthogonality of \( \lgdb \).
Now \eqref{klfu87fjhr57765645r} and \eqref{juifue36fhg65e3yehop} follow from Lemma~\ref{Lsemiorto}.
\end{proof}

For a general smooth function \( \lgd(\targv,\nuiv) \) satisfying \nameref{LLsS3ref}, 
we expect that these identities are fulfilled
in a local vicinity of \( \prmtvs \) up to the error of quadratic approximation.
The next lemma quantifies this guess.

\begin{proposition}
\label{PonepptLL}
Let \( \lgd(\prmtv) \) satisfy \nameref{IFLref}, \nameref{LLsS3ref} with \( \rr = \rru \), and
\begin{EQA}
	\dltwbss
	&=&
	\CONSTIFT \, \hmax_{3} \, \rru n^{-1/2} 
	\leq 
	1/3 .
\label{8uiiihkkcrdrdteggcjcjdan}
\end{EQA} 
Then the result \eqref{36gfijh94ejdvtwekoi} of Proposition~\ref{LsemiAv} continues to apply.
\end{proposition}

\begin{proof}
The linear transform \eqref{0vcucvf653nhftdreyenwe} reduces the statement 
to the case with one-point orthogonality. 
Condition \eqref{8uiiihkkcrdrdteggcjcjdan} ensures \eqref{8uiiihkkcrdrdteggcjcjd} after
this transform. 
Conditions \nameref{LLseS3ref}, \nameref{LLseS4ref}, \nameref{LLsoS3ref}, and \nameref{LLpS3ref}
follow from \nameref{LLsS3ref}.
\end{proof}

\Subsection{Composite nuisance variable}
For some situations, the nuisance variable \( \nuiv \) is by itself a composition of few other subvectors.
Checking condition \nameref{IFLref} in the scope of variables can be involved.
However, it can be reduced to a check for each subvectors.
We only consider the case of two variables \( \nuiv = (\zv,\nuov) \), that is,
\begin{EQA}
	\lgd(\prmtv)
	&=&
	\lgd(\targv,\nuiv)
	=
	\lgd(\targv,\zv,\nuov) .
\label{576vner6734nf9h4ede}
\end{EQA}
Denote by \( \IFT_{\targv,\targv} \), \( \IFT_{\targv,\zv} \), \( \IFT_{\targv,\nuov} \), 
\( \IFT_{\zv,\zv} \), \( \IFT_{\nuov,\nuov} \) the corresponding blocks of \( \IFT \).

\begin{lemma}
\label{LblocksIFT}
For the matrix \( \IFT = - \nabla^{2} \lgd(\targv,\zv,\nuov) \), suppose that 
\begin{EQA}
\label{vcjcfvedt7wdwhesqgLcomp}
	\| 
	\IFT_{\targv\targv}^{-1/2} \IFT_{\targv\zv} \, \IFT_{\zv\zv}^{-1} \, \IFT_{\zv\targv} \, \IFT_{\targv\targv}^{-1/2} 
	\|
	& \leq &
	\rhoIF_{\zv}
	< 
	1 ,
	\\
	\| 
	\IFT_{\targv\targv}^{-1/2} \IFT_{\targv\nuov} \, \IFT_{\nuov\nuov}^{-1} \, \IFT_{\nuov\targv} \, \IFT_{\targv\targv}^{-1/2} 
	\|
	& \leq &
	\rhoIF_{\nuov}
	< 
	1 .
\label{vcjcfvedt7wdwhesqgLcompn}
\end{EQA}
Then \nameref{IFLref} is fulfilled for \( \nuiv = (\zv,\nuov) \) with \( \rhoIF = \rhoIF_{\zv} + \rhoIF_{\nuov} \).
\end{lemma}

\begin{proof}
Define the scalar product 
\begin{EQA}
	\bigl\langle \prmtv_{1},\prmtv_{2} \bigr\rangle
	&=&
	\prmtv_{1}^{\T} \IFT \, \prmtv_{2} \, .
\label{8cney7f76ryerdeiek}
\end{EQA}
Let us fix any vectors \( \targv \), \( \zv \), \( \nuov \) with the natural embedding in the \( \prmtv \)-space and define
\begin{EQA}
	\targv_{1}
	& \eqdef &
	\targv - \frac{\langle \targv,\zv \rangle}{\langle \zv,\zv \rangle} \zv ,
	\qquad
	\nuov_{1}
	\eqdef 
	\nuov - \frac{\langle \nuov,\zv \rangle}{\langle \zv,\zv \rangle} \zv.
\label{3wfduf6w3yhesdicdurete}
\end{EQA}
Then \( \langle \targv_{1},\zv \rangle = 0 \), \( \langle \nuov_{1},\zv \rangle = 0 \).
Condition \eqref{vcjcfvedt7wdwhesqgLcomp} yields
\begin{EQA}
	\| \targv_{1} \|
	& \geq &
	(1 - \rhoIF_{\zv}) \| \targv \| .
\label{gc7e3jwka8dfyey3ehwk}
\end{EQA}
Similarly define
\begin{EQA}
	\targv_{2}
	& \eqdef &
	\targv_{1} - \frac{\langle \targv_{1},\nuov_{1} \rangle}{\langle \nuov_{1},\nuov_{1} \rangle} \nuov_{1}
\label{ydcf7w3ufhjvcfre6eyw39xj}
\end{EQA}
Then \( \langle \targv_{2},\zv \rangle = 0 \) and \( \langle \targv_{2},\nuov \rangle = 0 \) and by \eqref{vcjcfvedt7wdwhesqgLcompn}
\begin{EQA}
	\| \targv_{2} \|
	& \geq &
	(1 - \rhoIF_{\nuov}) \| \targv_{1} \| 
	\geq 
	(1 - \rhoIF_{\zv}) (1 - \rhoIF_{\nuov}) \| \targv \| 
	\geq 
	(1 - \rhoIF_{\zv} - \rhoIF_{\nuov}) \| \targv \|
\label{gc7e3jwka8dfyey3ehwk}
\end{EQA}
and \nameref{IFLref} follows.
\end{proof}


\Chapter{Dimension free bounds for Laplace approximation}
\label{STaylor}
Here we present several issues related to Laplace approximation.
Section~\ref{SboundsLapl} states general results about the accuracy of Laplace approximation for finite samples.
\ifinexact{Section~\ref{SLaplinexact} discusses inexact approximation and the use of posterior mean.}{}
Technical assertions and proofs are collected in Section~\ref{SLapltools}.

\Section{Setup and conditions}
\label{SsetupLapl}
Let \( \lgd(\xv) \) be a function in a high-dimensional Euclidean space \( \R^{\dimp} \) such that
\( \int \ex^{\lgd(\xv)} \, d\xv = \CONST < \infty \),
where the integral sign \( \int \) without limits means the integral over the whole space \( \R^{\dimp} \).
Then \( \lgd \) determines a distribution \( \PfL \) with the density
\( \CONST^{-1} \ex^{\lgd(\xv)} \).
Let \( \xvs \) be a point of maximum:
\begin{EQA}
	\lgd(\xvs)
	&=&
	\sup_{\uv \in \R^{\dimp}} \lgd(\xvs + \uv) .
\label{scdygw7ytd7wqqsquuqydtdtd}
\end{EQA}
We also assume that \( \lgd(\cdot) \) is at least three time differentiable. 
Introduce the negative Hessian \( \IFL = - \nabla^{2} \lgd(\xvs) \) and assume \( \IFL \) strictly positive definite.
We aim at approximating the measure \( \PfL \) by a Gaussian measure \( \ND(\xvs,\IFL^{-1}) \).
Given a function \( g(\cdot) \), define its expectation w.r.t. \( \PfL \) after centering at \( \xvs \):
\begin{EQA}
	\II(g)
	& \eqdef &
	\frac{\int g(\uv) \, \ex^{\lgd(\xvs + \uv)} \, d\uv}{\int \ex^{\lgd(\xvs + \uv)} \, d\uv} \, .
\label{IIgfigefxududu}
\end{EQA}
A Gaussian approximation \( \II_{\IFL}(g) \) for \( \II(g) \) is defined as
\begin{EQA}
	\II_{\IFL}(g)
	& \eqdef &
	\frac{\int g(\uv) \, \ex^{- \| \IFL^{1/2} \uv \|^{2}/2} \, d\uv}{\int \ex^{- \| \IFL^{1/2} \uv \|^{2}/2} \, d\uv} 
	=
	\E g(\gaussv_{\IFL}) ,
	\qquad
	\gaussv_{\IFL} \sim \ND(0,\IFL^{-1}) \, .
\label{IgiLapgeHvdu}
\end{EQA}
The choice of the distance between \( \PfL \) and \( \ND(\xvs,\IFL^{-1}) \) specifies the considered class 
of functions \( g \).
The most strong total variation distance can be obtained as 
the supremum of \( |\II(g) - \II_{\IFL}(g)| \) over all measurable functions \( g(\cdot) \) with 
\( |g(\uv)| \leq 1 \):
\begin{EQA}
	\TV\bigl( \PfL,\ND(\xvs,\IFL^{-1}) \bigr)
	&=&
	\sup_{\|g\|_{\infty} \leq 1} \bigl| \II(g) - \II_{\IFL}(g) \bigr| \, .
\label{ghdvcgftdftgegdftefd3434545}
\end{EQA}
The results can be substantially improved if only centrally symmetric functions \( g(\cdot) \)
with \( g(\xv) = g(-\xv) \) are considered.
Obviously, for any \( g(\cdot) \)
\begin{EQA}
	\II(g)
	& = &
	\frac{\int g(\uv) \, \ex^{\lgd(\xvs + \uv) - \lgd(\xvs)} \, d\uv}{\int \ex^{\lgd(\xvs + \uv) - \lgd(\xvs)} \, d\uv} \, .
\label{IIgfigefxuudufxudu}
\end{EQA}
Moreover, as \( \xvs = \argmax_{\xv} \lgd(\xv) \), it holds \( \nabla \lgd(\xvs) = 0 \) and 
\begin{EQA}
	\II(g)
	& = &
	\frac{\int g(\uv) \, \ex^{\lgd(\xvs;\uv)} \, d\uv}{\int \ex^{\lgd(\xvs;\uv)} \, d\uv} \, ,
\label{IIgfifxutudufxutdu}
\end{EQA}
where \( \lgd(\xv;\uv) \) is the Bregman divergence 
\begin{EQA}
	\lgd(\xv;\uv)
	&=&
	\lgd(\xv + \uv) - \lgd(\xv) - \bigl\langle \nabla \lgd(\xv), \uv \bigr\rangle .
\label{fxufxpufxfpxu}
\end{EQA}
Implicitly we assume that the negative Hessian \( \IFL = - \nabla^{2} \lgd(\xvs) \) is sufficiently large
in the sense that the Gaussian measure \( \ND(0,\IFL^{-1}) \) concentrates on a small local set \( \UVL \).
This allows to use a local Taylor expansion for 
\( \lgd(\xvs;\uv) \approx - \| \IFL^{1/2} \uv \|^{2}/2 \) in \( \uv \) on \( \UVL \).
If \( \lgd(\cdot) \) is also strongly concave, then the \( \PfL \)-mass  
of the complement of \( \UVL \) is exponentially small yielding the desirable Laplace approximation.

Our setup is motivated by Bayesian inference.
Assume for a moment that
\begin{EQA}
	\lgd(\xv)
	&=&
	\lgdL(\xv) - \| \GP (\xv - \xv_{0}) \|^{2} / 2
\label{jhdctrdfred4322edt7y}
\end{EQA}
for some \( \xv_{0} \) and a symmetric \( \dimp \)-matrix \( \GP^{2} \geq 0 \).
Here \( \lgdL(\cdot) \) stands for a log-likelihood function while the quadratic penalty 
\( \| \GP (\xv - \xv_{0}) \|^{2} / 2 \) corresponds to a Gaussian prior \( \ND(\xv_{0},\GP^{-2}) \).
Let also \( \lgdL(\cdot) \) be concave with \( \DVL^{2} \eqdef - \nabla^{2} \lgdL(\xvs) > 0 \).
Then  
\begin{EQA}
	- \nabla^{2} \lgd(\xvs)
	&=&
	- \nabla^{2} \lgdL(\xvs) + \GP^{2} 
	=
	\DVL^{2} + \GP^{2} .
\label{hydsf42wdsdtrdstdg}
\end{EQA} 
In typical asymptotic setups, the log-likelihood function \( \lgdL(\xv) \) scales with the sample size \( n \)
or inverse noise variance while the prior is kept fixed; see, e.g. \cite{SSW:2020,HeKr2022}.
We allow \( \GP^{2} \) depend on \( n \) as well, this is important for obtaining the dimension free results
and to obtain optimal contraction rate; see Section~\ref{Ssmoothprior}.

\Subsection{Concavity}
Below we implicitly assume decomposition \eqref{jhdctrdfred4322edt7y} with 
a \emph{weak\-ly concave} function \( \lgdL(\cdot) \).
More specifically, we assume the following condition.

\begin{description}
    \item[\label{LLf0ref} \( \bb{(\mathcal{C}_{0})} \)]
      \textit{ There exists an operator \( \GP^{2} \geq 0 \) in \( \R^{\dimp} \) such that 
      \( \GP^{2} \leq - \nabla^{2} \lgd(\xvs) \) and  
\begin{EQA}
	\lgdL(\xvs + \uv)
	& \eqdef &
	\lgd(\xvs + \uv) + \| \GP \uv \|^{2} / 2
\label{fTvufupumDTpDfx}
\end{EQA}
is a concave function. 
      }
\end{description}

If \( \lgdL(\cdot) \) in decomposition \eqref{jhdctrdfred4322edt7y} is concave then this condition is obviously fulfilled. 
More generally, if \( \lgdL(\cdot) \) in \eqref{jhdctrdfred4322edt7y} is weakly concave, so that 
\( \lgdL(\xvs + \uv) - \| \GP_{0} \uv \|^{2}/2 \) is concave in \( \uv \) with \( \GP_{0}^{2} \leq \GP^{2} \), then 
\nameref{LLf0ref} is fulfilled with \( \GP^{2} - \GP_{0}^{2} \) in place of \( \GP^{2} \).

The operator \( \DVL^{2} \) plays an important role in our conditions and results:
\begin{EQA}
	\DVL^{2} 
	&=&
	- \nabla^{2} \lgd(\xvs) - \GP^{2}
	\quad
	\bigl( = - \nabla^{2} \lgdL(\xvs) \text{ under \eqref{jhdctrdfred4322edt7y}} \bigr) .
\label{kc7c322dgdf43djhd}
\end{EQA}

\begin{remark}
The condition of strong concavity of \( \lgd \) on the whole space \( \R^{\dimp} \)
can be too restrictive%
\ifapp{; see Section~\ref{Slaplnonlin} for an example}{}.
This condition can be replaced by its local version: \( \lgd \) is concave on a set \( \XXL_{0} \) 
such that the Gaussian prior \( \ND(\xv_{0},\GP^{-2}) \) concentrates on \( \XXL_{0} \)
and the maximizer \( \xvs_{\GP} \) belongs to \( \XXL_{0} \)%
\ifapp{; see Section~\ref{Slaplnonlin}}{}.
In all the results, the integral over \( \R^{\dimp} \) has to be replaced by the integral over \( \XXL_{0} \).
\end{remark}

\Subsection{Laplace effective dimension}
\label{SeffdimLa}
With decomposition \eqref{kc7c322dgdf43djhd} in mind, we use a decomposition for \( \IFL = - \nabla^{2} \lgd(\xvs) \):
\begin{EQA}
	\IFL 
	&=& 
	- \nabla^{2} \lgd(\xvs)
	=
	\DVL^{2} + \GP^{2} .
\label{dscytf5w2edte5rw4e24gyd}
\end{EQA}
The \emph{Laplace effective dimension} \( \dimLL \) is given by
\begin{EQA}
	\dimLL
	& \eqdef & 
	\tr\bigl( \DVL^{2} \, \IFL^{-1} \bigr) 
	=
	\tr \bigl\{ \DVL^{2} \, (\DVL^{2} + \GP^{2})^{-1} \bigr\}.
\label{dAdetrH02Hm2}
\end{EQA}
Of course, \( \dimLL \leq \dimp \) but a proper choice of the penalty \( \GP^{2} \) in \eqref{jhdctrdfred4322edt7y}
allows to avoid the ``curse of dimensionality'' issue and ensure a small effective dimension \( \dimLL \) even for \( \dimp \) large or infinite; see 
\ifapp{Section~\ref{Ssmoothprior}}{\cite{SpPa2019}} 
for more rigorous discussion.

Later we write \( \dimL \) instead of \( \dimLL \) without risk of confusion because the parameter dimension \( \dimp \)
does not show up anymore.
The value \( \dimL \) helps to describe a local vicinity \( \UVL \) around \( \xvs \) such that the most of mass of \( \PfL \)
concentrates on \( \UVL \); see Section~\ref{StailLa}.
Namely, let us fix some \( \amax < 1 \), e.g. \( \amax = 2/3 \), and some \( \xx > 0 \) ensuring that \( \ex^{-\xx} \)
is our significance level.
Define
\begin{EQ}[rcccl]
	\rrL
	& = &
	2 \sqrt{\dimL} + \sqrt{2 \xx} ,
	\qquad
	\UVL 
	&=& 
	\bigl\{ \uv \colon \| \DVL \uv \| \leq \amax^{-1} \rrL \bigr\} .
\label{UvTDunm12spT}
\end{EQ}

\Subsection{Local smoothness conditions}
Let \( \dimp \leq \infty \) and
let \( \lgd(\cdot) \) be a three times continuously differentiable function on \( \R^{\dimp} \).
We fix a reference point \( \xv \) and local region around \( \xv \) given by the local set 
\( \UVL \subset \R^{\dimp} \) from \eqref{UvTDunm12spT}. 
Consider the remainder of the second and third order Taylor approximation 
\begin{EQ}[rcl]
	\dltw_{3}(\xv,\uv)
	&=&
	\lgd(\xv;\uv) - 
	\bigl\langle \nabla^{2} \lgd(\xv) , \uv^{\otimes 2} \bigr\rangle/2 ,
	\\
	\dltw_{4}(\xv,\uv)
	&=&
	\lgd(\xv;\uv) - 
	\bigl\langle \nabla^{2} \lgd(\xv) , \uv^{\otimes 2} \bigr\rangle/2 
	- 
	\bigl\langle \nabla^{3} \lgd(\xv), \uv^{\otimes 3} \bigr\rangle / 6
\label{d4fuv1216303}
\end{EQ}
with \( \lgd(\xv;\uv) \) from \eqref{fxufxpufxfpxu}.
The use of the Taylor formula allows to bound
\begin{EQ}[rcl]
	\bigl| \dltw_{k}(\xv,\uv) \bigr|
	& \leq &
	\sup_{t \in [0,1]}
	\frac{1}{k!} 
	\Bigl| \bigl\langle \nabla^{k} \lgd(\xv + t \uv), \uv^{\otimes k} \bigr\rangle \Bigr|,
	\quad
	k\geq 3. 
\label{k3t01f12n3fvtu}
\end{EQ}
Note that the quadratic penalty \( - \| \GP (\xv - \xv_{0}) \|^{2}/2 \) in \( \lgd \) does not affect 
the remainders \( \dltw_{3}(\xv,\uv) \) and \( \dltw_{4}(\xv,\uv) \).
Indeed, 
with \( \lgd(\xv) = \lgdL(\xv) - \| \GP (\xv - \xv_{0}) \|^{2}/2 \), it holds 
\begin{EQA}
	\lgd(\xv;\uv) 
	& \eqdef & 
	\lgd(\xv + \uv) - \lgd(\xv) - \bigl\langle \nabla \lgd(\xv), \uv \bigr\rangle 
	=
	\lgdL(\xv;\uv) - \| \GP \uv \|^{2}/2 
\label{D2fuvGu2fG}
\end{EQA}
and the quadratic term in definition of the values \( \dltw_{k}(\xv,\uv) \) cancels, \( k \geq 3 \). 
Local smoothness of \( \lgd(\cdot) \) or, equivalently, of \( \lgdL(\cdot) \), at \( \xv \) will be measured by
the value \( \dltwb(\xv) \):
\begin{EQA}
	\dltwb(\xv)
	& \eqdef &
	\sup_{\uv \in \UVL} \frac{1}{\| \DVL \uv \|^{2}/2} \bigl| \dltw_{3}(\xv,\uv) \bigr| ;
\label{om3esuU1H02d3}
\end{EQA}
cf. \eqref{dtb3u1DG2d3GPg}.
We also denote \( \dltwb \eqdef \dltwb(\xvs) \).
Our results apply if \( \dltwb \ll 1 \).
Local concentration of the measure \( \PfL \) requires \( \dltwb \leq 1/3 \);
see Proposition~\ref{PlocconLa}.
The results about Gaussian approximation are valid under a stronger condition 
\( \dltwb \, \dimL \leq 2/3 \) with 
\( \dimL \) from \eqref{dAdetrH02Hm2}.


\Section{Error bounds for Laplace approximation}
\label{SboundsLapl}
Our first result describes the quality of approximation of the measure \( \PfL \) by the Gaussian measure \( \ND(\xvs,\IFL^{-1}) \)
with mean \( \xvs \) and the covariance \( \IFL^{-1} \) in total variation distance.
In all our result, the value \( \xx \) is fixed to ensure that \( \ex^{-\xx} \) is negligible.
First we present the general results which will be specified later under the self-concordance condition.


\begin{theorem}
\label{TLaplaceTV}
Suppose \nameref{LLf0ref}.
Let also \( \dimL \) be defined by \eqref{dAdetrH02Hm2} and \( \rrL \) and \( \UVL \) by \eqref{UvTDunm12spT}.
If \( \dltwb \) from \eqref{om3esuU1H02d3} satisfies \( \dltwb \leq 1/3 \), then
\begin{EQA}
	\PfL(\Xv - \xvs \not\in \UVL)
	& \leq &
	\ex^{-\xx} .
\label{poybf3679jd532ff2}
\end{EQA}
If \( \dltwb \, \dimL \leq 2/3 \),
then for any \( g(\cdot) \) with \( |g(\uv)| \leq 1 \), it holds for \( \II(g) \) from \eqref{IIgfigefxududu}
\begin{EQA}
	\bigl| \II(g) - \II_{\GP}(g) \bigr|
	& \leq &
	\frac{2 (\err + \ex^{-\xx})}{1 - \err - \ex^{-\xx}}
	\leq 
	4(\err + \ex^{-\xx}) 
\label{ufgdt6df5dtgededsxd23gjg}
\end{EQA}
with 
\begin{EQA}
	\err
	&=&
	\err_{2} 
	= 
	\frac{0.75 \, \dltwb \, \dimL}{1 - \dltwb} \, .
\label{juytr90f2dzaryjhfyf}
\end{EQA}
\end{theorem}

This section presents more advanced bounds on the error of Laplace approximation
under conditions \nameref{LL3tref} and \nameref{LL4tref} with \( \upsv = \xvs \) 
or \nameref{LLtS3ref} and \nameref{LLtS4ref} with \( \upsv = \xvs \) and \( \HL(\xvs) = n^{-1/2} \, \DVL \);
see Section~\ref{Slocalsmooth}.

\begin{theorem}
\label{TLaplaceTV34}
Suppose \nameref{LLf0ref} and \nameref{LL3tref} 
and let \( \dltwu_{3} \, \amax^{-1} \rrL \leq 3/4 \) for \( \rrL \) from \eqref{UvTDunm12spT}.
Then the concentration bound \eqref{poybf3679jd532ff2} holds.
Moreover, if 
\begin{EQA}
	\dltwu_{3} \, \amax^{-1} \rrL \, \dimL
	& \leq &
	2 ,
\label{0hcde4dft3igthg94yr5twew}
\end{EQA}
then the accuracy bound \eqref{ufgdt6df5dtgededsxd23gjg} applies with 
\( \normG = \| \DVL \, \IFL^{-1} \DVL \| \leq 1 \)
\begin{EQA}
	\err
	&=&
	\err_{3}
	\eqdef
	\frac{\dltwu_{3} (\dimL + \normG)^{3/2}}{4 (1 - \dltwb)^{3/2}} 
	\leq 
	\frac{\dltwu_{3} \, (\dimL + \normG)^{3/2}}{2} \, ,
\label{5qw7dyf4e4354coefw9dufih}
\end{EQA}
where \( \dltwb \eqdef \dltwu_{3} \, \amax^{-1} \rrL / 3 \leq 1/4 \).
Furthermore, under \nameref{LL4tref}, 
for any symmetric function \( g(\uv) = g(-\uv) \), \( |g(\uv)| \leq 1 \),
the accuracy bound \eqref{ufgdt6df5dtgededsxd23gjg} applies with 
\begin{EQA}
	\err
	&=&
	\err_{4}
	=
	\frac{\dltwu_{3}^{2} \, (\dimL + 2\normG)^{3} + 2 \dltwu_{4} (\dimL + \normG)^{2}}{16 (1 - \dltwb)^{2}} 	
	\leq 
	\frac{\dltwu_{3}^{2} \, (\dimL + 2\normG)^{3} + 2 \dltwu_{4} (\dimL + \normG)^{2}}{8} \, \, .
\label{hg25t6mxwhydseg3hhfdr}
\end{EQA}
Under \nameref{LLtS3ref} and \nameref{LLtS4ref} instead of \nameref{LL3tref} and \nameref{LL4tref}, the results apply with 
\( \dltwu_{3} = \hmax_{3} \, n^{-1/2} \) and \( \dltwu_{4} = \hmax_{4} \, n^{-1} \).
\end{theorem}

Let \( \BBB(\R^{\dimp}) \) be the \( \sigma \)-field of all Borel sets in \( \R^{\dimp} \),
while \( \BBB_{s}(\R^{\dimp}) \) stands for all centrally symmetric sets from \( \BBB(\R^{\dimp}) \).
By \( \Xv \) we denote a random element with the distribution \( \PfL \), while \( \gaussv_{\IFL} \sim \ND(0,\IFL^{-1}) \).

\begin{corollary}
\label{CTLaplaceTV}
Under the conditions of Theorem~\ref{TLaplaceTV34}, it holds for \( \Xv \sim \PfL \)
\begin{EQA}
	\sup_{A \in \BBB(\R^{\dimp})} \bigl| \PfL(\Xv - \xvs \in A) - \P(\gaussv_{\IFL} \in A) \bigr|
	& \leq &
	4(\err_{3} + \ex^{-\xx}),
	\\
	\sup_{A \in \BBB_{s}(\R^{\dimp})} \bigl| \PfL(\Xv - \xvs \in A) - \P(\gaussv_{\IFL} \in A) \bigr|
	& \leq &
	4(\err_{4} + \ex^{-\xx}) .
\label{cbc5dfedrdwewwgerg}
\end{EQA}
\end{corollary}

\Subsection{Critical dimension}
Here we briefly discuss the important issue of \emph{critical dimension} meaning a relation between 
\( \dimL \) and \( n \) sufficient for our results.
Theorem~\ref{TLaplaceTV34} states concentration of \( \PfL \) under the condition \( \dltwu_{3} \rrL \leq 1 \).
Under \nameref{LLtS3ref}, we can use \( \dltwu_{3} = \hmax_{3} \, n^{-1/2} \).
Together with \( \rrL \approx \sqrt{\dimL} \), this yields the condition
\( \dimL \ll n \).
Gaussian approximation applies under \( \hmax_{3} \, \amax^{-1} \rrL \, \dimL \, n^{-1/2} \leq 2 \); see \eqref{0hcde4dft3igthg94yr5twew},
yielding \( \dimL^{3} \ll n \).
We see that there is a gap between these conditions. 
We guess that in the region \( n^{1/3} \lesssim \dimL \lesssim n \), 
non-Gaussian approximation of the posterior is possible; cf. 
\cite{bochkina2014}.

\ifKL{
\Subsection{A bound for Kullback-Leibler divergence}
Theorem~\ref{TLaplaceTV} through \ref{TLaplaceTV34} quantify the approximation \( \PfL \approx \ND(\xvs,\IFL^{-1}) \)
in the total variation distance.
Another useful characteristic could be the Kullback-Leibler (KL) divergence between \( \PfL \) and  
\( \ND(\xvs,\IFL^{-1}) \).
The KL divergence \( \kullb(\P_{1},\P_{2}) = \E_{1} \log(d\P_{1}/d\P_{2}) \) is asymmetric, 
\( \kullb(\P_{1},\P_{2}) \neq \kullb(\P_{2},\P_{1}) \) with few exceptions 
like the case of Gaussian measures \( \P_{1} \) and \( \P_{2} \).
Moreover, \( \kullb(\P_{1},\P_{2}) \) can explode if \( \P_{1} \) is not absolutely continuous w.r.t. \( \P_{2} \). 
We present two bounds for each ordering.
For ease of presentation, we limit ourselves to the case when either \nameref{LL3tref} or \nameref{LLtS3ref} meets.

\begin{theorem}
\label{TLaplaceKL}
Suppose \nameref{LLf0ref} and \nameref{LL3tref} and let \( \dltwb \, \dimL \leq 2/3 \).
Then 
\begin{EQA}
	\kullb(\PfL,\ND(\xvs,\IFL^{-1}))
	& \leq &
	4 \err_{3} + 4 \ex^{-\xx} 
	\leq 
	\frac{\E \dltwu_{3}(\gaussv_{\IFL})}{(1 - \dltwb)^{3/2}} + 4 \ex^{-\xx} .
\label{vlgvi8ugu7tr4ry43et31}
\end{EQA}
Moreover, under \nameref{LLtS3ref}
\begin{EQA}
	\kullb(\PfL,\ND(\xvs,\IFL^{-1}))
	& \leq &
	2 \hmax_{3} \, \sqrt{\frac{(\dimL + \normG)^{3}}{n}} + 4 \ex^{-\xx} .
\label{cheiufheurhrtdgwb3wesrtd}
\end{EQA}
\end{theorem}

Now we briefly discuss the value \( \kullb(\ND(\xvs,\IFL^{-1}),\PfL) \).
We already know that \( \PfL \) concentrates on \( \UVL \) and can be well approximated by 
\( \ND(\xvs,\IFL^{-1}) \) on \( \UVL \).
However, this does not guarantee a small value of \( \kullb(\ND(\xvs,\IFL^{-1}),\PfL) \).
It can even explode if e.g. \( \PfL \) has a compact support.
In fact, 
the log-density of \( \ND(\xvs,\IFL^{-1}) \) w.r.t. \( \PfL \) reads
\begin{EQA}
	\log \frac{d\ND(\xvs,\IFL^{-1})}{d\PfL}(\xv)
	&=&
	- \lgd(\xv) - \frac{1}{2} \| \IFL^{1/2} (\xv - \xvs) \|^{2} - \CONSTi_{\IFL}
\label{fduhefuifhew2tdw6ww3}
\end{EQA}
for some constant \( \CONSTi_{\IFL} \),
and an upper bound on \( \kullb(\ND(\xvs,\IFL^{-1}),\PfL) \) requires that the integral of \( \lgd(\xv) \) w.r.t. the measure \( \ND(\xvs,\IFL^{-1}) \) is finite.

\begin{theorem}
\label{TLaplaceKLi}
Suppose \nameref{LLf0ref} and \nameref{LLtS3ref} and let \( \dltwb \, \dimL \leq 2/3 \).
Let also \( \rho = 2\xx /\rrL^{2} \); see \eqref{UvTDunm12spT}.
If \( \lgdL(\xvs;\uv) = \lgdL(\xvs + \uv) - \lgdL(\xvs) - \langle \nabla \lgdL(\xvs), \uv \rangle \) fulfills
with \( \IFL = \DVL^{2} +\GP^{2} \)
\begin{EQA}
	\int \bigl| \lgdL(\xvs;\uv) \bigr| \, \exp \bigl\{ - \| \IFL^{1/2} \uv \|^{2}/2 + \rho \| \DVL \uv \|^{2}/2 \bigr\} \, d\uv
	& \leq &
	\CONSTi_{\lgdL} \, 
\label{dchbhwdhwwdgscsn2efty2162}
\end{EQA} 
for some fixed constant \( \CONSTi_{\lgdL} \) then
\begin{EQA}
	\kullb(\ND(\xvs,\IFL^{-1}),\PfL)
	& \leq &
	\hmax_{3} \, \sqrt{\frac{(\dimL + \normG)^{3}}{n}} + (2 + \CONSTi_{\lgdL}) \ex^{-\xx} .
\label{cheiufheurhrtdgwb3wesrtd}
\end{EQA}
\end{theorem}
}{}

\ifinexact{
\Subsection{Mean and MAP}
Here we present the bound on \( |\II(g) - \II_{\GP}(g)| \) for the case of a linear vector function \( g(\uv) = \QP \uv \) 
with \( \QP \colon \R^{\dimp} \to \R^{\dimq} \), \( \dimq \geq 1 \). 
A special case of \( \QP = \Id_{\dimp} \) corresponds to the mean value \( \xvb \) of \( \PfL \).
The next result presents an upper bound for \( \QP (\xvb - \xvs) \) under the conditions of Theorem~\ref{TLaplaceTV34}
including \nameref{LLtS3ref}.

\begin{theorem}
\label{TpostmeanLa}
Assume the conditions of Theorem~\ref{TLaplaceTV34} and let \( \QP^{\T} \QP \leq \DVL^{2} \).
Then 
\begin{EQA}
	\| \QP (\xvb - \xvs) \|
	& \leq &
	2.4 \, \hmax_{3} \, \| \QP \, \IFL^{-1} \QP^{\T} \|^{1/2} \,
	\sqrt{(\dimL + \normG)^{3}/n} + 4 \ex^{-\xx} .
\label{hcdtrdtdehfdewdrfrhgyjufger}
\end{EQA}
\end{theorem}

Now we specify the result for the special choice \( \QP = \DVL \).

\begin{corollary}
\label{CTpostmeanLa}
Assume the conditions of Theorem~\ref{TLaplaceTV34}.
Then 
\begin{EQA}
	\| \DVL (\xvb - \xvs) \|
	& \leq &
	2.4 \, \hmax_{3} \, \sqrt{(\dimL + \normG)^{3}/{n}} + 4 \ex^{-\xx} \, .
\label{klu8gitfdgregfkhj7yt}
\end{EQA}
\end{corollary}

\begin{remark}
\label{RTpostmeanc}
An interesting question is whether the result of Theorem~\ref{TpostmeanLa} or Corollary~\ref{CTpostmeanLa} applies 
with \( \QP = \IFL^{1/2} \).
This issue is important in connection to inexact Laplace approximation; see the next section.
The answer is negative. 
The problem is related to the last term \( 4 \ex^{-\xx} \) in the right hand-side of \eqref{hcdtrdtdehfdewdrfrhgyjufger} which is responsible for the tail integral in posterior mean.
This tail integral explodes for \( \QP = \IFL^{1/2} \) and \( \dimp = \infty \).
\end{remark}
}{}

\ifKL{

\Subsection{Posterior covariance}
Posterior covariance \( \IFCov \) is given by
\begin{EQA}
	\IFCov
	& \eqdef &
	\EfL (\Xv - \xvb) (\Xv - \xvb)^{\T} .
\label{hgyetdgdgwtwdfttetste5}
\end{EQA}
Also define
\begin{EQA}
	\IFCovs
	& \eqdef &
	\EfL (\Xv - \xvs) (\Xv - \xvs)^{\T} .
\label{hgyetdgdgwtwdfttetstes}
\end{EQA}
We use an obvious representation
\begin{EQA}
	\IFCov
	&=&
	\EfL (\Xv - \xvs) (\Xv - \xvs)^{\T} - (\xvs - \xvb) (\xvs - \xvb)^{\T} 
	=
	\IFCovs - (\xvs - \xvb) (\xvs - \xvb)^{\T} .
\label{p8txdwydnvivdydhdj}
\end{EQA}
Laplace approximation suggests that \( \IFCov \approx \IFL^{-1} \).
Here we intend to measure the accuracy of this approximation.

\begin{theorem}
\label{TpostvarL}
Assume the conditions of Theorem~\ref{TLaplaceTV34}.
Then
\begin{EQA}
	\| \QP (\IFCov - \IFL^{-1}) \QP^{\T} \|
	& \leq &
	\| \QP \, \IFL^{-1} \QP^{\T} \| ( 9 \hmax_{3}^{2} + \hmax_{4} ) \frac{(\dimA + 3)^{3}}{n}
	+ 8 \ex^{- \xx} .
\label{5hjkkoiioi1223hhhth}
\end{EQA}
\end{theorem}
}{}

\ifinexact{
\Section{Inexact approximation and the use of posterior mean}
\label{SLaplinexact}

Now we change the setup.
Namely, we suppose that the true maximizer \( \xvs \) of the function \( \lgd \) is not available, but 
\( \xv \) is somehow close to the point of maximum \( \xvs \).
Similarly, the negative Hessian \( \IFL = \IFL(\xvs) = - \nabla^{2} \lgd(\xvs) \) is hard to obtain 
and we use a proxy \( \IFLap \).
We already know that \( \PfL \) can be well approximated by \( \ND(\xvs,\IFL^{-1}) \).
This section addresses the question whether \( \ND(\xv,\IFLap^{-1}) \) can be used instead.
Here we may greatly benefit from the fact that Theorem~\ref{TLaplaceTV} provides a bound in the total variation distance 
between \( \PfL \) and \( \ND(\xvs,\IFL^{-1}) \) yielding
\begin{EQA}
	\TV\bigl( \PfL, \ND(\xv,\IFLap^{-1}) \bigr)
	& \leq &
	\TV\bigl( \PfL, \ND(\xvs,\IFL^{-1}) \bigr)
	+ \TV\bigl( \ND(\xv,\IFLap^{-1}), \ND(\xvs,\IFL^{-1}) \bigr) .
	\qquad
\label{iuhvfuggeyfgeyjhefhgdy}
\end{EQA}
Therefore, it suffices to bound the TV-distance between the Gaussian distribution 
\( \ND(\xvs,\IFL^{-1}) \) naturally arising in Laplace approximation, 
and the one used instead.
Pinsker's inequality provides an upper bound: for any two measures \( P,Q \)
\begin{EQA}
	\TV(P,Q)
	& \leq &
	\sqrt{\kullb(P,Q)/2} ,
\label{gffwdygweufeuyfeygfdwyd}
\end{EQA}
where \( \kullb(P,Q) \) is the Kullback-Leibler divergence between \( P \) and \( Q \).
The KL-divergence between two Gaussians has a closed form:
\begin{EQA}
	\kullb\bigl( \ND(\xvs,\IFL^{-1}), \ND(\xv,\IFLap^{-1}) \bigr)
	&=&
	\frac{1}{2} \bigl\{ \| \IFL^{1/2} (\xv - \xvs) \|^{2} + \tr (\IFLap^{-1} \IFL - \Id_{\dimp}) + \log \det (\IFLap^{-1} \IFL) \bigr\} .
\label{h0bede4bweve7yjdyeghy}
\end{EQA}
Moreover, if the matrix \( \BB = \IFLap^{-1/2} \, \IFL \, \IFLap^{-1/2} - \Id_{\dimp} \) satisfies 
\( \| \BB \| \leq 2/3 \) then
\begin{EQA}
	\TV\bigl( \ND(\xvs,\IFL^{-1}), \ND(\xv,\IFLap^{-1}) \bigr)
	& \leq &
	\frac{1}{2} \Bigl( \| \IFL^{1/2} (\xv - \xvs) \| + \sqrt{\tr \BB^{2}} \Bigr) .
\label{dferrfwbvf6nhdnghfkeiry}
\end{EQA}
However, Pinsker's inequality is only a general upper bound which is applied to any two distributions \( P \) and \( Q \).
If \( P \) and \( Q \) are Gaussian, it might be too rough.
Particularly the use of \( \tr \BB^{2} \) is disappointing, this quantity is full dimensional even 
if each of \( \IFL^{-1} \) and \( \IFLap^{-1} \) has a bounded trace.
Also dependence on \( \| \IFL^{1/2} (\xv - \xvs) \| \) is very discouraging.
\cite{Devroy2022} provides much sharper results, however, limited to the case of the same mean.
Even stronger results have been obtained in \cite{GNSUl2017} 
after restricting to the class \( \BBB_{el}(\R^{\dimp}) \) of elliptic sets \( A \) in \( \R^{\dimp} \) of the form
\begin{EQA}
	A
	&=&
	\bigl\{ \uv \in \R^{\dimp} \colon \| \QP (\uv - \xv) \| \leq \rr \bigr\}
\label{vhg4dfe5w3tfdf54wteg}
\end{EQA}
for some linear mapping \( \QP \colon \R^{\dimp} \to \R^{\dimq} \), \( \xv \in \R^{\dimp} \), and \( \rr > 0 \).
We refer to Section~\ref{SmainresGC} for exact formulations and further references.

\begin{theorem}
\label{TLaplaceTVin}
Assume the conditions of Theorem~\ref{TLaplaceTV} with \( \xvs \) being the maximizer of \( \lgd \) and 
\( \IFL = - \nabla^{2} \lgd(\xvs) \).
For any \( \xv \) and \( \IFLap \),
it holds with \( \gaussv_{\IFLap} \sim \ND(0,\IFLap^{-1}) \)
\begin{EQA}
	&& \nquad
	\sup_{A \in \BBB(\R^{\dimp})} \bigl| \PfL(\Xv - \xv \in A) - \P(\gaussv_{\IFLap} \in A) \bigr|\; 
	\\
	& \leq &
	4(\err + \ex^{-\xx}) + \TV\bigl( \ND(\xv,\IFLap^{-1}), \ND(\xvs,\IFL^{-1}) \bigr) ,
\label{nhjrwsdgdehgdftregwbdctsh}
\end{EQA}
where \( \err = \err_{2} \), see \eqref{juytr90f2dzaryjhfyf}, or \( \err = \err_{3} \), see \eqref{5qw7dyf4e4354coefw9dufih}.

Furthermore, for \( \Xv \sim \PfL \) and \( \gaussv \sim \ND(0,\Id_{\dimp}) \), any linear mapping 
\( \QP \colon \R^{\dimp} \to \R^{\dimq} \),
it holds 
under \( 3 \| \QP \, \IFL^{-1} \, \QP^{\T} \|^{2} \leq \| \QP \, \IFL^{-1} \, \QP^{\T} \|_{\Fr}^{2} \)
\begin{EQA}
	&& \nquad
	\sup_{\rr > 0}
	\left| \PfL\bigl( \| \QP (\Xv - \xv) \| \leq \rr \bigr) 
	- \P\bigl( \| \QP \, \IFLap^{-1/2} \gaussv \| \leq \rr \bigr) 
	\right|
	\\
	& \leq &
	4(\err_{3} + \ex^{-\xx}) + \frac{\CONST}{\| \QP \, \IFL^{-1} \, \QP^{\T} \|_{\Fr}}
	\left( 
		\| \QP (\IFL^{-1} - \IFLap^{-1}) \QP^{\T} \|_{1} + \| \QP (\xv - \xvs) \|^{2} 
	\right) .
\label{jrwguvyr23jbviufdsfsdgf6w}
\end{EQA}
\end{theorem}

\begin{proof}
The first bound follows from \eqref{iuhvfuggeyfgeyjhefhgdy} and Theorem~\ref{TLaplaceTV}.
For the second bound, we use Corollary~\ref{CTgaussiancomparisonS}.
\end{proof}

As a special case, consider the use of the posterior mean \( \xvb \) instead of \( \xvs \): 
\begin{EQA}
	\xvb
	& \eqdef &
	\frac{\int \xv \, \ex^{\lgd(\xv)} \, d\xv}{\int \ex^{\lgd(\xv)} \, d\xv} \, .
\label{054rfgofe3rdgbrfty6ry}
\end{EQA}

\begin{theorem}
\label{TpostmeanLan}
Assume the conditions of Theorem~\ref{TpostmeanLa} and Theorem~\ref{TLaplaceTVin}. 
Then it holds for any linear mapping \( \QP \colon \R^{\dimp} \to \R^{\dimq} \) with \( \QP^{\T} \QP \leq \DVL^{2} \)
\begin{EQA}
	\sup_{\rr > 0}
	\left| \PfL\bigl( \| \QP (\Xv - \xvb) \| \leq \rr \bigr) 
	- \P\bigl( \| \QP \gaussv_{\IFL} \| \leq \rr \bigr) 
	\right|
	& \leq & 
	4(\err_{3} + \ex^{-\xx}) + \frac{\CONST \| \QP (\xvb - \xvs) \|^{2}}{\| \QP \, \IFL^{-1} \, \QP^{\T} \|_{\Fr}} \, ,
	\qquad
	\quad
\label{jrwguvyr23jbviufdsfsdgf6wm}
\end{EQA}
where \( \| \QP (\xvb - \xvs) \| \) follows \eqref{hcdtrdtdehfdewdrfrhgyjufger} and \eqref{klu8gitfdgregfkhj7yt}.
\end{theorem}

The case \( \QP = \DVL \) is particularly transparent.
In view of \eqref{klu8gitfdgregfkhj7yt} of Corollary~\ref{CTpostmeanLa} and 
\( \| \QP \, \IFL^{-1} \, \QP^{\T} \|_{\Fr}^{2} = \tr\bigl\{ (\DVL \, \IFL^{-1} \, \DVL^{\T})^{2} \bigr\} \asymp \dimL \),
the following result holds.

\begin{corollary}
\label{CTpostmeanLab}
Under the conditions of Corollary~\ref{CTpostmeanLa}, it holds for \( \Xv \sim \PfL \)
\begin{EQA}
	\sup_{\rr > 0}
	\Bigl| \PfL\bigl( \| \DVL (\Xv - \xvb) \| \leq \rr \bigr) 
	- \P\bigl( \| \DVL \, \gaussv_{\IFL} \| \leq \rr \bigr) 
	\Bigr|
	& \leq &
	\CONST \biggl( \sqrt{{\dimL^{3}}/{n}} + \ex^{-\xx} \biggr) \, .
\label{klu8gitfdgregfkhj7yt3n}
\end{EQA}
\end{corollary}

The same bound applies with \( \QP = n^{1/2} \Id_{\dimp} \) in place of \( \DVL \) provided that 
\( \DVL^{2} \geq \CONSTi_{0} \, n \, \Id_{\dimp} \) for some fixed \( \CONSTi_{0} > 0 \).
We may conclude that the use of posterior mean \( \xvb \) in place of the posterior mode \( \xvs \) is justified 
under the same condition on critical dimension \( \dimL^{3} \ll n \) as required for the main result about Gaussian 
approximation.
}{}

{
\Section{Tools and proofs}
\label{SLapltools}
\renewcommand{\Section}[1]{\subsubsection{#1}}

Here we collect the proofs of the main results and some useful technical statements 
about the error of Laplace approximation.
Below we write \( \xv \) instead of \( \xvs \).
After passing to representation \eqref{IIgfifxutudufxutdu}, 
many results below apply to any \( \xv \), not necessarily for \( \xv = \xvs \).
We only use \( \DVL_{\GP}^{2} = \IFL = - \nabla^{2} \lgd(\xv) \) and 
\( \dltwb \) instead of \( \dltwb(\xv) \).
Everywhere we assume the local set \( \UVL \) to be fixed by \eqref{UvTDunm12spT}.
We separately study the integrals over \( \UVL \) and over its complement. 
The local error of approximation is measured by
\begin{EQA}
	\err
	=
	\err(\UVL)
	& \eqdef &
	\biggl| 
	\frac{\int_{\UVL} \ex^{\lgd(\xv;\uv)} \, g(\uv) \, d\uv - \int_{\UVL} \ex^{- \| \DVL_{\GP} \uv \|^{2}/2} \, g(\uv) \, d\uv} 
		 {\int \ex^{- \| \DVL_{\GP} \uv \|^{2}/2} d\uv} 
	\biggr| \, .
\label{errdefdiUaHu}
\end{EQA}
As a special case with \( g(\uv) \equiv 1 \) we obtain an approximation of the denominator in \eqref{IIgfifxutudufxutdu}.
In addition, we have to bound the tail integrals
\begin{EQ}[rcccl]
	\rho
	&=&
	\rho(\UVL)
	& \eqdef &
	\frac{\int \Ind(\uv \not\in \UVL) \, \ex^{\lgd(\xv;\uv)} \, d\uv}{\int \ex^{- \| \DVL_{\GP} \uv \|^{2}/2} \, d\uv} \, ,
	\\
	\rho_{\GP}
	&=&
	\rho_{\GP}(\UVL)
	& \eqdef &
	\frac{\int \Ind(\uv \not\in \UVL) \, \ex^{- \| \DVL_{\GP} \uv \|^{2}/2} \, d\uv}{\int \ex^{- \| \DVL_{\GP} \uv \|^{2}/2} \, d\uv} \, .
\label{rhfiUaceHu22m}
\end{EQ}
Everywhere later \( \gaussv_{\GP} \sim \ND(0,\DVL_{\GP}^{-2}) \) is a Gaussian element in \( \R^{\dimp} \).
The analysis will be split into several steps.

\Section{Overall error of Laplace approximation}
First we show how to seam together the error \( \err \) of local approximation and the bounds for the tail integrals
\( \rho \) and \( \rho_{\GP} \); see \eqref{rhfiUaceHu22m}.

\begin{proposition}
\label{PunbintLapl}
Suppose that for a function \( g(\uv) \in [0,1] \) and some \( \err , \err_{g} \)
\begin{EQA}[rcl]
	\biggl| 
		\frac{\int_{\UVL} \ex^{\lgd(\xv;\uv)} \, d\uv - \int_{\UVL} \ex^{- \| \DVL_{\GP} \uv \|^{2}/2} \, d\uv}			 
			 {\int \ex^{- \| \DVL_{\GP} \uv \|^{2}/2} \, d\uv} 
	\biggr|
	& \leq &
	\err \, ,
\label{erifxudieHu22}
	\\
\label{erifxudieHu22g}
	\biggl| 
		\frac{\int_{\UVL} g(\uv) \, \ex^{\lgd(\xv;\uv)} \, d\uv - \int_{\UVL} g(\uv) \, \ex^{- \| \DVL_{\GP} \uv \|^{2}/2} \, d\uv}			 {\int \ex^{- \| \DVL_{\GP} \uv \|^{2}/2} \, d\uv} 
	\biggr|
	& \leq &
	\err_{g} \, .
\end{EQA}
Then with \( \rho \) and \( \rho_{\GP} \) from \eqref{rhfiUaceHu22m}
\begin{EQ}[rcl]
	\frac{\int g(\uv) \, \ex^{\lgd(\xv;\uv)} \, d\uv}{\int \ex^{\lgd(\xv;\uv)} \, d\uv} 
	& \leq &
	\frac{1}{1 - \rho_{\GP} - \err} \,\, 
	\frac{\int g(\uv) \, \ex^{- \| \DVL_{\GP} \uv \|^{2}/2} \, d\uv}{\int \ex^{- \| \DVL_{\GP} \uv \|^{2}/2} \, d\uv}
	+ \frac{\rho + \err_{g}}{1 - \rho_{\GP} - \err} \, ,
	\qquad
	\\
	\frac{\int g(\uv) \, \ex^{\lgd(\xv;\uv)} \, d\uv}{\int \ex^{\lgd(\xv;\uv)} \, d\uv} 
	& \geq &
	\frac{1}{1 + \rho + \err} \,\, 
	\frac{\int g(\uv) \, \ex^{- \| \DVL_{\GP} \uv \|^{2}/2} \, d\uv}{\int \ex^{- \| \DVL_{\GP} \uv \|^{2}/2} \, d\uv}
	- \frac{\rho_{\GP} + \err_{g}}{1 + \rho + \err} \, .
	\qquad
\label{igefxumiguexHu22}
\end{EQ}
\end{proposition}

\begin{proof}
It follows from \eqref{erifxudieHu22} 
\begin{EQA}
	\int \ex^{\lgd(\xv;\uv)} \, d\uv
	& \geq &
	\int_{\UVL} \ex^{\lgd(\xv;\uv)} \, d\uv
	\geq 
	\int_{\UVL} \ex^{- \| \DVL_{\GP} \uv \|^{2}/2} \, d\uv - \err \int \ex^{- \| \DVL_{\GP} \uv \|^{2}/2} \, d\uv
	\\
	& \geq &
	(1 - \err - \rho_{\GP}) \int \ex^{- \| \DVL_{\GP} \uv \|^{2}/2} \, d\uv ,
\label{1erriexmH22du22}
	\\
	\int \ex^{\lgd(\xv;\uv)} \, d\uv
	& \leq &
	\int_{\UVL} \ex^{\lgd(\xv;\uv)} \, d\uv + \rho \int \ex^{- \| \DVL_{\GP} \uv \|^{2}/2} \, d\uv
	\\
	& \leq &
	(1 + \err + \rho) \int \ex^{- \| \DVL_{\GP} \uv \|^{2}/2} \, d\uv .
	\qquad
\label{1erriexmH22du22u}
\end{EQA}
Similarly for \( g(\uv) \geq 0 \)
\begin{EQA}
	\int g(\uv) \, \ex^{\lgd(\xv;\uv)} \, d\uv
	& \geq &
	\int_{\UVL} g(\uv) \, \ex^{- \| \DVL_{\GP} \uv \|^{2}/2} \, d\uv - \err_{g} \int \ex^{- \| \DVL_{\GP} \uv \|^{2}/2} \, d\uv
	\\
	& \geq &
	\int g(\uv) \, \ex^{- \| \DVL_{\GP} \uv \|^{2}/2} \, d\uv
	- (\rho_{\GP} + \err_{g}) \int \ex^{- \| \DVL_{\GP} \uv \|^{2}/2} \, d\uv ,
\label{jhdjsmswnqwyts5423nd}
\end{EQA}	
\begin{EQA}
	\int g(\uv) \, \ex^{\lgd(\xv;\uv)} \, d\uv
	& \leq &
	\int_{\UVL} g(\uv) \, \ex^{\lgd(\xv;\uv)} \, d\uv + \rho \int \ex^{- \| \DVL_{\GP} \uv \|^{2}/2} \, d\uv
	\\
	& \leq &
	\int g(\uv) \, \ex^{- \| \DVL_{\GP} \uv \|^{2}/2} \, d\uv + (\rho + \err_{g}) \int \ex^{- \| \DVL_{\GP} \uv \|^{2}/2} \, d\uv \, .
\label{iguexmHu22dreieH}
\end{EQA}
Putting together all these bounds yields \eqref{igefxumiguexHu22}.
\end{proof}

The next corollary is straightforward.

\begin{corollary}
\label{CPunbintLapl}
Let \( \rho_{\GP} \leq \rhos \), \( \rho \leq \rhos \); see \eqref{rhfiUaceHu22m}.
Let also for a function \( g(\uv) \) with \( |g(\uv)| \leq 1 \), \eqref{erifxudieHu22}, \eqref{erifxudieHu22g} hold with
\( \err_{g} \leq \err \).
If \( \err + \rhos \leq 1/2 \) then 
\begin{EQA}
	\left| 
		\frac{\int g(\uv) \, \ex^{\lgd(\xv;\uv)} \, d\uv}{\int \ex^{\lgd(\xv;\uv)} \, d\uv} 
		- \frac{\int g(\uv) \, \ex^{- \| \DVL_{\GP} \uv \|^{2}/2} \, d\uv}{\int \ex^{- \| \DVL_{\GP} \uv \|^{2}/2} \, d\uv}
	\right|
	& \leq &
	\frac{2 (\rhos + \err)}{1 - \rhos - \err} 
	\leq 
	4 (\rhos + \err)	\, . 
\label{2rherIHg2re}
\end{EQA}
\end{corollary}


\ifinexact{
Sometimes we need an extension to the case of an unbounded function \( g \).
This particularly arises when evaluating the moment of the posterior; see Theorem~\ref{TpostmeanLa}.
The next result corresponds to estimation of posterior mean with a linear function \( g \) 
and posterior variance with \( g \) quadratic.

\begin{proposition}
\label{PunbintLapg}
Given a function \( g(\uv) \), assume \eqref{erifxudieHu22}, \eqref{erifxudieHu22g}, and define
\begin{EQA}
	\rho_{g}
	& \eqdef &
	\frac{\int \Ind(\uv \not\in \UVL) \, | g(\uv) | \, \ex^{\lgd(\xv;\uv)} \, d\uv}{\int \ex^{- \| \DVL_{\GP} \uv \|^{2}/2} \, d\uv} 
	 \, ,
	\\
	\rho_{\GP,g}
	& \eqdef &
	\frac{\int \Ind(\uv \not\in \UVL) \, | g(\uv) | \, \ex^{- \| \DVL_{\GP} \uv \|^{2}/2} \, d\uv}{\int \ex^{- \| \DVL_{\GP} \uv \|^{2}/2} \, d\uv} \, ,
\label{rhfiUaceHu22mg}
\end{EQA}
while \( \rho \) and \( \rho_{\GP} \) are given in \eqref{rhfiUaceHu22m}.
Then for \( \II_{\GP}(g) = \E g(\gaussv_{\GP}) \), \( \gaussv_{\GP} \sim \ND(0,\DVL_{\GP}^{-2}) \), 
\begin{EQA}
	&& \nquad
	\left| 
		\frac{\int g(\uv) \, \ex^{\lgd(\xv;\uv)} \, d\uv}{\int \ex^{\lgd(\xv;\uv)} \, d\uv} 
		- \frac{\int g(\uv) \, \ex^{- \| \DVL_{\GP} \uv \|^{2}/2} \, d\uv}{\int \ex^{- \| \DVL_{\GP} \uv \|^{2}/2} \, d\uv}
	\right|
	\leq 
	\frac{\rho_{g} + \rho_{\GP,g} + \err_{g}}{1 - \rho_{\GP} - \err}
	+ \frac{| \II_{\GP}(g) | \, (\rho + \err)}{1 - \rho_{\GP} - \err} \, . 
	\qquad
\label{2rherIHg2red}
\end{EQA}
In particular, if \( \int g(\uv) \, \ex^{-\| \DVL_{\GP} \uv \|^{2}/2} \, d\uv = 0 \) then
\begin{EQA}[rcl]
	\left| \frac{\int g(\uv) \, \ex^{\lgd(\xv;\uv)} \, d\uv}{\int \ex^{\lgd(\xv;\uv)} \, d\uv} \right| 
	& \leq &
	\frac{\rho_{g} + \rho_{\GP,g} + \err_{g}}{1 - \rho_{\GP} - \err} \, .
\label{igefxumiguexHu22g}
\end{EQA}
\end{proposition}

\begin{proof}
Suppose that \( \II_{\GP}(g) \geq 0 \).
Then 
\begin{EQA}
	&& \nquad
	\left| 
		\frac{\int g(\uv) \, \ex^{\lgd(\xv;\uv)} \, d\uv}{\int \ex^{\lgd(\xv;\uv)} \, d\uv} 
		- \frac{\int g(\uv) \, \ex^{- \| \DVL_{\GP} \uv \|^{2}/2} \, d\uv}{\int \ex^{- \| \DVL_{\GP} \uv \|^{2}/2} \, d\uv}
	\right|
	\\
	& \leq &
	\left| 
		\frac{\int g(\uv) \, \ex^{\lgd(\xv;\uv)} \, d\uv}{\int \ex^{\lgd(\xv;\uv)} \, d\uv} 
		- \frac{\int g(\uv) \, \ex^{- \| \DVL_{\GP} \uv \|^{2}/2} \, d\uv}{\int \ex^{\lgd(\xv;\uv)} \, d\uv}
	\right|
	+ \II_{\GP}(g)
	\left| 
		\frac{\int \ex^{- \| \DVL_{\GP} \uv \|^{2}/2} \, d\uv}{\int \ex^{\lgd(\xv;\uv)} \, d\uv} - 1
	\right| .
\label{igefxvv22igh}
\end{EQA}
By definitions
\begin{EQA}
	&& \nquad
	\left| 
		\int g(\uv) \, \ex^{\lgd(\xv;\uv)} \, d\uv - \int g(\uv) \, \ex^{-\| \DVL_{\GP} \uv \|^{2}/2} \, d\uv
	\right|
	\\
	& \leq &
	\left| 
		\int_{\UVL} g(\uv) \, \ex^{\lgd(\xv;\uv)} \, d\uv - \int_{\UVL} g(\uv) \, \ex^{-\| \DVL_{\GP} \uv \|^{2}/2} \, d\uv
	\right|
	\\
	&&
	+ \, \left| 
		\int \Ind(\uv \not\in \UVL) \, g(\uv) \, \ex^{\lgd(\xv;\uv)} \, d\uv 
	\right|
	+ \left| 
		\int \Ind(\uv \not\in \UVL) \, g(\uv) \, \ex^{-\| \DVL_{\GP} \uv \|^{2}/2} \, d\uv
	\right|	
	\\
	& \leq &
	\bigl( \rho_{g} + \rho_{\GP,g} + \err_{g} \bigr) \int \ex^{-\| \DVL_{\GP} \uv \|^{2}/2} \, d\uv
\label{lrguegfxumHu22}
\end{EQA}
and the assertion follows by \eqref{1erriexmH22du22} and \eqref{1erriexmH22du22u}.
\end{proof}
}{}

\ifunivariate{}{
\Section{Local approximation. Univariate case}
This section presents a bound of local approximation for the univariate case.
Later we provide an independent study of the multivariate case.

Define \( \lgd(x;u) = \lgd(x+u) - \lgd(x) - f'(x) u \), 
\begin{EQA}
	\dltw_{3}(x,u)
	& \eqdef &
	\lgd(x+u) - \lgd(x) - f'(x) u - f''(x) \frac{u^{2}}{2} \, .
\label{dl3xufxufxufppxu2}
\end{EQA}

\begin{proposition}
\label{Lintapproxdu}
Let a function \( \lgd(u) \) be two times continuously differentiable and satisfy 
\( f''(x) = - \DU^{2} < 0 \).
If for some \( \zq > 0 \) 
\begin{EQA}
	\sup_{|u| \leq \zq} \frac{1}{\DU^{2} u^{2}/2}\bigl| \dltw_{3}(x,u) \bigr|
	& \leq &
	\dltwb  
	\leq 
	1/3 ,
\label{123stzf3tH2}
\end{EQA}
then for any function \( g(u) \) with \( |g(u)| \leq 1 \) 
\begin{EQA}
	\frac{\bigl| \int_{-\zq}^{\zq}  \ex^{\lgd(x;u)} \, g(u) \, du 
			- \int_{-\zq}^{\zq} \ex^{- \DU^{2} u^{2}/2} \, g(u) \, du \bigr|}
		 {\int \ex^{- \DU^{2} u^{2}/2} \, du}
	& \leq &
	\err_{2}
	\eqdef
	\frac{\dltwb}{2 (1 - \dltwb)} \, .
\label{HimzzeftgttmH3}
\end{EQA}
Moreover, if, in addition to \eqref{123stzf3tH2}, for some \( \fba^{(3)} = \fba^{(3)}(x) \)
\begin{EQA}
	\bigl| \dltw_{3}(x,u) \bigr|
	& \leq &
	\frac{1}{6} \bigl| \fba^{(3)} u^{3} \bigr|,
	\qquad
	|u| \leq \zq,
\label{123stzf3tH3uzq}
\end{EQA}
then
\begin{EQA}
	\frac{\bigl| \int_{-\zq}^{\zq}  \ex^{\lgd(x;u)} \, g(u) \, du 
			- \int_{-\zq}^{\zq} \ex^{- \DU^{2} u^{2}/2} \, g(u) \, du \bigr|}
		 {\int \ex^{- \DU^{2} u^{2}/2} \, du}
	& \leq &
	\err_{3}
	\eqdef
	\frac{0.7 \bigl| \fba^{(3)} \bigr|}{\DU^{3}} \, .
\label{HimzzeftgttmH1}
\end{EQA}
\end{proposition}

\begin{proof}
By \eqref{123stzf3tH2} and \( \dltwb \leq 1/3 \), it holds 
\begin{EQA}
	&& \nquad
	\frac{\DU}{\sqrt{2\pi}}
	\biggl| 
		\int_{-\zq}^{\zq}  \ex^{\lgd(x;u)} \, g(u) \, du - \int_{-\zq}^{\zq} \ex^{- \DU^{2} u^{2}/2} \, g(u) \, du 
	\biggr|
	\\
	& \leq &
	\frac{\DU}{\sqrt{2\pi}}
	\biggl| 
		\int_{-\zq}^{\zq}  \ex^{- (1 - \dltwb) \DU^{2} u^{2}/2} \, du 
		- \int_{-\zq}^{\zq} \ex^{- \DU^{2} u^{2}/2} \, du 
	\biggr|
	\\
	& \leq &
	\frac{1}{\sqrt{2\pi}}
	\biggl| 
		\int  \ex^{- (1 - \dltwb) u^{2}/2} \, du 
		- \int \ex^{- u^{2}/2} \, du 
	\biggr|
	\leq 
	(1 - \dltwb)^{-1/2} - 1 .
\label{Himzzeftgt}
\end{EQA}
We now use that \( (1 - \dltwb)^{-1/2} - 1 \leq  0.5 \dltwb/(1 - \dltwb) \) for \( \dltwb < 1 \),
and the result \eqref{HimzzeftgttmH3} follows.

Now we show \eqref{HimzzeftgttmH1}. It holds
\begin{EQA}
	\int_{-\zq}^{\zq} \ex^{\lgd(x;u)} \, g(u) \, du
	& = &
	\int_{-\zq}^{\zq} \ex^{- \DU^{2} u^{2}/2 + \dltw_{3}(x,u)} \, g(u) \, du .
\label{iUefgu22guu}
\end{EQA}
Define for \( t \geq 0 \)
\begin{EQA}
	\Rem(t)
	&=&
	\int_{-\zq}^{\zq} \ex^{- \DU^{2} u^{2}/2 + t \dltw_{3}(x,u)} \, g(u) \, du .
\label{PhtiUexmHugu}
\end{EQA}
Then for \( t \in [0,1] \) by \eqref{123stzf3tH2} and \eqref{123stzf3tH3uzq}
\begin{EQA}
	\bigl| \Rem'(t) \bigr|
	&=&
	\left| \int_{-\zq}^{\zq} \dltw_{3}(x,u) \ex^{- \DU^{2} u^{2}/2 + t \dltw_{3}(x,u)} \, g(u) \, du \right|
	\leq 
	\int_{-\zq}^{\zq} \bigl| \dltw_{3}(x,u) \bigr| \ex^{- (1 - \dltwb) \DU^{2} u^{2}/2} \, du 
	\\
	& \leq &
	\frac{1}{6} \int_{-\zq}^{\zq} \bigl| \fba^{(3)} u^{3} \bigr| \ex^{- (1 - \dltwb) \DU^{2} u^{2}/2} \, du
	\leq 
	\frac{\bigl| \fba^{(3)} \bigr| \, \DU^{-4}}{6 (1 - \dltwb)^{2}} \int |u|^{3} \ex^{- u^{2}/2} \, du
	= 
	\frac{2\bigl| \fba^{(3)} \bigr| \, \DU^{-4}}{3 (1 - \dltwb)^{2}}
\label{Phit12u2H02u32}
\end{EQA}
and it holds
\begin{EQA}
	\bigl| \Rem(1) - \Rem(0) \bigr|
	& \leq &
	\sup_{t \in [0,1]} \bigl| \Rem'(t) \bigr|
	\leq 
	\frac{2\bigl| \fba^{(3)} \bigr| \, \DU^{-4}}{3 (1 - \dltwb)^{2}} \, 
\label{Pd2d3H0wem22}
\end{EQA}
yielding for \( \dltwb \leq 1/3 \)
\begin{EQA}
	\frac{\bigl| \Rem(1) - \Rem(0) \bigr|}{\int \ex^{- \DU^{2} u^{2}/2} du}
	& \leq &
	\frac{2\bigl| \fba^{(3)} \bigr| \, \DU^{-4}}{3 (1 - \dltwb)^{2} \sqrt{2\pi} \, \DU^{-1}}
	\leq 
	\frac{0.27 \bigl| \fba^{(3)} \bigr|}{(1 - \dltwb)^{2} \DU^{3}} \, 
\label{Phd3Ph0iUmHu22}
\end{EQA}
and \eqref{HimzzeftgttmH1} follows.
\end{proof}

Now we present some results based on a higher order Taylor approximation.
Define 
\begin{EQA}
	\dltw_{4}(x,u)
	& \eqdef &
	\lgd(x+u) - \lgd(x) - f'(x) u + \frac{\DU^{2} u^{2}}{2} - \frac{f^{(3)}(x) u^{3}}{6} 
	\\
	&=&
	\lgd(x,u) + \frac{\DU^{2} u^{2}}{2} - \frac{f^{(3)}(x) u^{3}}{6} \, .
\label{124stzf4tH24}
\end{EQA}
With \( \dltw_{3}(x,u) = \lgd(x,u) + \DU^{2} u^{2}/2 \), it obviously holds
\begin{EQA}
	\dltw_{3}(x,u)
	& = &
	\dltw_{4}(x,u) + {f^{(3)}(x) u^{3}}/{6} .
\label{d3d4uf30u36}
\end{EQA}

\begin{proposition}
\label{Lintapproxdu4}
Let \( \lgd(u) \) be three times continuously differentiable and satisfy \eqref{123stzf3tH2}.
Define \( \dltwm_{3} = f^{(3)}(x) / \DU^{2} \).
If \( \dltw_{4}(x,u) \) from \eqref{124stzf4tH24} satisfies
\begin{EQA}
	\sup_{|u| \leq \zq} \frac{24}{\DU^{2} u^{4}} \bigl| \dltw_{4}(x,u) \bigr|
	& \leq &
	\dltwm_{4}  
\label{124stzf4tH2}
\end{EQA}
then for any function \( g(u) \) with \( |g(u)| \leq 1 \) and \( g(u) = g(-u) \)
\begin{EQA}
	\frac{
	\bigl| \int_{-\zq}^{\zq}  \ex^{\lgd(x;u)} \, g(u) \, du - \int_{-\zq}^{\zq} \ex^{- \DU^{2} u^{2}/2} \, g(u) \, du \bigr|}
	{\int \ex^{- \DU^{2} u^{2}/2} \, du}
	& \leq &
	\err_{4} \, .
\label{HimzzeftgttmH2}
\end{EQA}
with 
\begin{EQA}
	\err_{4}
	& \eqdef &
	\left\{ \frac{5 \dltwm_{3}^{2}}{12 (1 - \dltwb)^{7/2} }
	+ \frac{\dltwm_{4}}{4 (1 - \dltwb)^{5/2}} \right\} \DU^{-2} .
\label{e4l1d372441m3}
\end{EQA}
\end{proposition}

\begin{proof}
We write \( \dltw_{3}(u) \) in place of \( \dltw_{3}(x,u) \) and \( f^{(3)} \) in place of \( f^{(3)}(x) \).
It holds 
\begin{EQA}
	\int_{-\zq}^{\zq} \ex^{\lgd(x;u)} \, g(u) \, du
	& = &
	\int_{-\zq}^{\zq} \exp\Bigl\{ - \frac{\DU^{2} u^{2}}{2} + \dltw_{3}(u) \Bigr\} \, g(u) \, du .
\label{P3tHmzzmt22t3}
\end{EQA}
Define for \( t \in [0,1] \)
\begin{EQA}
	\Rem(t)
	& \eqdef &
	\int_{-\zq}^{\zq} \exp\Bigl\{ - \frac{\DU^{2} u^{2}}{2} + t \dltw_{3}(u) \Bigr\} \, g(u) \, du .
\label{P3tHmzzmt22t}
\end{EQA}
Symmetricity \( g(u) = g(-u) \) implies that 
\begin{EQA}
	\Rem'(0)
	&=&
	\frac{1}{2} \int_{-\zq}^{\zq} 
	\exp\Bigl( - \frac{\DU^{2} u^{2}}{2} \Bigr) \bigl\{ \dltw_{3}(u) + \dltw_{3}(-u) \bigr\} \, g(u) \, du
	=
	\int_{-\zq}^{\zq} \exp\Bigl( - \frac{\DU^{2} u^{2}}{2} \Bigr) \dltw_{4}(u) \, g(u) \, du \, .
\label{Pp012zzmH2u212}
\end{EQA}
Moreover, as  \( |\dltw_{3}(u)| \leq \dltwb \, \DU^{2} u^{2}/2 \),
it holds for \( t \in [0,1] \) 
\begin{EQA}
	|\Rem''(t)|
	&=&
	\left| \int_{-\zq}^{\zq} \dltw_{3}^{2}(u) 
		\exp\Bigl\{ - \frac{\DU^{2} u^{2}}{2} + t \dltw_{3}(u) \Bigr\} \, g(u) \, du 
	\right|
	\leq 
	\int_{-\zq}^{\zq} \dltw_{3}^{2}(u) 
		\exp\Bigl( - \frac{(1 - \dltwb) \DU^{2} u^{2}}{2} \Bigr) \, du \, .
\label{F3pptHm2t6}
\end{EQA}
As \( \dltw_{3}(u) = f^{(3)} u^{3}/6 + \dltw_{4}(u) \) and \( |\dltw_{4}(u)| \leq 1 \), it holds 
\begin{EQA}
	|\Rem''(t)| 
	& \leq & 
	2 \int_{-\zq}^{\zq} \bigl\{ \dltw_{4}^{2}(u) + \bigl| f^{(3)} u^{3}/6 \bigr|^{2} \bigr\} 
		\exp\Bigl( - \frac{(1 - \dltwb)\DU^{2} u^{2}}{2} \Bigr) \, du
	\\
	& \leq & 
	2 \int_{-\zq}^{\zq} \bigl\{ |\dltw_{4}(u)| + \bigl| f^{(3)} \bigr|^{2} u^{6}/36 \bigr\} 
	\exp\Bigl( - \frac{(1 - \dltwb) \DU^{2} u^{2}}{2} \Bigr) \, du .
\label{u221md3u636}
\end{EQA}
The use of \( \bigl| f^{(3)} \bigr| = \dltwm_{3} \, \DU^{2} \) and \eqref{124stzf4tH2} yields
in view of \( \dltwb \leq 1/3 \) and \( |\dltw_{4}(u)| \leq \dltwm_{4} \, \DU^{2} u^{4} \)
\begin{EQA}
	&& \nquad
	\bigl| \Rem(1) - \Rem(0) \bigr|
	\leq 
	\bigl| \Rem'(0) \bigr| + \frac{1}{2} \sup_{t \leq 1} |\Rem''(t)|
	\\
	& \leq &
	2 \int_{-\zq}^{\zq} |\dltw_{4}(u)| \, \exp\Bigl( - \frac{(1 - \dltwb) \DU^{2} u^{2}}{2} \Bigr) \, du
	+
	\frac{\dltwm_{3}^{2} \, \DU^{4}}{36} \int_{-\zq}^{\zq} u^{6} \, 
		\exp\Bigl( - \frac{(1 - \dltwb)\DU^{2} u^{2}}{2} \Bigr) \, du 
	\\
	& \leq &
	\frac{\dltwm_{4} \, \DU^{2}}{12} 
	\int_{-\zq}^{\zq} u^{4} \, \exp\Bigl( - \frac{(1 - \dltwb) \DU^{2} u^{2}}{2} \Bigr) \, du
	+ \frac{\dltwm_{3}^{2} \, \DU^{4}}{36} 
	\int_{-\zq}^{\zq} u^{6} \, \exp\Bigl( - \frac{(1 - \dltwb)\DU^{2} u^{2}}{2} \Bigr) \, du 
	\\
	& \leq &
	\frac{\dltwm_{4} (1 - \dltwb)^{-5/2}}{12 \, \DU^{3}}	
	\int u^{4} \, \exp\Bigl( - \frac{\DU^{2} u^{2}}{2} \Bigr) \, du
	+ \frac{\dltwm_{3}^{2} (1 - \dltwb)^{-7/2}}{36 \, \DU^{3}} 
	\int u^{6} \, \exp\Bigl( - \frac{\DU^{2} u^{2}}{2} \Bigr) \, du 
	\\
	&=&
	\frac{3 \sqrt{2\pi} \, \dltwm_{4} (1 - \dltwb)^{-5/2} }{12 \, \DU^{3}}	
	+ \frac{15 \sqrt{2\pi} \, \dltwm_{3}^{2} (1 - \dltwb)^{-7/2}}{36 \, \DU^{3}}
\label{C4Hm2Fm4tH2t4}
\end{EQA}
yielding \eqref{HimzzeftgttmH2} in view of 
\( \int \ex^{- \DU^{2} u^{2}/2} \, du = \sqrt{2\pi} \, \DU^{-1} \)
and \( \dltwm_{3} \leq 1/3 \).
\end{proof}

\begin{remark}
Suppose that \( \lgd(\cdot) \) is 
three times continuously differentiable and \( |f^{(3)}(x+u)| \leq \dltwu_{3} \DU^{2} \)
for \( |u| \leq \zq \).
Then the expansion \eqref{123stzf3tH2} applies with \( \dltwb \leq \dltwu_{3} \zq/3 \).
The use of \( \zq = \xx/\DU \) yields
\begin{EQA}
	\dltwb \leq \dltwu_{3} \zq/ 3
	& \leq &
	\CONST \, \dltwu_{3} \, \xx \, \DU^{-1} .
\label{Cw5x5Hm3120}
\end{EQA}
Similarly, if \( \lgd(\cdot) \) is four times continuously differentiable and \( |f^{(4)}(u)| \leq \dltwm_{4} \DU^{2} \)
for \( |u| \leq \zq \), then \eqref{124stzf4tH2} holds as well.
One can see that the use of a higher order approximation of \( \lgd \) allows to improve 
the accuracy of approximation  
from \( \DU^{-1} \) as in \eqref{HimzzeftgttmH1} of Proposition~\ref{Lintapproxdu} to \( \DU^{-2} \) as in \eqref{HimzzeftgttmH2} of Proposition~\ref{Lintapproxdu4}.

\end{remark}
}

\Section{Lower and upper Gaussian measures}
This section introduces the lower and upper Gaussian measure which locally sandwich the measure \( \PfL \)
using the decomposition from condition \nameref{LLf0ref}. 
Denote \( - \nabla^{2} \lgd(\xv) = \DVL_{\GP}^{2} \).
Definition \eqref{om3esuU1H02d3} enables us to bound with \( \dltwb = \dltwb(\xv) \)
\begin{EQA}
	\frac{1}{2} (\| \DVL_{\GP} \uv \|^{2} - \dltwb \| \DVL \uv \|^{2})
	& \leq &
	\lgd(\xv;\uv)
	\leq 
	\frac{1}{2} (\| \DVL_{\GP} \uv \|^{2} + \dltwb \| \DVL \uv \|^{2}) \, 
\label{ysdtydsfrtswedftedswfhd}
\end{EQA}
yielding two Gaussian measures which bounds \( \PfL \) locally from above and from below.
The next technical result provides sufficient conditions for their contiguity.

\begin{proposition}
\label{LlocintgrL}
Let \( \dltwb \) from \eqref{om3esuU1H02d3} satisfy \( \dltwb \leq 1/3 \).
Then with \( \dimL \) from \eqref{dAdetrH02Hm2}
\begin{EQA}
\label{12Ip2t3t4Hp3}
	\det\bigl( \Id + \dltwb \DVL_{\GP}^{-1} \DVL^{2} \, \DVL_{\GP}^{-1} \bigr)
	& \leq &
	\exp (\dltwb \, \dimL) \, ,
	\\
	\det\bigl( \Id - \dltwb \DVL_{\GP}^{-1} \DVL^{2} \, \DVL_{\GP}^{-1} \bigr)^{-1/2}
	& \leq &
	\exp \bigl\{ 3/2 \log (3/2) \, \dltwb \, \dimL \bigr\} .
	\qquad
\label{12Im2t3t4Hm30}
\end{EQA}
\end{proposition}

\begin{proof}
W.l.o.g. assume that \( \DVL_{\GP}^{-1} \DVL^{2} \, \DVL_{\GP}^{-1} \) is diagonal with 
eigenvalues \( \lambda_{j} \in [0,1] \).
As \( - x^{-1}\log(1 - x) \leq 3 \log(3/2) \) for \( x \in [0,1/3] \), it holds by \eqref{om3esuU1H02d3}
\begin{EQA}
	&& \nquad
	\log \det\bigl( \Id - \dltwb \, \DVL_{\GP}^{-1} \DVL^{2} \, \DVL_{\GP}^{-1} \bigr)^{-1}
	=
	- \sum_{j=1}^{\dimp} \log\bigl( 1 - \dltwb \lambda_{j} \bigr)
	\leq 
	3 \log (3/2) \sum_{j=1}^{\dimp} \dltwb \lambda_{j}
	\\
	&=&
	3 \log (3/2) \, \dltwb \, \tr\bigl( \DVL_{\GP}^{-1} \DVL^{2} \, \DVL_{\GP}^{-1} \bigr)
	=
	3 \log (3/2) \, \dltwb \, \dimL \, 
\label{125d3p0ldm1}
\end{EQA}
yielding \eqref{12Im2t3t4Hm30}.
The proof of \eqref{12Ip2t3t4Hp3} is similar using \( \log(1 + x) \leq x \) for \( x \geq 0 \).
\end{proof}

\Subsection{Gaussian moments}
The presented bounds involve the Gaussian moments \( \E \| \DVL \gaussv_{\GP} \|^{k} \) for \( k=3,4,6 \)
and \( \gaussv_{\GP} \sim \ND(0,\DVL_{\GP}^{-2}) \).
We make use of the following lemma.

\begin{lemma}
\label{LdltwLa}
It holds for \( \gaussv_{\GP} \sim \ND(0,\DVL_{\GP}^{-2}) \) and \( \normG = \| \DVL \, \DVL_{\GP}^{-2} \DVL \| \)
\begin{EQA}
	\E \| \DVL \, \gaussv_{\GP} \|^{3}
	& \leq &
	(\dimL + \normG)^{3/2} \, ,
	\\
	\E \| \DVL \, \gaussv_{\GP} \|^{4}
	& \leq &
	(\dimL + \normG)^{2} \, ,
	\\
	\E \| \DVL \, \gaussv_{\GP} \|^{6}
	& \leq &
	(\dimL + 2 \normG)^{3} \, ,
	\\
	\E \| \DVL \, \gaussv_{\GP} \|^{8}
	& \leq &
	(\dimL + 3 \normG)^{4} \, .
\label{eciu8hef8h8we3hy87u3y7y3fe7}
\end{EQA}
\end{lemma}

\begin{proof}
Represent \( \| \DVL \, \gaussv_{\GP} \|^{2} = \| \DVL \, \DVL_{\GP}^{-1} \gaussv \|^{2} 
= \langle \BB_{\GP} \gaussv, \gaussv \rangle \)
with \( \BB_{\GP} = \DVL \, \DVL_{\GP}^{-2} \DVL \leq \Id_{\dimp} \) and \( \gaussv \sim \ND(0,\Id_{\dimp}) \).
By Lemma~\ref{Gaussmoments}, for \( m=1,2,3 \)
\begin{EQA}
	\E \| \DVL \, \gaussv_{\GP} \|^{2m+2}
	&=&
	\E \bigl\langle \BB_{\GP} \gaussv, \gaussv \bigr\rangle^{m+1}
	\leq 
	\bigl\{ \tr (\BB_{\GP}) + m \, \normG \bigr\}^{m+1}
	=
	(\dimL + m \, \normG)^{m+1} ,
\label{lkvu7rerycfrrwt3ghdtndf}
\end{EQA}
and \( \E \| \DVL \, \gaussv_{\GP} \|^{3} \leq \E^{3/4} \| \DVL \, \gaussv_{\GP} \|^{4}	\leq (\dimL + \normG)^{3/2} \).
\end{proof}

\Section{Local approximation}
\label{SLaplapprmu}

This section presents the bounds on the error \( \err \) of local approximation \eqref{errdefdiUaHu}.
The first result only uses \( \dltwb \, \dimL \leq 2/3 \).
More advanced bounds also assume \nameref{LL3tref} and \nameref{LL4tref} with \( \upsv = \xv \). 
We also present some extensions for the moments of \( \PfL \).

\begin{proposition}
\label{Lintfxupp3}
Let \( \dltwb = \dltwb(\xv) \) from \eqref{om3esuU1H02d3} and \( \dimL \) from \eqref{dAdetrH02Hm2} satisfy 
\begin{EQA}
	\dltwb \, \dimL
	& \leq &
	2/3 \, .
\label{w3p0le13fx}
\end{EQA}
Then for any function \( g(\uv) \) with \( |g(\uv)| \leq 1 \) 
\begin{EQA}
	\biggl| 
	\frac{\int_{\UVL} \ex^{\lgd(\xv;\uv)} \, g(\uv) \, d\uv - \int_{\UVL} \ex^{- \| \DVL_{\GP} \uv \|^{2}/2} \, g(\uv) \, d\uv} 
		 {\int \ex^{- \| \DVL_{\GP} \uv \|^{2}/2} d\uv} 
	\biggr|
	& \leq &
	\err \, 
\label{ed3le2d3pGd}
\end{EQA}
with
\begin{EQA}
	\err
	&=&
	\err_{2} 
	= 
	\frac{0.75 \, \dltwb \, \dimL}{1 - \dltwb} \, .
\label{8fjre34et6dehgweyt6wsh}
\end{EQA}
\end{proposition}

\begin{proof}
The condition \( \dltwb \, \dimL \leq 2/3 \) from \eqref{w3p0le13fx} and bound \eqref{12Im2t3t4Hm30} imply 
\begin{EQA}
	\det\bigl( \Id - \dltwb \DVL_{\GP}^{-1} \DVL^{2} \, \DVL_{\GP}^{-1} \bigr)^{-1/2}
	& \leq &
	\exp \bigl\{ 3/2 \log (3/2) \, \dltwb \, \dimL \bigr\} 
	\leq 
	3/2 \, .
	\qquad
\label{12Im2t3t4Hm3}
\end{EQA}
Define for \( t \geq 0 \)
\begin{EQA}
	\Rem(t)
	&=&
	\int_{\UVL} \ex^{- \| \DVL_{\GP} \uv \|^{2}/2 + t \dltw_{3}(\xv,\uv)} \, g(\uv) \, d\uv .
\label{PhtiUexmHu2mtH0u2g}
\end{EQA}
Then for \( t \in [0,1] \) by \eqref{om3esuU1H02d3} 
\begin{EQA}
	\bigl| \Rem'(t) \bigr|
	&=&
	\left| \int_{\UVL} \dltw_{3}(\xv,\uv) \ex^{- \| \DVL_{\GP} \uv \|^{2}/2 + t \dltw_{3}(\xv,\uv)} \, g(\uv) \, d\uv \right|
	\\
	& \leq &
	\int_{\UVL} \bigl| \dltw_{3}(\xv,\uv) \bigr| \ex^{- (\| \DVL_{\GP} \uv \|^{2} - \dltwb \| \DVL \uv \|^{2})/2} \, d\uv .
\label{Pht12u212H02u22}
\end{EQA}
Now we make change of variable \( \Idd \uv \) to \( \uv \) with
\( \Idd^{2} = \Id - \dltwb \, \DVL_{\GP}^{-1} \DVL^{2} \, \DVL_{\GP}^{-1} \).
By \eqref{12Im2t3t4Hm3} \( \det \Idd^{-1} \leq 3/2 \) and also 
\( \| \Idd^{-1} \| \leq (1 - \dltwb)^{-1/2} \).
By \eqref{om3esuU1H02d3} and \eqref{Pht12u212H02u22}
\begin{EQA}
	&& \nquad
	\bigl| \Rem(1) - \Rem(0) \bigr|
	\leq 
	\sup_{t \in [0,1]} \bigl| \Rem'(t) \bigr|
	\leq 
	\frac{\dltwb}{2}
	\int_{\UVL} \| \DVL \uv \|^{2} \ex^{- (\| \DVL_{\GP} \uv \|^{2} - \dltwb \| \DVL \uv \|^{2})/2} \, d\uv 
	\\
	& \leq &
	\frac{3 \dltwb}{4}
	\int \| \DVL \Idd^{-1} \uv \|^{2} \ex^{- \| \DVL_{\GP} \uv \|^{2} /2} \, d\uv 
	\leq 
	\frac{3 \dltwb}{4(1 - \dltwb)}
	\int \| \DVL \uv \|^{2} \ex^{- \| \DVL_{\GP} \uv \|^{2} /2} \, d\uv 	.
\label{ppt2iUH0w2edw}
\end{EQA}
In view of \( \E \| \DVL \, \gaussv_{\GP} \|^{2} = \tr\bigl( \DVL^{2} \, \DVL_{\GP}^{-2} \bigr) \) for a standard normal
\( \gaussv \), we derive 
\begin{EQA}
	\frac{\bigl| \Rem(1) - \Rem(0) \bigr|}{\int \ex^{- \| \DVL_{\GP} \uv \|^{2}/2} d\uv}
	& \leq &
	\frac{3 \dltwb }{4 (1 - \dltwb)}
	\frac{\int \| \DVL \uv \|^{2} \ex^{- \| \DVL_{\GP} \uv \|^{2} /2} \, d\uv}
		{\int \ex^{- \| \DVL_{\GP} \uv \|^{2} /2} \, d\uv}
	\leq 
	\frac{3 \dltwb \, \dimL}{4 (1 - \dltwb)}
\label{Phd3Ph0iUmHu22}
\end{EQA}
and \eqref{ed3le2d3pGd} follows.
\end{proof}

\begin{proposition}
\label{Lintfxupp3T}
Assume \nameref{LL3tref} with \( \upsv = \xv \) 
and let \( \dltwu_{3} \, \amax^{-1} \rrL \, \dimL \leq 2 \). 
Define \( \dltwb \eqdef \dltwu_{3} \, \amax^{-1} \rrL / 3 \).
Then bound \eqref{ed3le2d3pGd} applies with 
\begin{EQA}
	\err
	&=&
	\err_{3}
	\eqdef
	\frac{\dltwu_{3} \E \| \DVL \gaussv_{\GP} \|^{3}}{4 (1 - \dltwb)^{3/2}} 
	\leq 
	\frac{\dltwu_{3} \, (\dimL + \normG)^{3/2}}{4 (1 - \dltwb)^{3/2}} \, .
\label{5qw7dyf4e4354co9dufih}
\end{EQA}
\end{proposition}

\begin{proof}
The proof follows the same line as for Proposition~\ref{Lintfxupp3}.
Under \nameref{LL3tref}, it holds \( |\dltw_{3}(\xv,\uv)| \leq \dltwu_{3} \| \DVL \uv \|^{3} / 6 \) for \( \uv \in \UVL \) and
\begin{EQA}
	\bigl| \Rem(1) - \Rem(0) \bigr|
	& \leq &
	\frac{\dltwu_{3} \, \det (\Idd^{-1})}{6} 
	\int \| \DVL \, \Idd^{-1} \uv \|^{3} \, \ex^{- \| \DVL_{\GP} \uv \|^{2} /2} \, d\uv 
	\\
	& \leq &
	\frac{\dltwu_{3}}{4(1 - \dltwb)^{3/2}}
	\int \| \DVL \uv \|^{3} \, \ex^{- \| \DVL_{\GP} \uv \|^{2} /2} \, d\uv 
\label{Pd3p2d3H0wem22g1}
\end{EQA}
yielding the statement in view of \( \E \| \DVL \, \gaussv_{\GP} \|^{3} \leq (\dimL + \normG)^{3/2} \); see Lemma~\ref{LdltwLa}.
It remains to note that 
by Lemma~\ref{LdltwLa3t} \( \dltwb \leq \dltwu_{3} \, \amax^{-1} \rrL \, /3\).
Hence, 
\( \dltwu_{3} \, \amax^{-1} \rrL \, \dimL \leq 2 \) implies \( \dltwb \dimL \leq 2/3 \).
\end{proof}

\noindent
The result can be extended to the case of a \( m \)-homogeneous function \( g(\uv) \).

\begin{proposition}
\label{Pd3t3hot3d3u}
Suppose the conditions of Proposition~\ref{Lintfxupp3} and \nameref{LL3tref}.
Then for \( m \geq 1 \) and any \( m \)-homogeneous function \( g(\cdot) \) with \( g(t\uv) = t^{m} g(\uv) \)
\begin{EQA}
	\biggl| 
	\frac{\int_{\UVL} \ex^{\lgd(\xv;\uv)} \, g(\uv) \, d\uv 
	- \int_{\UVL} \ex^{- \| \DVL_{\GP} \uv \|^{2}/2} \, g(\uv) \, d\uv} 
		 {\int \ex^{- \| \DVL_{\GP} \uv \|^{2}/2} d\uv} 
	\biggr|
	& \leq &
	\frac{\dltwu_{3} \E \bigl\{ |g(\gaussv_{\GP})| \, \| \DVL \gaussv_{\GP} \|^{3} \bigr\}}{4 (1 - \dltwb)^{(m+3)/2}} \, .
	\qquad
\label{e3dEf3xGHo3}
\end{EQA}
\end{proposition}

\begin{proof}
As in the proof of Proposition~\ref{Lintfxupp3}, under \nameref{LL3tref}, it holds for \( \Rem(t) \) from
\eqref{PhtiUexmHu2mtH0u2g}
\begin{EQA}
	\bigl| \Rem(1) - \Rem(0) \bigr|
	& \leq &
	\frac{\dltwu_{3} \, \det (\Idd^{-1})}{6} 
	\int \| \DVL \, \Idd^{-1} \uv \|^{3} \, \bigl| g(\Idd^{-1} \uv) \bigr| \, \ex^{- \| \DVL_{\GP} \uv \|^{2} /2} \, d\uv 
	\\
	& \leq &
	\frac{\dltwu_{3}}{4(1 - \dltwb)^{(m+3)/2}}
	\int \| \DVL \uv \|^{3} \, \bigl| g(\uv) \bigr| \, \ex^{- \| \DVL_{\GP} \uv \|^{2} /2} \, d\uv 
\label{Pd3p2d3H0wem22}
\end{EQA}
yielding similarly to \eqref{Pd3p2d3H0wem22g1}
\begin{EQA}
	\frac{\bigl| \Rem(1) - \Rem(0) \bigr|}{\int \ex^{- \| \DVL_{\GP} \uv \|^{2}/2} d\uv}
	& \leq &
	\frac{\dltwu_{3}}{4(1 - \dltwb)^{(m+3)/2}}
	\frac{\int \| \DVL \uv \|^{3} \, \bigl| g(\uv) \bigr| \, \ex^{- \| \DVL_{\GP} \uv \|^{2} /2} \, d\uv}
		{\int \ex^{- \| \DVL_{\GP} \uv \|^{2} /2} \, d\uv} \, .
\label{Phd3Ph0iUmHu22}
\end{EQA}
This yields \eqref{e3dEf3xGHo3}.
\end{proof}

Important special cases correspond to \( m=1 \).

\begin{proposition}
\label{Perrbound3}
Suppose the conditions of Proposition~\ref{Pd3t3hot3d3u}.
Then it holds 
for any linear mapping \( \QP \colon \R^{\dimp} \to \R^{\dimq} \) and any unit vector \( \av \in \R^{\dimq} \)
\begin{EQA}
	&& \nquad
	\frac{\bigl| 
		\int_{\UVL} \ex^{\lgd(\xv;\uv)} \, \langle \QP \uv,\av \rangle \, d\uv 
		- \int_{\UVL} \ex^{- \| \DVL_{\GP} \uv \|^{2}/2} \, \langle \QP \uv,\av \rangle \, d\uv 	 
		\bigr|} 
		 {\int \ex^{- \| \DVL_{\GP} \uv \|^{2}/2} d\uv} 
	\\
	& \leq &
	0.6 \, \dltwu_{3} \, (\dimL + \normG)^{3/2} \, \| \QP \, \DVL_{\GP}^{-2} \QP^{\T} \|^{1/2} .
	\qquad
\label{iUafxvdvH2vk3}
\end{EQA}
\end{proposition}

\begin{proof}
Proposition~\ref{Pd3t3hot3d3u} yields
\begin{EQA}
	\frac{\bigl| 
		\int_{\UVL} \ex^{\lgd(\xv;\uv)} \, \langle \QP \uv,\av \rangle \, d\uv 
		- \int_{\UVL} \ex^{- \| \DVL_{\GP} \uv \|^{2}/2} \, \langle \QP \uv,\av \rangle \, d\uv 	 
		\bigr|} 
		 {\int \ex^{- \| \DVL_{\GP} \uv \|^{2}/2} d\uv} 
	& \leq &
	\frac{\dltwu_{3} \, \E \bigl\{ \bigl| \langle \QP \gaussv_{\GP},\av \rangle \bigr| \, \| \DP \gaussv_{\GP} \|^{3} \bigr\}}
		 {4 (1 - \dltwb)^{2}} \, .
	\qquad 
\label{iUafxvdvH2vk3pr}
\end{EQA}
By Lemma~\ref{LdltwLa}
\begin{EQA}
	&& \nquad
	\E \bigl\{ | \langle \QP \, \gaussv_{\GP},\av \rangle | \,\, \| \DVL \, \gaussv_{\GP} \|^{3} \, \bigr\} 
	\\
	& \leq &
	\bigl( \E \| \DVL \, \gaussv_{\GP} \|^{4} \bigr)^{3/4} \,\, 
	\bigl( \E \langle \QP \gaussv_{\GP},\av \rangle^{4} \bigr)^{1/4}
	\leq 
	3^{1/4} \, (\dimL + \normG)^{3/2} \sqrt{ \av^{\T} \QP \, \DVL_{\GP}^{-2} \QP^{\T} \av } .
\label{sfdhysdfyf11ewds3wsded}
\end{EQA}
Here we used that 
\( \E \langle \QP \gaussv_{\GP},\av \rangle^{4} = \E \langle \gaussv, \DVL_{\GP}^{-1} \QP^{\T} \av \rangle^{4}
= 3 (\av^{\T} \QP \, \DVL_{\GP}^{-2} \QP^{\T} \av)^{2} \).
Now \eqref{iUafxvdvH2vk3} follows from 
\begin{EQA}
	\sup_{\av \in \R^{\dimq} \colon \| \av \| = 1} \av^{\T} \QP \, \DVL_{\GP}^{-2} \QP^{\T} \av
	&=&
	\| \QP \, \DVL_{\GP}^{-2} \QP^{\T} \| \, .
\label{hdrtesw5sdghww3r4d5tdy6gh}
\end{EQA}
and \( 3^{1/4} (1 - \dltwb)^{-2} \leq 2.4 \).
\end{proof}
Now we state a sharper result based on \nameref{LL4tref}.

\begin{proposition}
\label{Lintfxupp2}
Suppose the conditions of Proposition~\ref{Lintfxupp3} and \nameref{LL4tref}.
Then for any function \( g(\uv) \) with \( |g(\uv)| \leq 1 \) and 
\( g(\uv) = g(-\uv) \)
\begin{EQA}
	\biggl| 
	\frac{\int_{\UVL} \ex^{\lgd(\xv;\uv)} \, g(\uv) \, d\uv - \int_{\UVL} \ex^{- \| \DVL_{\GP} \uv \|^{2}/2} \, g(\uv) \, d\uv} 
		 {\int \ex^{- \| \DVL_{\GP} \uv \|^{2}/2} d\uv} 
	\biggr|
	& \leq &
	\err_{4} \, 
\label{efuiguduUem}
\end{EQA}
with
\begin{EQA}
	\err_{4}
	& \eqdef &
	\frac{1}{16 (1 - \dltwb)^{2}} \Bigl\{ \E \bigl\langle \nabla^{3} \lgd(\xv) , \gaussv_{\GP}^{\otimes 3} \bigr\rangle^{2} 
	+ 2 \dltwu_{4} \E \| \DVL \gaussv_{\GP} \|^{4} \Bigr\}
	\\
	& \leq &
	\frac{1}{16 (1 - \dltwb)^{2}} \Bigl\{ \dltwu_{3}^{2} \, (\dimL + 2 \normG)^{3} + 2 \dltwu_{4} (\dimL + \normG)^{2} \Bigr\}
	\, .
\label{errdef3322Hm2}
\end{EQA}
If the function \( g(\cdot) \) is not bounded by one but it is symmetric and \( 2m \)-homogeneous,
i.e. \( g(t \uv) = t^{2m} g(\uv) \), then
\eqref{efuiguduUem} still applies with 
\begin{EQA}
	\err_{4}
	& \eqdef &
	\frac{1}{16 (1 - \dltwb)^{2+m}} \E \Bigl\{ 
		\bigl\langle \nabla^{3} \lgd(\xv) , \gaussv_{\GP}^{\otimes 3} \bigr\rangle^{2} \, g(\gaussv_{\GP})
	+ 2 \dltwu_{4} \, \| \DVL \gaussv_{\GP} \|^{4} \, g(\gaussv_{\GP}) \Bigr\}
	\, .
\label{errdef3322Hm2m}
\end{EQA}
\end{proposition}

\begin{proof}
We write \( f^{(3)} \) and \( \dltw_{k}(\uv) \) in place of \( \nabla^{3} \lgd(\xv) \) and
\( \dltw_{k}(\xv,\uv) \), \( k=3,4 \).
It holds 
\begin{EQA}
	\int_{\UVL} \ex^{\lgd(\xv;\uv)} \, g(\uv) \, d\uv
	& = &
	\int_{\UVL} \exp\Bigl\{ - \frac{\| \DVL_{\GP} \uv \|^{2}}{2} + \dltw_{3}(\uv) \Bigr\} \, g(\uv) \, d\uv .
\label{P3tHmzzmt22t3v}
\end{EQA}
Define for \( t \in [0,1] \)
\begin{EQA}
	\Rem(t)
	& \eqdef &
	\int_{\UVL} \exp\Bigl\{ - \frac{\| \DVL_{\GP} \uv \|^{2}}{2} + t \dltw_{3}(\uv) \Bigr\} \, g(\uv) \, d\uv .
\label{P3tHmzzmt22t}
\end{EQA}
Symmetricity of \( \UVL \) and \( g(\uv) = g(-\uv) \) implies that 
\begin{EQA}
	\Rem'(0)
	&=&
	\frac{1}{2} \int_{\UVL} 
	\exp\Bigl( - \frac{\| \DVL_{\GP} \uv \|^{2}}{2} \Bigr) \bigl\{ \dltw_{3}(\uv) + \dltw_{3}(-\uv) \bigr\} \, g(\uv) \, d\uv
	\\
	&=&
	\int_{\UVL} \exp\Bigl( - \frac{\| \DVL_{\GP} \uv \|^{2}}{2} \Bigr) \, \bar{\dltw}_{4}(\uv) \, g(\uv) \, d\uv \, 
\label{Pp012zzmH2u212v}
\end{EQA}
with \( \bar{\dltw}_{4}(\uv) = \bigl\{ \dltw_{4}(\uv) + \dltw_{4}(-\uv) \bigr\}/2 \).
Moreover, as  \( |\dltw_{3}(\uv)| \leq \dltwb \| \DVL \uv \|^{2}/2 \),
it holds for \( t \in [0,1] \) 
\begin{EQA}
	|\Rem''(t)|
	& \leq &
	\int_{\UVL} \dltw_{3}^{2}(\uv) 
		\exp\Bigl\{ - \frac{\| \DVL_{\GP} \uv \|^{2}}{2} + t \dltw_{3}(\uv) \Bigr\} \, |g(\uv)| \, d\uv
	\\
	& \leq &
	\int_{\UVL} \dltw_{3}^{2}(\uv) 
		\exp\Bigl( - \frac{\| \DVL_{\GP} \uv \|^{2} - \dltwb \| \DVL \uv \|^{2}}{2} \Bigr) \, d\uv \, .
\label{F3pptHm2t6v}
\end{EQA}
As \( \dltw_{3}(\uv) = \bigl\langle f^{(3)}, \uv^{\otimes 3} \bigr\rangle/6 + \dltw_{4}(\uv) \) and \( |\dltw_{4}(\uv)| \leq 1 \), one can bound for \( t \in [0,1] \)
\begin{EQA}
	|\Rem''(t)| 
	& \leq & 
	2 \int_{\UVL} \bigl\{ \bar{\dltw}_{4}^{2}(\uv) + \bigl| \bigl\langle f^{(3)}, \uv^{\otimes 3} \bigr\rangle/6 \bigr|^{2} \bigr\} 
		\exp\Bigl( - \frac{\| \DVL_{\GP} \uv \|^{2} - \dltwb \| \DVL \uv \|^{2}}{2} \Bigr) \, d\uv
	\\
	& \leq & 
	2 \int_{\UVL} \bigl\{ |\bar{\dltw}_{4}(\uv)| + \bigl\langle f^{(3)}, \uv^{\otimes 3} \bigr\rangle^{2} /36 \bigr\} 
	\exp\Bigl( - \frac{\| \DVL_{\GP} \uv \|^{2} - \dltwb \| \DVL \uv \|^{2}}{2} \Bigr) \, d\uv .
\label{u221md3u636v}
\end{EQA}
This and \eqref{Pp012zzmH2u212v} yield
\begin{EQA}
	&& \nquad
	\bigl| \Rem(1) - \Rem(0) \bigr|
	\leq 
	\bigl| \Rem'(0) \bigr| + \frac{1}{2} \sup_{t \in [0,1]} |\Rem''(t)|
	\leq 
	2 \int_{\UVL} |\bar{\dltw}_{4}(\uv)| \, 
	\ex^{- (\| \DVL_{\GP} \uv \|^{2} - \dltwb \| \DVL \uv \|^{2})/{2} } \, d\uv
	\\
	&&
	+ \, \frac{1}{36} \int_{\UVL} \bigl\langle f^{(3)}, \uv^{\otimes 3} \bigr\rangle^{2} \, 
		\ex^{- (\| \DVL_{\GP} \uv \|^{2} - \dltwb \| \DVL \uv \|^{2})/{2} } \, d\uv .
\label{C4Hm2Fm4tH2tv}
\end{EQA}
Change of variable \( \bigl( \Id - \dltwb \, \DVL_{\GP}^{-1} \DVL^{2} \, \DVL_{\GP}^{-1} \bigr)^{1/2} \uv \) to \( \wv \) yields by 
\eqref{12Im2t3t4Hm3} in view of \( \dltwb \leq 1/3 \)
\begin{EQA}
	&& \nquad
	\frac{1}{36} \int_{\UVL} \bigl\langle f^{(3)}, \uv^{\otimes 3} \bigr\rangle^{2} \, 
		\exp\Bigl( - \frac{\| \DVL_{\GP} \uv \|^{2} - \dltwb \| \DVL \uv \|^{2}}{2} \Bigr) \, d\uv 
	\\
	& \leq &
	\frac{3/2}{36 (1 - \dltwb)^{3}} 
	\int \bigl\langle f^{(3)}, \wv^{\otimes 3} \bigr\rangle^{2} \, \exp\Bigl( - \frac{\| \DVL_{\GP} \wv \|^{2}}{2} \Bigr) \, d\wv .
\label{C4Hm2Fm4tH2t3v}
\end{EQA}
Similarly by \nameref{LL4tref}
\begin{EQA}
	\int_{\UVL} |\bar{\dltw}_{4}(\uv)| \,\exp\Bigl( - \frac{\| \DVL_{\GP} \uv \|^{2} - \dltwb \| \DVL \uv \|^{2}}{2} \Bigr) \, d\uv
	& \leq &
	\frac{3/2}{24(1 - \dltwb)^{2}} 
	\int  
	\dltwu_{4} \| \DVL \wv \|^{4} \, \exp\Bigl( - \frac{\| \DVL_{\GP} \wv \|^{2}}{2} \Bigr) \, d\wv .
\label{C4Hm2Fm4tH2t4v}
\end{EQA}
The use of \( \dltwb \leq 1/3 \) implies that
\begin{EQA}
	\frac{\bigl| \Rem(1) - \Rem(0) \bigr|}{\int_{\UVL} \ex^{- \| \DVL_{\GP} \uv \|^{2}/2} d\uv}
	& \leq &
	\frac{3/2}{24 (1 - \dltwb)^{2}}
	\Bigl\{  
	\E \bigl\langle f^{(3)} , \gaussv_{\GP}^{\otimes 3} \bigr\rangle^{2}
		+ 2 \dltwu_{4} \E \| \DVL \gaussv_{\GP} \|^{4}
	\Bigr\} 
	\leq 
	\err_{4} \, 
\label{dwHw22n3n432}
\end{EQA}
and \eqref{efuiguduUem} follows.
Further, \nameref{LL3tref} yields
\( \bigl\langle \nabla^{3} \lgd(\xv) , \uv^{\otimes 3} \bigr\rangle^{2} \leq \dltwu_{3}^{2} \| \DVL \uv \|^{6} \).
Now \eqref{errdef3322Hm2} follows from Lemma~\ref{LdltwLa}.
The proof of \eqref{errdef3322Hm2m} is similar.
\end{proof}

Again, an important special cases correspond to \( m=1 \) and Laplace covariance approximation.

\begin{proposition}
\label{Perrboundm2}
Suppose the conditions of Proposition~\ref{Lintfxupp2} with \( \dltwb \leq 1/3 \).
Then 
for any linear mapping \( \QP \colon \R^{\dimp} \to \R^{\dimq} \) 
with \( \QP \QP^{\T} \leq \DVL^{2} \)
and any unit vector \( \av \in \R^{\dimq} \)
\begin{EQA}
	&& \nquad
	\frac{\bigl| 
		\int_{\UVL} \ex^{\lgd(\xv;\uv)} \, \langle \QP \uv,\av \rangle^{2} \, d\uv 
		- \int_{\UVL} \ex^{- \| \DVL_{\GP} \uv \|^{2}/2} \, \langle \QP \uv,\av \rangle^{2} \, d\uv 	 
		\bigr|} 
		 {\int \ex^{- \| \DVL_{\GP} \uv \|^{2}/2} d\uv} 
	\\
	& \leq &
	\frac{\| \QP \, \DVL^{-2} \QP^{\T} \|}{16 (1 - \dltwb)^{3}} \Bigl\{ 
		10.25 \dltwu_{3}^{2} \, (\dimL + 3\normG)^{3} + 3.5 \dltwu_{4} (\dimL + 3\normG)^{2} \Bigr\} 
	\\
	& \leq &
	\| \QP \, \DVL^{-2} \QP^{\T} \| \Bigl\{ 
		3 \dltwu_{3}^{2} \, (\dimL + 3\normG)^{3} + \dltwu_{4} (\dimL + 3\normG)^{2} \Bigr\}.
\label{iUafxvdvH2vk3m2}
\end{EQA}
\end{proposition}

\begin{proof}
By Proposition~\ref{Lintfxupp2}, it holds
\begin{EQA}
	&& \nquad
	\frac{\bigl| 
		\int_{\UVL} \ex^{\lgd(\xv;\uv)} \, \langle \QP \uv,\av \rangle^{2} \, d\uv 
		- \int_{\UVL} \ex^{- \| \DVL_{\GP} \uv \|^{2}/2} \, \langle \QP \uv,\av \rangle^{2} \, d\uv 	 
		\bigr|} 
		 {\int \ex^{- \| \DVL_{\GP} \uv \|^{2}/2} d\uv} 
	\\
	& \leq &
	\frac{1}{16 (1 - \dltwb)^{3}} \E \Bigl\{ 
		\dltwu_{3}^{2} \| \DVL \gaussv_{\GP} \|^{6} \, \langle \QP \gaussv_{\GP},\av \rangle^{2}
	+ 2 \dltwu_{4} \, \| \DVL \gaussv_{\GP} \|^{4} \, \langle \QP \gaussv_{\GP},\av \rangle^{2} \Bigr\} .
\label{iUafxvdvH2vk3m21}
\end{EQA}
Now we use that for any unit vector \( \av \in \R^{\dimq} \), by Lemma~\ref{LdltwLa}
\begin{EQA}
	\E \bigl\{ \| \DVL \gaussv_{\GP} \|^{6} \, \langle \QP \gaussv_{\GP},\av \rangle^{2} \bigr\}
	& \leq &
	\bigl( \E \| \DVL \gaussv_{\GP} \|^{8} \bigr)^{3/4} \,\, 
	\E \bigl( \langle \QP \gaussv_{\GP},\av \rangle^{8} \bigr)^{1/4}
	\\
	& \leq &
	10.25 (\dimL + 3\normG)^{3} \| \DVL_{\GP}^{-1} \QP^{\T} \av \|^{2}
	\leq 
	10.25 (\dimL + 3\normG)^{3} \| \QP \, \DVL_{\GP}^{-2} \QP^{\T} \| ,
	\\
	\E \bigl\{ \| \DVL \gaussv_{\GP} \|^{4} \, \langle \QP \gaussv_{\GP},\av \rangle^{2} \bigr\}
	& \leq &
	\bigl\{ \E \| \DVL \gaussv_{\GP} \|^{8} \bigr\}^{1/2} \, 
	\bigl\{ \E \langle \QP \gaussv_{\GP},\av \rangle^{4} \bigr\}^{1/2}
	\\
	& \leq &
	\sqrt{3} (\dimL + 3\normG)^{2} \| \DVL_{\GP}^{-1} \QP^{\T} \av \|^{2}
	\leq 
	\sqrt{3} \, (\dimL + 3\normG)^{2} \| \QP \, \DVL_{\GP}^{-2} \QP^{\T} \| 
\label{hcuf6yw3hfycvtegehctgggt}
\end{EQA}
and the result follows.
\end{proof}

\Section{Tail integrals}
\label{StailLa}
In this section we also write \( \xv \) in place of \( \xvs \).
Below we evaluate \( \rho \) from \eqref{rhfiUaceHu22m} which bounds the integral of \( \ex^{\lgd(\xv;\uv)} \) 
over the complement of the local set \( \UVL \) of a special form 
\( \UVL = \bigl\{ \uv \colon \| \DVL \uv \| \leq \amax^{-1} \rrL \bigr\} \) for \( \DVL \) from \nameref{LLf0ref}.
Our results help to understand how the radius \( \rrL \) should be fixed to ensure \( \rho \) sufficiently small.
The main tools for the analysis are deviation probability bounds for Gaussian quadratic forms; see Section~\ref{SdevboundGauss}.

\begin{proposition}
\label{PUVstarLLfx}
Suppose \nameref{LLf0ref}.
Given \( \amax < 1 \) and \( \xx > 0 \), 
let \( \UVL \) and \( \rrL \) be defined by \eqref{UvTDunm12spT}.
Let also \( \dltwb \) from \eqref{om3esuU1H02d3} satisfy \( \dltwb \leq 1 - \amax \).
Then 
\begin{EQA}
\label{fiinIHUniUef0}
	\frac{\int \Ind\bigl( \uv \not\in \UVL \bigr) \, \ex^{\lgd(\xv;\uv)} \, d\uv}
		 {\int \ex^{ - \| \DVL_{\GP} \uv \|^{2}/2} \, d \uv}
	& \leq &
	4 \ex^{-\xx - \dimL/2} \, ,
	\\
	\frac{\int \Ind\bigl( \uv \not\in \UVL \bigr) \, \ex^{ - \| \DVL_{\GP} \uv \|^{2}/2} \, d\uv}
		 {\int \ex^{ - \| \DVL_{\GP} \uv \|^{2}/2} \, d \uv}
	& \leq &
	\ex^{-\xx - \dimL/2} \, .
\label{fiinIHUniUef}
\end{EQA}
\end{proposition}

\begin{proof}
Let \( \uv \not \in \UVL \), i.e. \( \| \DVL \uv \| > \rr  \) with \( \rr = \amax^{-1} \rrL \).
Define \( \uvc = \rr \| \DVL \uv \|^{-1} \uv \) yielding \( \| \DVL \uvc \| = \rr \).
We also write \( \uv = (1 + \tau) \uvc \) for \( \tau > 0 \).
By \eqref{om3esuU1H02d3} 
and \( \nabla^{2} \lgdL(0) = - \DVL^{2} \)
\begin{EQ}[rcl]
	\lgdL(\uvc) - \lgdL(0) - \bigl\langle \nabla \lgdL(0), \uvc \bigr\rangle 
	& \leq &
	- (1 - \dltwb) \| \DVL \uvc \|^{2}/2 ,
	\\
	\bigl\langle \nabla \lgdL(\uvc) - \nabla \lgdL(0), \uv - \uvc \bigr\rangle
	& \leq &
	- (1 - \dltwb) \bigl\langle \DVL^{2} \uvc, \uv - \uvc \bigr\rangle .
\label{2ud1mud3Hv1dfx}
\end{EQ}
Concavity of \( \lgdL(\uv) \) implies for \( \uv = (1 + \tau) \uvc \),
\begin{EQA}
	\lgdL(\uv) 
	& \leq &
	\lgdL(\uvc) + \bigl\langle \nabla \lgdL(\uvc), \uv - \uvc \bigr\rangle 
\label{fudfpudumudwuld}
\end{EQA}
yielding by \eqref{2ud1mud3Hv1dfx} in view of 
\( \bigl\langle \DVL \uvc, \DVL \uv \bigr\rangle = \| \DVL \uvc \| \, \| \DVL \uv \| \)
\begin{EQA}
	&& \nquad
	\lgdL(\uv) - \lgdL(0) - \bigl\langle \nabla \lgdL(0), \uv \bigr\rangle
	=
	\lgdL(\uv) - \lgdL(\uvc) - \bigl\langle \nabla \lgdL(\uvc), \uv - \uvc \bigr\rangle
	\\
	&&
	\qquad
	\qquad
	\qquad
	\qquad
	+ \, \lgdL(\uvc) - \lgdL(0) - \bigl\langle \nabla \lgdL(0), \uvc \bigr\rangle
	+ \bigl\langle \nabla \lgdL(\uvc) - \nabla \lgdL(0), \uv - \uvc \bigr\rangle
	\\
	& \leq &
	(1 - \dltwb) \| \DVL \uvc \|^{2}/2
	- (1 - \dltwb) \bigl\langle \DVL \uvc, \DVL \uv \bigr\rangle 
	\leq 
	- (1 - \dltwb) \| \DVL \uvc \| \, \| \DVL \uv \| /2.
\label{tnfudnfHm1Hm1d}
\end{EQA}
We now use that 
\( \| \DVL \uvc \| = \rr \),  
\( \uvc = \uv/(1 + \tau) \), and thus,
\begin{EQA}
	&& \nquad
	\lgd(\xv + \uv) - \lgd(0) - \bigl\langle \nabla \lgd(\xv), \uv \bigr\rangle
	=
	\lgdL(\uv) - \lgdL(0) - \bigl\langle \nabla \lgdL(0), \uv \bigr\rangle
	- \| \DVL_{\GP} \uv \|^{2}/2 + \| \DVL \uv \|^{2}/2 
	\\
	& \leq &
	- (1 - \dltwb) \rr \| \DVL \uv \|/2 
	- \| \DVL_{\GP} \uv \|^{2}/2 + \| \DVL \uv \|^{2}/2.
\label{DGa22DTa221mw3rT22}
\end{EQA}
This yields by \( \rrL = \amax \, \rr \leq (1 - \dltwb) \rr \) with \( \Tau = \DVL \DVL_{\GP}^{-1} \)
\begin{EQA}
	&& \nquad
	\frac{\int \Ind\bigl( \uv \not\in \UVL \bigr) 
		\exp\bigl\{  \lgd(\xv + \uv) - \lgd(\xv) - \langle \nabla \lgd(\xv), \uv \rangle \bigr\} \, d\uv}
		{\int \exp \bigl( - \| \DVL_{\GP} \uv \|^{2}/{2} \bigr) \, d \uv}
	\\
	& \leq &
	\frac{\int \Ind\bigl( \| \DVL \uv \| > \rr \bigr)
		\exp \bigl\{ - {(1 - \dltwb) \rr} \| \DVL \uv \|/2 - \| \DVL_{\GP} \uv \|^{2}/2 + \| \DVL \uv \|^{2}/2 \bigr\} \, d\uv}
		{\int \exp \bigl( - \| \DVL_{\GP} \uv \|^{2}/{2} \bigr) \, d \uv}
	\\
	& \leq &
	\E \exp \bigl\{ 
		- \rrL \| \Tau \gaussv \| / 2 
		+ \| \Tau \gaussv \|^{2} / 2 
		\bigr\} \Ind\bigl( \| \Tau \gaussv \| > \rrL \bigr) 
\label{2C2emx1212IganiU}
\end{EQA}
with \( \gaussv \) standard normal in \( \R^{\dimp} \).
Next, define
\begin{EQA}
	\riskt_{0}(\rr)
	& \eqdef &
	\E \exp\bigl( - \rr \| \DVL \gaussv \| / 2 + \| \DVL \gaussv \|^{2}/2 \bigr) 
		\Ind(\| \DVL \gaussv \| > \rr) .
\label{mp2mxr2g22T}
\end{EQA}
Integration by parts allows to represent the last integral as
\begin{EQA}
	\riskt_{0}(\rr)
	&=&
	- \int_{\rr}^{\infty} \exp\Bigl( - { \rr \, \zq}/{2} + {\zq^{2}}/{2} \Bigr) \, 
		d\P\bigl( \| \Tau \gaussv \| > \zq \bigr)
	\\
	&=&
	\P\bigl( \| \Tau \gaussv \| > \rr \bigr) 
	+ \int_{\rr}^{\infty} (\zq - \rr/2) \exp\Bigl( - {\rr \zq}/{2} + {\zq^{2}}/{2} \Bigr) \, 
		\P\bigl( \| \Tau \gaussv \| > \zq \bigr) \, d\zq \, .
\label{dzzTPzr2rT}
\end{EQA}
By Theorem~\ref{TexpbLGA}, for any \( \zq \geq \sqrt{\dimL} \) for \( \dimL = \tr(\Tau \, \Tau^{\T}) = \tr (\DVL^{2} \, \DVL_{\GP}^{-2}) \)
\begin{EQA}
	\P\bigl( \| \Tau \gaussv \| > \zq \bigr)
	& \leq & 
	\exp\bigl\{ - (\zq - \sqrt{\dimL})^{2}/2 \bigr\}
\label{2emxPTgasps2d}
\end{EQA}
yielding for \( \zq \geq \rrL = 2 \sqrt{\dimL} + \sqrt{2\xx} \)
\begin{EQA}
	\P\bigl( \| \Tau \gaussv \| > \zq \bigr)
	& \leq &
	\exp\bigl\{ - (\zq - \sqrt{\dimL})^{2}/2 \bigr\}
	\leq 
	\ex^{-\xx - \dimL/2}
\label{2exmp22ezdp2}
\end{EQA}
and for \( \rr \geq 2 \sqrt{\dimL} + \sqrt{2\xx} \) and \( \xx \geq 2 \)
\begin{EQA}
	\riskt_{0}(\rr)
	& \leq &
	\ex^{-\xx - \dimL/2}
	+ \int_{\rr}^{\infty} (\zq - \rr/2) 
	\exp\Bigl\{ - \frac{\rr \zq}{2} + \frac{\zq^{2}}{2} - \frac{(\zq - \sqrt{\dimL})^{2}}{2} \Bigr\} \, d\zq
	\\
	& \leq &
	\ex^{-\xx - \dimL/2}
	+ \exp\Bigl( - \frac{(\rr - \sqrt{\dimL})^{2}}{2} \Bigr) \int_{0}^{\infty} \Bigl(\zq + \frac{\rr}{2} \Bigr) 
	\exp\Bigl\{ - \frac{(\rr - 2 \sqrt{\dimL}) \zq}{2} \Bigr\} \, d\zq
	\\
	& \leq &
	2 \ex^{-\xx - \dimL/2} .
\label{Cexmxmp222Tr22fx}
\end{EQA}
This completes the proof of the result \eqref{fiinIHUniUef0}.
Statement \eqref{fiinIHUniUef} is about Gaussian probability 
\( \P\bigl( \| \Tau \gaussv \| \geq \rrL \bigr) \) for a standard normal element \( \gaussv \), and we derive
\begin{EQA}
	\P\bigl( \| \Tau \gaussv \| \geq 2 \sqrt{\dimL} + \sqrt{2\xx} \bigr)
	& \leq &
	\exp\bigl\{ - \bigl( \sqrt{\dimL} + \sqrt{2\xx} \bigr)^{2}/2 \bigr\}
	\leq 
	\exp\bigl( -\xx - \dimL/2 \bigr) 
\label{fiinIHUniUefi}
\end{EQA}
and \eqref{fiinIHUniUef} follows.
\end{proof}

\medskip

The next result extends \eqref{fiinIHUniUef0}.
\begin{proposition}
\label{PUVstarLLgenlifx}
Assume the conditions of Proposition~\ref{PUVstarLLfx}
with
\begin{EQA}
	\rrL
	& \geq &
	2 \sqrt{\dimL} + \sqrt{2 \xx} + m
\label{1w3rT2spTs2xlifx}
\end{EQA}
for some \( m \geq 0 \).
Then 
\eqref{fiinIHUniUef0} can be extended to
\begin{EQA}
\label{fiinIHUniUef0m}
	\frac{\int \Ind\bigl( \uv \not\in \UVL \bigr) \, \| \DVL \uv \|^{m} \, \ex^{\lgd(\xv;\uv)} \, d\uv}
		 {\int \ex^{ - \| \DVL_{\GP} \uv \|^{2}/2} \, d \uv}
	& \leq &
	4 \ex^{-\xx - \dimL/2} \, ,
	\\
	\frac{\int \Ind\bigl( \uv \not\in \UVL \bigr) \, \| \DVL \uv \|^{m} \, \ex^{ - \| \DVL_{\GP} \uv \|^{2}/2} \, d\uv}
		 {\int \ex^{ - \| \DVL_{\GP} \uv \|^{2}/2} \, d \uv}
	& \leq &
	\ex^{-\xx - \dimL/2} \, .
\label{fiinIHUniUefm}
\end{EQA}
\end{proposition}

\begin{proof}
The case \( m > 0 \) can be proved similarly to \( m=0 \) using 
\( m \log z \leq m z \).
\end{proof}

\Section{Local concentration}
Here we show that the measure \( \PfL \) well concentrates on the local set \( \UVL \) from \eqref{UvTDunm12spT}.
Again we fix \( \xv = \xvs \).

\begin{proposition}
\label{PlocconLa}
Assume \( \dltwb \leq 1/3 \).
Then
\begin{EQA}
	\int_{\UVL} \ex^{\lgd(\xv;\uv)} \, d\uv  
	& \geq &
	\ex^{ - \dltwb \, \dimL/2} \, \int_{\UVL} \ex^{- \| \DVL_{\GP} \uv \|^{2}/2} d\uv \, .
\label{emdw3dp0iUfxu}
\end{EQA}
Moreover,
\begin{EQA}
	\frac{\int_{\UVL^{c}} \ex^{\lgd(\xv;\uv)} \, d\uv}{\int \ex^{\lgd(\xv;\uv)} \, d\uv}  
	& \leq &
	4 \ex^{ - \xx - (1 - \dltwb) \, \dimL/2} 
	\leq 
	\ex^{-\xx} \, .
\label{emdw3dp0iUfxuL}
\end{EQA}
\end{proposition}

\begin{proof}
By \eqref{om3esuU1H02d3}
\begin{EQA}
	\int_{\UVL} \ex^{\lgd(\xv;\uv)} \, d\uv
	& = &
	\int_{\UVL} \ex^{- \| \DVL_{\GP} \uv \|^{2}/2 + \dltw_{3}(\xv,\uv)} \, d\uv 
	\geq 
	\int_{\UVL} \ex^{- \| \DVL_{\GP} \uv \|^{2}/2 - \dltwb \| \DVL \uv \|^{2}/2} \, d\uv \, .
\label{PhtiUexmHu222}
\end{EQA}
Change of variable \( \bigl( \Id + \dltwb \, \DVL_{\GP}^{-1} \DVL^{2} \, \DVL_{\GP}^{-1} \bigr)^{1/2} \uv \) to \( \wv \) yields by \eqref{12Ip2t3t4Hp3} 
\begin{EQA}
	\int_{\UVL} \ex^{\lgd(\xv;\uv)} \, d\uv
	& \geq &
	\det \bigl( \Id + \dltwb \, \DVL_{\GP}^{-1} \DVL^{2} \, \DVL_{\GP}^{-1} \bigr)^{-1/2} 
	\int_{\UVL} \ex^{- \| \DVL_{\GP} \wv \|^{2}/2 } \, d\wv
	\\
	& \geq &
	\ex^{- \dltwb \, \dimL / 2} \int_{\UVL} \ex^{- \| \DVL_{\GP} \wv \|^{2}/2 } \, d\wv ,
\label{iUefxudetId3H121}
\end{EQA}
and \eqref{emdw3dp0iUfxu} follows.
This and \eqref{fiinIHUniUef0}, \eqref{fiinIHUniUef} of Proposition~\ref{PUVstarLLfx} imply
\begin{EQA}
	&& \nquad
	\frac{\int_{\UVL^{c}} \ex^{\lgd(\xv;\uv)} \, d\uv}{\int \ex^{\lgd(\xv;\uv)} \, d\uv}
	= 
	\frac{\int_{\UVL^{c}} \ex^{\lgd(\xv;\uv)} \, d\uv}{\int_{\UVL} \ex^{\lgd(\xv;\uv)} \, d\uv + \int_{\UVL^{c}} \ex^{\lgd(\xv;\uv)} \, d\uv}
	\\
	& \leq &
	\frac{4 \ex^{-\xx - \dimL/2} \int \ex^{- \| \DVL_{\GP} \uv \|^{2}/2} d\uv}
		 {\ex^{ - \dltwb \, \dimL/2} \, \int_{\UVL} \ex^{- \| \DVL_{\GP} \uv \|^{2}/2} d\uv 
			+ 4 \ex^{-\xx - \dimL/2} \int \ex^{- \| \DVL_{\GP} \uv \|^{2}/2} d\uv}
	\leq 
	4 \ex^{-\xx - (1 - \dltwb) \dimL/2} 
\label{jyfftc434rdssdsaaa}
\end{EQA}
as required in \eqref{emdw3dp0iUfxuL}.
\end{proof}

\Section{Finalizing the proof of Theorem~\ref{TLaplaceTV} and \ref{TLaplaceTV34}}
These results are proved by compiling the previous technical statements.
Proposition~\ref{PUVstarLLfx} provides some upper bounds for the quantities \( \rho \) and \( \rho_{\GP} \),
while Proposition~\ref{Lintfxupp3}, Proposition~\ref{Lintfxupp3T}, and Proposition~\ref{Lintfxupp2} bound the local errors 
\( \err \) and \( \err_{g} \).
The final bound \eqref{ufgdt6df5dtgededsxd23gjg} follows from Corollary~\ref{CPunbintLapl}.

\ifKL{
\Section{Proof of Theorem~\ref{TLaplaceKL} and Theorem~\ref{TLaplaceKLi}}
Define
\begin{EQA}
	\CDG 
	& \eqdef &
	\log \int \ex^{- \| \DVL_{\GP} \uv \|^{2}/2} \, d\uv 
	- \log \int \ex^{\lgd(\xvs;\uv)} \, d\uv 
	\, . 
\label{sihedyfw3ytw3ydqwdbhwd}
\end{EQA}
For \( \uv = \xv - \xvs \), it holds as in \eqref{IIgfifxutudufxutdu}
\begin{EQA}
	\PfL
	& \sim &
	\frac{\ex^{\lgd(\xv)}}{\int \ex^{\lgd(\xv)} \, d\uv}
	=
	\frac{\ex^{\lgd(\xv) - \lgd(\xvs)}}{\int \ex^{\lgd(\xv) - \lgd(\xvs)} \, d\uv}
	=
	\frac{\ex^{\lgd(\xvs;\uv)}}{\int \ex^{\lgd(\xvs;\uv)} \, d\uv} 
	=
	\frac{\ex^{\lgd(\xvs;\uv) + \CDG}}{\int \ex^{- \| \DVL_{\GP} \uv \|^{2}/2} \, d\uv} \, .
\label{dfhuwffwt5wtcsbydyqyqh}
\end{EQA}
Further, with \( \P_{\GP} = \ND(\xvs,\DVL_{\GP}^{-2}) \)
\begin{EQA}
	\log \frac{d\PfL}{d\P_{\GP}}(\xv) 
	&=&
	\lgd(\xvs;\uv) - \| \DVL_{\GP} \uv \|^{2}/2 + \CDG
	=
	\dltw_{3}(\uv) + \CDG 
\label{sufuwdeuquqeygwyg22e24}
\end{EQA}
and
\begin{EQA}
	\kullb(\PfL,\P_{\GP})
	&=&
	\ex^{\CDG} \, \frac{\int \dltw_{3}(\uv) \, \ex^{\lgd(\xvs;\uv)} \, d\uv}
		 {\int \ex^{- \| \DVL_{\GP} \uv \|^{2}/2} \, d\uv}
	+ \CDG \, .
\label{dvjnue37e37e7weedhe}
\end{EQA}
Similarly to \eqref{Pht12u212H02u22} and \eqref{Pd3p2d3H0wem22g1}, we can bound
\begin{EQA}
	\left| \int_{\UVL} \dltw_{3}(\uv) \, \ex^{\lgd(\xvs;\uv)} \, d\uv \right|
	& \leq &
	\int_{\UVL} |\dltw_{3}(\uv)| \, \ex^{- \| \DVL_{\GP} \uv \|^{2}/2 - \dltwb \| \DVL \uv \|^{2}/2} \, d\uv 
	\\
	& \leq &
	\frac{1}{4(1 - \dltwb)^{3/2}}
	\int \dltwu_{3}(\uv) \, \ex^{- \| \DVL_{\GP} \uv \|^{2} /2} \, d\uv \, .
\label{ucuwudhuwhi2wefhw78w}
\end{EQA}
On the complement \( \uv \in \UVL^{c} \), we use that
\begin{EQA}
	\dltw_{3}(\uv)
	=
	\lgd(\xvs;\uv) + \frac{1}{2} \| \DVL_{\GP} \uv \|^{2}
	&=&
	\lgdL(\xvs;\uv) + \frac{1}{2} \| \DVL \uv \|^{2}
	\leq 
	\frac{1}{2} \| \DVL \uv \|^{2} .
\label{yuqt2sstwwtwbstdywhsyw}
\end{EQA}
The last inequality is based on concavity of \( \lgdL(\cdot) \) and local approximation 
\( \lgdL(\xvs;\uv) \approx - \| \DVL \uv \|^{2}/2 \) within \( \UVL \) yielding \( \lgdL(\xvs;\uv) < 0 \) for \( \uv \in \UVL^{c} \).
By Proposition~\ref{PUVstarLLgenlifx} with \( m = 2 \),
\begin{EQA}
	\int_{\UVL^{c}} \dltw_{3}(\uv) \, \ex^{\lgd(\xvs;\uv)} \, d\uv
	& \leq &
	\frac{1}{2} \int_{\UVL^{c}} \| \DVL \uv \|^{2} \, \ex^{\lgd(\xvs;\uv)} \, d\uv
	\leq 
	\ex^{-\xx} \int \ex^{- \| \DVL_{\GP} \uv \|^{2}/2} \, d\uv \, .
\label{sywdywydwbwtwgqb1qt213ba}
\end{EQA}
We conclude that
\begin{EQA}
	\frac{\int \dltw_{3}(\uv) \, \ex^{\lgd(\xvs;\uv)} \, d\uv}
		 {\int \ex^{- \| \DVL_{\GP} \uv \|^{2}/2} \, d\uv}
	& \leq &
	\frac{\E \dltwu_{3}(\gaussv_{\GP})}{4(1 - \dltwb)^{3/2}} + \ex^{- \xx} 
	\leq 
	\err_{3} + \ex^{-\xx} .
\label{0fmdfeyjdwhqkshyxiwck}
\end{EQA}
Similarly
\begin{EQA}
	&& \nquad
	\bigl| \ex^{\CDG} - 1 \bigr|
	=
	\left| \frac{\int \ex^{\lgd(\xvs;\uv)} \, d\uv - \int \ex^{- \| \DVL_{\GP} \uv \|^{2}/2} \, d\uv}
		 {\int \ex^{- \| \DVL_{\GP} \uv \|^{2}/2} \, d\uv} 
	\right|
	\\
	& \leq &
	\left| \frac{\int_{\UVL} \ex^{\lgd(\xvs;\uv)} \, d\uv - \int_{\UVL} \ex^{- \| \DVL_{\GP} \uv \|^{2}/2} \, d\uv}
		 {\int \ex^{- \| \DVL_{\GP} \uv \|^{2}/2} \, d\uv} 
	\right|
	+ \frac{\int_{\UVL^{c}} \ex^{\lgd(\xvs;\uv)} \, d\uv}{\int \ex^{- \| \DVL_{\GP} \uv \|^{2}/2} \, d\uv}
	+ \frac{\int_{\UVL^{c}} \ex^{- \| \DVL_{\GP} \uv \|^{2}/2} \, d\uv}
		 {\int \ex^{- \| \DVL_{\GP} \uv \|^{2}/2} \, d\uv}
	\\
	& \leq &
	\frac{\E \dltwu_{3}(\gaussv_{\GP})}{4(1 - \dltwb)^{3/2}} + 2 \ex^{- \xx} 
	\leq 
	\err_{3} + 2 \ex^{-\xx} .
\label{hsyc4w2tde6dyfdgyww23}
\end{EQA}
This yields \( \ex^{\CDG} \leq 1 + \err_{3} + 2 \ex^{-\xx} \) and \( \CDG \leq \err_{3} + 2 \ex^{-\xx} \).
Putting all together results in
\begin{EQA}
	\kullb(\PfL,\P_{\GP})
	& \leq &
	\bigl( 1 + \err_{3} + 2 \ex^{-\xx} \bigr) \bigl( \err_{3} + \ex^{-\xx} \bigr) + \err_{3} + 2 \ex^{-\xx}
	<  
	4 \err_{3} + 4 \ex^{-\xx} 
\label{gcusgu255eetfghbhjjkkook}
\end{EQA}
provided that \( 4 \err_{3} + 4 \ex^{-\xx} \leq 1 \),
and \eqref{vlgvi8ugu7tr4ry43et31} follows.

The proof of Theorem~\ref{TLaplaceKLi} is similar and even simpler except one special part, namely, 
the bound for the tail integral of \( - \dltw_{3}(\uv) \).
By definition \( \| \DVL \uv \| \geq \rr_{\GP} \) for \( \uv \in \UVL^{c} \).
This implies by \eqref{dchbhwdhwwdgscsn2efty2162} similarly to \eqref{yuqt2sstwwtwbstdywhsyw} 
\begin{EQA}
	&& \nquad
	- \int_{\UVL^{c}} \dltw_{3}(\uv) \, \ex^{- \| \DVL_{\GP} \uv \|^{2}/2} \, d\uv
	\leq
	\int_{\UVL^{c}} \bigl| \lgdL(\xvs;\uv) \bigr| \, \ex^{- \| \DVL_{\GP} \uv \|^{2}/2 + \rho_{\GP} \| \DVL \uv \|^{2}/2} \,
		\ex^{- \rho_{\GP} \| \DVL \uv \|^{2}/2} \, d\uv
	\\
	& \leq &
	\ex^{- \rho_{\GP} \rr_{\GP}^{2}/2} \, 
	\int \bigl| \lgdL(\xvs;\uv) \bigr| \, \ex^{- \| \DVL_{\GP} \uv \|^{2}/2 + \rho_{\GP} \| \DVL \uv \|^{2}/2} \,
		d\uv
	\leq \CONSTi_{\lgdL} \, \ex^{-\xx} ,
\label{sguswdw2e2j2wuhq2wdw}
\end{EQA}
and the result follows.
}{}

\ifKL{
\Section{Proof of Theorem~\ref{TpostmeanLa}, Corollary~\ref{CTpostmeanLa}, and Theorem~\ref{TpostvarL}}
For Theorem~\ref{TpostmeanLa}, we follow the same line as for Theorem~\ref{TLaplaceTV34}.
Note first that
\begin{EQA}
	\QP (\xvb - \xvs)
	=
	\frac{\int \QP(\xvs + \uv) \, \ex^{\lgd(\xvs + \uv)} \, d\uv}{\int \ex^{\lgd(\xvs + \uv)} \, d\uv} - \QP \xvs
	&=&
	\frac{\int \QP \uv \, \ex^{\lgd(\xvs;\uv)} \, d\uv}{\int \ex^{\lgd(\xvs;\uv)} \, d\uv} \, 
\label{dvhjt6efedfchsijcfte4ws}
\end{EQA}
and
\begin{EQA}
	\| \QP (\xvb - \xvs) \|
	&=&
	\sup_{\av \in \R^{\dimq} \colon \| \av \| = 1} \bigl| \langle \QP (\xvb - \xvs), \av \rangle \bigr|
	= 
	\sup_{\av \in \R^{\dimq} \colon \| \av \| = 1} 
		\left| \frac{\int \langle \QP \uv, \av \rangle \, \ex^{\lgd(\xvs;\uv)} \, d\uv}{\int \ex^{\lgd(\xvs;\uv)} \, d\uv} \right| \, .
\label{p3hb893dfvfdsiuw5whlh}
\end{EQA}
Now fix \( \av \in \R^{\dimq} \) with \( \| \av \| = 1 \) and \( g(\uv) = \langle \QP \uv, \av \rangle \).
Bound \eqref{igefxumiguexHu22g} implies
\begin{EQA}[rcl]
	\left| \frac{\int g(\uv) \, \ex^{\lgd(\xvs;\uv)} \, d\uv}{\int \ex^{\lgd(\xvs;\uv)} \, d\uv} \right| 
	& \leq &
	\frac{\rho_{g} + \rho_{\GP,g} + \err_{3,g}}{1 - \rho_{\GP} - \err_{3}} \, .
\label{igefxumiguexHu22gp}
\end{EQA}
The bound \( 1 - \err_{3} - \rho_{\GP} \geq 1/2 \) has been already checked.
Proposition~\ref{PUVstarLLgenlifx} for \( m=1 \) helps to bound the values \( \rho_{g} \) and \( \rho_{\GP,g} \) 
by \( \ex^{-\xx} \).
Next we bound \( \err_{3,g} \).
By Proposition~\ref{Perrbound3} 
\begin{EQA}
	4 \err_{3,g}
	& \leq &
	2.4 \, \dltwu_{3} \, (\dimL + \normG)^{3/2} \, \| \QP \, \DVL_{\GP}^{-2} \QP^{\T} \|^{1/2} .
\label{sfdhysdfyf11ewds3wsded}
\end{EQA}
This yields \eqref{jrwguvyr23jbviufdsfsdgf6wm} in view of \( \dltwu_{3} = \hmax_{3}/\sqrt{n} \).
With \( \QP = \DVL \), this also implies \eqref{klu8gitfdgregfkhj7yt}.

Now we proceed with the proof of Theorem~\ref{TpostvarL}.
For any linear mapping \( \QP \colon \R^{\dimp} \to \R^{\dimq} \), it holds 
\begin{EQA}
	\QP (\IFCovs - \IFL^{-1}) \QP^{\T}
	& \leq &
	\QP (\IFCov - \IFL^{-1}) \QP^{\T} + |\QP (\xvs - \xvb)|^{2}.
\label{cfhsxsuxgwugxuwgwgxrdrdh}
\end{EQA}
%
Now fix any unit vector \( \av \in \R^{\dimq} \) and consider the test function 
\( g(\uv) = \langle \QP \uv , \av \rangle^{2} \), so that \( \EfL \, g(\Xv - \xvs) = \av^{\T} \QP \IFCovs \QP^{\T} \av \).
Bound \eqref{2rherIHg2red} of Proposition~\ref{PunbintLapg} involves the already evaluated 
quantities \( \err = \err_{4} \), \( \rho \), and \( \rho_{\GP} \), 
as well as  the \( g \)-dependent quantities \( \rho_{g} \), \( \rho_{\GP,g} \), \( \II_{\GP}(g) \),
and \( \err_{g} \).
By Proposition~\ref{PUVstarLLgenlifx},
\( \rho_{g} + \rho_{\GP,g} \leq 2 \ex^{-\xx} \).
Proposition~\ref{Perrboundm2} states bound \eqref{iUafxvdvH2vk3m2} for \( \err_{g} \).
Further,
\begin{EQA}
	\II_{\GP}(g)
	&=&
	\E \langle \QP \gaussv_{\GP}, \av \rangle^{2}
	=
	\| \DVL_{\GP}^{-1} \QP^{\T} \av \|^{2}
	\leq 
	\| \QP \, \DVL_{\GP}^{-2} \QP^{\T} \| .
\label{0457uhg3rhdf9fjr}
\end{EQA}
We conclude that
\begin{EQA}
	&& \nquad
	\Bigl| \EfL \langle \QP (\Xv - \xvs), \av \rangle^{2} - \| \DVL_{\GP}^{-1} \QP^{\T} \av \|^{2} \Bigr|
	\\
	& \leq &
	\frac{\| \QP \, \DVL_{\GP}^{-2} \QP^{\T} \|}{n}
	\bigl( 3 \hmax_{3}^{2} + \hmax_{4} \bigr) (\dimA + 3)^{3} + |\QP (\xvs - \xvb)|^{2} + 8 \ex^{- \xx} .
\label{5hjkkoiioi1223hhh}
\end{EQA}
This and \eqref{hcdtrdtdehfdewdrfrhgyjufger} of Theorem~\ref{TpostmeanLa} yield the assertion of Theorem~\ref{TpostvarL}.
}{}

}


\Chapter{Examples of priors}
\label{Spriors}
This section presents two typical examples of priors and some properties 
including the bounds for \emph{effective dimension} and \emph{Laplace effective dimension}.

\Section{Truncation and smooth priors}
\label{Ssmoothprior}
Below we consider two non-trivial examples of Gaussian priors: 
truncation and smooth priors.
To make the presentation clear, we impose some assumptions on the considered setup.
Most of them are non-restrictive and can be extended to more general situations.
We assume to be given a growing sequence of nested linear approximation subspaces 
\( \VV_{1} \subset \VV_{2} \subset \ldots \subset \R^{\dimp} \) of dimension 
\( \dim(\VV_{\mm}) = \mm \).
Below \( \Proj_{\mm} \) is the projector on \( \VV_{\mm} \) and 
\( \VV_{\mm}^{c} \) is the orthogonal complement of \( \VV_{\mm} \).
A \emph{smooth prior} is described by a self-adjoint operator \( \GP \) such that
\( \| \GP \uv \| / \| \uv \| \) becomes large for \( \uv \in \VV_{\mm}^{c} \) and \( \mm \) large.
One can write this condition in the form
\begin{EQ}[rcl]
	\| \GP \uv \|^{2}
	& \leq &
	\gp_{\mm}^{2} \| \uv \|^{2} ,
	\qquad 
	\uv \in \VV_{\mm} \, ,
	\\
	\| \GP \uv \|^{2}
	& \geq &
	\gp_{\mm}^{2} \| \uv \|^{2} \, ,
	\qquad 
	\uv \in \VV_{\mm}^{c} \, .
\label{uniVmgm2u2uiVm}
\end{EQ}
Often one assumes that \( \VV_{\mm} \) is spanned by the eigenvectors
of \( \GP^{2} \) corresponding to its smallest eigenvalues 
\( \gp_{1}^{2} \leq \gp_{2}^{2} \leq \ldots \leq \gp_{\mm}^{2} \).
We only need \eqref{uniVmgm2u2uiVm}.
A typical example is given by \( \GP^{2} = \diag(\gp_{j}^{2}) \) with \( \gp_{j}^{2} = \CGP^{-1} j^{2\smp} \) for \( \smp > 1/2 \)
and some window parameter \( \CGP \).
Below we refer to this case as \( (\smp,\CGP) \)-\emph{smooth prior}.

A \( \mm \)-\emph{truncation prior} assumes that the prior distribution is restricted to \( \VV_{\mm} \). 
This formally corresponds to a covariance operator \( \GP_{\mm}^{-2} \) with 
\( \GP_{\mm}^{-2} \bigl( \Id - \Proj_{\mm} \bigr) = 0 \).
Equivalently, we set \( \gp_{\mm+1} = \gp_{\mm+2} = \ldots = \infty \) in \eqref{uniVmgm2u2uiVm}.

\Section{Effective dimension}
\label{Seffdima}
This section explains how the \emph{effective dimension} and \emph{Laplace effective dimension}
can be evaluated for some typical situations.  

Let \( \IF \) be a generic information matrix while \( \GP^{2} \) a penalizing matrix. 
With \( \IF_{\GP} \eqdef \IF + \GP^{2} \), define a sub-projector \( \proj_{\GP} \) in \( \R^{\dimp} \) by
\begin{EQA}[c]
	\proj_{\GP}
	\eqdef
	\IF_{\GP}^{-1} \IF
	=
	(\IF + \GP^{2})^{-1} \IF .
\label{9bvgu7fh54rfgy7byrfuer4}
\end{EQA}
Also define the \emph{Laplace effective dimension}
\begin{EQA}
	\dimA(\GP)
	& \eqdef &
	\tr \proj_{\GP}
	=
	\tr \bigl\{ (\IF + \GP^{2})^{-1} \IF \bigr\} .
\label{pjnv6yfkjwwyfvjvikekm}
\end{EQA}
The penalizing matrix \( \GP^{2} \) will be supposed diagonal, 
\( \GP^{2} = \diag\{ \gp_{1}^{2},\ldots,\gp_{\dimp}^{2} \} \).
Moreover, we implicitly assume that the values \( \gp_{j}^{2} \) \emph{grow} with \( j \) at some rate, 
\emph{polynomial} or exponential, yielding for all \( \mm \geq 1 \)
\begin{EQA}
	\sum_{j > \mm} \gp_{j}^{-2} 
	& \leq &
	\CONSTgp \, \mm \, \gp_{\mm}^{-2} \, .
\label{sumjJgjm2C}
\end{EQA}
Our leading example is given by \( \gp_{j}^{2} = \CGP^{-1} j^{2\smp} \) for \( \smp > 1/2 \).
Then \eqref{sumjJgjm2C} holds with \( \CONSTgp = (2 \smp -1)^{-1} \).

Concerning the matrix \( \IF \), we assume 
\begin{EQ}[rcccl]
	\CONSTIF^{-1} \, \nsize \| \uv \|^{2} 
	& \leq &
	\bigl\langle \IF \uv, \uv \bigr\rangle
	& \leq &
	\CONSTIF \, \nsize \| \uv \|^{2},
	\qquad
	\uv \in \R^{\dimp} \, ,
\label{uVmGu2C1C2gen}
\end{EQ}
for some \( \CONSTIF \geq 1 \).
It appears that the value \( \dimA(\GP) \) is closely related to the index \( \mm \) for which 
\( \gp_{\mm}^{2} \approx \nsize \).

\begin{lemma}
\label{LeffdimG}
Let \( \IF \) satisfy \eqref{uVmGu2C1C2gen}. 
Let also \( \GP^{2} = \diag\{ \gp_{1}^{2},\ldots,\gp_{\dimp}^{2} \} \) with \( \gp_{j}^{2} \) satisfying
\eqref{sumjJgjm2C}.
Define the index \( \mm \) as the smallest index \( j \) with \( \gp_{j}^{2} \geq \nsize \):
\begin{EQA}
	\mm
	& = &
	\mm(\GP)
	\eqdef
	\min \{ j \colon \gp_{j}^{2} \geq n \} .
\label{gklghi43mjd7skje8vhr}
\end{EQA}
Then 
\begin{EQA}
	\frac{1}{\CONSTIF + 1}
	& \leq &
	\frac{\dimA(\GP)}{\mm}
	\leq 
	1 + \CONSTIF \, \CONSTgp \, .
\label{4hhyc6vfjhfgyyrheudj}
\end{EQA}
\end{lemma}

\begin{proof}
By \eqref{sumjJgjm2C}
\begin{EQA}
	\tr (\proj_{\GP})
	& \leq &
	\sum_{j \geq 1} \frac{\CONSTIF \, \nsize}{\CONSTIF \, \nsize + \gp_{j}^{2}}
	\leq 
	\mm + \sum_{j > \mm} \frac{\CONSTIF \, \nsize}{\CONSTIF \, \nsize + \gp_{j}^{2}}
	\leq 
	\mm + \CONSTIF \, \nsize \sum_{j > \mm} \gp_{j}^{-2}
	\\
	& \leq &
	\mm + \CONSTIF \, \nsize \, \CONSTgp \, \mm \, \gp_{\mm}^{-2}
	\leq 
	\mm (1 + \CONSTIF \, \CONSTgp) .
\label{pendf6ycvfhewkevtye}
\end{EQA}
Similarly
\begin{EQA}
	\tr (\proj_{\GP})
	& \geq &
	\sum_{j=1}^{\mm} \frac{\CONSTIF^{-1} \, \nsize}{\CONSTIF^{-1} \, \nsize + \gp_{j}^{2}}
	\geq 
	\mm \, \frac{\CONSTIF^{-1} \, \nsize}{\CONSTIF^{-1} \, \nsize + \gp_{\mm}^{2}}
	\geq 
	\mm \, \frac{\CONSTIF^{-1}}{\CONSTIF^{-1} + 1} \, ,
\label{2fdvfg9rfmvvhywhsh}
\end{EQA}
and the assertion follows.
\end{proof}

This result yields an immediate corollary.

\begin{corollary}
\label{CLeffdimG}
Let \( \GP_{1}^{2} \) and \( \GP_{2}^{2} \) be two different penalizing matrices 
satisfying \eqref{sumjJgjm2C} and such that 
\( \mm(\GP_{1}) = \mm(\GP_{2}) \); see \eqref{gklghi43mjd7skje8vhr}.
Then 
\begin{EQA}
	\frac{\dimA(\GP_{1})}{\dimA(\GP_{2})}
	& \leq &
	( 1 + \CONSTIF \, \CONSTgp ) (1 + \CONSTIF) .
\label{c98dfjme478fcdje3otwe}
\end{EQA}
\end{corollary}

Now we evaluate the \emph{effective dimension} \( \dimG = \tr(\IF_{\GP}^{-1} \VP^{2}) \), 
where the \emph{variance matrix} \( \VP^{2} \) satisfies the condition 
\begin{EQ}[rcccl]
	\CONSTV^{-1} \, \| \IF \uv \|^{2} 
	& \leq &
	\| \VP \uv \|^{2}
	& \leq &
	\CONSTV \, \| \IF \uv \|^{2},
	\qquad
	\uv \in \R^{\dimp} \, 
\label{uVmGu2C1C2genV}
\end{EQ}
with some constant \( \CONSTV \geq 1 \); cf. \eqref{uVmGu2C1C2gen}.

\begin{lemma}
\label{LeffdimGV}
Assume \eqref{uVmGu2C1C2gen} for \( \IF \) and \eqref{uVmGu2C1C2genV} for \( \VP^{2} \). 
Let also \( \GP^{2} = \diag\{ \gp_{1}^{2},\ldots,\gp_{\dimp}^{2} \} \) with \( \gp_{j}^{2} \) satisfying
\eqref{sumjJgjm2C} and let \( \mm \) be given by \eqref{gklghi43mjd7skje8vhr}.
Then \( \dimG = \tr(\IF_{\GP}^{-1} \VP^{2}) \) satisfies
\begin{EQA}
	\frac{\CONSTV^{-1}}{\CONSTIF + 1}
	& \leq &
	\frac{\dimG}{\mm}
	\leq 
	\CONSTV (1 + \CONSTIF \, \CONSTgp) \, .
\label{4hhyc6vfjhfgyyrheudjV}
\end{EQA}
\end{lemma}

\begin{proof}
It follows from \eqref{uVmGu2C1C2gen} and \eqref{uVmGu2C1C2genV} that 
\begin{EQA}
	\tr (\CONSTIF \Id_{\dimp} + \GP^{2})^{-1} \CONSTIF \CONSTV^{-1}
	& \leq &
	\tr ( \IF_{\GP}^{-1} \VP^{2} )
	\leq 
	\tr (\CONSTIF^{-1} \Id_{\dimp} + \GP^{2})^{-1} \CONSTIF^{-1} \CONSTV \, .
\label{iuhevopedfijwu8efjfw3}
\end{EQA}
Further we may proceed as in the proof of Lemma~\ref{LeffdimG}.
\end{proof}

\Section{Sobolev classes and smooth priors}
This section illustrates the introduced notions and results for a typical situation of 
a \( (\smp,\CGP) \)-smooth prior with \( \gp_{j}^{2} = \CGP j^{2\smp} \).

\Section{Properties of the sub-projector \( \proj_\GP \)}
\label{SsubprojPG}
For the sub-projector \( \proj_{\GP} \) from \eqref{9bvgu7fh54rfgy7byrfuer4}, this section analyzes 
the operator \( \Id_{\dimp} - \proj_{\GP} \) which naturally appears 
in the evaluation of the bias \( \upsvs_{\GP} - \upsvs \)%
\ifapp{; see Section~\ref{Ssmoothbias}.}{.}
It turns out that the main characteristic of \( \proj_{\GP} \) is the index \( \mm \) defined by \eqref{gklghi43mjd7skje8vhr}.
The sub-projector \( \proj_{\GP} \) approximates the projector on the space \( \VV_{\mm} \).
The quality of approximation is controlled by the growth rate of the eigenvalues \( \gp_{j}^{2} \):
the faster is this rate the better is the approximation \( \proj_{\GP} \approx \Proj_{\mm} \).
To illustrate this point, we consider the situation with two operators 
\( \Id_{\nsize} - \proj_{\GP} \) and \( \Id_{\nsize} - \proj_{\GP_{0}} \) for two different penalizing matrices 
\( \GP^{2} \) and \( \GP_{0}^{2} \) with the same characteristic \( \mm \).
We slightly change the notations and assume that 
\begin{EQA}
	\GP^{2}
	&=&
	\CGP^{-1} \diag\{ \gp_{1}^{2}, \ldots, \gp_{\dimp}^{2} \},
	\qquad
	\GP_{0}^{2}
	=
	\CGP_{0}^{-1} \diag\{ \gp_{1,0}^{2}, \ldots, \gp_{\dimp,0}^{2} \} ,
\label{946e7f4e346dyeehdy}
\end{EQA}
with some fixed constants \( \CGP \) and \( \CGP_{0} \)
and growing sequences \( (\gp_{j}^{2}) \) and \( (\gp_{j,0}^{2}) \)
satisfying 
\begin{EQA}
	\CGP^{-1} \gp_{\mm}^{2}
	& \approx &
	\CGP_{0}^{-1} \gp_{\mm,0}^{2}
	\approx 
	\nsize .
\label{ywhsywshxetwtsgdew2w}
\end{EQA}
To simplify the presentation we later assume that these relations in \eqref{ywhsywshxetwtsgdew2w} are precisely fulfilled. 
We also assume that 
\begin{EQA}
	\gp_{j,0}^{2}/\gp_{\mm,0}^{2}
	& \leq &
	\gp_{j}^{2}/\gp_{\mm}^{2} \, ,
	\qquad
	j \leq \mm ,
\label{ohcybd7cyr6duyhcyf}
\end{EQA}
meaning that \( \gp_{j,0}^{2} \) grows faster than \( \gp_{j}^{2} \).

\begin{lemma}
\label{Lsmoothbiase}
Let \( \IF \) satisfy \eqref{uVmGu2C1C2gen}, and let \( \GP^{2} \) and \( \GP_{0}^{2} \)
be diagonal penalizing matrices satisfying
\eqref{ywhsywshxetwtsgdew2w} and \eqref{ohcybd7cyr6duyhcyf} for some \( \mm \leq \dimp \).
Then
\begin{EQA}
	(\Id_{\dimp} - \proj_{\GP}) \Proj_{\mm}
	& \leq &
	\CONSTIF^{2} (\Id_{\dimp} - \proj_{\GP_{0}}) \Proj_{\mm} \, ,
	\\
	(\CONSTIF + 1)^{-1} \Proj_{\mm}^{c}
	& \leq &
	(\Id_{\dimp} - \proj_{\GP}) \Proj_{\mm}^{c}
	\leq 
	\Proj_{\mm}^{c} \, ,
\label{8h8h8g5g77hduejwyddyh}
\end{EQA}
yielding
\begin{EQA}
	\Id_{\dimp} - \proj_{\GP}
	& \leq &
	\CONST (\Id_{\dimp} - \proj_{\GP_{0}}) ,
	\qquad
	\CONST = \CONSTIF^{2} \vee (\CONSTIF + 1) .
\label{4t66he5heudshwyseuw}
\end{EQA}
\end{lemma}

\begin{proof}
The definition implies \( \Id_{\dimp} - \proj_{\GP} = (\IF + \GP^{2})^{-1} \GP^{2} \)
and by \eqref{uVmGu2C1C2gen}
\begin{EQA}
	(\CONSTIF \, n \, \Id_{\dimp} + \GP^{2})^{-1} \GP^{2}
	\leq 
	(\IF + \GP^{2})^{-1} \GP^{2}
	& \leq &
	(\CONSTIF^{-1} \, n \, \Id_{\dimp} + \GP^{2})^{-1} \GP^{2} .
\label{6dtetwdtdywywddtetw}
\end{EQA}
Further, for \( j \leq \mm \), \eqref{ywhsywshxetwtsgdew2w} and
\( \gp_{j}^{2}/\gp_{\mm}^{2} \leq \gp_{j,0}^{2}/\gp_{\mm,0}^{2} \) imply
\( \CGP^{-1} \gp_{j}^{2} \leq \CGP_{0}^{-1} \gp_{j,0}^{2} \), i.e. 
\begin{EQA}
	\GP^{2} \, \Proj_{\mm}
	& \leq &
	\GP_{0}^{2} \, \Proj_{\mm} \, .
\label{8efdvjdieewwedhfhdvh}
\end{EQA}
Therefore,
\begin{EQA}
	(\Id_{\dimp} - \proj_{\GP}) \Proj_{\mm}
	&=&
	\IF_{\GP}^{-1} \GP^{2} \, \Proj_{\mm}
	\leq 
	(\CONSTIF^{-1} \, n \, \Id_{\dimp} + \GP^{2})^{-1} \GP^{2} \, \Proj_{\mm}
	\\
	& \leq &
	(\CONSTIF^{-1} \, n \, \Id_{\dimp} + \GP_{0}^{2})^{-1} \GP_{0}^{2} \, \Proj_{\mm}
	\leq 
	\CONSTIF^{2} (\IF + \GP_{0}^{2})^{-1} \GP_{0}^{2} \Proj_{\mm} \, .
\label{bdgdigeifgwiegfeydycgs}
\end{EQA}
After restricting to the orthogonal complement \( \VV_{\mm}^{c} \), both operators 
\( \Id_{\dimp} - \proj_{\GP} \) and \( \Id_{\dimp} - \proj_{\GP_{0}} \) behave nearly as 
projectors:
in view of \( \gp_{j}^{2} \geq \nsize \) for \( j > \mm \)
\begin{EQA}
	(\CONSTIF + 1)^{-1} \Proj_{\mm}^{c}
	& \leq &
	(\Id_{\dimp} - \proj_{\GP}) \Proj_{\mm}^{c}
	\leq 
	\Proj_{\mm}^{c}
\label{jhh6dtd5eghdgdetegke5}
\end{EQA}
and similarly for \( \Id_{\dimp} - \proj_{\GP_{0}} \).
\end{proof}

Finally we evaluate the quantity \( \| \QP (\Id - \proj_{\GP}) \upsv \| \)
assuming \( \| \GP_{0} \upsv \| \) bounded.

\begin{lemma}
\label{Lsmoothbias}
It holds for any \( \QP \colon \R^{\dimp} \to \R^{\dimq} \) and any \( \GP^{2} \)
\begin{EQA}
	\| \QP (\Id - \proj_{\GP}) \upsv \|^{2}
	& \leq &
	\| \QP \IF_{\GP}^{-1} \QP^{\T} \| \, \| \GP \upsv \|^{2} .
\label{yc65e5eftw3dtwfefdtw}
\end{EQA}
Moreover, let \( \IF \) satisfy \eqref{uVmGu2C1C2gen}, and let \( \GP^{2} \) and \( \GP_{0}^{2} \)
be diagonal penalizing matrices satisfying
\eqref{ywhsywshxetwtsgdew2w} and \eqref{ohcybd7cyr6duyhcyf} for some \( \mm \leq \dimp \).
Then
\begin{EQA}
	\| \QP (\Id - \proj_{\GP}) \upsv \|
	& \leq &
	\CONST \| \QP \IF_{\GP_{0}}^{-1} \QP^{\T} \|^{1/2} \, \| \GP_{0} \upsv \|
\label{6ew6td366cetyeyegey}
\end{EQA}
with \( \CONST \) from \eqref{4t66he5heudshwyseuw}.
\end{lemma}

\begin{proof}
It holds \( (\Id - \proj_{\GP}) \upsv = \IF_{\GP}^{-1} \GP^{2} \upsv \) and in view of \( \GP^{2} \leq \IF_{\GP} \)
\begin{EQA}
	\| \QP (\Id - \proj_{\GP}) \upsv \|
	& = &
	\| \QP \IF_{\GP}^{-1} \GP^{2} \upsv \|
	\leq 
	\| \QP \IF_{\GP}^{-1/2} \| \, \| \IF_{\GP}^{-1/2} \, \GP^{2} \upsv \|
	\leq 
	\| \QP \IF_{\GP}^{-1} \QP^{\T} \|^{1/2} \, \| \GP \upsv \| .
\label{uchfuehe6fr3ehdjd}
\end{EQA}
For the second statement, we apply \eqref{4t66he5heudshwyseuw} of Lemma~\ref{Lsmoothbiase}.
As \( \IF_{\GP}^{-1} \, \GP^{2} \leq \CONST \IF_{\GPa}^{-1} \, \GPa^{2} \) implies
\( \IF_{\GPa}^{1/2} \, \IF_{\GP}^{-1} \, \GP^{2} \leq \CONST \IF_{\GPa}^{-1/2} \, \GPa^{2} \), it follows in a similar way
\begin{EQA}
	\| \QP \IF_{\GP}^{-1} \GP^{2} \upsv \|
	& \leq &
	\| \QP \IF_{\GPa}^{-1/2} \| \, \| \IF_{\GPa}^{1/2} \, \IF_{\GP}^{-1} \, \GP^{2} \upsv \|
	\\
	& \leq &
	\CONST \| \QP \IF_{\GPa}^{-1/2} \| \, \| \IF_{\GPa}^{-1/2} \, \GPa^{2} \upsv \|
	\leq 
	\CONST \| \QP \IF_{\GPa}^{-1/2} \| \, \| \GPa \upsv \|
\label{uchfuehe6fr3ehdjda}
\end{EQA}
as required.
\end{proof}

For the case of \( \QP = \IF_{\GP_{0}}^{1/2} \), we obtain a corollary of \eqref{6ew6td366cetyeyegey}
\begin{EQA}
	\| \IF_{\GP_{0}}^{1/2} (\Id - \proj_{\GP}) \upsv \|
	& \leq &
	\CONST \| \GP_{0} \upsv \| .
\label{p0p0h9h6f3w33rfhjk}
\end{EQA}

\section{Some results for Gaussian quadratic forms}
\label{SGaussqf}
\def\xxe{\xx_{\ex}}

\Section{Moments of a Gaussian quadratic form}
\label{SmomentqfG}
Let \( \gaussv \) be standard normal in \( \R^{\dimp} \) for \( \dimp \leq \infty \).
Given a self-adjoint trace operator \( \BBH \), consider a quadratic form 
\( \bigl\langle \BBH \gaussv, \gaussv \bigr\rangle \).

\begin{lemma}
\label{Gaussmoments}
It holds
\begin{EQA}
	\E \bigl\langle \BBH \gaussv, \gaussv \bigr\rangle 
	&=& 
	\tr \BBH ,
	\\ 
	\Var \bigl\langle \BBH \gaussv, \gaussv \bigr\rangle 
	&=& 
	2 \tr \BBH^{2} .
\label{EAarAtrA2trA2}
\end{EQA}
Moreover, 
\begin{EQA}
	\E \bigl( \bigl\langle \BBH \gaussv, \gaussv \bigr\rangle - \tr \BBH \bigr)^{2}
	&=&
	2 \tr \BBH^{2}  ,
	\\
	\E \bigl( \bigl\langle \BBH \gaussv, \gaussv \bigr\rangle - \tr \BBH \bigr)^{3}
	&=&
	8 \tr \BBH^{3} ,
	\\
	\E \bigl( \bigl\langle \BBH \gaussv, \gaussv \bigr\rangle - \tr \BBH \bigr)^{4}
	&=&
	48 \tr \BBH^{4} + 12 (\tr \BBH^{2})^{2} ,
\label{2pG2trD2DGm22m2}
\end{EQA}
and
\begin{EQA}
	\E \bigl\langle \BBH \gaussv, \gaussv \bigr\rangle^{2}
	&=&
	(\tr \BBH)^{2} + 2 \tr \BBH^{2},
	\\
	\E \bigl\langle \BBH \gaussv, \gaussv \bigr\rangle^{3}
	& = &
	(\tr \BBH)^{3} + 6 \tr \BBH \,\, \tr \BBH^{2} + 8 \tr \BBH^{3} ,
	\\
	\E \bigl\langle \BBH \gaussv, \gaussv \bigr\rangle^{4}
	& = &
	(\tr \BBH)^{4} + 12 (\tr \BBH)^{2} \tr \BBH^{2}
	+ 32 (\tr \BBH) \tr \BBH^{3}
	+ 48 \tr \BBH^{4} + 12 (\tr \BBH^{2})^{2} ,
\label{2pG2trD2DGm22m2}
	\\
	\Var \bigl\langle \BBH \gaussv, \gaussv \bigr\rangle^{2}
	& = &
	8 (\tr \BBH)^{2} \tr \BBH^{2}
	+ 32 (\tr \BBH) \tr \BBH^{3}
	+ 48 \tr \BBH^{4} + 8 (\tr \BBH^{2})^{2} .
\label{2pG2trD2DGm22m4}
\end{EQA}
Moreover, if \( \BBH \leq \Id_{\dimp} \) and \( \dimA = \tr \BBH \), then \( \tr \BBH^{m} \leq \dimA \| \BBH \|^{m-1} \) for 
\( m \geq 1 \) and
\begin{EQA}[rcccl]
	\E \bigl\langle \BBH \gaussv, \gaussv \bigr\rangle^{2}
	& \leq &
	\dimA^{2} + 2 \dimA \| \BBH \|
	&\leq &
	(\dimA + \| \BBH \|)^{2},
	\\
	\E \bigl\langle \BBH \gaussv, \gaussv \bigr\rangle^{3}
	& \leq &
	\dimA^{3} + 6 \dimA^{2} \| \BBH \| + 8 \dimA \| \BBH \|^{2}
	&\leq &
	(\dimA + 2 \| \BBH \|)^{3},
	\\
	\E \bigl\langle \BBH \gaussv, \gaussv \bigr\rangle^{4}
	& \leq &
	\dimA^{4} + 12 \dimA^{3} \| \BBH \|
	+ 44 \dimA^{2} \| \BBH \|^{2}
	+ 48 \dimA \| \BBH \|^{3}
	&\leq &
	(\dimA + 3 \| \BBH \|)^{4},
\label{2pG2trD2DGm22m2}
	\\
	\Var \bigl\langle \BBH \gaussv, \gaussv \bigr\rangle^{2}
	& \leq &
	8 \dimA^{3} + 40 \dimA^{2} \| \BBH \| + 48 \dimA \| \BBH \|^{2}.
\label{2pG2trD2DGm22m4}
\end{EQA}
\end{lemma}

\begin{proof}
Let \( \chi = \gauss^{2} - 1 \) for \( \gauss \) standard normal.
Then \( \E \chi = 0 \), \( \E \chi^{2} = 2 \), \( \E \chi^{3} = 8 \), \( \E \chi^{4} = 60 \).
Without loss of generality assume \( \BBH \) diagonal: \( \BBH = \diag(\lambda_{1},\lambda_{2},\ldots,\lambda_{\dimp}) \).
Then 
\begin{EQA}
	\xi
	\eqdef
	\bigl\langle \BBH \gaussv, \gaussv \bigr\rangle - \tr \BBH
	&=&
	\sum_{j=1}^{\dimp} \lambda_{j} (\gauss_{j}^{2} - 1) ,
\label{j1ljgj2m1}
\end{EQA}
where \( \gauss_{j} \) are i.i.d. standard normal. 
This easily yields
\begin{EQA}
	\E \xi^{2}
	&=&
	\sum_{j=1}^{\dimp} \lambda_{j}^{2} \E (\gauss_{j}^{2} - 1)^{2}
	=
	\E \chi^{2} \, \tr \BBH^{2} 
	=
	2 \tr \BBH^{2}  ,
	\\
	\E \xi^{3}
	&=&
	\sum_{j=1}^{\dimp} \lambda_{j}^{3} \E (\gauss_{j}^{2} - 1)^{3}
	=
	\E \chi^{3} \, \tr \BBH^{3} 
	=
	8 \tr \BBH^{3} ,
	\\
	\E \xi^{4}
	&=&
	\sum_{j=1}^{\dimp} \lambda_{j}^{4} (\gauss_{j}^{2} - 1)^{4}
	+ \sum_{i\neq j} \lambda_{i}^{2} \lambda_{j}^{2} \E (\gauss_{i}^{2} - 1)^{2} \E (\gauss_{j}^{2} - 1)^{2}
	\\
	&=&
	\bigl( \E \chi^{4} - 3 (\E \chi^{2})^{2} \bigr) \tr \BBH^{4} + 3 (\E \chi^{2} \, \tr \BBH^{2})^{2}
	=
	48 \tr \BBH^{4} + 12 (\tr \BBH^{2})^{2} ,
\label{2pG2trD2DGm22m2}
\end{EQA}
ensuring
\begin{EQA}
	\E \bigl\langle \BBH \gaussv, \gaussv \bigr\rangle^{2}
	&=&
	\bigl( \E \bigl\langle \BBH \gaussv, \gaussv \bigr\rangle \bigr)^{2} 
	+ \E \xi^{2}
	= 
	(\tr \BBH)^{2} + 2 \tr \BBH^{2},
	\\
	\E \bigl\langle \BBH \gaussv, \gaussv \bigr\rangle^{3}
	& = &
	\E \bigl( \xi + \tr \BBH \bigr)^{3}
	=
	(\tr \BBH)^{3} + \E \xi^{3}
	+ 3 \tr \BBH \,\, \E \xi^{2}
	\\
	&=&
	(\tr \BBH)^{3} + 6 \tr \BBH \,\, \tr \BBH^{2} + 8 \tr \BBH^{3} ,
\label{2pG2trD2DGm22m2}
\end{EQA}
and 
\begin{EQA}
	\Var \bigl\langle \BBH \gaussv, \gaussv \bigr\rangle^{2}
	& = &
	\E \bigl( \xi + \tr \BBH \bigr)^{4}
	- \bigl( \E \bigl\langle \BBH \gaussv, \gaussv \bigr\rangle \bigr)^{2}
	\\
	&=&
	\bigl( \tr \BBH \bigr)^{4} + 6 (\tr \BBH)^{2} \E \xi^{2} + 4 \tr \BBH \, \E \xi^{3} + \E \xi^{4}
	- \bigl( (\tr \BBH)^{2} + 2 \tr \BBH^{2} \bigr)^{2}
	\\
	&=& 
	8 (\tr \BBH)^{2} \tr \BBH^{2}
	+ 32 (\tr \BBH) \tr \BBH^{3}
	+ 48 \tr \BBH^{4} + 8 (\tr \BBH^{2})^{2} .
\label{2pG2trD2DGm22m4}
\end{EQA}
This implies the results of the lemma.
\end{proof}

Now we compute the exponential moments of centered and non-centered quadratic forms.

\begin{lemma}
\label{Lqfexpmom}
Let \( \| \BBH \|_{\oper} = \lambda \) and \( \gaussv \sim \ND(0,\Id_{\dimp}) \).
Then for any \( \mu \in (0,\lambda^{-1}) \), 
\begin{EQA}
	\E \exp \Bigl\{ \frac{\mu}{2} \bigl( \langle \BBH \gaussv, \gaussv \rangle - \dimA \bigr) \Bigr\}
	&=&
	\det(\Id - \mu \BBH)^{-1/2} \, .
\label{m2v241m41m}
\end{EQA}
Moreover, with \( \dimA = \tr \BBH \) and \( \vH^{2} = \tr \BBH^{2} \)
\begin{EQA}
	\log \E \exp \Bigl\{ \frac{\mu}{2} \bigl( \langle \BBH \gaussv, \gaussv \rangle - \dimA \bigr) \Bigr\}
	& \leq &
	\frac{\mu^{2} \vH^{2}}{4 (1 - \lambda \mu)} \, .
\label{m2v241m41mb}
\end{EQA}
If \( \BBH \) is positive semidefinite, \( \lambda_{j} \geq 0 \), then 
\begin{EQA}
	\log \E \exp \Bigl\{ - \frac{\mu}{2} \bigl( \langle \BBH \gaussv, \gaussv \rangle - \dimA \bigr) \Bigr\}
	& \leq &
	\frac{\mu^{2} \vH^{2}}{4} \, .
\label{m2v241m41mbn}
\end{EQA}
\end{lemma}

\begin{proof}
W.l.o.g. assume \( \lambda = 1 \).
Let \( \lambda_{j} \) be the eigenvalues of \( \BBH \), 
\( |\lambda_{j}| \leq 1 \).
By an orthogonal transform, one can reduce the statement to the case of a diagonal matrix 
\( \BBH = \diag\bigl( \lambda_{j} \bigr) \). 
Then \( \langle \BBH \gaussv, \gaussv \rangle = \sum_{j=1}^{\dimp} \lambda_{j} \gauss_{j}^{2} \) and 
by independence of the \( \gauss_{j} \)'s
\begin{EQA}
	&& \nquad
	\E \Bigl\{ \frac{\mu}{2} \langle \BBH \gaussv, \gaussv \rangle  \Bigr\}
	=
	\prod_{j=1}^{\dimp} \E \exp \Bigl( \frac{\mu}{2} \lambda_{j} \eps_{j}^{2} \Bigr)
	=
	\prod_{j=1}^{\dimp} \frac{1}{\sqrt{1 - \mu \lambda_{j}}} 
	=
	\det \bigl( \Id - \mu \BBH \bigr)^{-1/2} .
\label{dOImuBm12EB}
\end{EQA}
Below we use the simple bound: 
\begin{EQ}[rcl]
\label{lo1uusk2iukkp}
	- \log(1 - u) - u
	&=&
	\sum_{k=2}^{\infty} \frac{u^{k}}{k}
	\leq 
	\frac{u^{2}}{2} \sum_{k=0}^{\infty} u^{k} 
	=
	\frac{u^{2}}{2 (1 - u)} \, ,
	\qquad 
	u \in (0,1),
	\\
	- \log(1 - u) + u
	&=&
	\sum_{k=2}^{\infty} \frac{u^{k}}{k}
	\leq 
	\frac{u^{2}}{2} \, ,
	\qquad \qquad
	u \in (-1,0).
\label{lo1uusk2iukk}
\end{EQ}
Now it holds 
\begin{EQA}
	&& \nquad
	\log \E \Bigl\{ \frac{\mu}{2} \bigl( \langle \BBH \gaussv, \gaussv \rangle - \dimA \bigr) \Bigr\}
	=
	\log \det(\Id - \mu \BBH)^{-1/2} - \frac{\mu \, \dimA}{2}
	\\
	&=&
	- \frac{1}{2} \sum_{j=1}^{\dimp} \bigl\{ \log(1 - \mu \lambda_{j}) + \mu \lambda_{j} \bigr\}
	\leq 
	\sum_{j=1}^{\dimp} \frac{\mu^{2} \lambda_{j}^{2}}{4 (1 - \mu)} 
	=
	\frac{\mu^{2} \vH^{2}}{4 (1 - \mu)} \, .
\label{m2v241m4mj1pd}
\end{EQA}
The last statement can be proved similarly.
\end{proof}

Now we consider the case of a non-centered quadratic form
\( \langle \BBH \gaussv,\gaussv \rangle/2 + \langle \Av,\gaussv \rangle \) for a fixed vector \( \Av \).

\begin{lemma}
\label{Lexpmomnoncen}
Let \( \lambda_{\max}(\BBH) < 1 \). 
Then for any \( \Av \)
\begin{EQA}
	\E \exp\Bigl\{ \frac{1}{2}\langle \BBH \gaussv,\gaussv \rangle + \langle \Av,\gaussv \rangle \Bigr\}
	&=&
	\exp\Bigl\{ \frac{\| (\Id - \BBH)^{-1/2} \Av \|^{2}}{2} \Bigr\} \, \det(\Id - \BBH)^{-1/2} .
\label{EeBf12BggA}
\end{EQA}
Moreover, for any \( \mu \in (0,1) \)
\begin{EQA}
	&& \nquad
	\log \E \exp\Bigl\{ 
		\frac{\mu}{2} \bigl( \langle \BBH \gaussv,\gaussv \rangle - \dimA \bigr) + \langle \Av,\gaussv \rangle 
	\Bigr\}
	\\
	&=&
	\frac{\| (\Id - \mu \BBH)^{-1/2} \Av \|^{2}}{2} + \log \det(\Id - \mu \BBH)^{-1/2} - \mu \dimA 
	\\
	& \leq &
	\frac{\| (\Id - \mu \BBH)^{-1/2} \Av \|^{2}}{2} + \frac{\mu^{2} \vH^{2}}{4 (1 - \mu)} \, .
\label{EeBf12BggAmu}
\end{EQA}
\end{lemma}

\begin{proof}
Denote \( \av = (\Id - \BBH)^{-1/2} \Av \). 
It holds by change of variables \( (\Id - \BBH)^{1/2} \xv = \uv \) for \( \CONSTi_{\dimp} = (2\pi)^{-\dimp/2} \)
\begin{EQA}
	&& \nquad
	\E \exp\Bigl\{ \frac{1}{2}\langle \BBH \gaussv,\gaussv \rangle + \langle \Av,\gaussv \rangle \Bigr\}
	=
	\CONSTi_{\dimp}
	\int \exp\Bigl\{ - \frac{1}{2}\langle (\Id - \BBH) \xv,\xv \rangle + \langle \Av,\xv \rangle \Bigr\} d\xv
	\\
	&=&
	\CONSTi_{\dimp}
	\det(\Id - \BBH)^{-1/2}
	\int \exp\Bigl\{ - \frac{1}{2} \| \uv \|^{2} + \langle \av,\uv \rangle \Bigr\} d\uv
	=
	\det(\Id - \BBH)^{-1/2} \, 	\ex^{\| \av \|^{2}/2}  	.
\label{EeBf12BggAp}
\end{EQA}
The last inequality \eqref{EeBf12BggAmu} follows by \eqref{m2v241m41mb}.
\end{proof}

\Section{Deviation bounds for Gaussian quadratic forms}
\label{SdevboundGauss}
The next result explains the concentration effect of \( \| \QP \xiv \|^{2} \)
for a centered Gaussian vector \( \xiv \sim \ND(0,\HVB^{2}) \) and a linear operator \( \QP \colon \R^{\dimp} \to \R^{\dimq} \),
\( \dimp,\dimq \leq \infty \).
We use a version from \cite{laurentmassart2000}.
For completeness, we present a simple proof of the upper bound.

\begin{theorem}
\label{TexpbLGA}
\label{Lxiv2LD}
\label{Cuvepsuv0}
Let \( \xiv \sim \ND(0,\HVB^{2}) \) be a Gaussian element in \( \R^{\dimp} \) and let
\( \QP \colon \R^{\dimp} \to \R^{\dimq} \) be such that \( \BBH = \QP \HVB^{2} \QP^{\T} \) 
is a trace operator in \( \R^{\dimq} \).
Then with \( \dimH = \tr(\BBH) \), \( \vH^{2} = \tr(\BBH^{2}) \), and 
\( \supH = \| \BBH \| \), it holds for each \( \xx \geq 0 \)
\begin{EQA}
\label{Pxiv2dimAvp12}
	\P\Bigl( \| \QP \xiv \|^{2} - \dimH > 2 \vH \, \sqrt{\xx} + 2 \supH \xx \Bigr)
	& \leq &
	\ex^{-\xx} ,
	\\
	\P\Bigl( \| \QP \xiv \|^{2} - \dimH \leq - 2 \vH \, \sqrt{\xx} \Bigr)
	& \leq &
	\ex^{-\xx} .
\label{Pxiv2dimAvp12m}
\end{EQA}
It also implies 
\begin{EQA}
	\P\bigl( \bigl| \| \QP \xiv \|^{2} - \dimH \bigr| > \zq_{2}(\BBH,\xx) \bigr)
	& \leq &
	2 \ex^{-\xx} ,
\label{PxivTBBdimA2vp}
\end{EQA}
with
\begin{EQA}
	\zq_{2}(\BBH,\xx)
	& \eqdef &
	2 \vH \, \sqrt{\xx} + 2 \supH \xx \,\, .
\label{zqdefGQF}
\end{EQA}
%
\end{theorem}

\begin{proof}
W.l.o.g. assume that \( \supH = \| \BBH \| = 1 \).
We use the identity \( \| \QP \xiv \|^{2} = \langle \BBH \gaussv, \gaussv \rangle \) with
 \( \gaussv \sim \ND(0,\Id_{\dimq}) \).
We apply the exponential Chebyshev inequality: with \( \mu > 0 \)
\begin{EQA}
	\P\Bigl( \langle \BBH \gaussv, \gaussv \rangle - \dimH > \zq_{2}(\BBH,\xx) \Bigr)
	& \leq &
	\E \exp \Bigl( \frac{\mu}{2} \bigl( \langle \BBH \gaussv, \gaussv \rangle - \dimH \bigr) - \frac{\mu \, \zq_{2}(\BBH,\xx)}{2} 
	\Bigr) \, .
\label{PBggiz2E2mz2}
\end{EQA}
Given \( \xx > 0 \), fix \( \mu < 1 \) by the equation
\begin{EQA}
	\frac{\mu}{1 - \mu} 
	&=&
	\frac{2 \sqrt{\xx}}{\vH} \, 
	\quad \text{ or } \quad
	\mu^{-1} 
	=
	1 + \frac{\vH}{2 \sqrt{\xx}} \, .
\label{1v2sxm12m1m}
\end{EQA}
Let \( \lambda_{j} \) be the eigenvalues of \( \BBH \), 
\( |\lambda_{j}| \leq 1 \).
It holds with \( \dimH = \tr \BBH \) in view of \eqref{m2v241m41mb}
\begin{EQA}
	&& \nquad
	\log \E \Bigl\{ \frac{\mu}{2} \bigl( \langle \BBH \gaussv, \gaussv \rangle - \dimH \bigr) \Bigr\}
	\leq 
	\frac{\mu^{2} \vH^{2}}{4 (1 - \mu)} \, .
\label{m2v241m4mj1p}
\end{EQA}
For \eqref{Pxiv2dimAvp12}, it remains to check that the choice \( \mu \) by \eqref{1v2sxm12m1m} yields
\begin{EQA}
	\frac{\mu^{2} \vH^{2}}{4 (1 - \mu)} - \frac{\mu \, \zq_{2}(\BBH,\xx)}{2}
	& = &
	\frac{\mu^{2} \vH^{2}}{4 (1 - \mu)} - \mu \bigl( \vH \sqrt{\xx} + \xx \bigr)
	=
	\mu \Bigl( \frac{\vH \sqrt{\xx}}{2} - \vH \sqrt{\xx} - \xx \Bigr)
	=
	- \xx .
	\qquad
	\qquad
\label{m2vA241muz2}
\end{EQA}
The bound \eqref{Pxiv2dimAvp12m} is obtained similarly by applying 
the exponential Chebyshev inequality to \( - \langle \BBH \gaussv, \gaussv \rangle + \dimH \) with \( \mu = 2 \vH^{-1} \sqrt{\xx} \).
The use of \eqref{m2v241m41mbn} yields
\begin{EQA}
	&& \nquad
	\P\Bigl( \langle \BBH \gaussv, \gaussv \rangle - \dimH < - 2 \vH \sqrt{\xx} \Bigr)
	\leq
	\E \exp \Bigl\{ \frac{\mu}{2} \bigl( - \langle \BBH \gaussv, \gaussv \rangle + \dimH \bigr) - \mu \, \vH \sqrt{\xx} 
	\Bigr\}
	\\
	& \leq &
	\exp \Bigl( \frac{\mu^{2} \vH^{2}}{4} - \mu \, \vH \sqrt{\xx} \Bigr) 
	=
	\ex^{-\xx} \, 
\label{PBggiz2E2mz2}
\end{EQA}
as required.
\end{proof}

\begin{corollary}
\label{CTexpbLGAd}
Assume the conditions of Theorem~\ref{TexpbLGA}.
Then for \( \zq > \vH \)
\begin{EQA}
	\P\bigl( \bigl| \| \QP \xiv \|^{2} - \dimH \bigr| \ge \zq \bigr)
	& \leq &
	2 \exp\biggl\{ - \frac{\zq^{2}}{\bigl( \vH + \sqrt{\vH^{2} + 2 \supH \zq} \bigr)^{2}} \biggr\}
	\leq 
	2 \exp\biggl( - \frac{\zq^{2}}{4\vH^{2} + 4 \supH \zq} \biggr) .
	\qquad
	\qquad
\label{3z2spsp2z3z2}
\end{EQA}
\end{corollary}

\begin{proof}
Given \( \zq \), define \( \xx \) by 
\( 2 \vH \sqrt{\xx} + 2 \supH \xx = \zq \) or 
\( 2 \supH \sqrt{\xx} = \sqrt{\vH^{2} + 2 \supH \zq} - \vH \).
Then
\begin{EQA}
	\P\bigl( \| \QP \xiv \|^{2} - \dimH \ge \zq \bigr)
	& \leq &
	\ex^{-\xx} 
	=
	\exp\biggl\{ - \frac{\bigl( \sqrt{\vH^{2} + 2 \supH \zq} - \vH \bigr)^{2}}{4 \supH^{2}} \biggr\}
	=
	\exp\biggl\{ - \frac{\zq^{2}}{\bigl( \vH + \sqrt{\vH^{2} + 2 \supH \zq} \bigr)^{2}} \biggr\}.
\label{3emzmsp22z2c}
\end{EQA}
This yields \eqref{3z2spsp2z3z2} by direct calculus.
\end{proof}

Of course, bound \eqref{3z2spsp2z3z2} is sensible only if \( \zq \gg \vH \).

\begin{corollary}
\label{RsochpHsA}
Assume the conditions of Theorem~\ref{TexpbLGA}.
If also \( \BBH \geq 0 \), then 
\begin{EQA}
\label{Pxiv2dimAxx12}
	\P\Bigl( \| \QP \xiv \|^{2} \geq \zq^{2}(\BBH,\xx) \Bigr)
	& \leq &
	\ex^{-\xx} 
\end{EQA}
with 
\begin{EQA}
	\zq^{2}(\BBH,\xx)
	& \eqdef &
	\dimH + 2 \vH \, \sqrt{\xx} + 2 \supH \xx
	\leq 
	\bigl( \sqrt{\dimH} + \sqrt{2 \supH \xx} \bigr)^{2} \, .
\label{zzxxppdBlroBB}
\end{EQA}
Also
\begin{EQA}
	\P\Bigl( \| \QP \xiv \|^{2} - \dimH < - 2 \vH \, \sqrt{\xx} \Bigr)
	& \leq &
	\ex^{-\xx} .
\label{Pxiv2dimAvp12d}
\end{EQA}
\end{corollary}

\begin{proof}
The definition implies \( \vH^{2} \leq \dimH \supH \).
One can use a sub-optimal choice of the value 
\( \mu(\xx) = \bigl\{ 1 + 2 \sqrt{\supH \dimH/\xx} \bigr\}^{-1} \) yielding the statement of the corollary.
\end{proof}

As a special case, we present a bound for the chi-squared distribution 
corresponding to \( \QP = \HVB^{2} = \Id_{\dimp} \), \( \dimp < \infty \).
Then \( \BBH = \Id_{\dimp} \), \( \tr (\BBH) = \dimp \), \( \tr(\BBH^{2}) = \dimp \) and \( \supH(\BBH) = 1 \).

\begin{corollary}
\label{Cchi2p}
Let \( \gaussv \) be a standard normal vector in \( \R^{\dimp} \).
Then for any \( \xx > 0 \)
\begin{EQA}[ccl]
\label{Pxi2pm2px}
	\P\bigl( \| \gaussv \|^{2} \geq \dimp + 2 \sqrt{\dimp \, \xx} + 2 \xx \bigr)
	& \leq &
	\ex^{-\xx},
	\\
	\P\bigl( \| \gaussv \| \,\,  \geq \sqrt{\dimp} + \sqrt{2 \xx} \bigr)
	& \leq &
	\ex^{-\xx} ,
\label{Pxi2pm2px12}
	\\
	\P\bigl( \| \gaussv \|^{2} \leq \dimp - 2 \sqrt{\dimp \, \xx} \bigr)
	& \leq &
	\ex^{-\xx}	.
\label{Pxi2pm2px22}
\end{EQA}
\end{corollary}

The bound of Theorem~\ref{TexpbLGA} 
can be represented as a usual deviation bound.

\begin{theorem}
\label{CTexpbLGA}
Assume the conditions of Theorem~\ref{TexpbLGA}.
For \( \yy > 0 \), define
\begin{EQA}
	\xx(\yy)
	& \eqdef &
	\frac{(\sqrt{\yy + \dimH} - \sqrt{\dimH})^{2}}{4 \supH} \, .
\label{iuvfiiow3kboieheuf}
\end{EQA}
Then
\begin{EQA}
	\P\bigl( \| \QP \xiv \|^{2} \ge \dimH + \yy \bigr)
	& \leq &
	\ex^{- \xx(\yy)} ,
\label{3emzmsp22z2}
	\\
	\E \bigl\{ (\| \QP \xiv \|^{2} - \dimH) \Ind\bigl( \| \QP \xiv \|^{2} \ge \dimH + \yy \bigr) \bigr\}
	& \leq &
	2 \Bigl( \frac{\yy + \dimH}{\supH \, \xx(\yy)} \Bigr)^{1/2} \, \, 
	\ex^{- \xx(\yy)} \, .
	\qquad
	\quad
\label{3emzmsp22z2e}
\end{EQA}
Moreover, let \( \muH > 0 \) fulfill \( \rexH =  \muH \supH + \muH \sqrt{\supH \dimH / \xx(\yy)} < 1 \). 
Then 
\begin{EQA}
	\E \bigl\{ \ex^{\muH (\| \QP \xiv \|^{2} - \dimH)/2} \Ind( \| \QP \xiv \|^{2} \ge \dimH + \yy) \bigr\}
	& \leq &
	\frac{1}{1 - \rexH} \, \exp\{ - (1 - \rexH) \xx(\yy) \} \, .
	\qquad
\label{llkknbononjm9hig4e}
\end{EQA}
\end{theorem}

\begin{proof}
Normalizing by \( \supH \) reduces the statements to the case with \( \supH = 1 \).
Define \( \eta = \| \QP \xiv \|^{2} - \dimH \) 
and
\begin{EQA}
	\zq(\xx)
	&=&
	2 \sqrt{\dimH \, \xx} + 2 \xx .
\label{0kmuy765433udgswhhh}
\end{EQA}
Then by \eqref{Pxiv2dimAvp12} \( \P(\eta \geq \zq(\xx)) \leq \ex^{-\xx} \).
Inverting the relation \eqref{0kmuy765433udgswhhh} yields
\begin{EQA}
	\xx(\zq)
	&=&
	\frac{1}{4} \bigl( \sqrt{\zq + \dimH} - \sqrt{\dimH} \bigr)^{2}
\label{jkv78fdjryfgsdfghgj}
\end{EQA}
and \eqref{3emzmsp22z2} follows by applying \( \zq = \yy \).
Further, 
\begin{EQA}
	\E \bigl\{ \eta \Ind(\eta \geq \yy) \bigr\}
	&=&
	\int_{\yy}^{\infty} \P(\eta \geq \zq) \, d\zq
	\leq 
	\int_{\yy}^{\infty} \ex^{ - \xx(\zq) } \, d\zq
	= 
	\int_{\xx(\yy)}^{\infty} \ex^{-\xx} \, \zq'(\xx) \, d\xx \, .
\label{zEe2Iezz2c2H23}
\end{EQA} 
As \( \zq'(\xx) = 2 + \sqrt{\dimH/\xx} \) monotonously decreases with \( \xx \), we derive
\begin{EQA}
	\E \bigl\{ \eta \Ind(\eta \geq \yy) \bigr\}
	& \leq &
	\zq'(\xx(\yy)) \ex^{-\xx(\yy)}
	=
	\frac{1}{\xx'(\yy)} \, \ex^{- \xx(\yy)}
	=
	\frac{4 \sqrt{\yy + \dimH}}{\sqrt{\yy + \dimH} - \sqrt{\dimH}} \, \ex^{- \xx(\yy)}
\label{e7ygv76bgughytuj}
\end{EQA}
and \eqref{3emzmsp22z2e} follows.

In a similar way, define \( \zqe(\xx) \) from the relation
\( \muH^{-1} \log \zqe(\xx) = \sqrt{\dimH \, \xx} + \xx \) yielding
\begin{EQA}
	\zqe(\xx)
	&=&
	\exp \bigl( \muH \sqrt{\dimH \, \xx} + \muH \, \xx \bigr) .
\label{jvcjjuvue37r6gtur4r}
\end{EQA}
The inverse relation reads
\begin{EQA}
	\xxe(\zqe)
	&=&
	\bigl( \sqrt{\muH^{-1} \log \zqe + \dimH/4} - \sqrt{\dimH/4} \bigr)^{2} .
\label{jkv78fdjryfgsdfghgjex}
\end{EQA}
Then with \( \xx(\yy) = \xxe(\ex^{\muH \yy/2}) = \bigl( \sqrt{\yy + \dimH} - \sqrt{\dimH} \bigr)^{2}/4 \)
\begin{EQA}
	\E \bigl\{ \ex^{\muH \eta/2} \Ind(\eta \geq \yy) \bigr\}
	&=&
	\int_{\ex^{\muH \yy/2}}^{\infty} \P(\ex^{\muH \eta/2} \geq \zqe) \, d\zqe
	=
	\int_{\ex^{\muH \yy/2}}^{\infty} \P(\eta \geq 2\muH^{-1} \log \zqe) \, d\zqe
	\\
	& \leq &
	\int_{\ex^{\muH \yy/2}}^{\infty} \ex^{ - \xxe(\zqe) } \, d\zqe
	= 
	\int_{\xx(\yy)}^{\infty} \ex^{-\xx} \, \zqe'(\xx) \, d\xx .
\label{zEe2Iezz2c2H23}
\end{EQA} 
Further, in view of \( \muH + 0.5 \,\muH \sqrt{\dimH/\xx} < \muH + \muH \sqrt{\dimH / \xx(\yy)} = \rexH < 1 \) for 
\( \xx \geq \xx(\yy) \), it holds
\begin{EQA}
	\zqe'(\xx)
	&=&
	\bigl( \muH + 0.5 \, \muH \sqrt{\dimH/\xx} \bigr) \exp \bigl( \muH \sqrt{\dimH \, \xx} + \muH \, \xx \bigr) 
	\leq 
	\exp \bigl( \muH \, \xx \sqrt{\dimH / \xx(\yy)} + \muH \, \xx \bigr)
	=
	\exp (\rexH \, \xx) 
\label{jcuyu3ww3jbkihjitwedk}
\end{EQA}
and  
\begin{EQA}
	\E \bigl\{ \ex^{\muH \eta/2} \Ind(\eta \geq \yy) \bigr\}
	& \leq &
	\int_{\xx(\yy)}^{\infty} \ex^{-(1 - \rexH)\xx} \, d\xx 
	=
	\frac{1}{1 - \rexH} \, \ex^{-(1 - \rexH)\xx(\yy)} \, 
\label{zEe2Iezz2c2H23}
\end{EQA} 
and \eqref{llkknbononjm9hig4e} follows.
\end{proof}



\Chapter{Gaussian comparison}
\label{SmainresGC}

This section collects some recent results on Ga\-us\-sian comparison from \cite{GNSUl2017}.
The reader is referred to that paper for an overview on the existing literature on this topic.
Throughout this section, the following notation are used. 
We write \( a \lesssim b \) (\( a \gtrsim b \)) if there exists some absolute constant \( C \) 
such that \( a \le C \, b \) (\( a \geq C \, b \) resp.). 
Similarly, \( a \asymp b \) means that there exist \( c, C \) such that \( c \, a \le b \le C \, a \).  
We assume that all random variables are defined on common probability space \( (\Omega, \mathcal{F}, \P) \) and take values in a real separable Hilbert space \(\HS \) with a scalar product \( \langle\cdot, \cdot\rangle \)  and  norm \( \| \cdot \| \). If dimension of \(\HS \) is finite and equals \(p \), we shall write \(\R^{p} \) instead of \( \HS \).  
We also denote by \(\mathcal{B}(\HS)\) the Borel \(\sigma\)-algebra. 

For a self-adjoint operator  \( \Sigma \) with eigenvalues \(\lambda_{k}(\Sigma) \), \( k \geq 1 \), 
define its operator norm \( \| \Sigma \| \), nuclear (Schatten-one) norm \( \| \Sigma \|_{1} \), and Frobenius norm \( \| \Sigma \|_{\Fr} \) by
\begin{EQA}[rcccl]
	\| \Sigma \|
	& \eqdef &
	\sup_{\|x\|=1} \| \Sigma x \|
	&=&
	\max_{k \geq 1} |\lambda_{k}(\Sigma)| \, ,
	\\
	\| \Sigma \|_{1}
	& \eqdef &
	\tr | \Sigma |
	&=& 
	\sum_{k=1}^{\infty} |\lambda_{k}(\Sigma)| \, . 
	\\
	\| \Sigma \|_{\Fr}^{2}
	& \eqdef &
	\tr \Sigma^{2} 
	&=& 
	\sum_{k=1}^{\infty} \lambda_{k}^{2}(\Sigma) \, . 
\end{EQA}  
We suppose below that \( \Sigma \) is a nuclear and \( \| \Sigma \|_{1} < \infty \).

Let \(\Sigma_{\xiv}\) be a covariance operator  of an arbitrary Gaussian random element \( \xiv \) in \(\HS\). 
By \(\{\lambda_{k\xiv}\}_{k \geq 1}\) we denote the set of its eigenvalues arranged in the non-increasing order, i.e. 
\(\lambda_{1\xiv} \geq \lambda_{2\xiv} \geq \ldots  \), and let 
\( \lambdav_{\xiv} \eqdef \diag(\lambda_{j \xiv})_{j=1}^{\infty} \). 
Note that  \(\sum_{j=1}^{\infty} \lambda_{j \xiv} < \infty\).  
Introduce the following quantities
\begin{EQA}[c]
	\Frobg_{k\xiv}^{2} 
	\eqdef 
	\sum_{j=k}^{\infty} \lambda_{j \xiv}^{2}, \quad k = 1,2,
	\label{Lambda def}
\end{EQA} 
and
\begin{EQA}
	\CONSTdlt(\Sigma_{\xiv})
	& = &
	\begin{cases}
		\Frobg_{1\xiv}^{-1} \, ,
		& \text{if } 
		3  \lambda_{1,\xiv}^{2} \le \Frobg_{1\xiv}^{2} \, ,  
		\\
		(\lambda_{1\xiv}\Frobg_{2\xiv})^{-1/2}  , 
		& \text{if } 3  \lambda_{1\xiv}^{2} > \Frobg_{1\xiv}^{2},\,\,
		3  \lambda_{2\xiv}^{2} \leq \Frobg_{2\xiv}^{2}, 
		\\
		(\lambda_{1\xiv} \lambda_{2\xiv})^{-1/2}  , 
		& \text{if } 3  \lambda_{1\xiv}^{2} > \Frobg_{1\xiv}^{2},\,\,
		3  \lambda_{2\xiv}^{2} > \Frobg_{2\xiv}^{2} .
	\end{cases}
	\label{CdelSxiSeta}
\end{EQA}
It is easy to see that \( \| \Sigma_{\xiv}\|_{\Fr} = \Frobg_{1\xiv} \). Moreover, it is straightforward to check that 
\begin{EQA}[c]
	\frac{0.9}{(\Frobg_{1\xiv}\Frobg_{2\xiv})^{1/2}}
	\leq 
	\CONSTdlt(\Sigma_{\xiv}) 
	\leq  
	\frac{1.8}{(\Frobg_{1\xiv}\Frobg_{2\xiv})^{1/2}}.
	\label{equal definitions}
\end{EQA}
Hence, \( 	\CONSTdlt(\Sigma_{\xiv})  \asymp (\Frobg_{1\xiv}\Frobg_{2\xiv})^{-1/2} \) and therefore equivalent results can be formulated in terms of any of the quantities introduced.
The following theorem is the main result of \cite{GNSUl2017}.

\begin{theorem}
	\label{l: explicit gaussian comparison}
	Let \(\xiv\) and \(\etav\) be Gaussian elements in \(\HS\) with zero mean and covariance operators 
	\(\Sigma_{\xiv}\) and \(\Sigma_{\etav}\) respectively. 
	For any \( \av \in \HS \) 
	\begin{EQA}[rcl]
		&& \nquad
		\sup_{x > 0} \left|\P( \| \xiv - \av \| \leq x) - \P( \| \etav  \| \leq x) \right| 
		\\ 
		&&
		\lesssim  
		\Bigl\{ \CONSTdlt(\Sigma_{\xiv}) + \CONSTdlt(\Sigma_{\etav}) \Bigr\}
		\bigg( \| \lambdav_{\xiv} - \lambdav_{\etav}  \|_{1} + \| \av\|^{2}\bigg).
		\label{expl_gauss}
	\end{EQA}
\end{theorem}

We see that the obtained bounds are expressed in terms of the specific 
characteristics of the matrices \( \Sigma_{\xiv} \) and \( \Sigma_{\etav} \)
such as their operator and the Frobenius norms rather than the dimension \( p \).
Another nice feature of the obtained bounds is that they do not involve 
the inverse of \( \Sigma_{\xiv} \) or \( \Sigma_{\etav} \).
In other words, small or vanishing eigenvalues of \( \Sigma_{\xiv} \) or
\( \Sigma_{\etav} \) 
do not affect the obtained bounds in the contrary to the Pinsker bound.
Similarly, only the squared norm \( \| \av \|^{2} \) of the shift \( \av \)  
shows up in the results, while the Pinsker bound involves 
\( \| \Sigma_{\xiv}^{-1/2} \av \| \) which can be very large or infinite if 
\( \Sigma_{\xiv} \) is not well conditioned.

Let us consider \(\CONSTdlt(\Sigma_{\xiv})\) in the first factor on the r.h.s of \eqref{expl_gauss}: \(\CONSTdlt(\Sigma_{\xiv}) + \CONSTdlt(\Sigma_{\etav})\). The representation \eqref{CdelSxiSeta} mimics well the three typical situations:
in the ``large-dimensional case'' with three or more significant eigenvalues
\( \lambda_{j\xiv} \), one can take 
\( \CONSTdlt(\Sigma_{\xiv}) = \| \Sigma_{\xiv} \|_{\Fr}^{-1} = \lambda_{1\xiv}^{-1} \).
In the ``two dimensional'' case, when the sum \( \Frobg_{k\xiv}^{2} \) is 
of the order \( \lambda_{k\xiv}^{2} \) for \( k=1,2 \), 
we have that \(\CONSTdlt(\Sigma_{\xiv})\)  behaves as the product \( (\lambda_{1\xiv} \lambda_{2\xiv})^{-1/2} \).
In the intermediate case of a spike model with one large eigenvalue 
\( \lambda_{1\xiv} \) and many small eigenvalues \( \lambda_{j\xiv}, j \geq 2 \),
we have that \(\CONSTdlt(\Sigma_{\xiv})\)  behaves as \( (\lambda_{1\xiv} \Frobg_{2\xiv})^{-1/2} \).

As mentioned earlier (see~\eqref{equal definitions}), the result of Theorem~\ref{l: explicit gaussian comparison} may be equivalently formulated in a ``unified'' way in terms of  \( (\Frobg_{1\xiv}\Frobg_{2\xiv})^{-1/2} \) and \( (\Frobg_{1\etav}\Frobg_{2\etav})^{-1/2} \). 
Moreover, we specify the bound \eqref{expl_gauss}  
in the ``high-dimensional'' case, \(3 \| \Sigma_{\xiv}\|^{2} \le \| \Sigma_{\xiv}\|_{\Fr}^{2}, 
3 \| \Sigma_{\etav}\|^{2} \le \| \Sigma_{\etav}\|_{\Fr}^{2}\), which means at least three significantly 
positive eigenvalues of the matrices \( \Sigma_{\xiv} \) and \( \Sigma_{\etav}\). In this case \( \Lambda_{2\xiv}^2 \geq 2 \Lambda_{1\xiv}^2/3, \Lambda_{2\eta}^2 \geq 2 \Lambda_{1\eta}^2/3\) and we get the following corollary.

\begin{corollary}
	\label{Tgaussiancomparison3}
	Let \(\xiv\) and \(\etav\) be Gaussian elements in \(\HS\) with zero mean and covariance operators 
	\(\Sigma_{\xiv}\) and \(\Sigma_{\etav}\) respectively.
	Then for any \( \av \in \HS \)
	\begin{EQA}[rcl]
		&& \nquad
		\sup_{x > 0} \left|\P( \| \xiv - \av \| \leq x) - \P( \| \etav  \| \leq x) \right| 
		\\ 
		&&
		\lesssim  
		\bigg(\frac{1}{(\Frobg_{1\xiv}\Frobg_{2\xiv})^{1/2}} + \frac{1}{(\Frobg_{1\etav}\Frobg_{2\etav})^{1/2}}\bigg) 
		\bigg( \| \lambdav_{\xiv} - \lambdav_{\etav}  \|_{1} + \| \av\|^{2}\bigg).
		\label{expl_gauss22}
	\end{EQA}
\end{corollary}

By the Weilandt--Hoffman inequality, 
\( \| \lambdav_{\xiv} - \lambdav_{\etav} \|_{1} \le \| \Sigma_{\xiv} - \Sigma_{\etav}\|_{1} \), see e.g.~\cite{MarkusEng}.
This yields the bound in terms of the nuclear norm of 
the difference \( \Sigma_{\xiv} - \Sigma_{\etav} \), which may be more useful in a number of applications.   

\begin{corollary}\label{CTgaussiancomparisonS} 
	Let \(\xiv\) and \(\etav\) be Gaussian elements in \(\HS\) with zero mean and covariance operators 
	\(\Sigma_{\xiv}\) and \(\Sigma_{\etav}\) respectively.
	Moreover, assume that
	\begin{EQA}[c]
		3 \| \Sigma_{\xiv}\|^{2} \le \| \Sigma_{\xiv}\|_{\Fr}^{2} 
		\quad \text{ and } 	\quad 
		3 \| \Sigma_{\etav}\|^{2} \le \| \Sigma_{\etav}\|_{\Fr}^{2} \, .
		\label{2plusdelta}
	\end{EQA} 
	Then for any \( \av \in \HS \) 
	\begin{EQA}
		&&
		\sup_{x > 0} \left|\P( \| \xiv - \av \| \leq x) - \P( \| \etav  \| \leq x) \right| 
		\\
		&& \qquad \lesssim  
		\bigg(
		\frac{1}{\| \Sigma_{\xiv}\|_{\Fr}} 
		+ \frac{1}{\| \Sigma_{\etav}\|_{\Fr}}
		\bigg) 
		\bigg( \| \Sigma_{\xiv} - \Sigma_{\etav} \|_{1} + \| \av\|^{2}\bigg).
		\label{expl_gauss 2}
	\end{EQA}
\end{corollary}

Since the right-hand-side of~~\eqref{expl_gauss} does not change if we exchange \( \xiv \) and \( \etav \), Theorem~\ref{l: explicit gaussian comparison} and its Corollaries hold for the balls with the same shift \( \av \). 
In particular, the following corollary is true.

\begin{corollary}\label{cor 10}
	Under conditions of Theorem~\ref{l: explicit gaussian comparison} we have 
	\begin{EQA}
		\label{expl_gauss_same_a}
		\sup_{x > 0} \Bigl|\P( \| \xiv - \av \| \leq x) - \P( \| \etav - \av \| \leq x) \Bigr| 
		& \lesssim &
		\Bigl\{ \CONSTdlt(\Sigma_{\xiv}) + \CONSTdlt(\Sigma_{\etav}) \Bigr\} 
		\Bigl( \| \lambdav_{\xiv} - \lambdav_{\etav}  \|_{1} + \| \av\|^{2}\Bigr).
	\end{EQA}	
\end{corollary}

The result of Theorem~\ref{l: explicit gaussian comparison} may be also rewritten in terms of the operator norm
\begin{EQA}[c]
	\| \Sigma_{\xiv}^{-1/2} \Sigma_{\etav} \Sigma_{\xiv}^{-1/2} - \Id \|.
\end{EQA}
Indeed, the inequality \(\| \Sigma_{\xiv} \Sigma_{\etav}\|_{1} \le \| \Sigma_{\xiv} \|_{1} \| \Sigma_{\etav} \| \) yields the following corollary.
\begin{corollary}\label{cor 1}
	Under conditions of Theorem~\ref{l: explicit gaussian comparison} we have
	\begin{EQA}
		&&	\sup_{x > 0} 
		\left|\P( \| \xiv - \av \| \leq x) - \P( \| \etav  \| \leq x) \right| 
		\\ 
		&&\qquad\qquad\qquad 
		\lesssim  
		\Bigl\{ \CONSTdlt(\Sigma_{\xiv}) + \CONSTdlt(\Sigma_{\etav}) \Bigr\} 
		\bigg( \tr \bigl( \Sigma_{\xiv} \bigr) \, 
		\| \Sigma_{\xiv}^{-1/2} \Sigma_{\etav} \Sigma_{\xiv}^{-1/2} - \Id \| 
		+ \| \av \|^{2}
		\bigg).
	\end{EQA}
\end{corollary}

We now discuss the origin of the value \( \CONSTdlt(\Sigma_{\xiv}) \) which appears in the main theorem and its corollaries. 
Analysing the proof of Theorem~\ref{l: explicit gaussian comparison} one may find out that it is necessary to get an upper bound for a probability density function (p.d.f.) \(p_{\xiv}(x)\) (resp. \(p_{\etav}(x)\)) of  \(\| \xiv\|^{2}\) (resp. \(\| \etav \|^{2}\)) and the more general p.d.f. \(p_{\xiv}(x, \av)\) of \(\| \xiv - \av\|^{2}\) for all  \(\av\in \HS \). 
The same arguments remain true for \(p_{\etav}(x)\).  The following theorem provides uniform bounds.
\begin{theorem}\label{l: density est 2}
	Let \(\xiv\) be a Gaussian element in \(\HS\) with zero mean and covariance operator \(\Sigma_{\xiv}\). 
	Then it holds for any \(\av\) that
	\begin{EQA}[c]
		\label{density bound 0}
		\sup_{x \geq 0} p_{\xiv}(x, \av) 
		\lesssim  
		\CONSTdlt(\Sigma_{\xiv})
	\end{EQA}
	with \( \CONSTdlt(\Sigma_{\xiv}) \) from \eqref{CdelSxiSeta}.
	In particular, \( \CONSTdlt(\Sigma_{\xiv}) \lesssim (\Frobg_{1\xiv}\Frobg_{2\xiv})^{-1/2} \).
\end{theorem}

Since \( \xiv \eqd \sum_{j=1}^{\infty} \sqrt{\lambda_{j\xiv}} Z_{j} \ev_{j\xiv}\), we obtain that \( \| \xiv\|^{2} \eqd \sum_{j=1}^{\infty} \lambda_{j \xiv} Z_{j}^{2} \). 
Here and in what follows \(\{\ev_{j\xiv}\ \}_{j=1}^{\infty}\) is the orthonormal basis formed by the eigenvectors of \(\Sigma_{\xiv}\) corresponding to \(\{\lambda_{j\xiv}\}_{j=1}^{\infty}\). 
In the case \(\HS = \R^{\dimp}\), \( \av = 0, \Sigma_{\xiv} \asymp \Id\) one has that 
the distribution of \( \| \xiv\|^{2} \) is close to standard \(\chi^{2}\) with 
\( \dimp \) degrees of freedom and 
\begin{EQA}[c]\label{identity case}
	\sup_{x \geq 0} p_{\xiv}(x, 0) \asymp p^{-1/2}.
\end{EQA} 
Hence, the bound~\eqref{density bound 0} gives the right dependence on \(\dimp\) because \(\CONSTdlt(\Sigma_{\xiv}) \asymp p^{-1/2}\).   
However, a lower bound for \(\sup_{x \geq 0} p_{\xiv}(x, \av) \) in the general case is still an open question. 

A direct corollary of Theorem~\ref{l: density est 2} is the following theorem which states for a rather general situation a dimension-free anti-concentration inequality for the squared norm of a Gaussian element \( \xiv \).
In the ``high dimensional situation'', this anti-concentration bound 
only involves the Frobenius norm of \(\Sigma_{\xiv}\). 

\begin{theorem}[\( \eps\)-band of the squared norm of a Gaussian element]
	\label{band of GE}
	Let \(\xiv\) be a Gaussian element in \(\HS\) with zero mean and a covariance operator \(\Sigma_{\xiv}\). 
	Then for arbitrary \(\eps > 0\), one has
	\begin{EQA}[c]
		\label{band of Gaussian2}
		\sup_{x  > 0} \P(x < \| \xiv - \av \|^{2} < x + \eps)  
		\lesssim  
		\CONSTdlt(\Sigma_{\xiv}) \, \eps
	\end{EQA}	
	with \( \CONSTdlt(\Sigma_{\xiv}) \) from \eqref{CdelSxiSeta}.
	In particular, \(\CONSTdlt(\Sigma_{\xiv}) \) can be replaced by 
	\( (\Frobg_{1\xiv} \, \Frobg_{2\xiv})^{-1/2} \).
\end{theorem}

\bibliography{exp_ts,listpubm-with-url}

\end{document}